
\input amstex
\documentstyle{amsppt}
\loadbold
\loadmsbm

\define\sigmat#1{{\sigma\kern .2em
\tilde{}_{\kern -.2em #1}}}

\def\qrfac#1#2{\left(#1\right)_{#2}} 

\def\bk{{\boldkey K}} 
\def\be{{\boldkey E}} 

\def\onehalf{\frac{1}{2}}

\def\sn{\operatorname{sn}}		
\def\cn{\operatorname{cn}}
\def\dn{\operatorname{dn}}
\def\sc{\operatorname{sc}}
\def\nc{\operatorname{nc}}
\def\dc{\operatorname{dc}}
\def\sd{\operatorname{sd}}
\def\cd{\operatorname{cd}}
\def\nd{\operatorname{nd}}
\def\ns{\operatorname{ns}}
\def\cs{\operatorname{cs}}

\def\sech{\operatorname{sech}}

\font\huge=cmb10 scaled\magstep 4
\def\contfracsymbol{\mathop{\hbox{\huge K}}}	
\def\contfrac#1#2{\lower.45em
\hbox{$\contfracsymbol\limits_{#1}^{#2}$}}

\setbox0=\hbox{$+$}
\newdimen\plusheight
\plusheight=\ht0
\def\+{\;\lower\plusheight\hbox{$+$}\;}

\setbox0=\hbox{$-$}
\newdimen\minusheight
\minusheight=\ht0
\def\-{\;\lower\minusheight\hbox{$-$}\;}

\setbox0=\hbox{$\cdots$}
\newdimen\cdotsheight     
\cdotsheight=\plusheight
\define\cds{\lower\cdotsheight\hbox{$\cdots$}}

\def\Im{\hbox{\rm Im}}
\def\Re{\hbox{\rm Re}}

\overfullrule=0pt
\magnification =1100
\vsize 9.25 truein
\hsize 6.5 truein
\voffset -.30in
\TagsOnRight

\document

\topmatter
\title Infinite families of exact sums of squares formulas, 
Jacobi elliptic functions, 
continued fractions, and Schur functions
\endtitle 
\author  Stephen C. Milne \endauthor
\thanks  S. C. Milne was partially supported by 
National Security Agency grants MDA 904-93-H-3032, 
MDA 904-97-1-0019 and MDA904--99--1--0003
\endthanks
\affil The Ohio State University \endaffil
\address Department of Mathematics, The Ohio State University, 
Columbus, Ohio, 43210
\endaddress
\email milne\@math.ohio-state.edu 
\endemail
\date 5 August 2000 (revised 31 December 2000)\enddate
\dedicatory Dedicated to the memory of Gian-Carlo Rota who
encouraged me to write this paper in the present style
\enddedicatory 
\leftheadtext{STEPHEN C. MILNE}
\rightheadtext{Exact sums of squares formulas and Jacobi elliptic
functions}
\abstract
In this paper we derive many infinite families of explicit exact
formulas involving either squares or triangular numbers, two of
which generalize Jacobi's $4$ and $8$ squares identities to
$4n^2$ or $4n(n+1)$ squares, respectively, without using cusp forms.
In fact, we similarly generalize to infinite families all of
Jacobi's explicitly stated degree $2$, $4$, $6$, $8$ Lambert series
expansions of classical theta functions. In addition, we extend
Jacobi's special analysis of $2$ squares, $2$ triangles, $6$
squares, $6$ triangles to $12$ squares, $12$ triangles, $20$
squares, $20$ triangles, respectively.  
Our $24$ squares identity leads to a different formula for
Ramanujan's tau function $\tau (n)$, when $n$ is odd.  These 
results, depending on new expansions for powers of various
products of classical theta functions, arise in the setting of Jacobi
elliptic functions, associated continued fractions, regular
C-fractions, Hankel or Tur\'anian determinants, Fourier series,
Lambert series, inclusion/exclusion, Laplace expansion formula for
determinants, and Schur functions.  The Schur function form of these
infinite families of identities are analogous to the $\eta$-function
identities  of Macdonald. Moreover, the powers $4n(n+1)$, $2n^2+n$,
$2n^2-n$ that appear in Macdonald's work also arise at appropriate
places in our analysis.  A special case of our general methods yields
a proof of the two Kac--Wakimoto conjectured identities
involving representing a positive integer by sums of $4n^2$ or
$4n(n+1)$ triangular numbers, respectively. Our $16$ and $24$ squares
identities were originally obtained via multiple basic
hypergeometric series, Gustafson's $C_{\ell}$ nonterminating 
${}_6\phi_5$ summation theorem,  and Andrews' basic hypergeometric
series proof of Jacobi's $2$, $4$, $6$, and $8$ squares identities. 
We have (elsewhere) applied symmetry and Schur function techniques
to this original approach to prove the existence of similar infinite
families of sums of squares identities for $n^2$ and $n(n+1)$
squares. Our sums of more than $8$ squares identities are not the
same as the formulas of
Mathews (1895), Glaisher (1907),  Sierpinski (1907), 
Uspensky (1913, 1925, 1928), 
Bulygin (1914, 1915), Ramanujan (1916),
Mordell (1917, 1919),  Hardy (1918, 1920), 
Bell (1919), Estermann (1936), 
Rankin (1945, 1962), Lomadze (1948), Walton (1949),  
Walfisz (1952), Ananda-Rau (1954),  
van der Pol (1954), Kr\"atzel (1961, 1962), 
Bhaskaran (1969),
Gundlach (1978),
Kac and Wakimoto (1994),
and, Liu (2000).
\endabstract 
\endtopmatter

\head 1. Introduction\endhead

In this paper we derive many infinite families of explicit exact
formulas involving either squares or triangular numbers, two of
which generalize Jacobi's \cite{117} $4$ and $8$ squares
identities to $4n^2$ or $4n(n+1)$ squares, respectively,
without using cusp forms.  In fact, we similarly generalize to 
infinite families all of Jacobi's \cite{117} explicitly stated degree
$2$, $4$, $6$, $8$ Lambert series expansions of classical theta
functions. In addition, we extend Jacobi's special analysis of $2$
squares, $2$ triangles, $6$ squares, $6$ triangles to $12$
squares, $12$ triangles, $20$ squares, $20$ triangles,
respectively.  We also utilize a special case of our general 
analysis, outlined in \cite{162}, to prove the two Kac--Wakimoto \cite{120, pp.
452} conjectured identities involving representing a positive integer by sums of
$4n^2$ or $4n(n+1)$ triangular numbers, respectively.  
Zagier in \cite{253} has also recently independently
proven these two identities.  In addition, he proved the more general 
Conjecture 7.2 of Kac--Wakimoto \cite{120, pp. 451}.  
The $n=1$ cases of these two Kac--Wakimoto conjectured identities 
for triangular numbers are equivalent to the classical identities of Legendre
\cite{140}, \cite{21, Eqns. (ii) and (iii), pp. 139}.  Our $24$ squares identity
leads to a different formula for Ramanujan's \cite{193} tau
function $\tau (n)$, when $n$ is odd.  

All of this work depends upon new expansions for powers of 
various products of the classical theta functions 
$$\spreadlines{6 pt}\allowdisplaybreaks\align 
\vartheta_3(0,q):=\ & 
\sum\limits_{j=-\infty}^{\infty}
q^{j^2},\tag 1.1\cr 
\vartheta_2 (0,q):=\ & 
\sum\limits_{j=-\infty}^{\infty}
q^{(j+1/2)^2},\tag 1.2\cr
\endalign$$
and $\vartheta_4(0,q):=\vartheta_3(0,-q)$, 
where $\vartheta_3(0,q)$, $\vartheta_2 (0,q)$, and 
$\vartheta_4 (0,q)$ are the $z=0$ cases of the theta functions
$\vartheta_3 (z,q)$, $\vartheta_2 (z,q)$, and 
$\vartheta_4 (z,q)$ in \cite{249, pp. 463--464}. We first give a
single determinant form of these expansions, then (where
applicable) a sum of determinants form, and finally a Schur function
form.  The single determinant form of our identities includes
expressing the quadratic powers $\vartheta_3(0,-q)^{4n^2}$, 
$\vartheta_3(0,-q)^{4n(n+1)}$, $\vartheta_2(0,q)^{4n^2}$, 
and $\vartheta_2(0,q^{1/2})^{4n(n+1)}$ of classical theta
functions as a constant multiple of an $n\times n$ Hankel
determinant whose entries are either certain Lambert series, or
Lambert series plus nonzero constants, for $\vartheta_2$ and
$\vartheta_3$, respectively.   Our sums of squares identities
arise from the sum of determinants form of our expansions by 
writing $\vartheta_3(0,-q)^{4n^2}$ and 
$\vartheta_3(0,-q)^{4n(n+1)}$ as explicit polynomials of
degree $n$ in $2n-1$ Lambert series similar to those in Jacobi's
$4$ and $8$ squares identities.  The Schur function form of our
analogous identities for $\vartheta_2(0,q)^{4n^2}$ and
$\vartheta_2(0,q^{1/2})^{4n(n+1)}$ completes our proof of
the Kac--Wakimoto conjectured identities for triangular
numbers in \cite{120, pp. 452}.  Depending on the analysis, we
either use $\vartheta_3(0,-q)$ or $\vartheta_4(0,q)$  in many of
our identities.  Some of the above work has already been
announced in \cite{162}.  We present a number of additional
related results in \cite{164--166}.

Our derivation of the above infinite families of identities entails
the analysis and combinatorics of Jacobi elliptic functions,
associated continued fractions, regular C-fractions, Hankel or
Tur\'anian determinants, Fourier series, Lambert series,
inclusion/exclusion, the Laplace expansion formula for
determinants, and Schur functions.  This background material is
contained in \cite{2, 15, 16, 21, 23, 43, 63, 79, 92, 93, 103, 
104, 111, 119,  125, 134, 149, 153, 184, 203, 206, 214, 247,
249, 258}.  
In order to make
the connection with divisor sums, we follow Jacobi in emphasizing
Lambert series $\sum_{r=1}^{\infty}f(r)q^r/(1-q^r)$, as defined
in \cite{8, Ex. 14, pp. 24}, as much as possible.  (For Glaisher's
historical note of where these series first appeared in Lambert's
writings see \cite{84, Section 13, pp. 163}.)
At the end of Section 2 we point out the relationship between the
Lambert series $U_{2m-1}(q)$, $G_{2m+1}(q)$, $C_{2m-1}(q)$,
$D_{2m+1}(q)$ appearing in our sums of squares and
sums of triangles identities and the Fourier expansions of the
classical Eisenstein series $E_n(\tau)$, with
$q:=\exp(2\pi i\tau)$ and $n$ an even positive integer, as
given by \cite{20, pp. 318} and \cite{203, pp. 194--195}.   
For the convenience of the reader we often reference or write down 
the $n=1$ and $n=2$ cases of our general infinite families of identities. 
These special cases can also be given short modular forms verification proofs.

The problem of representing an integer as a sum of squares of
integers is one of the chestnuts of number theory, where new
significant contributions occur infrequently over the centuries.  
The long and interesting history of this topic is surveyed in 
\cite{11, 12, 26, 43, 94, 101, 102, 188, 200, 204, 229, 241} 
and chapters 6-9 of \cite{56}.  The review article \cite{222}
presents many questions connected with representations of
integers as sums of squares.   Direct applications of sums of
squares to lattice point problems and crystallography models in
theoretical chemistry can be found in \cite{28, 90}.  One such
example is the computation of the constant $Z_N$ that occurs in
the evaluation of a certain Epstein zeta function, needed in the
study of the stability of rare gas crystals, and in that of the
so-called {\sl Madelung constants} of ionic salts.  More theoretical
applications to ``theta series'' appear in \cite{51}. 

The $s$ squares problem is to count the number 
$r_s(n)$ of integer solutions $(x_1,\dots ,x_s)$ 
of the Diophantine equation 
$$x_1^2+\cdots +x_s^2=n,\tag 1.3$$
in which changing the sign or order of the $x_i$'s give
distinct solutions.  

Diophantus (325-409 A.D.) knew that no integer of the
form $4n-1$ is a sum of two squares.  Girard conjectured
in 1632 that $n$ is a sum of two squares if and only if all
prime divisors $q$ of $n$ with $q\equiv 3\pmod{4}$ occur
in $n$ to an even power.  Fermat in 1641 gave an
``irrefutable proof'' of this conjecture.  Euler gave the
first known proof in 1749.  Early explicit formulas for 
$r_2(n)$ were given by Legendre in 1798 and Gau{\ss} in
1801.   It appears that Diophantus was aware that all 
positive integers are sums of four integral squares.  
Bachet conjectured this result in 1621, and Lagrange
gave the first proof in 1770.  

Jacobi in his famous {\sl Fundamenta Nova}
\cite{117} of 1829 introduced elliptic and theta
functions, and utilized them as tools in the study of (1.3). 
Motivated by Euler's work on $4$ squares, Jacobi observed
that the number $r_s(n)$ of integer solutions of (1.3) is
also determined by   
$$\vartheta_3(0,q)^s:=
1+\sum\limits_{n=1}^{\infty}r_s(n)q^n,\tag 1.4 $$
where $\vartheta_3(0,q)$ is the classical theta function in 
(1.1).  

Jacobi then used his theory of elliptic and theta
functions to derive remarkable identities for the 
$s=2,4,6,8$ cases of $\vartheta_3(0,q)^s$.  He
immediately obtained elegant explicit formulas for 
$r_s(n)$, where $s=2,4,6,8$.   We find it more convenient to 
work with Jacobi's equivalent identities for the $s=2,4,6,8$ cases
of $\vartheta_3(0,-q)^s$.  Dealing with powers of 
$\vartheta_3(0,-q)$ simplifies somewhat many of our identities
in Sections $5$, $7$, and $8$, especially the multiple power
series of Section $7$.  The signs and infinite products  in the
classical formulas for $\vartheta_3(0,-q)^{16}$ and 
$\vartheta_3(0,-q)^{24}$ in Theorems 1.2 and 1.3 below are
simpler.  Finally, $\vartheta_3(0,-q)^s$ arises more naturally in
the basic hypergeometric series or multiple basic hypergeometric
series analysis in \cite{4} and \cite{6, pp. 506--508}, or  
\cite{167}, respectively.  

We recall Jacobi's identities from \cite{117, Eqn. (10.), Section 40;
Eqn. (36.), Section 40} and \cite{117, Eqn. (7.), Section 42} for
the $s=4$ and $8$ cases of $\vartheta_3(0,-q)^s$ in the
following theorem.      
\proclaim{Theorem 1.1 (Jacobi)}
$$\spreadlines{6 pt}\allowdisplaybreaks\align 
\vartheta_3(0,-q)^4 =\ & 
1-8\sum\limits_{r=1}^{\infty}(-1)^{r-1}{rq^r \over 
{1+q^r}}
=1+8\sum\limits_{n=1}^{\infty}(-1)^n
\bigl[\kern-.5em 
\sum\limits_{d\mid n,d>0 \atop 4\nmid d}
\kern-.5em d\,\bigr]q^n,\tag 1.5\cr 
\kern -5.5 em\text{and}\kern 5.5 em
\vartheta_3(0,-q)^8 =\ & 
1+16\sum\limits_{r=1}^{\infty}(-1)^r{r^3q^r \over 
{1-q^r}}
=1+16\sum\limits_{n=1}^{\infty}
\bigl[\sum\limits_{d\mid n,d>0}
\kern-.5em (-1)^dd^3\,\bigr]q^n.\tag 1.6\cr
\endalign$$
Consequently, we have respectively,  
$$r_4(n)=8\kern-.5em
\sum\limits_{d\mid n,d>0 \atop 4\nmid d}
\kern-.5em d 
\qquad \hbox {\rm and} \qquad 
r_8(n)=16\kern-.5em
\sum\limits_{d\mid n,d>0}
\kern-.5em (-1)^{n+d}d^3.\tag 1.7
$$
\endproclaim 
Note that Jacobi's identities for the $s=4$ and $8$
cases of $\vartheta_3(0,q)^s$ appear in 
\cite{117, Eqn. (8.), Section 40; Eqn. (34.), Section 40}
and \cite{117, Eqn. (8.), Section 42}, respectively.

In general it is true that 
$$r_{2s}(n)=\delta_{2s}(n)+e_{2s}(n),\tag 1.8$$
where $\delta_{2s}(n)$ is a divisor function and 
$e_{2s}(n)$ is a function of order substantially lower than
that of $\delta_{2s}(n)$.  If $2s=2,4,6,8$, then 
$e_{2s}(n)=0$, and (1.8) becomes Jacobi's formulas for 
$r_{2s}(n)$, including (1.7).  On the other hand \cite{202,
203}, if $2s>8$ then $e_{2s}(n)$ is never $0$.  The
function $e_{2s}(n)$ is the coefficient of $q^n$ in a
suitable  ``cusp form''.   The difficulties of computing (1.8),
and especially the ``nondominate'' term $e_{2s}(n)$, 
increase rapidly with $2s$.   The modular function approach
to (1.8) and the  cusp form $e_{2s}(n)$ is discussed in
\cite{123}, \cite{202}, \cite{203, pp. 241--244}.  For $2s>8$
modular function  methods such as those in \cite{95, 99, 100,
147, 172, 201}, or the more classical elliptic function approach of
\cite{14, 25, 33, 34, 39, 126, 127}, \cite{155, pp. 140--143}, 
\cite{241, pp. 204--211},  are used to determine general
formulas for $\delta_{2s}(n)$  and $e_{2s}(n)$ in (1.8). 
Explicit, exact examples of (1.8)  have been worked out for
$2\leq 2s\leq 32$.   Similarly, explicit formulas for $r_{s}(n)$
have been  found for (odd) $s<32$.  Alternate, elementary
approaches to sums of squares formulas can be found in
\cite{154, 211, 226--229}. 

We next consider classical analogs of (1.5) and (1.6)
corresponding to the $s=8$ and $12$ cases of (1.8).  

Glaisher \cite{87, pp. 210} utilized elliptic function
methods, rather than modular functions, to prove
\proclaim{Theorem 1.2 (Glaisher)}
$$\spreadlines{6 pt}\allowdisplaybreaks\align 
\vartheta_3(0,-q)^{16}
=1&+{\tfrac{32}{17}}\kern-.9em
\sum\limits_{y_1,m_1\geq 1}
\kern-.6em(-1)^{m_1}m_1^7q^{m_1y_1}\tag 1.9a\cr
&-\tfrac{512}{17}\,
q\qrfac{q;q}{\infty}^{8} \qrfac{q^2;q^2}{\infty}^{8}
\tag 1.9b\cr
\kern -6.9 em\text{where we have }\kern 6.9 em
\qrfac{q;q}{\infty}:= 
\prod\limits_{r\geq 1}(1-&q^r),\kern 3 em
\text{with}\qquad 0<|q|<1.\tag 1.10\cr
\endalign$$
\endproclaim

Glaisher took the coefficient of $q^n$ to obtain
$r_{16}(n)$.  The same formula appears in 
\cite{203, Eqn. (7.4.32), pp. 242}.  

In order to find $r_{24}(n)$, Ramanujan 
\cite{193, Entry 7, table VI}, see also 
\cite{203, Eqn. (7.4.37), pp. 243} and \cite{89}, first
proved  
\proclaim{Theorem 1.3 (Ramanujan)}
Let $\qrfac{q;q}{\infty}$ be defined by 
\hbox {\rm (1.10)}.  Then 
$$\spreadlines{6 pt}\allowdisplaybreaks\align 
\vartheta_3(0,-q)^{24}=
1&+{\tfrac{16}{691}}\kern-.9em
\sum\limits_{y_1,m_1\geq 1}
\kern-.6em(-1)^{m_1}m_1^{11}q^{m_1y_1}\tag 1.11a\cr
&-\tfrac{33152}{691}\,
q\qrfac{q;q}{\infty}^{24}
-\tfrac{65536}{691}\,
q^2\qrfac{q^2;q^2}{\infty}^{24}.\tag 1.11b\cr
\endalign$$
\endproclaim

An analysis of (1.11b) depends upon Ramanujan's 
\cite{193} tau function $\tau (n)$ defined by 
$$q\qrfac{q;q}{\infty}^{24}:=
\sum\limits_{n=1}^{\infty}\tau (n)q^n.\tag 1.12$$
For example, $\tau (1)=1$, $\tau (2)=-24$, 
$\tau (3)=252$, $\tau (4)=-1472$, 
$\tau (5)=4830$, $\tau (6)=-6048$, and  
$\tau (7)=-16744$.   Ramanujan \cite{193, 
Eqn. (103)} conjectured, and Mordell \cite{171} 
proved that $\tau (n)$ is multiplicative.   

Taking the coefficient of $q^n$ in (1.11) yields the 
classical formula \cite{107, pp. 216}, \cite{193}, 
\cite{203, Eqn. (7.4.37), pp. 243} 
for $r_{24}(n)$ given by 
$$r_{24}(n)=\, {\tfrac{16}{691}}(-1)^n
\sigma_{11}^{\dag}(n)+
{\tfrac{128}{691}}
\{(-1)^{(n-1)}259\tau (n)-512
\tau (\tfrac{n}{2})\},\tag 1.13$$
where 
$$\sigma_r^{\dag}(n):=\sum\limits_{d\mid n,d>0}
\kern-.5em (-1)^{d}d^r,\tag 1.14$$
and $\tau (x) =0$ if $x$ is not an integer.

The classical formula (1.13) can be rewritten in terms of 
divisor functions by appealing to the formula for 
$\tau (n)$ from \cite{130, Eqn. (24), pp. 34; and, Eqn.
(11.1), pp.36}, \cite{140, Eqn. (9), pp. 111}, \cite{8, Ex.
10, pp. 140} given by  
$$\tau (n)=\,
{\tfrac{65}{756}}
\sigma_{11}(n)+{\tfrac{691}{756}}
\sigma_5(n)-{\tfrac{691}{3}}
\sum\limits_{m=1}^{n-1}
\sigma_5(m)\sigma_5(n-m),\tag 1.15$$
where
$$\sigma_r(n):=\sum\limits_{d\mid n,d>0}
\kern-.5em d^r. \tag 1.16$$
Another useful formula for $\tau (n)$ is 
\cite{140, Eqn. (10), pp. 111}.   Many similar identities
attributed to Ramanujan appear in \cite{128--132}, and an 
exposition of classical results on $\tau (n)$ can be found in
\cite{27}. 

For convenience, we work with (1.15) in the form  
$$q\qrfac{q;q}{\infty}^{24}=\,
{\tfrac{65}{756}}V_{11}(q)+
{\tfrac{691}{756}}V_{5}(q)-
{\tfrac{691}{3}}V_{5}^2(q),\tag 1.17$$
where
$$V_{s}\equiv V_{s}(q) :=
\sum\limits_{r=1}^{\infty}{r^sq^r \over 
{1-q^r}}
=\sum\limits_{n=1}^{\infty}
\bigl[\sum\limits_{d\mid n,d>0}
\kern-.5em d^s\,\bigr]q^n
=\sum\limits_{y_1,m_1\geq 1}
\kern-.6em m_1^{s}q^{m_1y_1}.\tag 1.18$$

The generating function version of applying (1.15) to 
(1.13) now becomes
\proclaim{Theorem 1.4}
Let $G_s(q)$ and $V_{s}(q)$ be defined by 
\hbox {\rm (1.23)} and \hbox {\rm (1.18)}, respectively.  
Then, 
$$\spreadlines{8 pt}\allowdisplaybreaks\align 
\vartheta_3(0,-q)^{24}=
1&+{\tfrac{16}{691}}G_{11}(q)-
{\tfrac{538720}{130599}}V_{11}(q)-
{\tfrac{1064960}{130599}}V_{11}(q^2)\cr
&-{\tfrac{8288}{189}}V_{5}(q)-
{\tfrac{16384}{189}}V_{5}(q^2)+
{\tfrac{33152}{3}}V_{5}^2(q)+
{\tfrac{65536}{3}}V_{5}^2(q^2).\tag 1.19\cr
\endalign$$
\endproclaim

If all we wanted to do is write $\vartheta_3(0,-q)^{24}$ as a
sum of products of at most $6$ or $3$ Lambert series, we
just take either the $6$th or $3$rd power of both sides of
(1.5) or (1.6), respectively.  A slightly more interesting
expansion of $\vartheta_3(0,-q)^{24}$ as a sum of products of
at most $3$ Lambert series results from applying the
formula for $\tau (n)$ in \cite{8, Ex. 13, pp. 24} to
(1.11)-(1.14). Equation (1.19) which expresses 
$\vartheta_3(0,-q)^{24}$ as a sum of products of at most two
Lambert series lies deeper.  

One of the main motivations for this paper was to 
generalize Theorem 1.1 to $4n^2$ or $4n(n+1)$
squares, respectively, without using cusp forms 
such as (1.9b) and (1.11b), while still utilizing just polynomials
of degree $n$ in $2n-1$ Lambert series similar to either (1.5) or
(1.6), respectively.   This condition on the maximal number  
$n$ of Lambert series factors in each term is consistent with the
$\eta$-function identities in Appendix I of  Macdonald \cite{151},
the two Kac--Wakimoto conjectured identities for triangular
numbers in \cite{120, pp. 452}, and the
above discussion of (1.19).  Essentially, expansions of
$\vartheta_3(0,-q)^{N}$, for $N$ large,  into a polynomial of at
most degree $M$ in Lambert series is ``trivial'' if $M=\alpha N$,
and lies deeper if $M=\alpha \sqrt {N}$, where
$\alpha >0$ is a constant.  
   
We carry out the above program in Theorems 5.4 and
5.6 below.   Here, we state the $n=2$ cases, which
determine different formulas for $16$ and $24$ squares. 
\proclaim{Theorem 1.5}
$$\spreadlines{6 pt}\allowdisplaybreaks\align 
\vartheta_3(0,-q)^{16}= \ & 
1-{\tfrac{32}{3}}\left(U_1+U_3+U_5\right)
+{\tfrac{256}{3}}\left(U_1U_5-U_3^2\right),
\tag 1.20\cr
\kern -6.3 em\text{where}\kern 6.3 em     
U_{s}\equiv U_{s}(q) := \ &
\sum\limits_{r=1}^{\infty}(-1)^{r-1}{r^sq^r \over 
{1+q^r}}
=\sum\limits_{n=1}^{\infty}
\bigl[\sum\limits_{d\mid n,d>0}
\kern-.5em (-1)^{d+n/d}d^s\,\bigr]q^n\cr
= \ & \sum\limits_{y_1,m_1\geq 1}
\kern-.6em(-1)^{y_1+m_1}m_1^{s}q^{m_1y_1}.
\tag 1.21\cr
\kern -12.6 em \text{and}\kern 12.6 em\quad &\cr
\vartheta_3(0,-q)^{24}= \ & 
1+{\tfrac{16}{9}}\left(17G_3+8G_5+2G_7\right)
+{\tfrac{512}{9}}\left(G_3G_7-G_5^2\right),
\tag 1.22\cr
\kern -6.05 em\text{where}\kern 6.05 em     
G_{s}\equiv G_{s}(q) := \ &
\sum\limits_{r=1}^{\infty}(-1)^{r}{r^sq^r \over{1-q^r}}
=\sum\limits_{n=1}^{\infty}
\bigl[\sum\limits_{d\mid n,d>0}
\kern -.5 em (-1)^{d}d^s\,\bigr]q^n\cr
= \ &\sum\limits_{y_1,m_1\geq 1}
\kern -.6 em(-1)^{m_1}m_1^{s}q^{m_1y_1}.
\tag 1.23\cr
\endalign$$
\endproclaim

A more compact way of writing Theorem 1.5 as a single
determinant is provided by 
\proclaim{Theorem 1.6}
$$\spreadlines{6 pt}\allowdisplaybreaks\align 
\vartheta_3(0,-q)^{16}= \ & 
{\tfrac{1}{3}}\det\vmatrix  16U_1-2 & 16U_3+1 \\
  &    \\
16U_3+1 & 16U_5-2  
\endvmatrix,\tag 1.24\cr
\kern -14.75 em \text{and}\kern 14.75 em\quad &\cr
\vartheta_3(0,-q)^{24}= \ & 
{\tfrac{1}{3^2}}\det\vmatrix  16G_3+1 & 16G_5-2 \\
  &    \\
32G_5-4 & 32G_7+17  
\endvmatrix,\tag 1.25\cr
\endalign$$
where $U_{s}\equiv U_{s}(q)$ and  
$G_{s}\equiv G_{s}(q)$ are given by 
\hbox {\rm (1.21)} and \hbox {\rm (1.23)}, respectively. 
\endproclaim
 
\noindent The general case of Theorem 1.6 for 
$\vartheta_3(0,-q)^{4n^2}$ and 
$\vartheta_3(0,-q)^{4n(n+1)}$ is given by Theorems 5.3 
and 5.5.

The formulas for $r_{16}(n)$ and $r_{24}(n)$ 
corresponding to Theorem 1.5 are given by 
\proclaim{Theorem 1.7} Let $n$ be any positive integer.   
Then 
$$\spreadlines{6 pt}\allowdisplaybreaks\align 
r_{16}(n)=\ &
-(-1)^n\cdot {\tfrac{32}{3}}
\left[\sigmat{1}(n)+\sigmat{3}(n)+
\sigmat{5}(n)\right]\tag 1.26a\cr
+&(-1)^n\cdot {\tfrac{256}{3}}
\sum\limits_{m=1}^{n-1}
\left[\sigmat{1}(m)\sigmat{5}(n-m)-
\sigmat{3}(m)\sigmat{3}(n-m)
\right],\tag 1.26b\cr
\kern -9.15 em\text{where}\kern 9.15 em  
 \ &\sigmat{r}(n):=\sum\limits_{d\mid n,d>0}
\kern-.5em (-1)^{d+n/d}d^r,\tag 1.27\cr
\kern -6.5 em\text{and}\kern 6.5 em
r_{24}(n)=\ &
(-1)^n\cdot {\tfrac{16}{9}}
\left[17\cdot\sigma_3^{\dag}(n)+
8\cdot\sigma_5^{\dag}(n)+
2\cdot\sigma_7^{\dag}(n)\right]\tag 1.28a\cr
+&(-1)^n\cdot {\tfrac{512}{9}}
\sum\limits_{m=1}^{n-1}
\left[\sigma_3^{\dag}(m)\sigma_7^{\dag}(n-m)-
\sigma_5^{\dag}(m)\sigma_5^{\dag}(n-m)
\right],\tag 1.28b\cr
\endalign$$
where $\sigma_r^{\dag}(n)$ is defined by 
\hbox {\rm (1.14)}.  
\endproclaim

The elementary analysis in \cite{94, pp. 122--123, 125} applied
to (1.26) and (1.28) immediately implies that the dominate
terms  (1.26b) and (1.28b) for $r_{16}(n)$ and $r_{24}(n)$
have orders of magnitude $n^7$ and $n^{11}$, respectively.  
Furthermore, the ``remainder terms'' (1.26a) and (1.28a) 
have lower orders of magnitude $n^5$ and $n^7$,
respectively.  The dominate term estimates are consistent
with \cite{94, Eqn. (9.20), pp. 122}.   We provide a more
detailed analyis of these estimates for (1.26) and (1.28) at the
end of Section 5.

In Sections 5 and 7 we present the infinite families of
explicit exact formulas that generalize Theorems 1.1, 1.5,
and 1.6.  

In the case where $n$ is an odd integer \cite{in particular
an odd prime}, equating (1.11) and (1.22) yields two  
formulas for $\tau (n)$, which we have presented in 
\cite{162}.  The first one is somewhat similar
to (1.15), and both are different from Dyson's 
\cite{61} formula.  We first obtain
\proclaim{Theorem 1.8} Let $\tau (n)$ be defined by 
 \hbox {\rm (1.12)} and let $n$ be odd.  Then 
$$\spreadlines{6 pt}\allowdisplaybreaks\align 
259\tau (n)=\,&
{\tfrac{1}{2^3\cdot 3^2}}
\left[17\cdot 691\sigma_3(n)+8\cdot
691\sigma_5(n)+2\cdot 691\sigma_7(n)-9
\sigma_{11}(n)\right]\cr
&-{\tfrac{691\cdot 2^2}{3^2}}
\sum\limits_{m=1}^{n-1}
\left[\sigma_3^{\dag}(m)\sigma_7^{\dag }(n-m)-
\sigma_5^{\dag }(m)\sigma_5^{\dag}(n-m)
\right],\tag 1.29\cr
\endalign$$
where $\sigma_r(n)$ and $\sigma_r^{\dag}(n)$ are 
defined by \hbox {\rm (1.16)} and \hbox {\rm (1.14)}, 
respectively.
\endproclaim

\remark{Remark} We can use (1.29) to compute $\tau (n)$
in $\leq 6n\ln n$ steps when $n$ is an odd integer, and
(1.15) to compute $\tau (n)$ in $\leq 3n\ln n$ steps when
$n$ is any positive integer.  On the other hand, this may also
be done in $n^{2+\epsilon}$ steps by appealing to Euler's
infinite-product-representation algorithm (EIPRA) 
\cite{5, pp.104 } applied to $\qrfac{q;q}{\infty}^{24}$ in
(1.12).
\endremark

A different simplification involving a power series 
formulation of (1.22) leads to 
\proclaim{Theorem 1.9} Let $\tau (n)$ be defined by 
 \hbox {\rm (1.12)} and let $n\geq 3$ be odd.  Then 
$$\spreadlines{6 pt}\allowdisplaybreaks\align 
259\tau (n)=\,&
-{\tfrac{1}{2^3}}\kern-.5em 
\sum\limits_{d\mid n,d>0}
\kern-.5em d^{11}\,
+{\tfrac{691}{2^3\cdot 3^2}}\kern-.5em 
\sum\limits_{d\mid n,d>0}
\kern-.4em d^{3}
\left(17+8d^2+2d^4\right)\tag 1.30a\cr
&-{\tfrac{691\cdot 2^2}{3^2}}\kern-1.3em 
\sum\limits_{{{m_1>m_2\geq 1}\atop 
{m_1+m_2\leq n}}\atop
{gcd(m_1,m_2)\mid n}}\kern-1em 
(-1)^{m_1+m_2}(m_1m_2)^3(m_1^2-m_2^2)^2
\kern-1.3em \sum\limits_{{y_1,y_2\geq 1}\atop
{m_1y_1+m_2y_2=n}}\kern-1.2em 1.
\tag 1.30b\cr
\endalign$$
\endproclaim

\remark{Remark} The inner sum in (1.30b) counts the
number of solutions $(y_1,y_2)$ of the classical linear
Diophantine equation $m_1y_1+m_2y_2=n$.  This relates
(1.30) to the combinatorics in sections 4.6 and 4.7 of 
\cite{214}. 
\endremark

Equation (1.15) and equation (11.3) of \cite{130, pp. 36} 
in the form 
$$\tau (n)=\,
{\tfrac{691}{1800}}\sigma_{3}(n)+
{\tfrac{691}{900}}\sigma_{7}(n)-
{\tfrac{91}{600}}\sigma_{11}(n)+
{\tfrac{2764}{15}}\sum\limits_{m=1}^{n-1}
\sigma_3(m)\sigma_7(n-m),\tag 1.31$$
yield formulas for $\tau (n)$ that are similar to (1.29)
and (1.30).  That is, motivated by Theorems 1.8 and 1.9, a
linear combination of (1.15) and (1.31) leads to 
\proclaim{Theorem 1.10} Let $n$ be any positive integer, 
and let $\tau (n)$ and $\sigma_r(n)$ be defined by 
\hbox {\rm (1.12)} and \hbox {\rm (1.16)}, respectively. 
Then
$$\spreadlines{6 pt}\allowdisplaybreaks\align 
\tau (n)=\,&
{\tfrac{1}{2^3\cdot 3^5\cdot 5\cdot 7}}
\left[3\cdot 7\cdot 691\sigma_3(n)+
2^3\cdot 5\cdot 691\sigma_5(n)\right.\cr
&\left.\hskip 5em +2\cdot 3\cdot 7\cdot 691
\sigma_7(n)- 
13\cdot 241\sigma_{11}(n)\right]\cr 
&+{\tfrac{691\cdot 2^2}{3^3}}
\sum\limits_{m=1}^{n-1}
\left[\sigma_3(m)\sigma_7(n-m)-
\sigma_5(m)\sigma_5(n-m)
\right],\tag 1.32\cr
\kern -9 em \text{and}\kern 9 em\quad &\cr 
\tau (n)=\,&
{\tfrac{1}{2^3\cdot 3^5\cdot 5\cdot 7}}\kern-.5em 
\sum\limits_{d\mid n,d>0}
\kern-.4em d^{3}
\left(3\cdot 7\cdot 691+
2^3\cdot 5\cdot 691d^2\right.\cr
&\left.\hskip 7em +2\cdot 3\cdot 7\cdot 691d^4-
13\cdot 241d^8\right)\cr
&+{\tfrac{691\cdot 2^2}{3^3}}\kern-1.3em 
\sum\limits_{{{m_1>m_2\geq 1}\atop 
{m_1+m_2\leq n}}\atop
{gcd(m_1,m_2)\mid n}}\kern-1em 
(m_1m_2)^3(m_1^2-m_2^2)^2
\kern-1.3em \sum\limits_{{y_1,y_2\geq 1}\atop
{m_1y_1+m_2y_2=n}}\kern-1.2em 1.
\tag 1.33\cr
\endalign$$
\endproclaim

Equating (1.29) and (1.32) immediately gives
\proclaim{Corollary 1.11} Let $n$ be odd, and let 
$\sigma_r^{\dag}(n)$ and $\sigma_r(n)$ be 
defined by \hbox {\rm (1.14)} and \hbox {\rm (1.16)}, 
respectively. Then 
$$\spreadlines{6 pt}\allowdisplaybreaks\align
\,&2160\sum\limits_{m=1}^{n-1}
\left[\sigma_3^{\dag}(m)\sigma_7^{\dag }(n-m)-
\sigma_5^{\dag }(m)\sigma_5^{\dag}(n-m)
\right]\cr
&+186480\sum\limits_{m=1}^{n-1}
\left[\sigma_3(m)\sigma_7(n-m)-
\sigma_5(m)\sigma_5(n-m)
\right]\cr
=\,&759\sigma_3(n)-200\sigma_5(n)-
642\sigma_7(n)+83\sigma_{11}(n).\tag 1.34\cr
\endalign$$
\endproclaim

We next give a brief description of some of the essential 
ingredients in our derivation of the sums of squares and related
identities in this paper.  

The first is a classical result of Heilermann \cite{103, 104}, more
recently presented in \cite{119, Theorem 7.14, pp. 244--246}, in
which Hankel determinants whose entries are the coefficients in a
formal power series $L$ can be expressed as a certain product of
the ``numerator'' coefficients of the associated continued
fraction $J$ corresponding to $L$, provided $J$ exists.  A similar
result holds for the related $\chi$ determinants.    We apply
Heilermann's product formulas for both $n\times n$ Hankel and 
the related $\chi$ determinants to Rogers' \cite{206}, 
Stieltjes' \cite{217--219}, and Ismail and Masson's \cite{111} 
associated continued fraction and/or regular C-fraction 
expansions of the Laplace transform of a small number of Jacobi
elliptic functions such as $\sn(u,k)$, $\cn(u,k)$, $\dn(u,k)$, 
$\sn^2(u,k)$, $\sn(u,k)\cn(u,k)/\dn(u,k)$, and 
$\sn(u,k)\cn(u,k)$.  Modular transformations,  row and column
operations, and Heilermann's \cite{103, 104}, \cite{119, Theorem
7.14, pp. 244--246; Theorem 7.2, pp. 223--224}
correspondence theorems between formal power series and
both types of continued fractions enables us to carry out a
similar product formula analysis for our associated continued 
fraction and/or regular C-fraction expansions of the
{\sl formal} Laplace transform of a number of ratios of Jacobi
elliptic functions in which $\cn(u,k)$ is in the denominator.  
These include both $\sn(u,k)\dn(u,k)/\cn(u,k)$ and 
$\sn^2(u,k)\dn^2(u,k)/\cn^2(u,k)$.       
Al@-Salam and Carlitz \cite{2, pp. 97--99} have already 
applied Heilermann's product formulas to Rogers' and Stieltjes'
continued fraction expansions of the Laplace transform of
$\sn(u,k)$, $\cn(u,k)$, $\dn(u,k)$, and  $\sn^2(u,k)$ to obtain
the product formulas for the corresponding Hankel determinants.
The direct relationship between orthogonal polynomials, J-fraction
expansions, and some of these continued fraction expansions of
the Laplace transform of Jacobi elliptic functions is surveyed in
\cite{19, 38, 44, 46, 47, 49, 50, 69, 92, 109, 110, 115,
116, 119, 148, 230--234}.      

All of the above analysis produces a large number of product
formulas for $n\times n$ Hankel or related $\chi$
determinants whose entries (polynomials in $k^2$) are the
coefficients in the formal Laplace transform of the 
Maclaurin series expansion (about $u=0$) of certain ratios 
$f_1(u,k)$ of Jacobi elliptic functions, with $k$ the
modulus.  These formulas express our determinants as a
constant multiple of a simple polynomial in $k^2$.  
By an analysis similar to that in \cite{15, 16, 21, 23, 38, 258}  we
compare the Maclaurin series expansion of $f_1(u,k)$ to its Fourier
series.  This immediately leads to a formula which factors the
entries in our determinants into the product of a Lambert
series plus a constant, a simple function of $k$, and a
suitable negative integral power of the Jacobi elliptic function 
parameter $z:={}_2F_1(1/2,1/2;1;k^2)=
2\bk(k)/\pi \equiv 2\bk/\pi$, with $\bk(k)\equiv \bk$ 
the complete elliptic integral of the first kind in 
\cite{134, Eqn. (3.1.3), pp. 51}, and $k$ the modulus.
We next solve for the resulting power of $z$ in our product
formula for the $n\times n$ Hankel or related $\chi$
determinants.  Appealing to well-known equalities such as 
$z=\vartheta_3(0,q)^2$ and $zk=\vartheta_2(0,q)^2$, where 
$q:=exp(-\pi  {\bk}(\sqrt{1-k^2})/{\bk}(k))$, and simplifying,
we are able to establish the single determinant form of our
identities.   Our proofs of the single determinant form of the
sums of squares identities in Theorems 5.3 and 5.5 utilizes the
ratios of Jacobi elliptic functions $f_1(u,k):=\sc(u,k)\dn(u,k)$
and $f_1(u,k):=\sc^2(u,k)\dn^2(u,k)$, respectively.  Similarly, 
the corresponding analysis for Theorem 5.11 in our proof of the
Kac@--Wakimoto conjectures utilizes
$f_1(u,k):=\sd(u,k)\cn(u,k)$ and $f_1(u,k):=\sn^2(u,k)$.      

An inclusion/exclusion argument applied to the above single
determinant form of our identities yields a sum of determinants
formulation.  We transform these latter identities into a multiple
power series Schur function form by employing classical
properties of Schur functions \cite{153},  symmetry and
skew-symmetry arguments, row and column operations, and the
Laplace expansion formula \cite{118, pp. 396-397} for a
determinant.  

We organize our paper as follows.  In Section 2 we
compare Fourier and Maclaurin series expansions of 
various ratios $f_1(u,k)$ of Jacobi elliptic functions to
derive the Lambert series formulas we need.  Section 3
contains a summary of the associated continued fraction
and regular C-fraction expansions of the Laplace transform of
these ratios $f_1(u,k)$.  This includes a sketch of Rogers'
\cite{206} integration-by-parts derivation of several of these
expansions, followed by his application of Landen's 
transformation \cite{134, Eqn. (3.9.15), (3.9.16), (3.9.17),
pp. 78--79} to one of them.  Here, we have continued fraction
expansions of both Laplace transforms and formal Laplace
transforms.  Section 4 utilizes Theorems 7.14 and 7.2 of
\cite{119, pp. 244--246; pp. 223--224}, row and column
operations, and modular transformations to deduce our Hankel
and $\chi$ determinant evaluations from the continued
fraction expansions of Section 3.  In Section 5 we first establish
an inclusion-exclusion lemma, recall necessary elliptic function
parameter relations,  and then obtain the single determinant and
sum of determinants form of our infinite families of sums of squares
and related identities.  We also obtain (see Theorem 5.19) our 
generalization to infinite families all $21$ of Jacobi's \cite{117,
Sections 40, 41, 42} explicitly stated degree $2, 4, 6, 8$ Lambert
series expansions of classical theta functions.  An elegant
generating function involving  the number of ways of writing $N$ as
a sum of $2n$ squares and $(2n)^2$ triangular numbers appears
in Corollary 5.15.  Section 6  contains the derivation of two key
theorems which expand certain general $n\times n$ determinants,
whose entries are either constants or Lambert series, into a
multiple power series whose terms include classical Schur functions
\cite{153} as factors.   Section 7 applies the key determinant
expansion formulas in Section 6 to most of the main identities in
Section 5 to obtain the Schur function form of our infinite
families of sums of squares and related identities.   These include
the two Kac--Wakimoto
\cite{120, pp. 452} conjectured identities involving representing a positive
integer by sums of $4n^2$ or $4n(n+1)$ triangular numbers,
respectively.  The $n=1$ case is equivalent to the classical
identities of Legendre \cite{139}, \cite{21, Eqns. (ii) and (iii), pp.
139}.   In addition, we obtain our analog of these Kac--Wakimoto
identities which involves representing a positive integer by sums of
$2n$ squares and $(2n)^2$ triangular numbers.  
Motivated by the analysis in \cite{41, 42}, we use 
Jacobi's transformation of the theta functions
$\vartheta_4$ and $\vartheta_2$ to derive a direct 
connection between our identities involving $4n^2$ or
$4n(n+1)$ squares, and the identities involving $4n^2$
or $4n(n+1)$ triangular numbers.    In a different direction
we also apply the classical techniques in \cite{90, pp. 96--101}
and \cite{28, pp. 288--305} to several of the theorems in this
section to obtain the corresponding infinite families of lattice sum
transformations. Our explicit multiple power series formulas for
$16$, $24$, $36$, and $48$ squares appear in Section 8.  

The Schur function form of our infinite families of
identities are analogous to the $\eta$-function identities in
Appendix I of  Macdonald \cite{151}.  Moreover, the powers 
$4n(n+1)$, $2n^2+n$, $2n^2-n$ that appear in Macdonald
\cite{151} also arise at appropriate places in our analysis. 
An important part of our approach to the infinite families of
identities in this paper is based upon a limiting process  (computing
associated continued fraction expansions), followed by a sieving
procedure (inclusion/exclusion).  On the other hand, the
derivation of the $\eta$-function and related identities in
\cite{76, 77, 141, 142, 151} relies on first a sieving procedure,
followed by taking limits.    

Theorems 1.5 and 1.6 were originally obtained via multiple
basic hypergeometric series \cite{143, 158--161, 163, 168, 169}
and Gustafson's \cite{96} $C_{\ell}$ nonterminating 
$_6\phi_5$ summation theorem combined with Andrews'
\cite{4}, \cite{6, pp. 506--508}, \cite{79, pp. 223--226} basic
hypergeometric series proof of Jacobi's $2$, $4$, $6$, and $8$
squares identities, and computer algebra \cite{250}.  We have
in \cite{167} applied symmetry and Schur function techniques
to this original approach to prove the existence of similar infinite
families of sums of squares identities for $n^2$ or $n(n+1)$
squares, respectively. 

Our sums of more than $8$ squares identities are not the
same as the formulas of Mathews \cite{154}, 
Glaisher \cite{85--88}, 
Sierpinski \cite{211}, 
Uspensky \cite{226--228},
Bulygin \cite{33, 34}, 
Ramanujan \cite{193}, 
Mordell \cite{172, 173},  
Hardy \cite{99, 100}, 
Bell \cite{14},
Estermann \cite{64}, 
Rankin \cite{200, 201}, 
Lomadze \cite{147}, 
Walton \cite{248}, 
Walfisz \cite{246}, 
Ananda-Rau \cite{3}, 
van der Pol \cite{186}, 
Kr\"atzel \cite{126, 127}, 
Bhaskaran \cite{25},
Gundlach \cite{95}, 
Kac and Wakimoto \cite{120}, 
and Liu \cite{146}.

We have found in \cite{164--166} a number of additional new results involving or
inspired by Hankel determinants.  In \cite{166} we
apply the Hankel determinant evaluations in the present paper to
the analysis in \cite{119, pp. 244--250} to yield a large number of
more complex $\chi$ determinant evaluations.  All of these
determinant evaluations lead to new determinantal formulas for
powers of classical theta functions.  The paper \cite{164} presents some new
evaluations of Hankel determinants of Eisenstein series which
generalize the classical formula for the modular discriminant $\Delta$ in
\cite{20, Entry 12(i), pp. 326}, \cite{203, Eqn. (6.1.14), pp.
197}, and \cite{210, eqn (42), pp. 95}.  The work in \cite{164}
was motivated by F. Garvan's comments and conjectured formula
for $\Delta^2$ in \cite{78} as a $3$ by $3$ Hankel determinant
of classical Eisenstein series after seeing an earlier version of the
present paper.  Finally, the paper \cite{165} derives a new formula for the 
modular discriminant $\Delta$ and a corresponding new formula that 
expresses Ramanujan's tau function as the difference of two positive
integers.  These integers arise naturally in areas as diverse as 
affine Lie super-algebras, higher-dimensional unimodular lattices, 
combinatorics, and number theory.  The analysis in \cite{165} utilizes special
cases of the methods in Section 5 to extend part of the work of Jacobi in
\cite{117}.  

\head 2. Fourier Expansions and Lambert Series\endhead

In this section we compare Fourier and Maclaurin series
expansions of various ratios $f_1(u,k)$ of Jacobi elliptic
functions to derive the Lambert series formulas we need. 
We first recall the classical Fourier expansions from 
\cite{93, 117, 134, 249}, write them as a double sum as in 
\cite{15, 16, 21, 38, 258}, and finally equate coefficients
with the corresponding Maclaurin series to obtain the Lambert
series formulas.  At the end of this section we point out the 
relationship between the Lambert series $U_{2m-1}(q)$, 
$G_{2m+1}(q)$, $C_{2m-1}(q)$, $D_{2m+1}(q)$ appearing
in our sums of squares and sums of triangles identities and the
Fourier expansions of the classical Eisenstein series
$E_n(\tau)$, with $q:=\exp(2\pi i\tau)$ and $n$ an
even positive integer, as given by \cite{20, pp. 318} and
\cite{203, pp. 194--195}.   

We utilize the Jacobi elliptic function parameter
$$z:={}_2F_1\left[\left.\matrix 
\tfrac{1}{2}\ ,\ \tfrac{1}{2}\\ 1\endmatrix\ \right|
\ k^2\right]  = 2\bk(k)/\pi
\equiv 2\bk/\pi,\tag 2.1$$ with 
$$\bk(k)\equiv \bk := \int_0^1 \frac{dt}
{\sqrt{(1-t^2)(1-k^2 t^2)}}
		= \tfrac{\pi}{2}\ 
{}_2F_1\left[\left.\matrix \tfrac{1}{2}\ ,\ 
\tfrac{1}{2}\\
1\endmatrix\ \right|\ k^2\right]\tag 2.2$$   
the complete elliptic integral of the first kind in 
\cite{134, Eqn. (3.1.3), pp. 51}, and $k$ the modulus.  
We also sometimes use the complete elliptic integral of the
second kind 
$$\be(k)\equiv \be := \int_0^1 
\sqrt{\frac{1-k^2 t^2}{1-t^2}}\,dt
		= \tfrac{\pi}{2}\ 
{}_2F_1\left[\left.\matrix \tfrac{1}{2}\ ,\ 
-\tfrac{1}{2}\\
1\endmatrix\ \right|\ k^2\right],\tag 2.3$$
and the complementary modulus $k':=\sqrt{1-k^2}$.
Finally, we take 
$$q:=exp(-\pi  {\bk}(\sqrt{1-k^2})/{\bk}(k))\tag 2.4$$

The classical Fourier expansions from \cite{93,  pp.
911-912}, \cite{117, Sections 39, 41, and 42}, 
\cite{134, pp. 222-225}, and 
\cite{249, pp. 510-520} that we need
are summarized by the following theorem.  Essentially
all of these Fourier expansions originally appeared in 
\cite{117, Sections 39, 41, and 42}. 
\proclaim{Theorem 2.1} Let  $z:=2\bk/\pi$, as in 
\hbox{\rm(2.1)}, with $\bk$ and $\be$ given by  
\hbox{\rm(2.2)} and \hbox{\rm(2.3)}, respectively.  Also
take $k':=\sqrt{1-k^2}$ and $q$ as in \hbox{\rm(2.4)}. 
Then,  
$$\spreadlines{6 pt}\allowdisplaybreaks\align
\sn(u,k)	& = \frac{2\pi}{k\bk}
\sum\limits_{n=1}^{\infty}
		\frac{q^{n-\onehalf}}{1-q^{2n-1}} \sin
\tfrac{(2n-1)u}{z}
	\tag 2.5\cr
\cn(u,k)	& = \frac{2\pi}{k\bk}
\sum\limits_{n=1}^{\infty}
		\frac{q^{n-\onehalf}}{1+q^{2n-1}} \cos
\tfrac{(2n-1)u}{z}
	\tag 2.6\cr
\dn(u,k)	& = \frac{\pi}{2\bk} +
	\frac{2\pi}{\bk} \sum\limits_{n=1}^{\infty}
		\frac{q^n}{1+q^{2n}} \cos \tfrac{2nu}{z}
	\tag 2.7\cr
\sd(u,k)	& = \frac{2\pi}{kk'\bk}
\sum\limits_{n=1}^{\infty}
		\frac{(-1)^{n-1} q^{n-\onehalf}}{1+q^{2n-1}} \sin
\tfrac{(2n-1)u}{z}
	\tag 2.8\cr
\cd(u,k)	& = \frac{2\pi}{k\bk}
\sum\limits_{n=1}^{\infty}
		\frac{(-1)^{n-1} q^{n-\onehalf}}{1-q^{2n-1}} 
\cos\tfrac{(2n-1)u}{z}
	\tag 2.9\cr
\nd(u,k)	& = \frac{\pi}{2k'\bk} +
	\frac{2\pi}{k'\bk} \sum\limits_{n=1}^{\infty}
		\frac{(-1)^n q^{n}}{1+q^{2n}} \cos \tfrac{2nu}{z}
	\tag 2.10\cr
\sc(u,k)	& = \frac{\pi}{2k'\bk}\tan \tfrac{u}{z} +
	\frac{2\pi}{k'\bk} \sum\limits_{n=1}^{\infty}
		\frac{(-1)^n q^{2n}}{1+q^{2n}} 
\sin \tfrac{2nu}{z}\tag 2.11\cr
\dc(u,k)	& = \frac{\pi}{2\bk}\sec\tfrac{u}{z} +
	\frac{2\pi}{\bk} \sum\limits_{n=1}^{\infty}
		\frac{(-1)^{n-1} q^{2n-1}}{1-q^{2n-1}} 
\cos\tfrac{(2n-1)u}{z}
	\tag 2.12\cr
\nc(u,k)	& = \frac{\pi}{2k'\bk}\sec\tfrac{u}{z} -
	\frac{2\pi}{k'\bk} \sum\limits_{n=1}^{\infty}
		\frac{(-1)^{n-1} q^{2n-1}}{1+q^{2n-1}} 
\cos\tfrac{(2n-1)u}{z}
	\tag 2.13\cr
\sn^2(u,k)
	& = \frac{\bk-\be}{k^2 \bk} - \frac{2\pi^2}{k^2 \bk^2} 
\sum\limits_{n=1}^{\infty}
		\frac{n q^n}{1-q^{2n}} \cos \tfrac{2nu}{z}
	\tag 2.14\cr 
\sc^2(u,k) & = -\frac{1}{{k'}^2}\frac{\be}{\bk}
  +\frac{1}{z^2{k'}^2}\sec^2 \tfrac{u}{z} -
\frac{8}{z^2{k'}^2}
\sum\limits_{n=1}^{\infty}
		\frac{(-1)^n n q^{2n}}{1-q^{2n}} \cos \tfrac{2nu}{z}
	\tag 2.15\cr 
\sd^2(u,k) & =\frac{\be - {k'}^2\bk}{k^2{k'}^2\bk}+
 	\frac{8}{z^2k^2{k'}^2}
\sum\limits_{n=1}^{\infty}
		\frac{(-1)^n n q^n}{1-q^{2n}} \cos \tfrac{2nu}{z}
	\tag 2.16\cr 
\frac{\sn(u,k)\ \cn(u,k)}{\dn(u,k)}
	& = \frac{4\pi}{k^2 \bk} \sum\limits_{n=1}^{\infty}
			\frac{q^{2n-1}}{1-q^{4n-2}} \sin \tfrac{(4n-2)u}{z}
	\tag 2.17\cr
\frac{\sn(u,k)\ \dn(u,k)}{\cn(u,k)}
	& = \frac{\pi}{2\bk} \tan \tfrac{u}{z}
		+ \frac{2\pi}{\bk} \sum\limits_{n=1}^{\infty}
			\frac{q^n}{1+(-1)^n q^n} \sin \tfrac{2nu}{z}
	\tag 2.18\cr
\frac{\sn(u,k)}{\cn(u,k)\ \dn(u,k)}
	& = \frac{\pi}{2{k'}^2 \bk} \tan \tfrac{u}{z}
		+ \frac{2\pi}{{k'}^2 \bk} \sum\limits_{n=1}^{\infty}
			\frac{(-1)^n q^n}{1+q^n} \sin \tfrac{2nu}{z}
	\tag 2.19\cr
\frac{\sn^2(u,k)\ \cn^2(u,k)}{\dn^2(u,k)}	
& =\frac{1}{k^4} + \frac{{k'}^2}{k^4} - 
\frac{2}{k^4}\frac{\be}{\bk}-
\frac{32}{z^2k^4}\sum\limits_{n=1}^{\infty}
		\frac{n q^{2n}}{1-q^{4n}} \cos \tfrac{4nu}{z}
\tag 2.20\cr
\frac{\sn^2(u,k)\ \dn^2(u,k)}{\cn^2(u,k)}	
& = 1-2\frac{\be}{\bk}+
\frac{1}{z^2}\sec^2 \tfrac{u}{z}-
\frac{8}{z^2}\sum\limits_{n=1}^{\infty}
		\frac{n q^n}{1-(-1)^n q^n} \cos \tfrac{2nu}{z}
\tag 2.21\cr
\frac{\sn^2(u,k)}{\cn^2(u,k)\ \dn^2(u,k)}	
& = \frac{1}{{k'}^2}-\frac{2}{{k'}^4}\frac{\be}{\bk} + 
\frac{1}{z^2{k'}^4}\sec^2 \tfrac{u}{z}-
\frac{8}{z^2{k'}^4}\sum\limits_{n=1}^{\infty}
		\frac{(-1)^n n q^n}{1-q^n} \cos \tfrac{2nu}{z}
\tag 2.22\cr
\sn(u,k)\ \dn(u,k)	& = \frac{\pi^2}{k\bk^2}
\sum\limits_{n=1}^{\infty}
		\frac{(2n-1)q^{n-\onehalf}}{1+q^{2n-1}} \sin
\tfrac{(2n-1)u}{z}
	\tag 2.23\cr
\sn(u,k)\ \cn(u,k)	& = 
	\frac{2\pi^2}{k^2\bk^2}\sum\limits_{n=1}^{\infty}
		\frac{nq^n}{1+q^{2n}} \sin \tfrac{2nu}{z}
	\tag 2.24\cr
\frac{\sn(u,k)}{\dn^2(u,k)}	& = 
\frac{\pi^2}{k{k'}^2\bk^2}
\sum\limits_{n=1}^{\infty}
		\frac{(-1)^{n-1}(2n-1) q^{n-\onehalf}}{1-q^{2n-1}} 
\sin\tfrac{(2n-1)u}{z}
	\tag 2.25\cr
\frac{\sn(u,k)\ \cn(u,k)}{\dn^2(u,k)}	& = 
-	\frac{2\pi^2}{k^2k'\bk^2} 
\sum\limits_{n=1}^{\infty}
		\frac{(-1)^n nq^{n}}{1+q^{2n}} \sin \tfrac{2nu}{z}
	\tag 2.26\cr
\endalign$$
$$\spreadlines{6 pt}\allowdisplaybreaks\align
\hskip -3 em 
\frac{\sn(u,k)}{\cn^2(u,k)} &  = 
\frac{1}{z^2{k'}^2}
\sec\tfrac{u}{z}\tan\tfrac{u}{z} -
	\frac{\pi^2}{\bk^2} \sum\limits_{n=1}^{\infty}
		\frac{(-1)^{n-1} (2n-1)q^{2n-1}}{1-q^{2n-1}} 
\sin\tfrac{(2n-1)u}{z}
	\tag 2.27\cr
\hskip -3 em 
\frac{\sn(u,k)\ \dn(u,k)}{\cn^2(u,k)} & = 
\frac{\pi^2}{4k'\bk^2}
\sec\tfrac{u}{z}\tan\tfrac{u}{z} +
	\frac{\pi^2}{k'\bk^2} \sum\limits_{n=1}^{\infty}
		\frac{(-1)^{n-1} (2n-1)q^{2n-1}}{1+q^{2n-1}} 
\sin\tfrac{(2n-1)u}{z}
	\tag 2.28\cr
\endalign$$
\endproclaim

The Fourier series expansions in (2.18) and (2.21) are
important ingredients of our derivation of the $4n^2$ and 
$4n(n+1)$ squares identities, respectively, in Section 5.  
Both of these expansions are direct consequences of
simpler Fourier expansions.  

A derivation of (2.18) is outlined in \cite{134, Ex. 6,  pp.
243}.  The Fourier expansions for $\ns(2u,k)$ and
$\cs(2u,k)$ are substituted into 
$$\sn(u,k)\ \dc(u,k) = \ns(2u,k) - \cs(2u,k).\tag 2.29$$
The terms from the sum in (2.18) with $n$ odd come from
$\ns(2u,k)$, and with $n$ even from $\cs(2u,k)$.  The
$\tan\tfrac{u}{z}$ in (2.18) follows from 
$\csc\theta - \cot\theta = \tan\tfrac{\theta}{2}$.  

In order to derive (2.21) we first apply the Pythagorean
relations for $\dn^2(u,k)$ and $\sn^2(u,k)$ to obtain
$$\spreadlines{6 pt}\allowdisplaybreaks\align
\sc^2(u,k)\ \dn^2(u,k) & = 
\sc^2(u,k)\left[1-k^2\sn^2(u,k)\right]\cr
& = \sc^2(u,k)\left[1-k^2(1-\cn^2(u,k))\right]\cr
& = (1-k^2)\sc^2(u,k) + k^2\sn^2(u,k)
\tag 2.30\cr
\endalign$$
Next, substitute (2.14) and (2.15) into (2.30).  Finally,
combine the two resulting sums termwise while appealing to
the difference of squares 
$$(1-q^{2n}) = (1 + (-1)^nq^n)(1-(-1)^nq^n).\tag 2.31$$
This last step is motivated by \cite{94, pp. 119}.  

The Fourier series expansions in (2.23)--(2.28) 
are obtained by differentiating those in 
(2.6), (2.7), (2.9), (2.10), (2.12), and (2.13).  

The Fourier expansions in Theorem 2.1 may be written as a
double sum by first expanding the $\sin\tfrac{Nu}{z}$ and
$\cos\tfrac{Nu}{z}$ as Maclaurin series, interchanging
summation, and then simplifying.  We obtain
\proclaim{Theorem 2.2} Let  $z:=2\bk/\pi$, as in 
\hbox{\rm(2.1)}, with $\bk$ and $\be$ given by  
\hbox{\rm(2.2)} and \hbox{\rm(2.3)}, respectively.  Also
take $k':=\sqrt{1-k^2}$ and $q$ as in \hbox{\rm(2.4)}. 
Then,  
$$\spreadlines{6 pt}\allowdisplaybreaks\align
zk&\cdot\sn(u,k)	 =
4\kern-.2em\sum\limits_{m=1}^{\infty}
   \frac{(-1)^{m-1}}{z^{2m-1}}
    \kern-.25em\left[\sum\limits_{r=1}^{\infty}
		\frac{(2r-1)^{2m-1}
    q^{r-\onehalf}}{1-q^{2r-1}}\right]
  \kern-.25em\frac{u^{2m-1}}{(2m-1)!}\tag 2.32\cr
zk&\cdot\cn(u,k)	 =
4\kern-.2em\sum\limits_{m=0}^{\infty}
   \frac{(-1)^{m}}{z^{2m}}
    \kern-.25em\left[\sum\limits_{r=1}^{\infty}
		\frac{(2r-1)^{2m}
    q^{r-\onehalf}}{1+q^{2r-1}}\right]
  \kern-.25em\frac{u^{2m}}{(2m)!}\tag 2.33\cr
z&\cdot\dn(u,k)	 =1+
4\kern-.2em\sum\limits_{m=0}^{\infty}
   \frac{(-1)^{m}2^{2m}}{z^{2m}}
    \kern-.25em\left[\sum\limits_{r=1}^{\infty}
		\frac{r^{2m}q^{r}}{1+q^{2r}}\right]
  \kern-.25em\frac{u^{2m}}{(2m)!}\tag 2.34\cr
zkk'&\cdot\sd(u,k)	 =
4\kern-.2em\sum\limits_{m=1}^{\infty}
   \frac{(-1)^{m-1}}{z^{2m-1}}
    \kern-.25em\left[\sum\limits_{r=1}^{\infty}
		\frac{(-1)^{r+1}(2r-1)^{2m-1}
    q^{r-\onehalf}}{1+q^{2r-1}}\right]
\kern-.25em\frac{u^{2m-1}}{(2m-1)!}\tag 2.35\cr
zk&\cdot\cd(u,k)	 =
4\kern-.2em\sum\limits_{m=0}^{\infty}
   \frac{(-1)^{m}}{z^{2m}}
    \kern-.25em\left[\sum\limits_{r=1}^{\infty}
		\frac{(-1)^{r+1}(2r-1)^{2m}
    q^{r-\onehalf}}{1-q^{2r-1}}\right]
  \kern-.25em\frac{u^{2m}}{(2m)!}\tag 2.36\cr
zk'&\cdot\nd(u,k)	 =
1+4\kern-.2em\sum\limits_{m=0}^{\infty}
   \frac{(-1)^{m}2^{2m}}{z^{2m}}
    \kern-.25em\left[\sum\limits_{r=1}^{\infty}
		\frac{(-1)^{r}r^{2m} q^{r}}{1+q^{2r}}\right]
  \kern-.25em\frac{u^{2m}}{(2m)!}\tag 2.37\cr
zk'&\cdot\sc(u,k)	 = \tan \tfrac{u}{z} +
	4\kern-.2em\sum\limits_{m=1}^{\infty}
   \frac{(-1)^{m-1}2^{2m-1}}{z^{2m-1}}
    \kern-.25em\left[\sum\limits_{r=1}^{\infty}
		\frac{(-1)^{r}r^{2m-1} q^{2r}}{1+q^{2r}}\right]
  \kern-.25em\frac{u^{2m-1}}{(2m-1)!}\tag 2.38\cr
z&\cdot\dc(u,k)	 = \sec \tfrac{u}{z} 
	+4\kern-.2em\sum\limits_{m=0}^{\infty}
   \frac{(-1)^{m}}{z^{2m}}
    \kern-.25em\left[\sum\limits_{r=1}^{\infty}
		\frac{(-1)^{r+1}(2r-1)^{2m}
    q^{2r-1}}{1-q^{2r-1}}\right]
  \kern-.25em\frac{u^{2m}}{(2m)!}\tag 2.39\cr
zk'&\cdot\nc(u,k)	 = \sec \tfrac{u}{z} 
	-4\kern-.2em\sum\limits_{m=0}^{\infty}
   \frac{(-1)^{m}}{z^{2m}}
    \kern-.25em\left[\sum\limits_{r=1}^{\infty}
		\frac{(-1)^{r-1}(2r-1)^{2m}
    q^{2r-1}}{1+q^{2r-1}}\right]
  \kern-.25em\frac{u^{2m}}{(2m)!}\tag 2.40\cr
z^2k^2&\cdot\sn^2(u,k)	 = z^2-z^2\frac{\be}{\bk}
 	-8\kern-.2em\sum\limits_{m=0}^{\infty}
   \frac{(-1)^{m}2^{2m}}{z^{2m}}
    \kern-.25em\left[\sum\limits_{r=1}^{\infty}
		\frac{r^{2m+1} q^{r}}{1-q^{2r}}\right]
  \kern-.25em\frac{u^{2m}}{(2m)!}\tag 2.41\cr
z^2{k'}^2&\cdot\sc^2(u,k)	 = -z^2\frac{\be}{\bk}
  +\sec^2 \tfrac{u}{z} +
	8\kern-.2em\sum\limits_{m=0}^{\infty}
   \frac{(-1)^{m}2^{2m}}{z^{2m}}\cr
&\kern 15 em\times
    \kern-.25em\left[\sum\limits_{r=1}^{\infty}
		\frac{(-1)^{r-1}r^{2m+1} q^{2r}}{1-q^{2r}}\right]
  \kern-.25em\frac{u^{2m}}{(2m)!}\tag 2.42\cr
z^2k^2{k'}^2&\cdot\sd^2(u,k)	
 =-z^2{k'}^2+z^2\frac{\be}{\bk}
 	-8\kern-.2em\sum\limits_{m=0}^{\infty}
   \frac{(-1)^{m}2^{2m}}{z^{2m}}\cr
&\kern 15 em\times
    \kern-.25em\left[\sum\limits_{r=1}^{\infty}
		\frac{(-1)^{r+1}r^{2m+1} q^{r}}{1-q^{2r}}\right]
  \kern-.25em\frac{u^{2m}}{(2m)!}\tag 2.43\cr
zk^2&\cdot\frac{\sn(u,k)\ \cn(u,k)}{\dn(u,k)}
	 = 8\kern-.2em\sum\limits_{m=1}^{\infty}
   \frac{(-1)^{m-1}2^{2m-1}}{z^{2m-1}}\cr
&\kern 12 em\times
    \kern-.25em\left[\sum\limits_{r=1}^{\infty}
		\frac{(2r-1)^{2m-1} q^{2r-1}}{1-q^{4r-2}}\right]
  \kern-.25em\frac{u^{2m-1}}{(2m-1)!}\tag 2.44\cr
z&\cdot\frac{\sn(u,k)\ \dn(u,k)}{\cn(u,k)}	
 = \tan \tfrac{u}{z} 
	+4\kern-.2em\sum\limits_{m=1}^{\infty}
   \frac{(-1)^{m-1}2^{2m-1}}{z^{2m-1}}\cr
&\kern 15 em\times
    \kern-.25em\left[\sum\limits_{r=1}^{\infty}
		\frac{r^{2m-1} q^{r}}{1+(-1)^{r}q^{r}}\right]
  \kern-.25em\frac{u^{2m-1}}{(2m-1)!}\tag 2.45\cr
z{k'}^2&\cdot\frac{\sn(u,k)}{\cn(u,k)\ \dn(u,k)}
	 =\tan \tfrac{u}{z}
		+4\kern-.2em\sum\limits_{m=1}^{\infty}
   \frac{(-1)^{m-1}2^{2m-1}}{z^{2m-1}}\cr
&\kern15 em\times    
\kern-.25em\left[\sum\limits_{r=1}^{\infty}
		\frac{(-1)^{r}r^{2m-1} q^{r}}{1+q^{r}}\right]
  \kern-.25em\frac{u^{2m-1}}{(2m-1)!}\tag 2.46\cr
z^2k^4&\cdot\frac{\sn^2(u,k)\ \cn^2(u,k)}{\dn^2(u,k)}	
 = z^2+z^2{k'}^2-2z^2\frac{\be}{\bk}\cr
 &\kern11em 	
-32\kern-.2em\sum\limits_{m=0}^{\infty}
   \frac{(-1)^{m}2^{4m}}{z^{2m}}
    \kern-.25em\left[\sum\limits_{r=1}^{\infty}
		\frac{r^{2m+1} q^{2r}}{1-q^{4r}}\right]
  \kern-.25em\frac{u^{2m}}{(2m)!}\tag 2.47\cr
z^2&\cdot\frac{\sn^2(u,k)\ \dn^2(u,k)}{\cn^2(u,k)}	
 = z^2-2z^2\frac{\be}{\bk}
  +\sec^2 \tfrac{u}{z}\cr 
&\kern11em	
-8\kern-.2em\sum\limits_{m=0}^{\infty}
   \frac{(-1)^{m}2^{2m}}{z^{2m}}
    \kern-.25em\left[\sum\limits_{r=1}^{\infty}
		\frac{r^{2m+1} q^{r}}{1-(-1)^{r}q^{r}}\right]
  \kern-.25em\frac{u^{2m}}{(2m)!}\tag 2.48\cr
z^2{k'}^4&\cdot\frac{\sn^2(u,k)}{\cn^2(u,k)\ \dn^2(u,k)}	
 = z^2{k'}^2-2z^2\frac{\be}{\bk}
  +\sec^2 \tfrac{u}{z}\cr 
	&\kern11em +8\kern-.2em\sum\limits_{m=0}^{\infty}
   \frac{(-1)^{m}2^{2m}}{z^{2m}}
    \kern-.25em\left[\sum\limits_{r=1}^{\infty}
		\frac{(-1)^{r-1}r^{2m+1} q^{r}}{1-q^{r}}\right]
  \kern-.25em\frac{u^{2m}}{(2m)!}\tag 2.49\cr
z^2k&\cdot\sn(u,k)\ \dn(u,k)		 =
4\kern-.2em\sum\limits_{m=1}^{\infty}
   \frac{(-1)^{m-1}}{z^{2m-1}}
    \kern-.25em\left[\sum\limits_{r=1}^{\infty}
		\frac{(2r-1)^{2m}
    q^{r-\onehalf}}{1+q^{2r-1}}\right]
  \kern-.25em\frac{u^{2m-1}}{(2m-1)!}\tag 2.50\cr
z^2k^2&\cdot\sn(u,k)\ \cn(u,k)	 =
8\kern-.2em\sum\limits_{m=1}^{\infty}
   \frac{(-1)^{m-1}2^{2m-1}}{z^{2m-1}}
    \kern-.25em\left[\sum\limits_{r=1}^{\infty}
		\frac{r^{2m}q^{r}}{1+q^{2r}}\right]
  \kern-.25em\frac{u^{2m-1}}{(2m-1)!}\tag 2.51\cr
z^2k{k'}^2&\cdot\frac{\sn(u,k)}{\dn^2(u,k)}	 =
4\kern-.2em\sum\limits_{m=1}^{\infty}
   \frac{(-1)^{m-1}}{z^{2m-1}}
    \kern-.25em\left[\sum\limits_{r=1}^{\infty}
		\frac{(-1)^{r+1}(2r-1)^{2m}
    q^{r-\onehalf}}{1-q^{2r-1}}\right]
  \kern-.25em\frac{u^{2m-1}}{(2m-1)!}\tag 2.52\cr
z^2k'k^2&\cdot
\frac{\sn(u,k)\ \cn(u,k)}{\dn^2(u,k)} =
-8\kern-.2em\sum\limits_{m=1}^{\infty}
   \frac{(-1)^{m-1}2^{2m-1}}{z^{2m-1}}
    \kern-.25em\left[\sum\limits_{r=1}^{\infty}
		\frac{(-1)^{r}r^{2m} q^{r}}{1+q^{2r}}\right]
  \kern-.25em\frac{u^{2m-1}}{(2m-1)!}\tag 2.53\cr
z^2{k'}^2&\cdot
\frac{\sn(u,k)}{\cn^2(u,k)} 	 = 
\sec \tfrac{u}{z}\tan \tfrac{u}{z}
	-4\kern-.2em\sum\limits_{m=1}^{\infty}
   \frac{(-1)^{m-1}}{z^{2m-1}}\cr
&\kern 10.65 em\times    
\kern-.25em\left[\sum\limits_{r=1}^{\infty}
		\frac{(-1)^{r+1}(2r-1)^{2m}
    q^{2r-1}}{1-q^{2r-1}}\right]
  \kern-.25em\frac{u^{2m-1}}{(2m-1)!}\tag 2.54\cr
z^2k'&\cdot
\frac{\sn(u,k)\ \dn(u,k)}{\cn^2(u,k)}  = 
\sec \tfrac{u}{z}\tan \tfrac{u}{z} 
+4\kern-.2em\sum\limits_{m=1}^{\infty}
   \frac{(-1)^{m-1}}{z^{2m-1}}\cr
&\kern 10.65 em\times
    \kern-.25em\left[\sum\limits_{r=1}^{\infty}
		\frac{(-1)^{r+1}(2r-1)^{2m}
    q^{2r-1}}{1+q^{2r-1}}\right]
  \kern-.25em\frac{u^{2m-1}}{(2m-1)!}\tag 2.55\cr
\endalign$$
\endproclaim

In order to equate coefficients of $u^N$ in Theorem 2.2 we
first need the expansions \cite{93, pp. 35} of 
$\tan\tfrac{u}{z}$, $\sec^2\tfrac{u}{z}$, 
$\sec\tfrac{u}{z}$, and  
$\sec\tfrac{u}{z}\tan\tfrac{u}{z}$ given by 
$$\spreadlines{6 pt}\allowdisplaybreaks\alignat 2
\tan\tfrac{u}{z}& =\sum\limits_{m=1}^{\infty}
\frac{2^{2m}(2^{2m}-1)|B_{2m}|}{(2m)!\,z^{2m-1}}
u^{2m-1},
&\qquad \text{ for}\  
\frac{u^2}{z^2}<\frac{\pi^2}{4},\tag 2.56\cr
\sec^2\tfrac{u}{z}& =\sum\limits_{m=0}^{\infty}
\frac{2^{2m+1}(2^{2m+2}-1)
|B_{2m+2}|}{(m+1)(2m)!\,z^{2m}}
u^{2m},&\tag 2.57\cr
\sec\tfrac{u}{z}& =\sum\limits_{m=0}^{\infty}
\frac{|E_{2m}|}{(2m)!\,z^{2m}}
u^{2m},&
\qquad \text{ for}\  
\frac{u^2}{z^2}<\frac{\pi^2}{4},\tag 2.58\cr
\sec\tfrac{u}{z}\tan\tfrac{u}{z}& =
\sum\limits_{m=1}^{\infty}
\frac{|E_{2m}|}{(2m-1)!\,z^{2m-1}}
u^{2m-1},&\tag 2.59\cr
\endalignat$$
where the Bernoulli numbers $B_n$ and Euler numbers
$E_n$ are defined  in \cite{48, pp. 48--49} by 
$${t\over {e^t-1}}:=\sum\limits_{n=0}^{\infty}
B_n{t^n\over n!}, 
\qquad \text{ for}\  |t|<2\pi,\tag 2.60$$
and 
$${2e^t\over {e^{2t}+1}}:=\sum\limits_{n=0}^{\infty}
E_n{t^n\over n!}, 
\qquad \text{ for}\  |t|<\pi/2.\tag 2.61$$
Convenient, explicit formulas for the Bernoulli numbers $B_n$ 
can be found in \cite{91}.  

We next write down the Maclaurin series expansions of the
ratios of Jacobi elliptic functions in Theorem 2.2.  The
coefficients of $u^N$ in these expansions are polynomials
in $k^2$, where $k$ is the modulus.  We have 
\definition{Definition 2.3  (Maclaurin series
expansion polynomials for Jacobi elliptic functions)}
Let the elliptic function polynomials 
$(\hbox{\rm elliptic})_m(k^2)$
of $k^2$, with $k$ the modulus, be determined by the
coefficients of the following Maclaurin series expansions:   
$$\spreadlines{6 pt}\allowdisplaybreaks\align
	\sn(u,k) & = \sum_{m=1}^\infty
		 (sn)_m(k^2)\frac{ u^{2m-1}} { (2m-1)! }
\text{;}\qquad 
	\cn(u,k)  = \sum_{m=0}^\infty
		 (cn)_m(k^2)\frac{ u^{2m}} { (2m)! },
\tag 2.62\cr
	\dn(u,k) & = \sum_{m=0}^\infty
		 (dn)_m(k^2)\frac{ u^{2m}} { (2m)! }
\text{;}\qquad 
	\sd(u,k)  = \sum_{m=1}^\infty
		 (s/d)_m(k^2)\frac{ u^{2m-1}} { (2m-1)! },
\tag 2.63\cr
	\cd(u,k) & = \sum_{m=0}^\infty
		 (c/d)_m(k^2)\frac{ u^{2m}} { (2m)! }
\text{;}\qquad 
	\nd(u,k)  = \sum_{m=0}^\infty
		 (nd)_m(k^2)\frac{ u^{2m}} { (2m)! },
\tag 2.64\cr
\sc(u,k) & = \sum_{m=1}^\infty
		 (s/c)_m(k^2)\frac{ u^{2m-1}} { (2m-1)! }
\text{;}\qquad 
	\dc(u,k)  = \sum_{m=0}^\infty
		 (d/c)_m(k^2)\frac{ u^{2m}} { (2m)! },
\tag 2.65\cr
\nc(u,k) & = \sum_{m=0}^\infty
		 (nc)_m(k^2)\frac{ u^{2m}} { (2m)! }
\text{;}\qquad 
	\sn^2(u,k)  = \sum_{m=1}^\infty
		 (sn^2)_m(k^2)\frac{ u^{2m}} { (2m)! },
\tag 2.66\cr
\sc^2(u,k) & = \sum_{m=1}^\infty
		 (s^2/c^2)_m(k^2)\frac{ u^{2m}} { (2m)! }
\text{;}\qquad 
	\sd^2(u,k)  = \sum_{m=1}^\infty
		 (s^2/d^2)_m(k^2)\frac{ u^{2m}} { (2m)! },
\tag 2.67\cr
\endalign$$

$$\spreadlines{6 pt}\allowdisplaybreaks\align
\frac{\sn(u,k)\ \cn(u,k)}{\dn(u,k)} & 
= \sum_{m=1}^\infty
		 (sc/d)_m(k^2)\frac{ u^{2m-1}} { (2m-1)! },
   \tag 2.68\cr
\frac{\sn(u,k)\ \dn(u,k)}{\cn(u,k)} & 
= \sum_{m=1}^\infty
		 (sd/c)_m(k^2)\frac{ u^{2m-1}} { (2m-1)! },	
   \tag 2.69\cr
\frac{\sn(u,k)}{\cn(u,k)\ \dn(u,k)} & 
= \sum_{m=1}^\infty
		 (s/cd)_m(k^2)\frac{ u^{2m-1}} { (2m-1)! },	
	\tag 2.70\cr
\frac{\sn^2(u,k)\ \cn^2(u,k)}{\dn^2(u,k)}	
& =\sum_{m=1}^\infty
		 (s^2c^2/d^2)_m(k^2)\frac{ u^{2m}} { (2m)! },
\tag 2.71\cr
\frac{\sn^2(u,k)\ \dn^2(u,k)}{\cn^2(u,k)}	
& = \sum_{m=1}^\infty
		 (s^2d^2/c^2)_m(k^2)\frac{ u^{2m}} { (2m)! },
\tag 2.72\cr
\frac{\sn^2(u,k)}{\cn^2(u,k)\ \dn^2(u,k)}	
& = \sum_{m=1}^\infty
		 (s^2/c^2d^2)_m(k^2)\frac{ u^{2m}} { (2m)! },
\tag 2.73\cr
\sn(u,k)\ \dn(u,k)	
& = \sum_{m=1}^\infty
		 (sd)_m(k^2)\frac{ u^{2m-1}} { (2m-1)! },
\tag 2.74\cr
\sn(u,k)\ \cn(u,k)	
& = \sum_{m=1}^\infty
		 (sc)_m(k^2)\frac{ u^{2m-1}} { (2m-1)! },
\tag 2.75\cr
\frac{\sn(u,k)}{\dn^2(u,k)}	
& = \sum_{m=1}^\infty
		 (s/d^2)_m(k^2)\frac{ u^{2m-1}} { (2m-1)! },
\tag 2.76\cr
\frac{\sn(u,k)\ \cn(u,k)}{\dn^2(u,k)}	
& = \sum_{m=1}^\infty
		 (sc/d^2)_m(k^2)\frac{ u^{2m-1}} { (2m-1)! },
\tag 2.77\cr
\frac{\sn(u,k)}{\cn^2(u,k)}	
& = \sum_{m=1}^\infty
		 (s/c^2)_m(k^2)\frac{ u^{2m-1}} { (2m-1)! },
\tag 2.78\cr
\frac{\sn(u,k)\ \dn(u,k)}{\cn^2(u,k)}	
& = \sum_{m=1}^\infty
		 (sd/c^2)_m(k^2)\frac{ u^{2m-1}} { (2m-1)! }.
\tag 2.79\cr
\endalign$$
\enddefinition  

Recursions, explicit computations, and tables of the 
polynomials $(\hbox{\rm elliptic})_m(k^2)$ in the simpler
expansions above such as for $\sn$, $\cn$, $\dn$, and 
$\sn^2$ in (2.62), (2.63), and (2.66) can be found in 
\cite{15--18, 21, 31, 49, 50, 58--60, 93, 117, 134, 170,
185, 187, 207, 208, 249, 251, 252, 258}. 
Corresponding applications to combinatorics appear in
\cite{21, 54, 58--60, 69--71, 207, 208, 223, 243, 
251, 252, 256}.  

Keeping in mind (2.56)-(2.61) and Definition 2.3, we find
that equating coefficients of $u^N$ in Theorem 2.2 yields
the Lambert series formulas in
\proclaim{Theorem 2.4} Let  $z:= 2\bk(k)/\pi
\equiv 2\bk/\pi$, as in \hbox{\rm(2.1)}, with $k$ the
modulus.  Let $q$ be as in \hbox{\rm(2.4)}, and   
let the Bernoulli numbers $B_n$ and Euler numbers
$E_n$ be defined by \hbox{\rm(2.60)} and
\hbox{\rm(2.61)},  respectively.  Take   
$(sd/c)_m(k^2)$, $(s^2d^2/c^2)_m(k^2)$, 
$(nc)_{m-1}(k^2)$, $(sc/d)_m(k^2)$, $(sn^2)_m(k^2)$, 
$(cn)_{m-1}(k^2)$, $(dn)_m(k^2)$, $(sd)_m(k^2)$, 
 $(sc)_m(k^2)$, and $(sd/c^2)_m(k^2)$ 
to be the  elliptic function polynomials of $k^2$
determined by \hbox{\rm Definition 2.3}.  Let
$m=1,2,3,\cdots$.  Then,  
$$\spreadlines{6 pt}\allowdisplaybreaks\align 
U_{2m-1}(-q):&=\sum\limits_{r=1}^{\infty}
{-r^{2m-1}q^r \over  {1+(-1)^rq^r}}\cr
&\kern-3.5em =(-1)^{m-1}{(2^{2m}-1)
\over {4m}}\cdot |B_{2m}|
+(-1)^m{z^{2m}\over {2^{2m+1}}}
\cdot (sd/c)_m(k^2),\tag 2.80\cr
G_{2m+1}(-q):&=\sum\limits_{r=1}^{\infty}
{r^{2m+1}q^r \over  {1-(-1)^{r}q^r}}\cr
&\kern-3.5em =(-1)^{m}{(2^{2m+2}-1)\over {4(m+1)}}
\cdot |B_{2m+2}|
+(-1)^{m-1}{z^{2m+2}\over {2^{2m+3}}}\cdot
(s^2d^2/c^2)_m(k^2),\tag 2.81\cr
R_{2m-2}(q):&=\sum\limits_{r=1}^{\infty}(-1)^{r+1}
{(2r-1)^{2m-2}q^{2r-1} \over  {1+q^{2r-1}}}\cr
&\kern-3.5em =(-1)^{m-1}\cdot {\tfrac {1}{4}}\cdot
|E_{2m-2}| +(-1)^m{z^{2m-1}\over 4}\,\sqrt{1-k^2}
\cdot (nc)_{m-1}(k^2),\tag 2.82\cr
C_{2m-1}(q):&=\sum\limits_{r=1}^{\infty}
{(2r-1)^{2m-1}q^{2r-1} \over  {1-q^{2(2r-1)}}}=
(-1)^{m-1}{z^{2m}k^2 \over
{2^{2m+2}}}\cdot (sc/d)_m(k^2),\tag 2.83\cr
D_{2m+1}(q):&=\sum\limits_{r=1}^{\infty}
{r^{2m+1}q^r \over  {1-q^{2r}}}=
(-1)^{m-1}{z^{2m+2}k^2 \over {2^{2m+3}}}\cdot
(sn^2)_m(k^2),\tag 2.84\cr
T_{2m-2}(q):&=\sum\limits_{r=1}^{\infty}
{(2r-1)^{2m-2}q^{r-\tfrac{1}{2}} \over {1+q^{2r-1}}}=
(-1)^{m-1}{z^{2m-1}k \over {4}}\cdot
(cn)_{m-1}(k^2),\tag 2.85\cr
N_{2m}(q):&=\sum\limits_{r=1}^{\infty}
{r^{2m}q^{r} \over {1+q^{2r}}}=
(-1)^m{z^{2m+1} \over {2^{2m+2}}}\cdot
(dn)_{m}(k^2),\tag 2.86\cr
\kern -5.35 em \text{and}\kern 5.35 em\quad 
N_{0}(q):&=\sum\limits_{r=1}^{\infty}
{q^{r} \over {1+q^{2r}}}=
-\tfrac{1}{4}+\frac{z}{4}\cdot (dn)_{0}(k^2)=
-\tfrac{1}{4}+\frac{z}{4},\tag 2.87\cr
T_{2m}(q)&=(-1)^{m+1}{z^{2m+1}k \over {4}}\cdot
(sd)_{m}(k^2),\tag 2.88\cr
N_{2m}(q)&=(-1)^{m+1}{z^{2m+1}k^2 \over {2^{2m+2}}}
\cdot(sc)_{m}(k^2),\tag 2.89\cr
R_{2m}(q)&=(-1)^{m}\cdot {\tfrac {1}{4}}\cdot
|E_{2m}| +(-1)^{m+1}{z^{2m+1}\over 4}\,\sqrt{1-k^2}
\cdot (sd/c^2)_{m}(k^2).\tag 2.90\cr
\endalign$$
\endproclaim

The Lambert series in (2.32) and (2.35) are constant
multiples of the one in (2.44), with $q\mapsto q^{1/2}$ and 
$q\mapsto i\sqrt{q}$, respectively.  The Lambert
series in (2.33), (2.34), (2.40), (2.41), (2.45), (2.48) are
transformed by $q\mapsto -q$ into the corresponding
Lambert series in (2.36), (2.37), (2.39), (2.43), (2.46),
(2.49), respectively.  The substitution $q\mapsto -q^2$ 
transforms the Lambert series in (2.45) and (2.48) into
those in (2.38) and (2.42), respectively.  Taking 
$q\mapsto q^2$ transforms the Lambert series in (2.41)
into the one in (2.47).  The relations in (2.88), (2.89),
and (2.90) also follow from combining (2.85), (2.86), and
(2.82) with the Maclaurin series expansions of the derivative
formulas for $\cn(u,k)$, $\dn(u,k)$, and $\nc(u,k)$.  
The above substitutions, combined with the corresponding
modular transformations,  is why we only need the $10$
Lambert series formulas in Theorem 2.4.  

The analysis leading to Theorem 2.4 is similar to that in
\cite{15, 16, 21, 23, 38, 258}.   The Lambert series identities
in (2.80)-(2.86) are equivalent to the identities in Tables
1(x), 1(ii), 1(xiv), 1(vii), 1(iii), 1(xv), 1(xi), respectively,
of \cite{258}.  In this paper we mainly work with the
above rational functions of $\sn(u,k)$, $\cn(u,k)$, and
$\dn(u,k)$ which {\sl do not} have $\sn(u,k)$ as a factor of
the denominator.  Our analysis in Sections $3$ and $4$ of
classical continued fraction expansions of the Laplace
transforms of these rational functions leads to our
elegant product formulas for Hankel and $\chi$ 
determinants that are crucial in Section 5.    

Table 1 of Zucker's paper \cite{258} also contains
formulas for the Lambert series determined by the Fourier
series expansions of rational functions, analogous to
several of those in Definition 2.3, of $\sn(u,k)$, $\cn(u,k)$, 
and $\dn(u,k)$, that {\sl do} have $\sn(u,k)$ as a factor of
the denominator.  This second class of rational functions 
(and their corresponding Lambert series) does not lead to
simple, ``closed'' product formulas for Hankel and $\chi$
determinants.    

The Lambert series $V_s(q)$ in (1.18) is one such
example.  Moreover, a formula for $V_s(q^2)$,
corresponding to $\ns^2(u,k)$, appears in Table 1(i) of
\cite{258}.  Even so, in \cite{164}, we discuss how the
Fourier series expansion of $\ns^2(u,k)$ has applications
to Hankel determinants of classical Eisenstein series.  

We conclude this section with a discussion of the 
relationship between the Lambert series $U_{2m-1}(q)$, 
$G_{2m+1}(q)$, $C_{2m-1}(q)$, $D_{2m+1}(q)$ and the
Fourier expansions of the classical Eisenstein series
$E_n(\tau)$ as given by \cite{20, pp. 318} and
\cite{203, pp. 194--195} in the following definition.
\definition{Definition 2.5} Let $q:=\exp(2\pi i\tau)$, 
where $\tau$ is in the upper half-plane $\Cal H$, and take 
$y:=\Im(\tau)>0$.  Let $n=1,2,3,\cdots$.  We then have 
$$\spreadlines{6 pt}\allowdisplaybreaks\align 
E_2(\tau)\equiv E_2(q):=&1-24
\sum\limits_{r=1}^{\infty}
{rq^r\over {1-q^r}} - {3\over {\pi y}},\tag 2.91\cr
\kern -13.35 em \text{and for } n\geq 2,
\kern 13.35 em\quad&\cr
E_{2n}(\tau)\equiv E_{2n}(q):=&1-{\frac{4n}{B_{2n}}}
\sum\limits_{r=1}^{\infty}
{r^{2n-1}q^r\over {1-q^r}},\tag 2.92\cr
\endalign$$
with the $B_{2n}$ the Bernoulli numbers in (2.60).
\enddefinition  

As examples, Ramanujan in \cite{193} studied the series 
$L:=E_2(\tau)+3/\pi y$, $M:=E_4(\tau)$, and 
$N:=E_6(\tau)$.  Note that we equate the $q$'s in (2.4) and 
Definition 2.5.  That is $2\tau = i {\bk}(\sqrt{1-k^2})/{\bk}(k)$. 

We now write our Lambert series $U_{2m-1}(q)$, 
$G_{2m+1}(q)$, $C_{2m-1}(q)$, and $D_{2m+1}(q)$ 
as linear combinations of the $E_{n}(\tau)$ in the following
lemma.
\proclaim{Lemma 2.6}Let $U_{2m-1}(q)$, 
$G_{2m+1}(q)$, $C_{2m-1}(q)$, and $D_{2m+1}(q)$ 
be determined by \hbox{\rm(2.80)}, \hbox{\rm(2.81)},
\hbox{\rm(2.83)}, and \hbox{\rm(2.84)}, respectively, with 
$q$ as in \hbox{\rm(2.4)}.  Take $E_{2n}(q)$ as in 
\hbox{\rm Definition 2.5}, with 
$2\tau = i {\bk}(\sqrt{1-k^2})/{\bk}(k)$.  Let
$m=1,2,3,\cdots$. We then have 
$$\spreadlines{6 pt}\allowdisplaybreaks\align 
U_{2m-1}(q)&={\frac{B_{2m}}{4m}}
\biggl\{(2^{2m}-1)-E_{2m}(q)\cr
&\kern 3.5 em+2(1+2^{2m-1})E_{2m}(q^2)-
2^{2m+1}E_{2m}(q^4)\biggr\},\tag 2.93\cr
G_{2m+1}(q)&={\frac{B_{2m+2}}{4(m+1)}}
\biggl\{(2^{2m+2}-1)+E_{2m+2}(q)
-2^{2m+2}E_{2m+2}(q^2)\biggr\},\tag 2.94\cr
C_{2m-1}(q)&={\frac{B_{2m}}{4m}}
\biggl\{-E_{2m}(q)+(1+2^{2m-1})E_{2m}(q^2)
-2^{2m-1}E_{2m}(q^4)\biggr\},\tag 2.95\cr
D_{2m+1}(q)&={\frac{B_{2m+2}}{4(m+1)}}
\biggl\{-E_{2m+2}(q)+E_{2m+2}(q^2)
\biggr\},\tag 2.96\cr
\endalign$$
with $B_{2n}$ the Bernoulli numbers defined by 
\hbox{\rm(2.60)}.
\endproclaim
\demo{Proof} To establish (2.93), consider the elementary
identity 
$$\sum\limits_{r=1}^{\infty}(-1)^rf(r)= 
2\sum\limits_{r=1}^{\infty}f(2r)-
\sum\limits_{r=1}^{\infty}f(r),\tag 2.97$$
with $f(r)$ determined by the $q\mapsto -q$ case of
(2.80).  Apply the trivial identity 
$${x\over {1+x}} =  {x\over {1-x}} - {2x^2\over {1-x^2}}
\tag 2.98$$
termwise to the two resulting sums.  Finally, use (2.91) or 
(2.92) to solve for each of the next four sums in terms of the 
$E_{2m}$, and simplify. For the case $m=1$, note that 
$y$, $2y$, and $4y$ correspond to $q$, $q^2$, and 
$q^4$, respectively, and that $(-3/\pi y)+(18/2\pi y)+
(-24/4\pi y)=0$.

Equation (2.94) is an immediate consequence of (2.97), with 
$f(r)$ determined by the $q\mapsto -q$ case of
(2.81), using (2.92) to solve for each of the two resulting
sums in terms of the $E_{2m+2}$, and simplifying.  

To establish (2.95), consider the elementary identity 
$$\sum\limits_{r=1}^{\infty}g(2r-1)= 
\sum\limits_{r=1}^{\infty}g(r)-
\sum\limits_{r=1}^{\infty}g(2r),\tag 2.99$$
with $g(r)$ determined by (2.83).  Apply the trivial identity 
$${x\over {1-x^2}} =  {x\over {1-x}} - {x^2\over {1-x^2}}
\tag 2.100$$
termwise to the two resulting sums.  Finally, use (2.91) or 
(2.92) to solve for each of the next four sums in terms of the 
$E_{2m}$, and simplify. For the case $m=1$, note that 
$(-3/\pi y)+(9/2\pi y)+(-6/4\pi y)=0$.

Equation (2.96) is an immediate consequence of applying
(2.100) termwise to the sum for $D_{2m+1}(q)$ in (2.84),
using (2.92) to solve for each of the two resulting sums in
terms of the $E_{2m+2}$, and simplifying.  
\qed\enddemo

The $m=1,2,3,4$ cases of the identities in Lemma 2.6 either
appear or are implicit in \cite{21}.  The $m=1,2,3$ cases of
(2.93) are given by \cite{21, Entries 14(i)--(iii) and their
proofs, pp. 129--131}, and the $m=4$ case is implicit in 
\cite{21, Entry 14(iv), pp. 130}.  The $m=1,2$ cases of
(2.94) are given by \cite{21, Entries 14(v)--(vi) and their
proofs, pp. 130--131}, and the $m=3,4$ cases are implicit in 
\cite{21, Entries 14(vii)--(viii), pp. 130}. The $m=1,2,3$ 
cases of (2.95) are given by \cite{21, Entries 15(ix)--(xi) and
their proofs, pp. 132--133}, and the $m=4$ case is implicit in 
\cite{21, Entry 15(xii), pp. 133}.  The $m=1,2$ cases of
(2.96) are given by \cite{21, Entries 15(i)--(ii) and their
proofs, pp. 132--133}, and the $m=3,4$ cases are implicit in 
\cite{21, Entries 15(iii)--(iv), pp. 132}.  

Equation (2.98) is utilized in series manipulations in \cite{21,
pp. 226, 260, 383}, and (2.100) is applied in the proof of 
\cite{21, Entry 15(ix), pp. 132--133}.  See also the comment
just after equation (5.180).  

The Lambert series $U_{2m-1}(q)$, $G_{2m+1}(q)$,
$C_{2m-1}(q)$, and $D_{2m+1}(q)$ in Lemma 2.6 appear
in the sums of squares and sums of triangles identities in
Theorems 5.3, 5.4, 5.5, 5.6, 5.11, and 5.12.

\head 3. Continued Fraction Expansions\endhead

In this section we derive associated continued fraction
and regular C-fraction expansions of the
Laplace transform and formal Laplace transform of
various ratios of Jacobi elliptic functions.  We first
survey Rogers' \cite{206} integration-by-parts proof
of the associated continued fraction expansions of the
Laplace transform of $\sn$, $\cn$, $\dn$, and $\sn^2$.  
We also provide a similar proof of the associated 
continued fraction expansions of the Laplace transform 
of $\sn\cn$ and $\sn\dn$.  The $\sn\cn$ case was first
obtained by Ismail and Masson in \cite{111} using a
more refined integration-by-parts analysis.  We next 
recall Rogers' application of Landen's 
transformation \cite{249, pp. 507}, \cite{134, Eqn.
(3.9.15), (3.9.16), (3.9.17), pp. 78--79} to his result
for $\sn$ to obtain the associated continued fraction
expansion of the Laplace transform of $\sn\cn/\dn$.  We
apply a similar modular transformation technique to obtain
serveral other continued fraction expansions. Our formal
Laplace transform continued fraction expansions are a
consequence of those for $\sn$, $\cn$, $\dn$, $\sn^2$, 
$\sn\cn$, $\sn\dn$, combined with Heilermann's \cite{103,
104} correspondence between formal power series and either
associated continued fractions or regular C-fractions
(see also \cite{119, Theorem 7.14, pp. 244--246; Theorem
7.2, pp. 223--224}), and modular transformations.  

Following Jones and Thron in \cite{119, pp. 18--19}, 
Lorentzen and Waadeland in \cite{149, pp. 5--8}, 
and Berndt in \cite{20, pp. 104--105},   
we adopt the following notation for continued fractions:

$$\contfrac{n=1}{\infty}
{\frac{a_n}{b_n}}
:=\cfrac{a_1}\\
{b_1 +\cfrac{a_2}\\
{b_2 +\tfrac{a_3}{b_3 +
{\lower.6em\hbox{$\,\ddots\,$}}}
\endcfrac}\endcfrac
= {\frac{a_1}{b_1}}\+
	 {\frac{a_2}{b_2}}\+
	 {\frac{a_3}{b_3}}\+\cds}\tag 3.1$$

The types of continued fractions that we need in this
paper are summarized in the following definition. 
\definition{Definition 3.1}
Let $\left\{\alpha_{\nu}\right\}_{\nu=1}^\infty$, 
$\left\{\beta_{\nu}\right\}_{\nu=1}^\infty$, and  
$\left\{\gamma_{\nu}\right\}_{\nu=1}^\infty$ be 
sequences in ${\Bbb C}^\times$ with 
$\alpha_{\nu}\gamma_{\nu}\neq 0$, and let $w$ and 
$\xi$ be indeterminate. The ``corresponding type'' 
continued fraction or regular C-fraction is given by 
$$1 + \contfrac{n=1}{\infty}
{\frac{\gamma_{n}w}{1}},\qquad\gamma_n\neq 0.
\tag 3.2$$  
The associated continued fraction is given by 
$$1 + \frac{\alpha_1 w}{1+\beta_1 w}\+
		\contfrac{n=2}{\infty}
			\frac{-\alpha_n w^2}
				{1+\beta_n w},\qquad\alpha_n\neq 0. 
     \tag 3.3$$
The Jacobi continued fraction or J-fraction is given by 
$$\frac{\alpha_1}{\beta_1 + \xi}\+
		\contfrac{n=2}{\infty}
			\frac{-\alpha_n}
				{\beta_n +\xi},\qquad\alpha_n\neq 0. 
     \tag 3.4$$
\enddefinition  

The following sources have been used for Definition 3.1. 
For regular C-fractions, see \cite{119, Eqn. (7.1.1), pp.
221}, \cite{149, pp. 252--253}, \cite{184, pp. 304}, 
\cite{247, Eqn. (99.2), pp. 399; Eqn. (54.2), pp. 208}. 
The associated continued fraction appears in 
\cite{119, Eqn. (7.2.1), pp. 241}, 
\cite{184, pp. 322, 324; Eqn (8), pp. 376}, 
\cite{247, Eqn. (54.1), pp. 208}.  Finally, the 
J-fraction can be found in \cite{119, Eqn. (7.2.35), pp.
249}, \cite{149, pp. 346}, \cite{184, Eqn. (9), pp. 376; 
Eqn. (1), pp. 390}, \cite{247, Eqn. (23.8), pp. 103; 
Eqn. (51.1), pp. 196}. J-fractions are so named because
the related quadratic form has long been called a J-form. 
For all three types of continued fractions, see 
\cite{119, pp. 128--129; Appendix A. pp. 386--394}.  

It is well-known (see \cite{119, pp. 129}, \cite{206, pp.
74}) that the even part \cite{119, Eqn. (2.4.24), pp. 42} 
of a regular C-fraction is an associated continued
fraction. Moreover, \cite{119, pp. 249}, if in the 
associated continued fraction (3.3) we let $w=1/\xi$,
omit the initial term $1$ and make an equivalence
transformation, we obtain the J-fraction in (3.4). See 
\cite{119, pp. 249--256} for the connection between
J-fractions and orthogonal polynomials. Most of our work
is with associated continued fractions, some is with
regular C-fractions, and we just mention J-fractions to
make the connection with orthogonal polynomials.  

Several recent authors used the term J-fraction for what
is really an associated continued fraction.  They include:
Flajolet \cite{69, pp. 130}, \cite{70, pp. 146}, 
Goulden and Jackson \cite{92, Definition 5.2.1, pp. 291}, 
and Zeng \cite{256, pp. 374}.  

We now have the following associated continued fraction
expansions of Laplace transforms of various ratios of
Jacobi elliptic functions.  
\proclaim{Theorem 3.2} Let the Jacobi elliptic functions 
have modulus $k$, and let $k':=\sqrt{1-k^2}$.  We then 
have the associated continued fraction expansions:
$$\spreadlines{6 pt}\allowdisplaybreaks\align
\int_0^\infty \sn u\ e^{-u/x}\,du
	&= \frac{x^2}{1+(1+k^2)x^2}\+
		\contfrac{n=2}{\infty}
			\frac{-(2n-1)(2n-2)^2 (2n-3) k^2 x^4}
				{1 + (2n-1)^2 (1+k^2) x^2} \tag 3.5\cr
\int_0^\infty \cn u\ e^{-u/x}\,du
	&= \frac{x}{1+x^2}\+
		\contfrac{n=2}{\infty}
			\frac{-(2n-2)^2 (2n-3)^2 k^2 x^4}
				{1 + ((2n-1)^2 + (2n-2)^2 k^2) x^2} \tag 3.6\cr
\int_0^\infty \dn u\ e^{-u/x}\,du
	&= \frac{x}{1+k^2 x^2}\+
		\contfrac{n=2}{\infty}
			\frac{-(2n-2)^2 (2n-3)^2 k^2 x^4}
				{1 + ((2n-1)^2 k^2 + (2n-2)^2) x^2} \tag 3.7\cr
\int_0^\infty \sn u\ \cn u\ e^{-u/x}\,du
	&= \frac{x^2}{1+(4+k^2)x^2}\+
		\contfrac{n=2}{\infty}
			\frac{-(2n-2)^2 (2n-1)^2 k^2 x^4}
				{1 + ((2n)^2 + (2n-1)^2 k^2) x^2} \tag 3.8\cr
\int_0^\infty \sn u\ \dn u\ e^{-u/x}\,du
	&= \frac{x^2}{1+(1+4k^2)x^2}\+
		\contfrac{n=2}{\infty}
			\frac{-(2n-2)^2 (2n-1)^2 k^2 x^4}
				{1 + ((2n-1)^2 + (2n)^2 k^2) x^2} \tag 3.9\cr
\int_0^\infty \sd u\ e^{-u/x}\,du
	&= \frac{x^2}{1+(1-2k^2)x^2}\+
		\contfrac{n=2}{\infty}
			\frac{(2n-1)(2n-2)^2 (2n-3) (kk')^2 x^4}
				{1 + (2n-1)^2 (1-2k^2) x^2} \tag 3.10\cr
\int_0^\infty \cd u\ e^{-u/x}\,du
	&= \frac{x}{1+{k'}^2 x^2}\+
		\contfrac{n=2}{\infty}
			\frac{(2n-2)^2 (2n-3)^2 (kk')^2 x^4}
				{1 + ((2n-1)^2 {k'}^2 - (2n-2)^2 k^2) x^2} 
\tag 3.11\cr
\int_0^\infty \nd u\ e^{-u/x}\,du
	&= \frac{x}{1-k^2 x^2}\+
		\contfrac{n=2}{\infty}
			\frac{(2n-2)^2 (2n-3)^2 (kk')^2 x^4}
				{1 + ((2n-2)^2 {k'}^2 - (2n-1)^2 k^2 ) x^2} 
\tag 3.12\cr
\int_0^\infty \frac{\sn u\ \cn u}{\dn^2 u}\,e^{-u/x}\,du
	&= \frac{x^2}{1+(4-5k^2) x^2}\+
		\contfrac{n=2}{\infty}
			\frac{(2n-2)^2 (2n-1)^2 (kk')^2 x^4}
				{1 + ((2n)^2 {k'}^2 - (2n-1)^2 k^2) x^2} 
\tag 3.13\cr
\int_0^\infty \frac{\sn u}{\dn^2 u}\,e^{-u/x}\,du
	&= \frac{x^2}{1+(1-5k^2) x^2}\+
		\contfrac{n=2}{\infty}
			\frac{(2n-2)^2 (2n-1)^2 (kk')^2 x^4}
				{1 + ((2n-1)^2 {k'}^2 - (2n)^2 k^2) x^2} 
\tag 3.14\cr
\int_0^\infty \sn^2 u\ e^{-u/x}\,du
	&= \frac{2x^3}{1+4(1+k^2)x^2}\+
		\contfrac{n=2}{\infty}
			\frac{-(2n)(2n-1)^2 (2n-2) k^2 x^4}
				{1 + (2n)^2 (1+k^2) x^2} \tag 3.15\cr
\int_0^\infty \sd^2 u\ e^{-u/x}\,du
	&= \frac{2x^3}{1+4(1-2k^2)x^2}\+
		\contfrac{n=2}{\infty}
			\frac{(2n) (2n-1)^2 (2n-2) (kk')^2 x^4}
				{1 + (2n)^2 (1-2k^2) x^2} \tag 3.16\cr
\int_0^\infty \frac{\sn u\ \cn u}{\dn u}
\,e^{-u/x}\,du
	&= \frac{x^2}{1 + (4-2k^2)x^2}\+
		\contfrac{n=2}{\infty}
			\frac{-(2n-1) (2n-2)^2 (2n-3) k^4 x^4}
				{1 + (2n-1)^2 (4-2k^2) x^2} \tag 3.17\cr
\int_0^\infty \frac{\sn^2 u\ \cn^2 u}{\dn^2 u}
\,e^{-u/x}\,du
	&= \frac{2x^3}{1 + 4(4-2k^2)x^2}\+
		\contfrac{n=2}{\infty}
			\frac{-(2n) (2n-1)^2 (2n-2) k^4 x^4}
				{1 + (2n)^2 (4-2k^2) x^2} \tag 3.18\cr
\int_0^\infty 
\frac{1-k\,\sn^2 u}{1+k\,\sn^2 u}
\,e^{-u/x}\,du
	&= \frac{x}{1+4k\,x^2}\+
		\contfrac{n=2}{\infty}
			\frac{-4(2n-2)^2 (2n-3)^2k(1+k)^2 x^4}
				{1 + ((2n-1)^24k + (2n-2)^2(1+ k)^2) x^2}
\tag 3.19\cr
\endalign$$
$$\allowdisplaybreaks\multline
\quad\kern .5em\int_0^\infty 
\frac{\sn u}{1+k\,\sn^2 u}\,e^{-u/x}\,du\\
	= \frac{x^2}{1+(1+6k+k^2)x^2}\+
		\contfrac{n=2}{\infty}
			\frac{-4(2n-1)(2n-2)^2 (2n-3)k(1+k)^2 x^4}
				{1 + (2n-1)^2 (1+6k+k^2) x^2}
\endmultline\tag 3.20$$
$$\allowdisplaybreaks\multline
\quad\kern .5em\int_0^\infty 
\frac{\cn u\ \dn u}{1+k\,\sn^2 u}\,e^{-u/x}\,du\\
	= \frac{x}{1+(1+k)^2x^2}\+
		\contfrac{n=2}{\infty}
			\frac{-4(2n-2)^2 (2n-3)^2k(1+k)^2 x^4}
				{1 + ((2n-1)^2(1+ k)^2 + (2n-2)^24k) x^2}
\endmultline\tag 3.21$$
\endproclaim
\demo{Proof}We first use Rogers' \cite{206}
integration-by-parts argument to establish equations
(3.5), (3.6), (3.7), and (3.15).   

We begin by defining $S_n$, $C_n$, $D_n$, with 
$n=0,1,2,\cdots$, by 
$$\spreadlines{6 pt}\allowdisplaybreaks\align
S_n:&=\int_0^\infty \sn^n(u,k)\,e^{-u/x}\,du,
\tag 3.22\cr
C_n:&=\int_0^\infty \sn^n(u,k)\cn(u,k)\,e^{-u/x}\,du,
\tag 3.23\cr
D_n:&=\int_0^\infty \sn^n(u,k)\dn(u,k)\,e^{-u/x}\,du.
\tag 3.24\cr
\endalign$$ 

We integrate by parts and use various differentiation
formulas and Pythagorean theorems for Jacobi elliptic
functions to show that:
$$\spreadlines{6 pt}\allowdisplaybreaks\align
S_1 &=x^2 - x^2(1+k^2) S_1 + 2k^2x^2 S_3,
\tag 3.25\cr
\kern -7.25 em \text{and}\kern 7.25 em\quad &\cr 
S_n &=n(n-1)x^2 S_{n-2} - n^2(1+k^2) x^2 S_n + 
n(n+1)k^2 x^2 S_{n+2},
\tag 3.26\cr
\endalign$$ 
for $n=2,3,4\cdots$.

\noindent
Solving for $S_1$ and ${S_n}/{S_{n-2}}$ in just the right 
way, we obtain:
$$\spreadlines{6 pt}\allowdisplaybreaks\align
	S_1 &= \frac{x^2}{1 + (1+k^2)x^2
			- 2k^2 x^2 S_3/S_1},
\tag 3.27\cr
	S_n/S_{n-2}
		&= \frac{n(n-1)x^2}{1 + n^2(1+k^2) x^2
			- n(n+1)k^2 x^2 S_{n+2}/S_n},
\tag 3.28\cr
\endalign$$
for $n=2,3,4\cdots$.

\noindent
After iterating and simplifying, we obtain the continued
fraction expansion in (3.5).

Using machinery from the proof of (3.5), and noting that
$S_0=x$, we have
$$S_2 = \frac{2x^3}{1 + 4(1+k^2) x^2 - 
6 k^2 x^2 S_4/S_2},\tag 3.29$$
with $S_n/S_{n-2}$ as in (3.28).  Again, we iterate and 
simplify to obtain the continued fraction expansion in
(3.15).

Proceeding as in the proof of (3.5), we obtain:
$$\spreadlines{6 pt}\allowdisplaybreaks\align
\kern -6 em C_0 &=x - x^2 C_0 + 2k^2x^2 C_2, 
\tag 3.30\cr
\kern -5.75 em \text{and}\kern 5.75 em \kern -2 em &\cr 
\kern -6 em C_n &=n(n-1)x^2 C_{n-2} - 
((n+1)^2+k^2n^2) x^2 C_n + 
(n+1)(n+2)k^2 x^2 C_{n+2},
\tag 3.31\cr
\endalign$$ 
for $n=2,3,4\cdots$.

\noindent
Solving for $C_0$ and ${C_n}/{C_{n-2}}$, we obtain:
$$\spreadlines{6 pt}\allowdisplaybreaks\align
	C_0 &= \frac{x}{1 + x^2	- 2k^2 x^2 C_2/C_0},
\tag 3.32\cr
	C_n/C_{n-2}
		&= \frac{n(n-1)x^2}{1 + ((n+1)^2+k^2n^2) x^2
			- (n+1)(n+2)k^2 x^2 C_{n+2}/C_n},
\tag 3.33\cr
\endalign$$
for $n=2,3,4\cdots$.

\noindent
After iterating and simplifying, we obtain the continued
fraction expansion in (3.6).

Next, proceeding as in the proof of (3.6), we obtain:
$$\spreadlines{6 pt}\allowdisplaybreaks\align
\kern -6 em D_0 &=x - k^2x^2 D_0 + 2k^2x^2 D_2,
\tag 3.34\cr
\kern -5.7 em \text{and}\kern 5.7 em \kern -2 em &\cr 
\kern -6 em D_n &=n(n-1)x^2 D_{n-2} - 
(n^2+(n+1)^2k^2) x^2 D_n + 
(n+1)(n+2)k^2 x^2 D_{n+2},
\tag 3.35\cr
\endalign$$ 
for $n=2,3,4\cdots$.

\noindent
Solving for $D_0$ and ${D_n}/{D_{n-2}}$, we obtain:
$$\spreadlines{6 pt}\allowdisplaybreaks\align
	D_0 &= \frac{x}{1 + k^2x^2	- 2k^2x^2 D_2/D_0},
\tag 3.36\cr
	D_n/D_{n-2}
		&= \frac{n(n-1)x^2}{1 + (n^2+(n+1)^2k^2) x^2
			- (n+1)(n+2)k^2 x^2 D_{n+2}/D_n},
\tag 3.37\cr
\endalign$$
for $n=2,3,4\cdots$.

\noindent
After iterating and simplifying, we obtain the continued
fraction expansion in (3.7).

We now use a slight variation of Rogers' integration by
parts arguments to obtain equations (3.8) and (3.9).  
The key idea is to utilize just one integration by parts
at a time, instead of two successive integrations, to
derive formulas for $C_1$ and $D_1$ that are analogous to
those for $C_0$ and $D_0$ in (3.32) and (3.36),
respectively.  The rest of the analysis then iterates
suitable cases of (3.33) and (3.37) as before.  

Starting with the $C_n$ and $D_n$ in (3.23) and (3.24),
it follows that one integration by parts and the
Pythagorean relation for $\cn^2(u,k)$ or $\dn^2(u,k)$,
respectively, gives the identities 
$$\spreadlines{6 pt}\allowdisplaybreaks\align
	C_n &= nxD_{n-1} -(n+1)xD_{n+1},
\tag 3.38\cr
D_n &= nxC_{n-1} -(n+1)k^2xC_{n+1},
\tag 3.39\cr
\endalign$$
for $n=1,2,3\cdots$.

We first obtain our formula for $C_1$.  Setting $n=1$ in
(3.38) gives
$$C_1 = xD_0 - 2xD_2.\tag 3.40$$
One integration by parts applied to $D_0$ gives 
$$D_0 = x - xk^2C_1.\tag 3.41$$
From the $n=2$ case of (3.39) we have 
$$D_2 = 2xC_1 - 3k^2xC_3.\tag 3.42$$
Substituting (3.41) and (3.42) into (3.40) immediately
gives 
$$C_1 = x^2 - x^2(4+k^2)C_1 + 6k^2x^2C_3.\tag 3.43$$
Solving for $C_1$ in (3.43) in just the right way yields 
$$C_1 = \frac{x^2}{1+(4+k^2)x^2-6k^2x^2C_3/C_1}.
\tag 3.44$$

We next obtain the formula for $D_1$.  Setting $n=1$ in 
(3.39) gives
$$D_1 = xC_0 - 2k^2xC_2.\tag 3.45$$
One integration by parts applied to $C_0$ gives 
$$C_0 = x - xD_1.\tag 3.46$$
From the $n=2$ case of (3.38) we have 
$$C_2 = 2xD_1 - 3xD_3.\tag 3.47$$
Substituting (3.46) and (3.47) into (3.45) immediately
gives 
$$D_1 = x^2 - x^2(1+4k^2)D_1 + 6k^2x^2D_3.\tag 3.48$$
Solving for $D_1$ in (3.48) in just the right way yields 
$$D_1 = \frac{x^2}{1+(1+4k^2)x^2-6k^2x^2D_3/D_1}.
\tag 3.49$$

Note that substituting suitable cases of (3.39) into
(3.38), or of (3.38) into (3.39) yields (3.31) and
(3.35), respectively.  These in turn lead to (3.33) and
(3.37).  

After iterating (3.44) and (3.33), and simplifying, we
obtain the continued fraction expansion in (3.8).  
Similarly, iterating (3.49) and (3.37) leads to
(3.9).  

A short alternate proof of (3.8) and (3.9) is
included just after the proof of Theorem 3.11 near the
end of this section.

We now use Rogers' \cite{206} modular transformation
technique to establish (3.10)--(3.14) and
(3.16)--(3.21).  

Rogers \cite{206} deduced (3.17) and (3.18) by applying
Landen's transformation \cite{249, pp. 507} in the form  
$$\frac{\sn(u,k)\ \cn(u,k)}{\dn(u,k)}
	= \frac{1}{1+k'}\,\sn\!\left((1+k')u,\frac{1-k'}{1+k'}
\right)\tag 3.50$$
to the integrands in (3.17) and (3.18), utilizing the
change of variables $v=(1+k')u$, appealing to (3.5) and
(3.15), and then simplifying.  Here, we use the fact that 
$$k_1^2(1+k')^4=k^4\qquad\text{and}\qquad 
(1+k_1^2)(1+k')^2=2(2-k^2),\tag 3.51$$
where $k':=\sqrt{1-k^2}$ and $k_1:=(1-k')/(1+k')$.  

In order to establish (3.10), (3.11), (3.12), and (3.16)
consider the modular transformations 
\cite{134, Ex. 33, pp. 90; and Eqn. (9.3.9), pp.
249--250}, \cite{249, sec. {\bf 22$\cdot$421}, pp. 508}
in 
$$\spreadlines{6 pt}\allowdisplaybreaks\align
\sd(u,k)&=\frac{1}{k'}\,\sn(k'u, i k/k'),
\tag 3.52a\cr
\cd(u,k)&=\cn(k'u, i k/k'),
\tag 3.52b\cr
\nd(u,k)&=\dn(k'u, i k/k').
\tag 3.52c\cr
\endalign$$
Next, apply (3.52a), (3.52b), (3.52c) as needed to the 
integrands in (3.10), (3.11), (3.12), (3.13), (3.14),
and (3.16).  Change variables by $v=k'u$, appeal to (3.5),
(3.6), (3.7), (3.8), (3.9), (3.15), and then simplify.

The continued fractions in (3.20), (3.21), and (3.19) are
immediate consequences of applying the Gau{\ss} modular
transformations \cite{134, pp. 80}, \cite{93, pp. 915}
in  
$$\spreadlines{6 pt}\allowdisplaybreaks\align
\frac{\sn(u,k)}{1+k\,\sn^2(u,k)}
&=\frac{1}{1+k}\,
\sn\!\left((1+k)u,\frac{2\sqrt{k}}{1+k}\right),
\tag 3.53a\cr
\frac{\cn(u,k)\ \dn(u,k)}{1+k\,\sn^2(u,k)}
&=\cn\!\left((1+k)u,\frac{2\sqrt{k}}{1+k}\right),
\tag 3.53b\cr
\frac{1-k\,\sn^2(u,k)}{1+k\,\sn^2(u,k)}
&=\dn\!\left((1+k)u,\frac{2\sqrt{k}}{1+k}\right),
\tag 3.53c\cr
\endalign$$
to the integrands in (3.20), (3.21), and (3.19),
respectively, changing variables by $v=(1+k)u$, 
appealing to (3.5), (3.6), and (3.7), and then 
simplifying.
\qed\enddemo

The integration-by-parts proof of (3.5), (3.6), (3.7),
and (3.15) first appeared in Rogers' \cite{206, 
pp. 76--77}.  A more recent discussion of these
calculations can be found in \cite{92, pp. 307--308, Ex.
5.2.8, pp. 517--519} and \cite{71}. The $\sn\cn$ case 
in (3.8) was first obtained by Ismail and Masson in
\cite{111} using a more refined integration-by-parts
analysis. Stieltjes \cite{217, 219} first derived (3.5),
(3.6), (3.7), and (3.15) by his addition theorem for
elliptic functions method. Rogers \cite{206} also
rediscovered the addition  theorem techniques of
Stieltjes.  Elegant combinatorial applications of (3.6)
and/or the addition theorem techniques of Rogers and
Stieltjes are studied in \cite{29, 30, 58, 59, 69--71,
97, 98, 195--199, 223, 255--257}.
Rogers \cite{206} was the first to derive the associated
continued fraction expansions in (3.17) and (3.18) by
applying Landen's transformation \cite{249, pp. 507},
\cite{135, Eqn. (3.9.15), (3.9.16), (3.9.17), pp. 78--79}
to (3.5), and (3.15), respectively. Ramanujan was the
first to obtain continued fraction expansions equivalent
to (3.10) and (3.11), in which the integral is written as a
hyperbolic series.  In particular, Ramanujan had the
associated continued fraction equivalent to (3.10), and
the  regular C-fraction equivalent to (3.107) below. 
Berndt \cite{21, pp. 165--167} obtains Ramanujan's
two results by applying the appropriate modular
transformations to (3.5) and (3.6). The rest of the
associated continued fraction expansions in Theorem 3.2
appear to be new. Subsequently, a simplified and more
symmetrical approach to Theorem 3.2 appears in
\cite{49}.   

In order to derive our formal Laplace transform
associated continued fraction expansions from the
continued fractions in Theorem 3.2 we first need
Heilermann's \cite{103, 104} correspondence between formal
power series and associated continued fractions. 
Furthermore, the derivation of our formal Laplace
transform C-fraction expansions in (3.109) and (3.110)
below requires Heilermann's \cite{103, 104} correspondence
between formal power series and C-fractions.    

With any formal power series 
$$L(w)=1+c_1w+c_2w^2+c_3w^3+\cdots,\tag 3.54$$
where $\left\{c_{\nu}\right\}_{\nu=1}^\infty$ is a
sequence in ${\Bbb C}^\times$, 
we associate the two sequences of determinants of 
$n\times n$ square matrices given by the following
definition. 
\definition{Definition 3.3}
Let $\left\{c_{\nu}\right\}_{\nu=1}^\infty$ be a
sequence in ${\Bbb C}^\times$, and let
$m,n=1,2,3,\cdots$.  We take $H_n^{(m)}$ and 
$\chi_n$ to be the determinants of 
$n\times n$ square matrices  
$$\spreadlines{8 pt}\allowdisplaybreaks\align
	H_n^{(m)}\equiv H_n^{(m)} (\{c_\nu\}) 
&:= \det \pmatrix
		c_m & c_{m+1} & \ldots & c_{m+n-2} & c_{m+n-1} \\
		c_{m+1} & c_{m+2} & \ldots & c_{m+n-1} &c_{m+n} \\
		\vdots & \vdots & \ddots & \vdots & \vdots \\
		c_{m+n-1} & c_{m+n} & \ldots & c_{m+2n-3} & c_{m+2n-2}
	\endpmatrix,\tag 3.55\cr
	\chi_n\equiv \chi_n (\{c_\nu\}) 
&:= \det \pmatrix
		c_1 & c_2 & \ldots & c_{n-1} & c_{n+1} \\
		c_2 & c_3 & \ldots & c_n &c_{n+2} \\
		\vdots & \vdots & \ddots & \vdots & \vdots \\
		c_n & c_{n+1} & \ldots & c_{2n-2} & c_{2n}
	\endpmatrix.\tag 3.56\cr
\endalign$$

\noindent
The matrix for $\chi_n$ is obtained from the 
matrix for $H_{n+1}^{(1)}$ by deleting the
next to last column and the last row.  In particular, for 
$n=1$ we have $H_1^{(1)} = c_1$, $H_1^{(2)} = c_2$, and
$\chi_1 = c_2$.  
\footnote{In Jones and Thron \cite{119, pp. 244--246}
Theorem 7.14, $\chi_1$ is incorrectly given as $c_1$.}
We also have $H_0^{(n)} = 1$ and $\chi_0 = 0$.  Note
that $H_n^{(m)} (\{c_\nu\})$ is not the Hankel function in
\cite{6, pp. 208}.
\enddefinition  
 
The following theorem \cite{103, 104}, \cite{119, Theorem
7.14, pp. 244--246} provides explicit necessary and
sufficient conditions for expanding a formal power
series into an associated continued fraction.
\proclaim{Theorem 3.4 (Heilermann)} If for a given
formal power series $L(w)$ in \hbox {\rm (3.54)} we
have 
$$\spreadlines{8 pt}\allowdisplaybreaks\align
\kern -5 em 1+ \sum_{m=1}^\infty c_m w^m &= 1+ 
\frac{\alpha_1 w}{1+\beta_1 w}\+
		\contfrac{n=2}{\infty}
			\frac{-\alpha_n w^2}
				{1+\beta_n w},\qquad\alpha_n\neq 0, 
     \tag 3.57\cr
\kern -8.15 em \text{then}\kern 8.15 em
H_n^{(1)} (\{c_\nu\})&\neq 0,\kern 4 em\text{for}\quad 
n=1,2,3,\cdots,\tag 3.58\cr
\endalign$$
where $H_n^{(1)}(\{c_\nu\})$ is the Hankel determinant in 
\hbox {\rm (3.55)} associated with $L(w)$.  Moreover, 
$$\spreadlines{6 pt}\allowdisplaybreaks\alignat 3
\kern -4 em\alpha_n&=
\frac{H_n^{(1)}\,H_{n-2}^{(1)}}
{\bigl(H_{n-1}^{(1)}\bigr)^2},
&\qquad\text{for}&\quad n=1,2,3,\cdots, 
&\quad&(H_{-1}^{(1)}=H_0^{(1)}=1),\cr  
\kern -4.48 em\text{and}\kern 4.48 em\qquad\beta_n&=
\frac{\chi_{n-1}}{H_{n-1}^{(1)}}-
\frac{\chi_{n}}{H_{n}^{(1)}},
&\qquad\text{for}&\quad n=1,2,3,\cdots, 
&\quad&\left(\chi_0=0,\quad\chi_1=c_2\right),
\tag 3.59\cr  
\endalignat$$
where $\chi_{n}$ is given by \hbox {\rm (3.56)}.

Conversely, suppose that \hbox {\rm (3.58)} holds.  Then, 
\hbox {\rm (3.57)} also holds with coefficients 
$\left\{\alpha_n\right\}$ and $\left\{\beta_n\right\}$
given by \hbox {\rm (3.59)}. In addition, we have 
$$\spreadlines{6 pt}\allowdisplaybreaks\alignat 2
\kern -4 em H_n^{(1)} (\{c_\nu\})&=
\prod\limits_{r=1}^{n}\alpha_r^{n+1-r},
&\qquad\text{for}&\quad n=1,2,3,\cdots, 
\tag 3.60\cr  
\kern -3.85 em\text{and}\kern 3.85 em\qquad
\chi_n (\{c_\nu\}) &=
-(\beta_1+\beta_2+\cdots +\beta_n)\,H_n^{(1)} (\{c_\nu\})
&\ &\ \cr
&=-(\beta_1+\beta_2+\cdots +\beta_n)\,
\prod\limits_{r=1}^{n}\alpha_r^{n+1-r},
&\qquad\text{for}&\quad n=1,2,3,\cdots. 
\tag 3.61\cr  
\endalignat$$
\endproclaim

The following theorem \cite{103, 104}, \cite{119, Theorem
7.2, pp. 223--224} provides explicit necessary and
sufficient conditions for expanding a formal power
series into a regular C-fraction.
\proclaim{Theorem 3.5 (Heilermann)} If for a given
formal power series $L(w)$ in \hbox {\rm (3.54)} we
have 
$$\kern -5 em 1+ \sum_{m=1}^\infty c_m w^m = 1+ 
\contfrac{n=1}{\infty}
{\frac{\gamma_{n}w}{1}},\qquad\gamma_n\neq 0, 
\tag 3.62$$
$$\kern -2.01 em\text{then}\kern 2.01 em\kern 3 em
H_n^{(1)}(\{c_\nu\})\neq 0\qquad\text{and}\qquad
H_n^{(2)}(\{c_\nu\})\neq 0
,\kern 3 em\text{for}\quad 
n=1,2,3,\cdots,\tag 3.63$$
where $H_n^{(1)}(\{c_\nu\})$ and $H_n^{(2)}(\{c_\nu\})$ are
the Hankel determinants in \hbox {\rm (3.55)} associated
with $L(w)$.  Moreover, 
$$\spreadlines{6 pt}\allowdisplaybreaks\alignat 3
\kern -6 em\gamma_{2m}&=
-\frac{H_{m-1}^{(1)}\,H_{m}^{(2)}}
{H_{m}^{(1)}\,H_{m-1}^{(2)}},
&\quad\text{for}&\quad m=1,2,3,\cdots, 
&\quad&(H_{0}^{(1)}=H_0^{(2)}=1),\cr  
\kern -4 em\text{and}\quad \gamma_1= H_{1}^{(1)}\ ;\quad
\gamma_{2m+1}&=
-\frac{H_{m+1}^{(1)}\,H_{m-1}^{(2)}}
{H_{m}^{(1)}\,H_{m}^{(2)}},
&\quad\text{for}&\quad m=1,2,3,\cdots, 
&\quad&\left(H_0^{(2)}=1\right).
\tag 3.64\cr  
\endalignat$$

Conversely, suppose that \hbox{\rm (3.63)} holds.  Then, 
\hbox{\rm (3.62)} also holds with coefficients 
$\left\{\gamma_n\right\}$ given by \hbox{\rm (3.64)}. In
addition, we have 
$$\spreadlines{6 pt}\allowdisplaybreaks\align
\kern -4 em H_n^{(2)}(\{c_\nu\})&=
(-1)^nH_n^{(1)}(\{c_\nu\})
\prod\limits_{r=1}^{n}\gamma_{2r}\cr  
&=(-1)^nH_{n+1}^{(1)}(\{c_\nu\})
\prod\limits_{r=0}^{n}\gamma_{2r+1}^{-1},
\qquad\text{for}\quad n=1,2,3,\cdots. 
\tag 3.65\cr  
\endalign$$
\endproclaim

In our applications of Theorems 3.4 and 3.5 in this paper,
we find the following elementary lemma very useful.
\proclaim{Lemma 3.6} Let $\left\{c_{\nu}
\right\}_{\nu=1}^\infty$ be a
sequence in ${\Bbb C}^\times$, let
$n=1,2,3,\cdots$, and take $H_n^{(m)}$ and 
$\chi_n$ as in \hbox{\rm Definition 3.3}.  If $x$ is a
constant we then have
$$\spreadlines{6 pt}\allowdisplaybreaks\align
H_n^{(1)}(\{x^{\nu}c_\nu\})&=
x^{n^2}\,H_n^{(1)}(\{c_\nu\}),\tag 3.66\cr
H_n^{(1)}(\{x^{\nu-1}c_{\nu-1}\})&=
x^{n(n-1)}\,H_n^{(1)}(\{c_{\nu-1}\})
=x^{2{n\choose 2}}\,H_n^{(1)}(\{c_{\nu-1}\}),\tag 3.67\cr  
\chi_n(\{x^{\nu}c_\nu\})&=
x^{1+n^2}\,\chi_n(\{c_\nu\}),\tag 3.68\cr
\chi_n(\{x^{\nu-1}c_{\nu-1}\})&=
x^{1+n(n-1)}\,\chi_n(\{c_{\nu-1}\})
=x^{1+2{n\choose 2}}\,\chi_n(\{c_{\nu-1}\}).\tag 3.69\cr  
H_n^{(2)}(\{x^{\nu}c_\nu\})&=
x^{n(n+1)}\,H_n^{(2)}(\{c_\nu\}),\tag 3.70\cr
H_n^{(2)}(\{x^{\nu-1}c_{\nu-1}\})&=
x^{n^2}\,H_n^{(2)}(\{c_{\nu-1}\}).\tag 3.71\cr
\endalign$$
\endproclaim

The relations in (3.66)--(3.71) follow immediately by
first factoring suitable powers of $x$ from the rows,
then the columns.

Applying (3.66) and (3.68) to (3.58) and (3.59), it is not
difficult to see that we have the following lemma.
\proclaim{Lemma 3.7} Suppose that the associated
continued fraction expansion in \hbox{\rm (3.57)} holds,
and that $A$ and $B$ are nonzero constants. We then have
$$\kern -5 em 1+ \sum_{m=1}^\infty AB^m c_m w^m = 1+ 
\frac{AB\alpha_1 w}{1+B\beta_1 w}\+
		\contfrac{n=2}{\infty}
			\frac{-B^2\alpha_n w^2}
				{1+B\beta_n w},\tag 3.72$$  
where $\left\{\alpha_n\right\}$ and 
$\left\{\beta_n\right\}$ are given by \hbox{\rm (3.59)}. 
That is, if $L(w)$ is the formal power series in
\hbox{\rm (3.54)}, and \hbox{\rm (3.57)} holds, then 
\hbox{\rm (3.72)} gives the associated continued fraction
expansion for $1+A(L(Bw)-1)$.
\endproclaim
\demo{Proof}The associated continued fraction expansion
in (3.72) is an immediate consequence of assuming
(3.57)--(3.59), and then using (3.66) and (3.68) to
simplify (3.58) and (3.59) where the sequence 
$\left\{c_{\nu}\right\}_{\nu=1}^\infty$ is replaced by 
$\left\{AB^{\nu}c_{\nu}\right\}_{\nu=1}^\infty$.  Just
note that 
$$\spreadlines{6 pt}\allowdisplaybreaks\align
\kern -5 em H_n^{(1)}(\{AB^{\nu}c_\nu\})&=
A^n B^{n^2}\,H_n^{(1)}(\{c_\nu\}),\tag 3.73\cr
\kern -12.6 em \text{and}\kern 12.6 em
\chi_n(\{AB^{\nu}c_\nu\})&=
A^n B^{1+n^2}\,\chi_n(\{c_\nu\}),\tag 3.74\cr
\endalign$$
and then substitute into the right-hand sides of (3.59).
 We find that $\alpha_1$, $\alpha_n$, and
$\beta_n$ become $AB\alpha_1$, $B^2\alpha_n$, and
$B\beta_n$, respectively.  We used the necessary and
then sufficient conditions in Theorem 3.4. 
\qed\enddemo

Next, applying (3.70) to (3.63) and (3.64), it
is not difficult to see that we have the following lemma.
\proclaim{Lemma 3.8} Suppose that the regular C-fraction 
fraction expansion in \hbox{\rm (3.62)} holds, and that
$A$ and $B$ are nonzero constants. We then have
$$\kern -5 em 1+ \sum_{m=1}^\infty AB^m c_m w^m = 1+ 
\frac{AB\gamma_1 w}{1}\+
		\contfrac{n=2}{\infty}
{\frac{B\gamma_{n}w}{1}},\tag 3.75$$  
where $\left\{\gamma_n\right\}$ is given by \hbox{\rm
(3.64)}.  That is, if $L(w)$ is the formal power series in
\hbox{\rm (3.54)}, and \hbox{\rm (3.62)} holds, then 
\hbox{\rm (3.75)} gives the regular C-fraction
expansion for $1+A(L(Bw)-1)$.
\endproclaim
\demo{Proof}The regular C-fraction expansion
in (3.75) is an immediate consequence of assuming
(3.62)--(3.64), and then using (3.66) and (3.70) to
simplify (3.63) and (3.64) where the sequence 
$\left\{c_{\nu}\right\}_{\nu=1}^\infty$ is replaced by 
$\left\{AB^{\nu}c_{\nu}\right\}_{\nu=1}^\infty$.  Just
note that 
$$\spreadlines{6 pt}\allowdisplaybreaks\align
\kern -5 em H_n^{(1)}(\{AB^{\nu}c_\nu\})&=
A^n B^{n^2}\,H_n^{(1)}(\{c_\nu\}),\tag 3.76\cr
\kern -11.5 em \text{and}\kern 11.5 em
H_n^{(2)}(\{AB^{\nu}c_\nu\})&=
A^n B^{n(n+1)}\,H_n^{(2)}(\{c_\nu\}),\tag 3.77\cr
\endalign$$
and then substitute into the right-hand sides of 
(3.64). We find that $\gamma_1$, $\gamma_{2m+1}$, and
$\gamma_{2m}$ become $AB\gamma_1$, $B\gamma_{2m+1}$, and
$B\gamma_{2m}$, respectively.  We used the necessary and
then sufficient conditions in Theorem 3.5. 
\qed\enddemo

In order to apply Lemmas 3.7 and 3.8 to Theorem 3.2 and
the first part of Theorem 3.11 to obtain additional
continued fraction expansions we first need the formal
Laplace transform in the following definition.  
\definition{Definition 3.9}
Given a Maclaurin series expansion
$$ f(u)=\sum_{m=0}^{\infty}\frac{a_m u^m}{m!},
\tag 3.78$$
we obtain a formal Laplace transform by integrating 
term by term:
$${\Cal L}(f,x^{-1}) :=
		\sum\limits_{m=0}^\infty \frac{a_m}{m!}
			\int_0^\infty u^m \; e^{-u/x} \; du
		=\sum\limits_{m=0}^\infty a_m  x^{m+1}.\tag 3.79$$
\enddefinition

We now have the following theorem.
\proclaim{Theorem 3.10} Let the Jacobi elliptic functions 
have modulus $k$, and let $k':=\sqrt{1-k^2}$.
Furthermore, take the formal Laplace transform in 
\hbox{\rm Definition 3.9} of the indicated Maclaurin
series expansions from \hbox{\rm Definition 2.3}.    
We then have the associated continued fraction expansions:
$$\spreadlines{6 pt}\allowdisplaybreaks\align
\int_0^\infty \sc u\ e^{-u/x}\,du
	&= \frac{x^2}{1+(k^2-2)x^2}\+
		\contfrac{n=2}{\infty}
			\frac{-(2n-1)(2n-2)^2 (2n-3) {k'}^2 x^4}
				{1 + (2n-1)^2 (k^2 - 2) x^2} \tag 3.80\cr
\int_0^\infty \dc u\ e^{-u/x}\,du
	&= \frac{x}{1-{k'}^2 x^2}\+
		\contfrac{n=2}{\infty}
			\frac{-(2n-2)^2 (2n-3)^2 {k'}^2 x^4}
				{1 - ((2n-1)^2 {k'}^2 + (2n-2)^2) x^2} \tag 3.81\cr
\int_0^\infty \nc u\ e^{-u/x}\,du
	&= \frac{x}{1-x^2}\+
		\contfrac{n=2}{\infty}
			\frac{-(2n-2)^2 (2n-3)^2 {k'}^2 x^4}
				{1 - ((2n-1)^2 + (2n-2)^2 {k'}^2) x^2} \tag 3.82\cr
\int_0^\infty \sc^2 u\ e^{-u/x}\,du
	&= \frac{2x^3}{1+4(k^2-2)x^2}\+
		\contfrac{n=2}{\infty}
			\frac{-(2n)(2n-1)^2 (2n-2) {k'}^2 x^4}
				{1 + (2n)^2 (k^2-2) x^2} \tag 3.83\cr
\int_0^\infty \frac{\sn u\ \dn u}{\cn^2 u}\,e^{-u/x}\,du
	&= \frac{x^2}{1-(1+4{k'}^2) x^2}\+
		\contfrac{n=2}{\infty}
			\frac{-(2n-2)^2 (2n-1)^2 {k'}^2 x^4}
				{1 - ((2n-1)^2 + (2n)^2 {k'}^2) x^2} 
\tag 3.84\cr
\int_0^\infty \frac{\sn u}{\cn^2 u}\,e^{-u/x}\,du
	&= \frac{x^2}{1-(4+{k'}^2) x^2}\+
		\contfrac{n=2}{\infty}
			\frac{-(2n-2)^2 (2n-1)^2 {k'}^2 x^4}
				{1 - ((2n)^2 + (2n-1)^2 {k'}^2) x^2} 
\tag 3.85\cr
\int_0^\infty \frac{\sn u\ \dn u}{\cn u}
\,e^{-u/x}\,du
	&= \frac{x^2}{1 + 2(2k^2-1)x^2}\+
		\contfrac{n=2}{\infty}
			\frac{-(2n-1) (2n-2)^2 (2n-3)x^4}
				{1 + 2(2n-1)^2 (2k^2-1) x^2}\tag 3.86\cr
\int_0^\infty \frac{\sn u}{\cn u\ \dn u}
\,e^{-u/x}\,du
	&= \frac{x^2}{1 - (2+2k^2)x^2}\+
		\contfrac{n=2}{\infty}
			\frac{-(2n-1) (2n-2)^2 (2n-3) {k'}^4 x^4}
				{1 - (2n-1)^2 (2+2k^2) x^2} \tag 3.87\cr
\int_0^\infty \frac{\sn^2 u\ \dn^2 u}{\cn^2 u}
\,e^{-u/x}\,du
	&= \frac{2x^3}{1 + 8(2k^2-1)x^2}\+
		\contfrac{n=2}{\infty}
			\frac{-(2n) (2n-1)^2 (2n-2)x^4}
         {1 + 2(2n)^2 (2k^2-1)x^2}\tag 3.88\cr
\int_0^\infty \frac{\sn^2 u}{\cn^2 u\ \dn^2 u}
\,e^{-u/x}\,du
	&= \frac{2x^3}{1 - 4(2+2k^2)x^2}\+
		\contfrac{n=2}{\infty}
			\frac{-(2n) (2n-1)^2 (2n-2) {k'}^4 x^4}
				{1 - (2n)^2 (2+2k^2) x^2} \tag 3.89\cr
\endalign$$
\endproclaim
\demo{Proof}We start with a modular transformation 
$$f(u,k)=A\,g(Bu,k_1),\tag 3.90$$
where $A$ and $B$ are nonzero constants, $k_1$ is a
function of the modulus $k$, and $f$ and $g$ are the
quotients of Jacobi elliptic functions in the integrands
of the (formal) Laplace transforms in Theorems 3.10 and
3.2, respectively.  

Let the Maclaurin series expansions of $f(u,k)$ and
$g(u,k)$ in Definition 2.3 be given by one of 
$$\spreadlines{6 pt}\allowdisplaybreaks\alignat 3
\kern -3 em f(u,k)&=
\sum_{m=0}^\infty
		 f_m(k^2)\frac{u^{2m}}{(2m)!}
&\qquad\text{and}&\quad
&\quad g(u,k)&=\sum_{m=0}^\infty
		 g_m(k^2)\frac{u^{2m}}{(2m)!},\tag 3.91\cr  
\kern -3 em f(u,k)&=
\sum_{m=1}^\infty
		 f_m(k^2)\frac{u^{2m-1}}{(2m-1)!}
&\qquad\text{and}&\quad
&\quad g(u,k)&=\sum_{m=1}^\infty
		 g_m(k^2)\frac{u^{2m-1}}{(2m-1)!},\tag 3.92\cr  
\kern -3 em f(u,k)&=
\sum_{m=1}^\infty
		 f_m(k^2)\frac{u^{2m}}{(2m)!}
&\qquad\text{and}&\quad
&\quad g(u,k)&=\sum_{m=1}^\infty
		 g_m(k^2)\frac{u^{2m}}{(2m)!}.\tag 3.93\cr  
\endalignat$$

Substitute (3.91), (3.92), or (3.93) into (3.90), take
the formal Laplace transform of both sides, multiply both
sides by $x$, $1$, or $x^{-1}$, respectively, and then
add $1$ to both sides. We obtain 
$$1+\sum_{m=1}^\infty\overline{f}_m(k^2) w^m=
1+\sum_{m=1}^\infty 
CB^{2m}\,\overline{g}_m(k_1^2) w^m,\tag 3.94$$
where $w=x^2$ and either $\overline{f}_m=f_{m-1}$, 
$\overline{g}_m=g_{m-1}$, $C=AB^{-2}$ in (3.91);  
$\overline{f}_m=f_{m}$, 
$\overline{g}_m=g_{m}$, $C=AB^{-1}$ in (3.92); and  
$\overline{f}_m=f_{m}$, 
$\overline{g}_m=g_{m}$, $C=A$ in (3.93). 

The formulas for the derivatives (with respect to $u$) of 
$\sn(u,k)$, $\cn(u,k)$, and $\dn(u,k)$ are the same for 
$k>1$ as they are for $k<1$.  Thus, we are able to use 
$g_{m-1}(k_1^2)$ and $g_{m}(k_1^2)$ in (3.94).  

The associated continued fraction expansions in
Theorem 3.10 are a direct consequence of applying 
Lemma 3.7 to the right hand side of the relations in
(3.94), corresponding to suitable cases of (3.90), while
keeping in mind the associated continued fraction 
expansions from Theorem 3.2 of the Laplace transform of
the $g(u,k)$.

We complete the proof of Theorem 3.10 by writing down the
necessary cases of (3.90).

For the associated continued fraction  expansions (3.80),
(3.81), and (3.82) consider Jacobi's imaginary modular
transformations \cite{134, Ex. 34, pp. 91; and Eqn.
(9.4.2), pp. 250}, \cite{249, sec. {\bf 22$\cdot$4},
pp. 505--506} in 
$$\spreadlines{6 pt}\allowdisplaybreaks\align
\sc(u,k)&=-\frac{i}{k'}\,\sn(i k'u,1/k'),
\tag 3.95a\cr
\dc(u,k)&=\cn(i k'u,1/k'),
\tag 3.95b\cr
\nc(u,k)&=\dn(i k'u,1/k').
\tag 3.95c\cr
\endalign$$
For (3.83), we just need the square of (3.95a) in 
$$\sc^2(u,k)=-\frac{1}{{k'}^2}\,\sn^2(i k'u,1/k').
\tag 3.96$$
For (3.84) and (3.85) we utilize the combinations of
(3.95a), (3.95b), (3.95c) given by 
$$\spreadlines{6 pt}\allowdisplaybreaks\align
\frac{\sn(u,k)\ \dn(u,k)}{\cn^2(u,k)} & = 
-\frac{i}{k'}\,\sn(i k'u,1/k')\ 
\cn(i k'u,1/k'),\tag 3.97\cr
\frac{\sn(u,k)}{\cn^2(u,k)} & = 
-\frac{i}{k'}\,\sn(i k'u,1/k')\ 
\dn(i k'u,1/k').\tag 3.98\cr
\endalign$$ 

The associated continued fraction expansion in (3.87)
depends on the modular transformation \cite{134, Ex.
35(i), pp. 91} in 
$$\frac{\sn(u,k)}{\cn(u,k)\ \dn(u,k)}
	= -\frac{i}{1+k}
\,\sn\!\left(i (1+k) u,\frac{1-k}{1+k}\right),
\tag 3.99$$
while (3.89) requires the square of (3.99) in 
$$\frac{\sn^2(u,k)}{\cn^2(u,k)\ \dn^2(u,k)}
	= -\frac{1}{(1+k)^2}
\,\sn^2\!\left(i (1+k) u,\frac{1-k}{1+k}\right).
\tag 3.100$$

To obtain (3.86) we observe from (3.95a)--(3.95c) that 
$$\kern -2 em\frac{\sn(u,k)\ \dn(u,k)}{\cn(u,k)}=
\frac{\sc(u,k)\ \dc(u,k)}{\nc(u,k)}=
-\frac{i}{k'}\,
\frac{\sn(i k'u,1/k')\ \cn(i k'u,1/k')}
{\dn(i k'u,1/k')}.\tag 3.101$$
Finally, (3.88) follows from the square of (3.101) in 
$$\kern -2 em\frac{\sn^2(u,k)\ \dn^2(u,k)}{\cn^2(u,k)}=
-\frac{1}{{k'}^2}\,
\frac{\sn^2(i k'u,1/k')\ \cn^2(i k'u,1/k')}
{\dn^2(i k'u,1/k')}.\tag 3.102$$
\qed\enddemo

As far as we know, the associated continued fraction
expansions in Theorem 3.10 are new. It would be
interesting to use suitable (modular) transformations to
extend Theorems 3.2 and 3.10 to the setting of the more
general types of continued fractions (associated with the
Lam\'e equation) discussed in \cite{46, pp. 28--31}. 

In order to obtain our regular C-fraction expansions we
first recall from \cite{119, pp. 129} that the even part 
\cite{119, Eqn. (2.4.24), pp. 42} of a regular C-fraction 
$$\contfrac{n=1}{\infty}
{\frac{\gamma_{n}w}{1}},\qquad\gamma_n\neq 0,
\tag 3.103$$  
is the associated continued fraction 
$$\frac{\gamma_1 w}{1+\gamma_2 w}\+
		\contfrac{n=1}{\infty}
			\frac{-\gamma_{2n}\gamma_{2n+1} w^2}
				{1+(\gamma_{2n+1}+\gamma_{2n+2})w}. 
     \tag 3.104$$
To see why (3.103) and (3.104) are equal, note Rogers' 
\cite{206, pp. 74} observation that the $2m$-th
convergent of (3.103) is identical with the $m$-th
convergent of (3.104).  

We now have the following regular C-fraction expansions 
of the Laplace transform and formal Laplace transform of
various ratios of Jacobi elliptic functions.  
\proclaim{Theorem 3.11} Let the Jacobi elliptic functions 
have modulus $k$, and let $k':=\sqrt{1-k^2}$.  We then 
have the regular C-fraction expansions:
$$\spreadlines{6 pt}\allowdisplaybreaks\alignat 3
&&\int_0^\infty \cn u\ e^{-u/x}\,du
	= \frac{x}{1}\+
		\contfrac{n=2}{\infty}
			\frac{\gamma_n x^2}{1},&&\tag 3.105a\cr
&\kern -2.4 em\text{where}\kern 2.4 em&\qquad
\gamma_{2m}=(2m-1)^2\quad\text{and}\quad 
\gamma_{2m+1}=(2m)^2k^2,&\quad
\text{for}&\quad m=1,2,3,\cdots.\qquad\tag 3.105b\cr
&&\int_0^\infty \dn u\ e^{-u/x}\,du
	= \frac{x}{1}\+
		\contfrac{n=2}{\infty}
			\frac{\gamma_n x^2}{1},&&\tag 3.106a\cr
&\kern -2.4 em\text{where}\kern 2.4 em&\qquad
\gamma_{2m}=(2m-1)^2k^2\quad\text{and}\quad 
\gamma_{2m+1}=(2m)^2,&\quad
\text{for}&\quad m=1,2,3,\cdots.\qquad\tag 3.106b\cr
&&\int_0^\infty \cd u\ e^{-u/x}\,du
	= \frac{x}{1}\+
		\contfrac{n=2}{\infty}
			\frac{\gamma_n x^2}{1},&&\tag 3.107a\cr
&\kern -2.4 em\text{where}\kern 2.4 em&\qquad
\gamma_{2m}=(2m-1)^2{k'}^2\quad\text{and}\quad 
\gamma_{2m+1}=-(2m)^2k^2,&\quad
\text{for}&\quad m=1,2,3,\cdots.\qquad\tag 3.107b\cr
&&\int_0^\infty \nd u\ e^{-u/x}\,du
	= \frac{x}{1}\+
		\contfrac{n=2}{\infty}
			\frac{\gamma_n x^2}{1},&&\tag 3.108a\cr
&\kern -2.4 em\text{where}\kern 2.4 em&\qquad
\gamma_{2m}=-(2m-1)^2k^2\quad\text{and}\quad 
\gamma_{2m+1}=(2m)^2{k'}^2,&\quad
\text{for}&\quad m=1,2,3,\cdots.\qquad\tag 3.108b\cr
\endalignat$$
Next, take the formal Laplace transform in 
\hbox{\rm Definition 3.9} of the indicated Maclaurin
series expansions from \hbox{\rm Definition 2.3}.    
We then have the regular C-fraction expansions:
$$\spreadlines{6 pt}\allowdisplaybreaks\alignat 3
&&\int_0^\infty \dc u\ e^{-u/x}\,du
	= \frac{x}{1}\+
		\contfrac{n=2}{\infty}
			\frac{\gamma_n x^2}{1},&&\tag 3.109a\cr
&\kern -2.45 em\text{where}\kern 2.45 em&\qquad
\gamma_{2m}=-(2m-1)^2{k'}^2\quad\text{and}\quad 
\gamma_{2m+1}=-(2m)^2,&\quad
\text{for}&\quad m=1,2,3,\cdots.\qquad\tag 3.109b\cr
&&\int_0^\infty \nc u\ e^{-u/x}\,du
	= \frac{x}{1}\+
		\contfrac{n=2}{\infty}
			\frac{\gamma_n x^2}{1},&&\tag 3.110a\cr
&\kern -2.45 em\text{where}\kern 2.45 em&\qquad
\gamma_{2m}=-(2m-1)^2\quad\text{and}\quad 
\gamma_{2m+1}=-(2m)^2{k'}^2,&\quad
\text{for}&\quad m=1,2,3,\cdots.\qquad\tag 3.110b\cr
\endalignat$$
\endproclaim
\demo{Proof}To obtain (3.105)--(3.108), we take (3.6),
(3.7), (3.11), (3.12), and multiply both sides by $x$, set 
$x^2=w$, and then identify by inspection what the
$\gamma_n$'s are in (3.103) and (3.104).  We then replace
(3.104) by (3.103), divide by $x$, and simplify. The
expansions in (3.109) and (3.110) follow in the same way
from (3.81) and (3.82), once we note that the left hand
sides are formal Laplace transforms.  

The regular C-fraction expansions in (3.107)--(3.110) are 
also a direct consequence of (3.105), (3.106), Lemma 3.8,
and the modular transformations in (3.52b), (3.52c),
(3.95b), and (3.95c), respectively.  Apply Lemma 3.8 to
the (3.91) case of the relations in (3.94) corresponding
to these four cases of (3.90), while keeping in mind the
regular C-fraction expansions in (3.105) and (3.106).  
\qed\enddemo

The regular C-fraction expansions (3.105) and (3.106)
appear in Rogers \cite{206, pp. 77}, and Ramanujan was the
first to obtain the regular C-fraction expansion
equivalent to (3.107).  Berndt \cite{21, pp. 165} obtains
(3.107) by applying the appropriate modular transformations
to (3.105). The rest of the regular C-fraction expansions
in Theorem 3.11 appear to be new.  

The associated continued fraction expansions in (3.8),
(3.9), (3.13), (3.14), (3.84), and (3.85) can also
be obtained directly from Theorem 3.11 in a simple
manner.  This alternate proof of (3.8) from (3.106) was
also first given by Ismail and Masson in \cite{111}.  
Their method of proof extends to the next three expansions
in (3.9), (3.13), (3.14).  The last two derivations 
of (3.84) and (3.85) from (3.116) and (3.117), 
respectively, are formal.

Keeping in mind the differentiation formulas for $\dn$,
$\cn$, $\nd$, $\cd$, $\nc$, $\dc$ it turns out that 
(3.8), (3.9), (3.13), (3.14), (3.84), (3.85) are
a consequence of one integration by parts, Theorem 3.11,
and the following identity for continued fractions that
can be derived from Lemmas I and II in Rogers
\cite{206}.  
$$1-\frac{1}{1}\+\contfrac{n=2}{\infty}
{\frac{c_{n}x^2}{1}} == 
\frac{c_2x^2}{1+(c_2 + c_3)x^2}\+
		\contfrac{n=2}{\infty}
			\frac{-c_{2n-1}c_{2n}x^4}
				{1+(c_{2n}+c_{2n+1})x^2}. 
     \tag 3.111$$
 
Just apply Theorem 3.11 and equation (3.111) to the
right-hand-sides of the following integration by parts
identities.
$$\spreadlines{6 pt}\allowdisplaybreaks\alignat 2
\int_0^\infty \sn (u,k)\ \cn (u,k)\ e^{-u/x}\,du
	& =\, & \frac{1}{k^2}&\left[1-\frac{1}{x} 
\int_0^\infty \dn (u,k)\ e^{-u/x}\,du\right]
\tag 3.112\cr
\int_0^\infty \sn (u,k)\ \dn (u,k)\ e^{-u/x}\,du
	& =\, &&\left[1-\frac{1}{x} 
\int_0^\infty \cn (u,k)\ e^{-u/x}\,du\right]
\tag 3.113\cr
\int_0^\infty \frac{\sn (u,k)\ \cn (u,k)}
{\dn^2 (u,k)}\,e^{-u/x}\,du
	& =\, &\frac{-1}{k^2}&\left[1-\frac{1}{x} 
\int_0^\infty \nd (u,k)\ e^{-u/x}\,du\right]
\tag 3.114\cr
\int_0^\infty \frac{\sn (u,k)}
{\dn^2 (u,k)}\,e^{-u/x}\,du
	& =\, & \frac{1}{{k'}^2}&\left[1-\frac{1}{x} 
\int_0^\infty \cd (u,k)\ e^{-u/x}\,du\right]
\tag 3.115\cr
\int_0^\infty \frac{\sn (u,k)\ \dn (u,k)}
{\cn^2 (u,k)}\,e^{-u/x}\,du
	& =\, & - &\left[1-\frac{1}{x} 
\int_0^\infty \nc (u,k)\ e^{-u/x}\,du\right]
\tag 3.116\cr
\int_0^\infty \frac{\sn (u,k)}
{\cn^2 (u,k)}\,e^{-u/x}\,du
	& =\, &\frac{-1}{{k'}^2}&\left[1-\frac{1}{x} 
\int_0^\infty \dc (u,k)\ e^{-u/x}\,du\right]
\tag 3.117\cr
\endalignat$$

Theorems 3.2 and 3.10 have interesting special cases when 
$k=0$ and $k=1$.  It is well-known \cite{135, pp. 26 and
39} that 
$$\spreadlines{6 pt}\allowdisplaybreaks\alignat 2
\kern -2 em \sn(u,0)&=\sin u,
&\qquad \cn(u,0)&=\cos u,\qquad \dn(u,0)=1,
\tag 3.118\cr  
\kern -6.2 em\text{and}\kern 6.2 em
\kern 4 em\sn(u,1)&=\tanh u,
&\qquad \cn(u,1)&=\dn(u,1)=\sech u.
\tag 3.119\cr  
\endalignat$$
That is, the Jacobi elliptic functions interpolate
between circular and hyperbolic functions. 

Equations (3.11) and (3.12) interpolate between the Laplace
transform of $\cos u$, $1$, and $\cosh u$ when $k=0$ and
$k=1$.  Equations (3.10) and (3.16) interpolate between the
Laplace transform of $\sin u$ and $\sinh u$, and 
$\sin^2 u$ and $\sinh^2 u$, respectively, when $k=0$ and
$k=1$.

The $k=1$ case of (3.5) and (3.86) both give the Laplace
transform of $\tanh u$.  Furthermore, (3.86) interpolates
between the formal Laplace transform of $\tan u$ and the
Laplace transform of $\tanh u$ when $k=0$ and $k=1$. 
That is, we have 
$$\spreadlines{6 pt}\allowdisplaybreaks\align
	\int_0^\infty \tan u\ e^{-u/x}\,du &=
		\frac{x^2}{1-2x^2}\+
		\contfrac{n=2}{\infty}
			\frac{-(2n-1) (2n-2)^2 (2n-3)x^4}
				{1 - 2(2n-1)^2 x^2},\tag 3.120\cr
	\int_0^\infty \tanh u\ e^{-u/x}\,du &=
		\frac{x^2}{1+2x^2}\+
		\contfrac{n=2}{\infty}
			\frac{-(2n-1) (2n-2)^2 (2n-3)x^4}
				{1 + 2(2n-1)^2 x^2}.\tag 3.121\cr
\endalign$$

When viewed as associated continued fraction expansions
of Laplace transforms, the convergence conditions that
would be required in Theorem 3.10 depend upon the methods
used to interpret the integrals and/or associated
continued fractions.  First, Laplace transforms
of a periodic function can be rewritten as an infinite sum
of integrals over a fixed minimal period.  In this case
we only have to deal with at most $2$ singularities of
the integrand.  Next, the Laplace transform integrals can
be discussed in the context of the Hadamard integral in
\cite{9, Section 5, pp. 45--46}.  Finally, we can consider
multisection of the associated continued fractions.  For
example, look at the numerators and/or denominators of the
even and/or odd partial quotients.  For additional
analytic work relating to associated continued fractions
and Jacobi elliptic functions, see \cite{9, 19, 21--23,
38, 44, 46, 47, 49, 57, 59, 60, 83, 109, 110, 112--116,
148, 217--219, 223, 230--234}, \cite{184, pp. 330}.  We
do not pursue these matters further here.

\head 4. Hankel and $\chi$ determinant
evaluations 
\endhead

In this section we utilize Theorems 7.14 and 7.2 of
\cite{119, pp. 244--246; pp. 223--224}, row and 
column operations, and modular transformations to
deduce our Hankel and $\chi$ determinant
evaluations from the continued fraction expansions of
Section 3.

We start with the Hankel determinant evaluations in
the following theorem. 
\proclaim{Theorem 4.1} Take $H_n^{(1)}(\{c_\nu\})$  
and $H_n^{(2)}(\{c_\nu\})$ to be the $n\times n$
Hankel determinants in \hbox{\rm(3.55)}.   Let the
elliptic function polynomials 
$(\hbox{\rm elliptic})_m(k^2)$ of $k^2$, with $k$
the modulus, be determined by the coefficients of the
Maclaurin series expansions in 
\hbox{\rm Definition 2.3}.  Let $n=1,2,3,\cdots$.  
Then,  
$$\spreadlines{6 pt}\allowdisplaybreaks\align
H_n^{(1)}(\{(sn)_{\nu}(k^2)\}) &=
(k^2)^{n\choose 2}
\prod\limits_{r=1}^{2n-1} r!, \tag 4.1\cr
H_n^{(1)}(\{(cn)_{\nu-1}(k^2)\}) &=
H_n^{(1)}(\{(dn)_{\nu-1}(k^2)\})=
(k^2)^{n\choose 2}
\prod\limits_{r=1}^{n-1}(2r)!^2, \tag 4.2\cr
H_n^{(1)}(\{(sc)_{\nu}(k^2)\}) &=
H_n^{(1)}(\{(sd)_{\nu}(k^2)\})=
(k^2)^{n\choose 2}
\prod\limits_{r=1}^{n}(2r-1)!^2, \tag 4.3\cr
H_n^{(1)}(\{(s/d)_{\nu}(k^2)\}) &=
(-1)^{n\choose 2}
[k^2(1-k^2)]^{n\choose 2}
\prod\limits_{r=1}^{2n-1} r!, \tag 4.4\cr
H_n^{(1)}(\{(c/d)_{\nu-1}(k^2)\}) &=
H_n^{(1)}(\{(nd)_{\nu-1}(k^2)\})\cr
&=(-1)^{n\choose 2}
[k^2(1-k^2)]^{n\choose 2}
\prod\limits_{r=1}^{n-1}(2r)!^2, \tag 4.5\cr
H_n^{(1)}(\{(sc/d^2)_{\nu}(k^2)\}) &=
H_n^{(1)}(\{(s/d^2)_{\nu}(k^2)\})\cr
&=(-1)^{n\choose 2}
[k^2(1-k^2)]^{n\choose 2}
\prod\limits_{r=1}^{n}(2r-1)!^2, \tag 4.6\cr
H_n^{(1)}(\{(sn^2)_{\nu}(k^2)\}) &=
(k^2)^{n\choose 2}
\prod\limits_{r=1}^{2n} r!, \tag 4.7\cr
H_n^{(1)}(\{(s^2/d^2)_{\nu}(k^2)\}) &=
(-1)^{n\choose 2}
[k^2(1-k^2)]^{n\choose 2}
\prod\limits_{r=1}^{2n} r!, \tag 4.8\cr
H_n^{(1)}(\{(sc/d)_{\nu}(k^2)\}) &=
(k^2)^{2{n\choose 2}}
\prod\limits_{r=1}^{2n-1} r!, \tag 4.9\cr
H_n^{(1)}(\{(s^2c^2/d^2)_{\nu}(k^2)\}) &=
(k^2)^{2{n\choose 2}}
\prod\limits_{r=1}^{2n} r!, \tag 4.10\cr
H_n^{(1)}(\{(s/c)_{\nu}(k^2)\}) &=
(1-k^2)^{n\choose 2}
\prod\limits_{r=1}^{2n-1} r!, \tag 4.11\cr
H_n^{(1)}(\{(d/c)_{\nu-1}(k^2)\}) &=
H_n^{(1)}(\{(nc)_{\nu-1}(k^2)\})=
(1-k^2)^{n\choose 2}
\prod\limits_{r=1}^{n-1}(2r)!^2, \tag 4.12\cr
H_n^{(1)}(\{(sd/c^2)_{\nu}(k^2)\})\kern 1.5 em 
&\kern -1.5 em=
H_n^{(1)}(\{(s/c^2)_{\nu}(k^2)\})=
(1-k^2)^{n\choose 2}
\prod\limits_{r=1}^{n}(2r-1)!^2, \tag 4.13\cr
H_n^{(1)}(\{(s^2/c^2)_{\nu}(k^2)\}) &=
(1-k^2)^{n\choose 2}
\prod\limits_{r=1}^{2n} r!, \tag 4.14\cr
H_n^{(1)}(\{(sd/c)_{\nu}(k^2)\}) &=
\prod\limits_{r=1}^{2n-1} r!, \tag 4.15\cr
H_n^{(1)}(\{(s/cd)_{\nu}(k^2)\}) &=
(1-k^2)^{2{n\choose 2}}
\prod\limits_{r=1}^{2n-1} r!, \tag 4.16\cr
H_n^{(1)}(\{(s^2d^2/c^2)_{\nu}(k^2)\}) &=
\prod\limits_{r=1}^{2n} r!, \tag 4.17\cr
H_n^{(1)}(\{(s^2/c^2d^2)_{\nu}(k^2)\}) &=
(1-k^2)^{2{n\choose 2}}
\prod\limits_{r=1}^{2n} r!, \tag 4.18\cr
H_n^{(2)}(\{(cn)_{\nu-1}(k^2)\}) &=
H_n^{(1)}(\{(cn)_{\nu}(k^2)\})\cr
&=(-1)^n(k^2)^{n\choose 2}
\prod\limits_{r=1}^{n}(2r-1)!^2, \tag 4.19\cr
H_n^{(2)}(\{(dn)_{\nu-1}(k^2)\}) &=
H_n^{(1)}(\{(dn)_{\nu}(k^2)\})\cr
&=(-1)^n(k^2)^{n+1\choose 2}
\prod\limits_{r=1}^{n}(2r-1)!^2, \tag 4.20\cr
H_n^{(2)}(\{(c/d)_{\nu-1}(k^2)\}) &=
H_n^{(1)}(\{(c/d)_{\nu}(k^2)\})\cr
&=(-1)^{n+1\choose 2}
(k^2)^{n\choose 2}(1-k^2)^{n+1\choose 2}
\prod\limits_{r=1}^{n}(2r-1)!^2, \tag 4.21\cr
H_n^{(2)}(\{(nd)_{\nu-1}(k^2)\}) &=
H_n^{(1)}(\{(nd)_{\nu}(k^2)\})\cr
&=(-1)^{n\choose 2}
(k^2)^{n+1\choose 2}(1-k^2)^{n\choose 2}
\prod\limits_{r=1}^{n}(2r-1)!^2, \tag 4.22\cr
H_n^{(2)}(\{(d/c)_{\nu-1}(k^2)\})\kern 1 em  
&\kern -1 em =
H_n^{(1)}(\{(d/c)_{\nu}(k^2)\})=
(1-k^2)^{n+1\choose 2}
\prod\limits_{r=1}^{n}(2r-1)!^2, \tag 4.23\cr
H_n^{(2)}(\{(nc)_{\nu-1}(k^2)\}) &=
H_n^{(1)}(\{(nc)_{\nu}(k^2)\})=
(1-k^2)^{n\choose 2}
\prod\limits_{r=1}^{n}(2r-1)!^2, \tag 4.24\cr
\endalign$$
\endproclaim
\demo{Proof}Substitute the Maclaurin series in
Definition 2.3 into the left hand sides of
(3.5)--(3.18),  and (3.80)--(3.89).  Take the formal
Laplace transform, multiply both sides by $x$, $1$, or
$x^{-1}$, where the Maclaurin series is of
the form (3.91), (3.92), or (3.93), respectively, and
then add $1$ to both sides.  The Hankel determinant
evaluations in (4.1)--(4.18) are now immediate
consequences of equation (3.60) of Theorem 3.4.  

An alternate approach is to first establish (4.1), (4.2),
(4.3),  and (4.7) as above, and then deduce the rest 
of (4.1)--(4.18) from these by appealing to suitable
modular transformations and equations (3.66) and
(3.67) of Lemma 3.6.  

We start with a modular transformation such as
(3.90) where $A$ and $B$ are nonzero constants,
$k_1$ is a function of the modulus $k$, $g(u,k)$ 
is one of the Jacobi elliptic functions $\sn(u,k)$, 
$\cn(u,k)$,  $\dn(u,k)$,  $\sn(u,k)\cn(u,k)$, 
$\sn(u,k)\dn(u,k)$, $\sn^2(u,k)$, and
$f(u,k)$ is one of the other quotients of Jacobi elliptic
functions in the integrands of the (formal) Laplace
transforms in (3.80)--(3.89), and (3.10)--(3.14),  
(3.16)--(3.18).  Let the Maclaurin series expansions
of $f(u,k)$ and $g(u,k)$ be given by one of    
(3.91), (3.92), and (3.93).  
Equating coefficients of powers of $u$ in the Maclaurin
series expansion of both sides of (3.90) gives 
$$\overline{f}_m(k^2)= CB^{2m}
\overline{g}_m(k_1^2),\tag 4.25$$
where $m=1,2,3,\cdots$, and either
$\overline{f}_m=f_{m-1}$, 
$\overline{g}_m=g_{m-1}$, $C=AB^{-2}$ in (3.91);
$\overline{f}_m=f_{m}$, 
$\overline{g}_m=g_{m}$, $C=AB^{-1}$ in (3.92); and 
$\overline{f}_m=f_{m}$, 
$\overline{g}_m=g_{m}$, $C=A$ in (3.93).  
The rest of our formulas in (4.1)--(4.18) are a direct
consequence of 
$$H_n^{(1)}(\{\overline{f}_{\nu}(k^2)\})=
C^nH_n^{(1)}(\{B^{2\nu}
\overline{g}_{\nu}(k_1^2)\}),\tag 4.26$$
equations (3.66) and (3.67) of Lemma 3.6, and
equations (4.1), (4.2), (4.3),  and (4.7).  

The necessary cases of (3.90) are given by (3.50),
(3.52a)--(3.52c), (3.95a)--(3.95c), (3.96), (3.97), 
(3.98), (3.99), (3.100), squaring both sides of (3.50) and
(3.52a), composing each of (3.101) and (3.102)
with (3.50) to give  
$$\kern -2 em\frac{\sn(u,k)\ \dn(u,k)}{\cn(u,k)}=
-(k+i k')\sn((i k'-k)u,1-2k(k+i k')),
\tag 4.27$$
$$\kern -2.5 em\text{and}\kern 2.5 em\qquad
\frac{\sn^2(u,k)\ \dn^2(u,k)}{\cn^2(u,k)}=
-[1-2k(k+i k')]
\sn^2((i k'-k)u,1-2k(k+i k')),
\tag 4.28$$
and finally, the combinations of (3.52a), (3.52b), (3.52c)
given by
$$\spreadlines{6 pt}\allowdisplaybreaks\align
\frac{\sn(u,k)\ \cn(u,k)}{\dn^2(u,k)} & = 
\frac{1}{k'}\,\sn(k'u,i k/k')\ 
\cn(k'u,i k/k'),\tag 4.29\cr
\frac{\sn(u,k)}{\dn^2(u,k)} & = 
\frac{1}{k'}\,\sn(k'u,i k/k')\ 
\dn(k'u,i k/k').\tag 4.30\cr
\endalign$$ 

To establish (4.19)--(4.24) from Theorem 3.11, we
recall the relevant Maclaurin series in
(2.62)--(2.66), apply (3.79) to the left hand
side of each regular C-fraction expansion in Theorem
3.11, multiply by $x$, set $x^2=w$, then add
$1$ to both sides.  We obtain
$$\kern -5 em 1+ \sum_{m=1}^\infty 
(elliptic)_{m-1}(k^2) w^m = 1+ 
\frac{w}{1}\+
		\contfrac{n=2}{\infty}
{\frac{\gamma_{n}w}{1}},\tag 4.31$$  
where $\left\{\gamma_n\right\}$ is determined by
Theorem 3.11.  Equations (4.19)--(4.24) now follow
immediately from (3.65) of Theorem 3.5 and the 
$H_n^{(1)}(\{c_{\nu-1}\})$ evaluations in (4.2), (4.5), 
and (4.12).   
\qed\enddemo

The Hankel determinant evaluations in (4.1), (4.2),
(4.7), (4.19), and (4.20) were also obtained earlier
from (3.60) in \cite{2, Eqns. (10.13), (10.17), (10.14),
(10.18), pp. 97--99}. The rest appear to be new.

Note that the product sides of (4.15) and (4.17) are
independent of $k$.  Analogous independence results
for classical orthogonal polynomials appear in \cite{80,
81}. For a more elementary type of independence
result of Sylvester for Hankel determinants, see
\cite{177, Vol. III, pp. 316--317}.  

The $H_n^{(1)}$ evaluations in (4.19)--(4.24) are not
utilized directly in the proof of the $\chi_n$ determinant
evaluations in Theorem 4.2.  However, they are equivalent
to the Hankel determinant evaluations in (4.3),
(4.6), and (4.13) which are used in the proof of
Theorem 4.2.  To see this equivalence first equate
coefficients of $u^N$ in both sides of the Maclaurin series
expansion of the differentiation formulas for $\dn$,
$\cn$, $\nd$, $\cd$, $\nc$, $\dc$.  The evaluations in 
(4.3), (4.6), and (4.13) then follow immediately
from factoring suitable constants from the resulting
determinants and appealing to (4.19)--(4.24).  Note that
the differentiation formula for $\nc$ thus explains why
the evaluations in (4.13) and (4.24) have the same
answer.   

It is not hard to see from Theorems 3.2, 3.10, 4.1,
and equation (3.61) of Theorem 3.4 that we have the
following theorem.  
\proclaim{Theorem 4.2} Take $\chi_n(\{c_\nu\})$ 
to be the $n\times n$ determinant in
\hbox{\rm(3.56)}.   Let the elliptic function
polynomials $(\hbox{\rm elliptic})_m(k^2)$ of $k^2$,
with $k$ the modulus, be determined by the
coefficients of the Maclaurin series expansions in 
\hbox{\rm Definition 2.3}.  Let $n=1,2,3,\cdots$.  
Then,  
$$\spreadlines{6 pt}\allowdisplaybreaks\align
\chi_n(\{(sn)_{\nu}(k^2)\})\kern 2.5 em  
&\kern -2.5 em =
-\tfrac{n(4n^2-1)}{3}(k^2)^{n\choose 2}(1+k^2)
\prod\limits_{r=1}^{2n-1} r!, \tag 4.32\cr
\chi_n(\{(cn)_{\nu-1}(k^2)\})\kern 2.5 em 
&\kern -2.5 em =
-\tfrac{n(2n-1)}{3}(k^2)^{n\choose 2}
[2n(1+k^2)+(1-2k^2)]
\prod\limits_{r=1}^{n-1}(2r)!^2, \tag 4.33\cr
\chi_n(\{(dn)_{\nu-1}(k^2)\})\kern 2.5 em 
&\kern -2.5 em =
-\tfrac{n(2n-1)}{3}(k^2)^{n\choose 2}
[2n(1+k^2) - (2-k^2)]
\prod\limits_{r=1}^{n-1}(2r)!^2, \tag 4.34\cr
\chi_n(\{(sc)_{\nu}(k^2)\})\kern 3.75 em 
&\kern -3.75 em =
-\tfrac{n(2n+1)}{3}(k^2)^{n\choose 2}
[2n(1+k^2) + (2-k^2)]
\prod\limits_{r=1}^{n}(2r-1)!^2, \tag 4.35\cr
\chi_n(\{(sd)_{\nu}(k^2)\})\kern 3.75 em 
&\kern -3.75 em =
-\tfrac{n(2n+1)}{3}(k^2)^{n\choose 2}
[2n(1+k^2)-(1-2k^2)]
\prod\limits_{r=1}^{n}(2r-1)!^2, \tag 4.36\cr
\chi_n(\{(s/d)_{\nu}(k^2)\})\kern 2.5 em  
&\kern -2.5 em =
-(-1)^{n\choose 2}\cdot\tfrac{n(4n^2-1)}{3}
[k^2(1-k^2)]^{n\choose 2}(1-2k^2)
\prod\limits_{r=1}^{2n-1} r!, \tag 4.37\cr
\chi_n(\{(c/d)_{\nu-1}(k^2)\}) &=
-(-1)^{n\choose 2}\cdot\tfrac{n(2n-1)}{3}
[k^2(1-k^2)]^{n\choose 2}\cr
&\qquad\times[2n(1-2k^2) + (1+k^2)]
\prod\limits_{r=1}^{n-1}(2r)!^2, \tag 4.38\cr
\chi_n(\{(nd)_{\nu-1}(k^2)\})&=
-(-1)^{n\choose 2}\cdot\tfrac{n(2n-1)}{3}
[k^2(1-k^2)]^{n\choose 2}\cr
&\qquad\times[2n(1-2k^2) - (2-k^2)]
\prod\limits_{r=1}^{n-1}(2r)!^2, \tag 4.39\cr
\chi_n(\{(sc/d^2)_{\nu}(k^2)\})&=
-(-1)^{n\choose 2}\cdot\tfrac{n(2n+1)}{3}
[k^2(1-k^2)]^{n\choose 2}\cr
&\qquad\times[2n(1-2k^2) + (2-k^2)]
\prod\limits_{r=1}^{n}(2r-1)!^2, \tag 4.40\cr
\chi_n(\{(s/d^2)_{\nu}(k^2)\}) &=
-(-1)^{n\choose 2}\cdot\tfrac{n(2n+1)}{3}
[k^2(1-k^2)]^{n\choose 2}\cr
&\qquad\times[2n(1-2k^2) - (1+k^2)]
\prod\limits_{r=1}^{n}(2r-1)!^2, \tag 4.41\cr
	\chi_n(\{(sn^2)_{\nu}(k^2)\}) &=
-\tfrac{2n(n+1)(2n+1)}{3}
(k^2)^{n\choose 2}(1+k^2)
\prod\limits_{r=1}^{2n} r!, \tag 4.42\cr
\chi_n(\{(s^2/d^2)_{\nu}(k^2)\}) &=
-(-1)^{n\choose 2}\cdot\tfrac{2n(n+1)(2n+1)}{3}\cr
&\qquad\times[k^2(1-k^2)]^{n\choose 2}(1-2k^2)
\prod\limits_{r=1}^{2n} r!, \tag 4.43\cr
\chi_n(\{(sc/d)_{\nu}(k^2)\}) &=
-\tfrac{2n(4n^2-1)}{3}
(k^2)^{2{n\choose 2}}(2-k^2)
\prod\limits_{r=1}^{2n-1} r!, \tag 4.44\cr
\chi_n(\{(s^2c^2/d^2)_{\nu}(k^2)\}) &=
-\tfrac{4n(n+1)(2n+1)}{3}
(k^2)^{2{n\choose 2}}(2-k^2)
\prod\limits_{r=1}^{2n} r!, \tag 4.45\cr
\chi_n(\{(s/c)_{\nu}(k^2)\}) &=
\tfrac{n(4n^2-1)}{3}(1-k^2)^{n\choose 2}
(2-k^2)\prod\limits_{r=1}^{2n-1} r!, \tag 4.46\cr
\chi_n(\{(d/c)_{\nu-1}(k^2)\}) &=
\tfrac{n(2n-1)}{3}
(1-k^2)^{n\choose 2}\cr
&\qquad\times[2n(2-k^2) - (1+k^2)]
\prod\limits_{r=1}^{n-1}(2r)!^2, \tag 4.47\cr
\chi_n(\{(nc)_{\nu-1}(k^2)\})&=
\tfrac{n(2n-1)}{3}
(1-k^2)^{n\choose 2}\cr
&\qquad\times[2n(2-k^2) - (1-2k^2)]
\prod\limits_{r=1}^{n-1}(2r)!^2, \tag 4.48\cr
\chi_n(\{(sd/c^2)_{\nu}(k^2)\})&=
\tfrac{n(2n+1)}{3}
(1-k^2)^{n\choose 2}\cr
&\qquad\times[2n(2-k^2) + (1-2k^2)]
\prod\limits_{r=1}^{n}(2r-1)!^2, \tag 4.49\cr
\chi_n(\{(s/c^2)_{\nu}(k^2)\}) &=
\tfrac{n(2n+1)}{3}
(1-k^2)^{n\choose 2}\cr
&\qquad\times[2n(2-k^2) + (1+k^2)]
\prod\limits_{r=1}^{n}(2r-1)!^2, \tag 4.50\cr
\chi_n(\{(s^2/c^2)_{\nu}(k^2)\}) &=
\tfrac{2n(n+1)(2n+1)}{3}
(1-k^2)^{n\choose 2}(2-k^2)
\prod\limits_{r=1}^{2n} r!, \tag 4.51\cr
\chi_n(\{(sd/c)_{\nu}(k^2)\}) &=
\tfrac{2n(4n^2-1)}{3}(1-2k^2)
\prod\limits_{r=1}^{2n-1} r!, \tag 4.52\cr
\chi_n(\{(s/cd)_{\nu}(k^2)\}) &=
\tfrac{2n(4n^2-1)}{3}
(1-k^2)^{2{n\choose 2}}(1 + k^2)
\prod\limits_{r=1}^{2n-1} r!, \tag 4.53\cr
\chi_n(\{(s^2d^2/c^2)_{\nu}(k^2)\}) &=
\tfrac{4n(n+1)(2n+1)}{3}(1-2k^2)
\prod\limits_{r=1}^{2n} r!, \tag 4.54\cr
\chi_n(\{(s^2/c^2d^2)_{\nu}(k^2)\}) &=
\tfrac{4n(n+1)(2n+1)}{3}
(1-k^2)^{2{n\choose 2}}(1 + k^2)
\prod\limits_{r=1}^{2n} r!.\tag 4.55\cr
\endalign$$
\endproclaim

The same Maclaurin series expansions of modular
transformations that were utilized in the second proof
of Theorem 4.1 can also be combined with equations
(3.68) and (3.69) of Lemma 3.6 to deduce the rest of
equations (4.32)--(4.55) from (4.32)--(4.36), and
(4.42).

All of the nonconstant factors in the products on the
right hand sides of the identities in Theorem 4.2, which are
not powers of $k^2$ or $(1-k^2)$, are linear
combinations of certain subsets of $(1+k^2)$, $(2-k^2)$,
and $(1-2k^2)$.  The $24$ identities in Theorem 4.2 also
have a number of elegant symmetries involving these last
three expressions.  In Section 5, the $(1+k^2)$, 
$(2-k^2)$, and $(1-2k^2)$ correspond to Lambert series.

Keeping in mind (3.118) and (3.119), we see that
Theorems 4.1 and 4.2 have interesting special
cases when $k=0$ and $k=1$. In this paper we need
the $k=0$ cases in the following theorem.
\proclaim{Theorem 4.3} Take $H_n^{(1)}(\{c_\nu\})$,  
$H_n^{(2)}(\{c_\nu\})$, and $\chi_n(\{c_\nu\})$ 
to be the $n\times n$ determinants 
in \hbox{\rm Definition 3.3}.   
Let the Bernoulli numbers $B_n$ and Euler numbers
$E_n$ be defined by \hbox{\rm(2.60)} and
\hbox{\rm(2.61)},  respectively.  Let
$n=1,2,3,\cdots$.   Then,  
$$\spreadlines{6 pt}\allowdisplaybreaks\align
H_n^{(1)}(\{(-1)^{\nu-1}{(2^{2\nu}-1)
\over {4\nu}}\cdot |B_{2\nu}|\}) &=
2^{-(2n^2 + n)}
\prod\limits_{r=1}^{2n-1} r!, \tag 4.56\cr
H_n^{(1)}(\{(-1)^{\nu}{(2^{2\nu+2}-1)
\over {4(\nu+1)}}\cdot |B_{2\nu+2}|\}) &=
(-1)^n2^{-(2n^2 + 3n)}
\prod\limits_{r=1}^{2n} r!, \tag 4.57\cr
H_n^{(1)}(\{(-1)^{\nu-1}\cdot {\tfrac {1}{4}}\cdot
|E_{2\nu-2}|\}) &=2^{-2n}
\prod\limits_{r=1}^{n-1}(2r)!^2, \tag 4.58\cr
H_n^{(2)}(\{(-1)^{\nu-1}\cdot {\tfrac {1}{4}}\cdot
|E_{2\nu-2}|\}) &=
H_n^{(1)}(\{(-1)^{\nu}\cdot {\tfrac {1}{4}}\cdot
|E_{2\nu}|\})\cr
&=(-1)^n2^{-2n}
\prod\limits_{r=1}^{n}(2r-1)!^2, \tag 4.59\cr
\chi_n(\{(-1)^{\nu-1}{(2^{2\nu}-1)
\over {4\nu}}\cdot |B_{2\nu}|\}) &=
-\tfrac{n(4n^2-1)}{3}2^{-(2n^2 + n +1)}
\prod\limits_{r=1}^{2n-1} r!, \tag 4.60\cr
\chi_n(\{(-1)^{\nu}{(2^{2\nu+2}-1)
\over {4(\nu+1)}}\cdot |B_{2\nu+2}|\}) &=
(-1)^{n-1}\cdot\tfrac{n(n+1)(2n+1)}{3}\cr
&\qquad\times2^{-(2n^2 + 3n)}
\prod\limits_{r=1}^{2n} r!, \tag 4.61\cr
\chi_n(\{(-1)^{\nu-1}\cdot {\tfrac {1}{4}}
\cdot |E_{2\nu-2}|\})&=
-\tfrac{n(2n-1)(4n-1)}{3}2^{-2n}
\prod\limits_{r=1}^{n-1}(2r)!^2, \tag 4.62\cr
\chi_n(\{(-1)^{\nu}\cdot {\tfrac {1}{4}}\cdot
|E_{2\nu}|\})\kern 2.5 em 
&\kern -2.5 em =
-(-1)^{n}\tfrac{n(2n+1)(4n+1)}{3}2^{-2n}
\prod\limits_{r=1}^{n}(2r-1)!^2. \tag 4.63\cr
\endalign$$
\endproclaim
\demo{Proof}Since $\tan u =\sc(u,0)$, we find after
equating powers of $u$ in the $u\mapsto uz$ case of
(2.56), and the $k=0$ case of (2.65) that 
$$(-1)^{m-1}{(2^{2m}-1)
\over {4m}}\cdot |B_{2m}|=
\frac{(-1)^{m-1}}{2^{2m+1}}\cdot 
(s/c)_m(0),\tag 4.64$$
for $m=1,2,3,\cdots$.  Equation (4.56) now follows
by (3.66) and the $k=0$ case of (4.11).  

Next, from $\sec^2 u = 1 + \sc^2(u,0)$, we obtain
after equating powers of $u$ in the $u\mapsto uz$
case of (2.57), and the $k=0$ case of (2.67) that 
$$(-1)^{m}{(2^{2m+2}-1)\over {4(m+1)}}
\cdot |B_{2m+2}|=
\frac{(-1)^{m}}{2^{2m+3}}\cdot 
(s^2/c^2)_m(0),\tag 4.65$$
for $m=1,2,3,\cdots$.  Equation (4.57) now follows
by (3.66) and the $k=0$ case of (4.14). 

By $\sec u = \nc(u,0)$, we find after equating powers
of $u$ in the $u\mapsto uz$ case of (2.58), and the
$k=0$ case of (2.66) that 
$$(-1)^{m-1}\cdot{\tfrac {1}{4}}\cdot |E_{2m-2}|=
(-1)^{m-1}\cdot{\tfrac {1}{4}}\cdot 
(nc)_{m-1}(0),\tag 4.66$$
for $m=1,2,3,\cdots$.  Equation (4.58) now follows
by (3.67) and the $k=0$ case of (4.12).  We establish 
(4.59) the same way we obtained (4.58).  Equation
(4.59) is an immediate consequence of (4.66),
(3.71), and the $k=0$ case of (4.24).  

Equation (4.60) is immediate from (4.64), (3.68),
and the $k=0$ case of (4.46).  Equation (4.61)
follows from (4.65), (3.68), and the $k=0$ case of
(4.51).  Equation (4.62) is immediate from 
(4.66), (3.69), and the $k=0$ case of (4.48).  

By $\sec u \tan u =\sc(u,0) \nc(u,0)$, it follows after
equating powers of $u$ in the $u\mapsto uz$ case of
(2.59), and the $k=0$ case of (2.78) that 
$$(-1)^{m}\cdot{\tfrac {1}{4}}\cdot |E_{2m}|=
(-1)^{m}\cdot{\tfrac {1}{4}}\cdot 
(s/c^2)_{m}(0),\tag 4.67$$
for $m=1,2,3,\cdots$.  Equation (4.63) is now immediate
from (4.67), (3.68), and the $k=0$ case of (4.50).  
\qed\enddemo

In addition to our proofs of Theorems 4.1, 4.2, 
and 4.3 we have used Mathematica \cite{250} to
directly verify a number of these determinant
evaluations up to $n=15$.  We computed the various
polynomials $(elliptic)_m(k^2)$ by utilizing Dumont's
recursions and symmetries of the ordinary and
symmetric Schett polynomials of $x$, $y$, and $z$
from \cite{58, Eqn. (2.1), pp. 4 ;Corollaire 2.2, pp. 6 ;
Eqn. (4.1), pp. 11 ;Proposition 4.1, pp. 11}.  The
general ordinary Schett polynomials were first
introduced by Schett \cite{207, 208} in a slightly
different form.  (The fourth Schett polynomial
$X_4$ appeared much earlier as 
$T_{\alpha}^{\text{(IV)}}(\nu)$ in 
\cite{185, pp. 84}).  
In 1926, Mitra \cite{170} employed classical methods to
compute explicit Taylor series expansions of 
$\sn(u,k)$, and $\cn(u,k)$ and $\dn(u,k)$ about
$u=0$, as far as $u^{21}$ and $u^{18}$,
respectively.  Bell, in \cite{15--18}, 
utilized a more direct and more practicable method to
explicitly compute these series and their reciprocals.  
Wrigge \cite{252} obtained Maclaurin series expansions 
of $\sn(u,k)$ in terms of $k^2$. Dumont's bimodular
functions in \cite{58} are very similar to but slightly
more symmetric than Abel's elliptic functions $\phi$,
$f$, and $F$ introduced in \cite{1}.  Abel's functions
are related to Dumont's functions as follows:
$$\phi(u;a,b) = \sn(u;i a,b),\qquad
f(u;a,b) = \cn(u;i a,b),\qquad
F(u;a,b) = \dn(u;i a,b).\tag 4.68$$	

The determinants $H_n^{(m)}(\{c_\nu\})$ in (3.55)
are known in the literature as Hankel, Tur\'anian,
or persymmetric determinants, and sometimes
as determinants of recurrent type.  They are also
called ``Tur\'anians'' by Karlin and Szeg{\H o} \cite{122,
pp. 5}.  An excellent survey of the classical literature
on these determinants can be found in Muir's books
and articles \cite{177--181}, and Krattenthaler's 
summary in \cite{125, pp. 20--23 ; pp. 46--48}. 
For more
details on \cite{177}, see \cite{177, Vol. I, pp.
485--487}, \cite{177, Vol. II, pp. 324--357 ; pp.
461--462}, \cite{177, Vol. III, pp. 309--326 ; pp.
469--473}, \cite{177, Vol. IV, pp. 312--331 ; pp.
464--465}.  For classical work on the related
determinants known as recurrents, see \cite{177, Vol.
II, pp. 210--211}, \cite{177, Vol. III, pp. 208--247},
\cite{177, Vol. IV, pp. 224--241}.  Sylvester \cite{177,
Vol. II, pp. 341} was the first to use the term
persymmetric.  We prefer the more common term of
Hankel determinant, in view of Hankel's very
important work described in \cite{177, Vol. III, pp.
312--316}.  Additional and/or more recent research
involving or directly related to Hankel determinants can be
found in \cite{2, 10, 13, 29, 30, 35, 36, 40, 44--47,
49, 50, 52--55, 62, 70, 72--75, 78, 80--82, 92,
97, 98, 105,  106, 109--116, 119, 122, 124, 125, 133,
135--138, 149, 153, 156, 157, 174--176, 181, 184,
189--192, 195--199, 205, 206, 209, 219, 220, 225,
235--240, 242, 244,  245, 247, 249, 254--257}.    

Heilermann's correspondences in Theorems 3.4 and
3.5, and similar results are well-known.  For example,
works that involved both associated continued fraction
and regular C-fraction correspondences include the
following: Rogers \cite{206, pp. 72--74}, Datta
\cite{54, pp. 109--110 ; pp. 127}, Perron \cite{184,
Satz 5, pp. 304--305 ; Satz 11, pp. 324--325}, Jones
and Thron \cite{119, pp. 244--246 ; pp. 223--224},
and D. \& G. Chudnovsky \cite{47, pp. 140--149}.  For
the associated continued fraction correspondence we
have: Sz\'asz \cite{220}, Beckenbach, Seidel and
Sz\'asz \cite{13, pp. 6}, Viennot \cite{244}, Flajolet
\cite{70, pp. 152}, Goulden and Jackson \cite{92, Ex.
5.2.21, pp. 312 ; pp. 528--529}, Hendriksen and Van
Rossum \cite{106, pp. 321}, and Zeng \cite{256, pp.
381}.  Finally, for the regular C-fraction correspondence
we note: Muir \cite{174--176}, Frobenius
\cite{72}, Stieltjes \cite{219}, Wall \cite{247, Eqn.
(44.7), pp. 172 ; see also pp. 399--413}, and
Lorentzen and Waadeland \cite{149, Eqn. (2.5.2), pp.
257}.  

Heilermann \cite{103, 104} seems to be the first to
discover and prove the explicit determinental formulas
in Theorems 3.4 and 3.5 for expanding a formal power
series into either an associated continued fraction or a
regular C-fraction.  However, special cases and/or a
discussion of the general methods he used to arrive at
his expansions (without his closed formulas) appeared
earlier.  For example, see page 7 of \cite{119} and the
survey in \cite{31, pp. 190--193} for these matters,
as well as references to some of the later work 
involving Theorems 3.4 and 3.5.  The papers
\cite{174--176} of Muir, which include a
rediscovery of many of Heilermann's results, also
provides additional information on the earlier work of
Euler, Gau{\ss}, and Stern.  The paper \cite{176} has
many interesting explicit regular C-fraction expansions. 
Stern \cite{215, pp. 245--259},
\cite{216} knew of these general methods and applied
them to particular series.  Gau{\ss} \cite{119, pp.
198--214} obtained a simple, direct special case of
\cite{103, 104} by finding explicit regular C-fraction
expansions for the ratio $F(a,b;c;z)/F(a,b+1;c+1;z)$
of hypergeometric functions, and of $F(a,1;c;z)$. 
Euler's general method in \cite{68, pp. 138} leads to
general T-fraction \cite{119, Eqn. (A.28), pp. 390}
expansions instead of regular C-fraction expansions. 
For additional related information on Euler's work on
continued fractions, see \cite{31, pp. 97--109},
\cite{65, 66, 68}, and \cite{67, chapter 18}. 
Further comments on Muir's papers \cite{174, 176}
appear in \cite{31, pp. 181--182}.  

The most common use of Theorem 3.4 is to first
obtain (3.57) and then appeal to (3.60) to compute
the Hankel determinants $H_n^{(1)}(\{c_\nu\})$. 
Similar calculations are carried out in \cite{2, 13,
70, 92, 106, 220, 244, 256}.  On the other hand,
continued fraction expansions analogous to (3.57)
and/or (3.62) are derived in \cite{47, 54} by first
directly evaluating suitable Hankel determinants
in order to obtain the required constants such as
(3.59) and/or (3.64) in the continued fraction.  
Further evaluations of Hankel determinants that do not
seem to follow from known results about continued
fractions or orthogonal polynomials appear in \cite{105,
Prop. 14} and \cite{75, Sec. 4}. The paper of Robbins \cite{205}
contains yet a different kind of Hankel determinant evaluation that is
related to those of Cauchy \cite{153, Ex. 6, pp. 67}, Frobenius
and Stickelberger \cite{73, 74}, 
Muir \cite{181, pp. 352; Ex. 29 and 31, pp. 357;  Ex. 43, pp.
360}, and D. V. and G. V. Chudnovsky \cite{47, pp. 143, 144, 147}.

The method of first obtaining (3.57) and then
appealing to (3.60) or (3.61),  
and/or deriving (3.62) and then using (3.65) is
responsible for Theorems 4.1, 4.2, and
(indirectly) Theorem 4.3.  The Hankel determinant 
evaluations in (4.1), (4.2), (4.7), (4.19), and (4.20)
were also obtained earlier from (3.60) in \cite{2, pp.
97--99}. Equations (4.56) and (4.58) are equivalent to
results which Flajolet described in \cite{70, pp. 149, 152}. 
Flajolet refers to applying a formula similar to (3.60) to the
associated continued fraction expansions for $\sec u$ and
$\tan u$ in \cite{70, Table 1, pp. 149}.  Flajolet's product
formulas in \cite{70, Table 1, pp. 149} for Hankel
determinants of Bell numbers and Hankel determinants of
the number of derangements also appear in \cite{55, 189,
191}.  Flajolet also points out that there is a similar
connection between his version of (3.60) and the classical
regular C-fraction expansion equivalent to (3.105) in
\cite{69, pp. 149}.  However, (3.64) is not mentioned. 
Equations (4.58) and (4.59) are equivalent to equations 
(3.52) and (3.53), respectively, in \cite{125, Theorem
52, pp. 46}.  In addition, equation (4.58) is treated in 
\cite{2, 138, 190--192}.  Equations (4.58) and
(4.59) are also equivalent to a (corrected) version, 
$$\det (c'_{r+s})_{0\leq r,s\leq n-1}
\equiv\det (|E_{r+s}|)_{0\leq r,s\leq n-1}
=\prod\limits_{r=0}^{n-1}(r!)^2, \tag 4.69$$
of the second formula on page 323 of \cite{157}.  
A correct version of the third formula on page 323 of
\cite{157} should be related to equations (4.56) and
(4.57). Very few
$\chi_n(\{c_\nu\})$ determinants such as (3.56) have
been explicitely evaluated in the literature for any
$\{c_\nu\}$.  One such example appears in
\cite{13, pp. 6} where $c_\nu:=P_{\nu -1}(x)$ is
the $(\nu -1)$-st Legendre polynomial.    

The determinant evaluations, other than (4.1), (4.2),
(4.7), (4.19), (4.20), (4.56)--(4.59) in
Theorems 4.1, 4.2, and 4.3 do not appear to have
been written down in the literature before.  

One of the deepest and most interesting
investigations in the literature of Hankel determinants
is the work of Karlin and Szeg{\H o} in \cite{122}.  In
order to generalize Tur\'an's \cite{221, 224}
inequality for Legendre polynomials they showed that
the $2m$ by $2m$ Hankel determinants, whose
entries are classical orthogonal polynomials (Legendre,
ultraspherical, Laguerre, Hermite) of $x$, has a
constant sign $(-1)^m$ for all $x$ in a suitable
interval.  For $2m+1$ by $2m+1$ Hankel
determinants of these polynomials they generalize the
well known oscillation property of the orthogonal
polynomials involved.  A crucial part of their analysis
was to express each of these $\ell$ by $\ell$ Hankel
determinants as a constant multiple of $(x^2-1)^N$,
$x^N$, or $1$ times the $n$ by $n$ Wronskian
determinant of certain orthogonal polynomials of
another class (subject to a change of variables and
normalization).  This reduced their analysis to that of
the sign of a Wronskian determinant, which is easier. 
These relations between Hankel determinants and
Wronskians are given by \cite{122, Eqn. (12.1), pp.
47 ; Eqn. (14.1), pp. 53 ; Eqn. (16.1), pp. 60--61 ;
Eqn. (18.2), pp. 69}.  Karlin and Szeg{\H o} describe in
\cite{122, pp. 4--5} how these relations are partially
motivated by the duality between $m$ and $x$ of the
classical orthogonal polynomials $\{Q_m(x)\}$.  By
setting $n=0$ in the above relations, Karlin and
Szeg{\H o} express certain of their $\ell$ by $\ell$
Hankel determinants as a constant multiple of 
$(x^2-1)^N$ or $x^N$, since a $0$ by $0$
Wronskian determinant equals $1$.  For example, see
\cite{122, Eqn. (12.3), pp. 47 ; Eqn. (14.3), pp. 53 ; 
Eqn. (16.5), pp. 61}.  Karlin and Szeg{\H o}'s results from
\cite{122} are discussed in \cite{156, Sections 9.1--9.2,
pp. 133--135}. 

Karlin and Szeg{\H o}'s proofs in \cite{122} of their
relations between Hankel determinants and
Wronskians do not clearly exhibit a general algebraic
transformation, which when specialized in the context
of suitable properties of the orthogonal polynomials
under consideration, yields their results.  Recently,
Leclerc \cite{137, 138} utilized Lascoux's \cite{133}
Schur function approach to orthogonal polynomials,
the `master identity' of Turnbull on minors of a
matrix \cite{225, pp. 48}, \cite{135}, and additional
identities from \cite{136} involving staircase Schur
functions, to put together an extremely elegant
derivation of a general algebraic transformation which
specializes to give Karlin and Szeg{\H o}'s relations
between Hankel determinants and Wronskians. 
Leclerc points out that this technique of realizing
identities between orthogonal polynomials
as specializations of identities for symmetric
functions is motivated by Littlewood's \cite{145,
chapter 7} important idea that any sequence
$\{a_m\}$, $m\ge 1$ of elements of a commutative
ring $R$ can be regarded as the sequence of complete
homogeneous symmetric functions of a fictitious set 
of variables $E$: $S_m(E) = a_m$.  Leclerc's
general algebraic identity states that a Wronskian
of orthogonal polynomials is proportional to a Hankel 
determinant whose elements form a new sequence of
polynomials.  In order to complete his derivation of
Karlin and Szeg{\H o}'s relations, Leclerc verifies that in
each case considered by Karlin and Szeg{\H o}, this new
sequence is transformed by a suitable change of
variables and normalization into another class of
classical orthogonal polynomials.  Leclerc reduces
these verifications to algebraic properties discovered
by Burchnall \cite{35} of the classical orthogonal
polynomials.  

It would be interesting to extend Karlin and
Szeg{\H o}'s  relations between determinants of
classical orthogonal polynomials of Hankel and
Wronskian type to a Jacobi elliptic function setting in
which the elements of the Hankel determinants are
functions of the Maclaurin series coefficients in
Definition 2.3.  This could provide more insights into
the identities and number theory discussed in this
paper. A hint in this direction is provided by such a
relation in \cite{122, Eqn. (12.1), pp. 47} involving
Legendre polynomials and the complicated explicit
formulas from \cite{251, pp. 556--557} which
express the Maclaurin series coefficients of
$\sn(x,k)$ and $\sn^2(x,k)$ as sums of products of
Legendre polynomials.  One direct approach to this
problem is to apply Theorem 1 of \cite{138} to the
orthogonal polynomials in \cite{38, 115}.  A second
approach starts by seeking a continued
fraction/lattice paths understanding of Theorem 1 of
\cite{138}.  This combinatorial analysis should be
directly related to \cite{69, pp. 148--150 ; pp.
152--154}, \cite{70, 71, 133, 244, 245}, and
\cite{82, pp. 304}.  In all the above analysis it is
important to keep in mind the `discrete Wronskian'
in \cite{122, Eqn. (1.3), pp. 5}.  Roughly speaking,
we want to use \cite{137, 138} as a bridge between
\cite{122}, this paper, and \cite{69}.  These
connections deserve much further study.  

In a different direction than \cite{137, 138}, it is
instructive to equate the evaluations of the Hankel
determinant of ultraspherical polynomials in equation
(20) of \cite{13, pp. 5} with that in equation
(14.3) of \cite{122, pp. 53}.  This immediately gives
an elegant evaluation of an integral equivalent to
the classical case of Selberg's multiple beta integral
\cite{152, pp. 992} in which the Vandermonde
determinant in the integrand occurs to the second
power.  It would be interesting to see what other
cases of Selberg's integral arise in a similar fashion
from \cite{122, 138}.

\head 5. The determinant form of sums of squares 
identities\endhead

In this section we combine Theorems 2.4, 4.1, 4.2, 4.3, some
row and column operations from Lemma 3.6, and well-known
relations between classical theta functions/Lambert series
and the elliptic function parameters $z$ and $k$, to derive
first the single determinant, and then sum of determinants,
form of our infinite families of sums of squares and related
identities.  We also provide (see Theorem 5.19 below) 
our generalization to infinite families all $21$ of Jacobi's
\cite{117, Sections 40, 41, 42} explicitly stated degree 
$2, 4, 6, 8$ Lambert series expansions of classical theta
functions.  An elegant generating function involving 
the number of ways of writing $N$ as a sum of $2n$ squares
and $(2n)^2$ triangular numbers appears in Corollary 5.15. 

Throughout this section we need the Bernoulli numbers 
$B_n$ and Euler numbers $E_n$, defined by (2.60) and 
(2.61), respectively.  We also use the notation 
$I_n:=\{1,2,\ldots,n\}$, $\Vert S\Vert$ is the cardinality of
the set $S$, and $\det (M)$ is the determinant of the
$n\times n$  matrix $M$.  

As background, we begin with two lemmas.  The first is an 
inclusion-exclusion computation for Hankel determinants, and
the second provides the necessary elliptic function 
parameter relations.

We first establish the inclusion-exclusion lemma.
\proclaim{Lemma 5.1} 
Let $v_1,\dots ,v_{2n-1}$ and $w_1,\dots ,w_{2n-1}$ be 
indeterminate, and let $n=1,2,3,\cdots$.  Suppose that
the $(i,j)$ entries of the $n\times n$ matrices 
$(w_{i+j-1})$ and $(v_{i+j-1})$ are 
$w_{i+j-1}$ and $v_{i+j-1}$, respectively, and that 
$M_{n,S}$ is the $n\times n$ matrix whose $i$-th row is 
$$\spreadlines{8 pt}\allowdisplaybreaks\alignat 2
\kern -4 em
&v_i+w_i,v_{i+1}+w_{i+1},\cdots,
v_{i+n-1}+w_{i+n-1},
&\kern -4 em\text{if}&\quad \ i\in S,\cr  
\kern -5.2 em\text{and}\kern 5.2 em
&v_i,v_{i+1},\cdots,v_{i+n-1},
&\kern -4 em\text{if}&\quad \ i\notin S. 
\tag 5.1\cr  
\kern -4.55 em \text{Then}\kern 4.55 em &
\kern -.5 em \det(w_{i+j-1})=
\sum\limits_{\emptyset \subseteq S
\subseteq I_n}\kern-.5em 
(-1)^{n-\Vert S\Vert }\det (M_{n,S})&\tag 5.2\cr
&\kern 4.7 em = (-1)^n\det(v_{i+j-1})+
\sum\limits_{p=1}^n(-1)^{n-p}
\kern-.5em
\sum\limits_{{\emptyset \subset S\subseteq I_n}
\atop {\Vert S\Vert =p}}\kern-.5em \det (M_{n,S}).
&\tag 5.3\cr
\endalignat$$
\endproclaim
\demo{Proof}We start with the sum in (5.2) and expand 
$\det (M_{n,S})$ along each row corresponding to 
$i\in S$.  For convenience, let $B_{n,T}$ denote the 
$n\times n$ matrix whose $i$-th row is 
$$w_i,w_{i+1},\cdots,w_{i+n-1},\quad 
\text{if}\quad \ i\in T\quad\text{and}\quad 
v_i,v_{i+1},\cdots,v_{i+n-1},\quad 
\text{if}\quad \ i\notin T. 
\tag 5.4$$
Expanding each row of $\det (M_{n,S})$ corresponding to 
$i\in S$ gives
$$\det (M_{n,S})=\sum\limits_{\emptyset \subseteq T
\subseteq S}\kern-.5em\det (B_{n,T}).\tag 5.5$$
Substituting (5.5) into the sum in (5.2) and then
interchanging summation yields
$$\sum\limits_{\emptyset \subseteq S
\subseteq I_n}\kern-.5em 
(-1)^{n-\Vert S\Vert }\det (M_{n,S}) = 
\sum\limits_{\emptyset \subseteq T
\subseteq I_n}\kern-.5em\det (B_{n,T}) 
\sum\limits_{T \subseteq S\subseteq I_n}\kern-.5em 
(-1)^{n-\Vert S\Vert }.\tag 5.6$$

We use the binomial theorem to simplify the inner sum in the
right-hand side of (5.6).  Let $\Vert T\Vert = p$, for 
$p=0,1,2,\cdots,n$.  We then have 
$$\spreadlines{6 pt}\allowdisplaybreaks\align
&\sum\limits_{T \subseteq S\subseteq I_n}\kern-.5em 
(-1)^{n-\Vert S\Vert } = 
\sum\limits_{\nu = p}^{n}(-1)^{n-\nu}
{n-p\choose \nu -p}\cr
&=(-1)^{n-p}\sum\limits_{\nu = 0}^{n-p}(-1)^{\nu}
{n-p\choose \nu} = (-1)^{n-p}(1-1)^{n-p}\cr
&=0,\quad\text{if}\quad n\neq p,\quad\text{and}\quad 
1, \quad\text{if}\quad n=p.\tag 5.7\cr
\endalign$$

It now follows that (5.6) becomes 
$$\spreadlines{6 pt}\allowdisplaybreaks\align
&\sum\limits_{\emptyset \subseteq S
\subseteq I_n}\kern-.5em 
(-1)^{n-\Vert S\Vert }\det (M_{n,S}) = 
\sum\limits_{{\emptyset \subseteq T
\subseteq I_n}\atop {\Vert T\Vert =n}}\kern-.5em 
\det (B_{n,T})\cr 
&= \det (B_{n,I_n}) = \det(w_{i+j-1}),\tag 5.8\cr
\endalign$$
which is (5.2).

The equality  of (5.2) and (5.3) is immediate.
\qed\enddemo

The elliptic function parameter relations are given by the
following lemma.
\proclaim{Lemma 5.2} Let  $z:=2\bk/\pi$, as in 
\hbox{\rm(2.1)}, with $\bk$ given by \hbox{\rm(2.2)} 
and $k$ the modulus.  Take $q$ as in \hbox{\rm(2.4)}.  
Let the classical theta functions $\vartheta_3 (0,q)$ and 
$\vartheta_2 (0,q)$ be defined by \hbox{\rm(1.1)} and 
\hbox{\rm(1.2)}, respectively. Then,
$$\spreadlines{6 pt}\allowdisplaybreaks\align
z	& = \vartheta_3(0,q)^2,\tag 5.9\cr
z\sqrt{1-k^2}	& = \vartheta_3(0,-q)^2
= \vartheta_4(0,q)^2,\tag 5.10\cr
zk	& = \vartheta_2 (0,q)^2,\tag 5.11\cr
4z^2k	& = \vartheta_2 (0,q^{1/2})^4,\tag 5.12\cr
z^2(1+k^2) & = 1+24\sum\limits_{r=1}^{\infty}
{rq^{r} \over {1+q^{r}}},\tag 5.13\cr
z^2(2-k^2) & = 2+24\sum\limits_{r=1}^{\infty}
{2rq^{2r} \over {1+q^{2r}}},\tag 5.14\cr
z^2(1-2k^2) & = 1-24\sum\limits_{r=1}^{\infty}
{(2r-1)q^{2r-1} \over {1+q^{2r-1}}}.\tag 5.15\cr
\endalign$$
\endproclaim
\demo{Proof} Essentially all of (5.9)--(5.15) can be found in
\cite{117, Eqn. (4.), Section 40; Eqn. (6.), Section 40; Eqn.
(5.), Section 40; Eqn. (7.), Section 40;}.  For convenience, we
give more recent references.   For the fundamental (5.9) see
either \cite{249, pp. 499--500}, \cite{23, Eqn. (7.10),  pp.
27},  or \cite{21, Entry 6,  pp. 101 (keep in mind Eqn. (3.1)
on pp. 98 and Eqn. (6.12) on pp. 102)}.  For (5.10) see 
\cite{23, Theorem 8.1 (Eqn. (ii)), pp. 27} or \cite{21, 
Entry 10(ii), pp. 122}.  Equation (5.11) is equivalent to the 
relation 
$$\vartriangle \kern-.3em (q^2)=
\tfrac{1}{2}z^{1/2}(k^2/q)^{1/4}\tag 5.16$$
in \cite{23, Theorem 8.2 (Eqn. (iii)), pp. 28} or \cite{21, 
Entry 11(iii), pp. 123}, where $\vartriangle \kern-.3em (q)$
from \cite{23, Eqn. (1.5), pp. 3} is defined by 
$$\vartriangle \kern-.3em (q):=
 \sum\limits_{j=0}^{\infty}
q^{\frac{1}{2}j(j+1)}.\tag 5.17$$
To see this, take the square root of both sides of (5.11), 
apply the relation 
$$\vartheta_2(0,q)=2q^{1/4}
\vartriangle \kern-.3em (q^2),\tag 5.18$$
and simplify.  
Similarly, equation (5.12) is equivalent to the 
relation 
$$\vartriangle \kern-.3em (q)=
2^{-{1/2}}z^{1/2}k^{1/4}q^{-{1/8}}\tag 5.19$$
in \cite{23, Theorem 8.2 (Eqn. (i)), pp. 28} or \cite{21, 
Entry 11(i), pp. 123}.  Equation (5.13) appears in \cite{21, 
Entry 13(viii), pp. 127 (keep in mind Eqn. (6.12) on pp.
102)}.   Equation (5.14) is immediate from multiplying 
 by $2$ both sides of the relation in \cite{21,  Entry 13(ix),
pp. 127 (keep in mind Eqn. (6.12) on pp. 102)}.  Finally,
equation (5.15) is an immediate consequence of subtracting
(5.13) from (5.14).
\qed\enddemo

We are now ready to obtain the main identities of this paper.
We begin with the single Hankel determinant form of the
$4n^2$ squares  identity in the following theorem. 
\proclaim{Theorem 5.3 } Let 
$\vartheta_3 (0,-q)$ be determined by 
\hbox{\rm(1.1)}, and let $n=1,2,3,\cdots$.   
We then have
$$\spreadlines{6 pt}\allowdisplaybreaks\align
\vartheta_3 (0,-q)^{4n^2}= &\ 
\left\{(-1)^n2^{2n^2+n}
\prod\limits_{r=1}^{2n-1}(r!)^{-1}\right\}\cdot 
\det (g_{r+s-1})_{1\leq r,s\leq n},\tag 5.20\cr
\kern -7.45 em \text{with}\kern 7.45 em\kern 2 em   
g_i:= &\  U_{2i-1}-c_i,\tag 5.21\cr 
\endalign$$
where $U_{2i-1}$ and $c_i$ are defined by 
$$\spreadlines{6 pt}\allowdisplaybreaks\alignat 2
U_{2i-1}\equiv &\  U_{2i-1}(q) 
 :=\sum\limits_{r=1}^{\infty}(-1)^{r-1}
{r^{2i-1}q^r \over  {1+q^r}},
&\qquad \text{ for}&\quad \ i=1,2,3,\cdots,
\tag 5.22\cr  
\kern -6.45 em\text{and}\kern 6.45 em
c_i:= &\ (-1)^{i-1}{(2^{2i}-1)\over {4i}}\cdot |B_{2i}|,
&\qquad \text{ for}&\quad \ i=1,2,3,\cdots,
\tag 5.23\cr  
\endalignat$$
with $B_{2i}$ the Bernoulli numbers defined by 
\hbox{\rm(2.60)}.
\endproclaim
\demo{Proof} Our analysis deals with 
$\sc(u,k)\dn(u,k)$.  Starting with (2.80), applying row
operations and (3.66), appealing to (4.15), and then
utilizing (5.9) we have the following computation.   
$$\spreadlines{10 pt}\allowdisplaybreaks\align
& H_n^{(1)}(\{U_{2\nu -1}(-q)
-(-1)^{\nu-1}{(2^{2\nu}-1)
\over {4\nu}}\cdot |B_{2\nu}|\})\tag 5.24\cr
& \kern -2 em  = H_n^{(1)}(\{(-1)^{\nu}{z^{2\nu}
\over 2^{2\nu + 1}}\cdot 
(sd/c)_{\nu}(k^2)\})\tag 5.25\cr
& \kern -2 em =  \left\{(-1)^n2^{-(2n^2+n)}
\prod\limits_{r=1}^{2n-1}r!\right\}
\cdot z^{2n^2}\tag 5.26\cr 
& \kern -2 em  =   \left\{(-1)^n2^{-(2n^2+n)}
\prod\limits_{r=1}^{2n-1}r!\right\}
\cdot \vartheta_3 (0,q)^{4n^2}.\tag 5.27\cr
\endalign$$ 
Replacing $q$ by $-q$ in (5.24) and (5.27) gives 
$$\spreadlines{10 pt}\allowdisplaybreaks\align
& H_n^{(1)}(\{U_{2\nu -1}(q)
-(-1)^{\nu-1}{(2^{2\nu}-1)
\over {4\nu}}\cdot |B_{2\nu}|\})  \cr
& \kern -2 em  =   \left\{(-1)^n2^{-(2n^2+n)}
\prod\limits_{r=1}^{2n-1}r!\right\}
\cdot \vartheta_3 (0,-q)^{4n^2}.\tag 5.28\cr
\endalign$$ 
Solving for $\vartheta_3 (0,-q)^{4n^2}$ in (5.28) yields
(5.20).
\qed\enddemo

Lemma 5.1 and Theorem 5.3 lead to the determinant sum
form of the $4n^2$ squares identity in the following
theorem.   
\proclaim{Theorem 5.4} 
Let $n=1,2,3,\cdots$.  Then 
$$\vartheta_3(0,-q)^{4n^2} = 1+\sum\limits_{p=1}^n
(-1)^p2^{2n^2+n}\prod\limits_{r=1}^{2n-1}(r!)^{-1}
\kern-.5em
\sum\limits_{{\emptyset \subset S\subseteq I_n}\atop 
{\Vert S\Vert =p}}\kern-.5em \det (M_{n,S}),
\tag 5.29$$
where $\vartheta_3(0,-q)$ is determined by 
\hbox{\rm(1.1)}, 
and $M_{n,S}$ is the 
$n\times n$ matrix whose $i$-th row is 
$$U_{2i-1},U_{2(i+1)-1},\cdots,
U_{2(i+n-1)-1},\quad 
\text{ if} \ i\in S \quad 
\text{ and}\quad
c_i,c_{i+1},\cdots,c_{i+n-1}, 
\quad \text{ if} \ i\notin S, \tag 5.30$$
where $U_{2i-1}$ and $c_i$ are defined by 
\hbox{\rm(5.22)} and \hbox{\rm(5.23)}, respectively, 
with $B_{2i}$ the Bernoulli numbers defined by 
\hbox{\rm(2.60)}.
\endproclaim
\demo{Proof} Specialize $v_i$ and $w_i$ as follows in
equation (5.3) of Lemma 5.1 and then utilize (2.80) to write
$v_i + w_i$ as a Lambert series.  
$$\spreadlines{6 pt}\allowdisplaybreaks\alignat 2
v_i= &\  (-1)^{i-1}{(2^{2i}-1)
\over {4i}}\cdot |B_{2i}|,&\qquad \text{for}\quad 
i=1,2,3,\cdots,\tag 5.31\cr
w_i= &\  (-1)^{i}{z^{2i}
\over 2^{2i+ 1}}\cdot 
(sd/c)_{i}(k^2),&\qquad \text{for}\quad 
i=1,2,3,\cdots,\tag 5.32\cr
v_i+w_i= &\  U_{2i -1}(-q)=
\sum\limits_{r=1}^{\infty}
{-r^{2i-1}q^r \over  {1+(-1)^rq^r}},
&\quad \text{for}\quad 
i=1,2,3,\cdots.\tag 5.33\cr
\endalignat$$

Equating (5.25) and (5.27) it is immediate that 
$$\det(w_{i+j-1}) =  \left\{(-1)^n2^{-(2n^2+n)}
\prod\limits_{r=1}^{2n-1}r!\right\}
\cdot \vartheta_3 (0,q)^{4n^2}.\tag 5.34$$
From (4.56) we have 
$$\det(v_{i+j-1}) =  2^{-(2n^2+n)}
\prod\limits_{r=1}^{2n-1}r!.\tag 5.35$$

Equation (5.29) is now a direct consequence of applying the
determinant evaluations in (5.34) and (5.35), replacing $q$
by $-q$, and then multiplying both sides of the resulting
transformation of (5.3) by 
$$(-1)^n2^{2n^2+n}
\prod\limits_{r=1}^{2n-1}(r!)^{-1},\tag 5.36$$
and simplifying.
\qed\enddemo

The single Hankel determinant form of the $4n(n+1)$ 
squares identity is given by the following theorem. 
\proclaim{Theorem 5.5 } Let 
$\vartheta_3 (0,-q)$ be determined by 
\hbox{\rm(1.1)}, and let $n=1,2,3,\cdots$.   
We then have
$$\spreadlines{6 pt}\allowdisplaybreaks\align
\vartheta_3 (0,-q)^{4n(n+1)}= &\ 
\left\{2^{2n^2+3n}
\prod\limits_{r=1}^{2n}(r!)^{-1}\right\}\cdot 
\det (g_{r+s-1})_{1\leq r,s\leq n},\tag 5.37\cr
\kern -9.65 em \text{with}\kern 9.65 em\kern 2 em   
g_i:= &\  G_{2i+1}-a_i,\tag 5.38\cr 
\endalign$$
where $G_{2i+1}$ and $a_i$ are defined by 
$$\spreadlines{6 pt}\allowdisplaybreaks\alignat 2
G_{2i+1}\equiv &G_{2i+1}(q) 
 :=\sum\limits_{r=1}^{\infty}(-1)^{r}
{r^{2i+1}q^r \over  {1-q^r}},
&\qquad \text{ for}&\quad \ i=1,2,3,\cdots,
\tag 5.39\cr  
\kern -7.15 em\text{and}\kern 7.15 em
a_i:= &\ (-1)^{i}{(2^{2i+2}-1)\over {4(i+1)}}
\cdot |B_{2i+2}|,
&\qquad \text{ for}&\quad \ i=1,2,3,\cdots,
\tag 5.40\cr  
\endalignat$$
with $B_{2i+2}$ the Bernoulli numbers defined by 
\hbox{\rm(2.60)}.
\endproclaim
\demo{Proof} Our analysis deals with 
$\sc^2(u,k)\dn^2(u,k)$.  Starting with (2.81), applying row
operations and (3.66), appealing to (4.17), and then
utilizing (5.9) we have the following computation.   
$$\spreadlines{10 pt}\allowdisplaybreaks\align
& H_n^{(1)}(\{G_{2\nu +1}(-q)
-(-1)^{\nu}{(2^{2\nu + 2}-1)
\over {4(\nu + 1)}}\cdot |B_{2\nu + 2}|\})\tag 5.41\cr
& \kern -2 em  = H_n^{(1)}(\{(-1)^{\nu - 1}{z^{2\nu + 2}
\over 2^{2\nu + 3}}\cdot 
(s^2d^2/c^2)_{\nu}(k^2)\})\tag 5.42\cr
& \kern -2 em =  \left\{2^{-(2n^2+3n)}
\prod\limits_{r=1}^{2n}r!\right\}
\cdot z^{2n(n+1)}\tag 5.43\cr 
& \kern -2 em  =   \left\{2^{-(2n^2+3n)}
\prod\limits_{r=1}^{2n}r!\right\}
\cdot \vartheta_3 (0,q)^{4n(n+1)}.\tag 5.44\cr
\endalign$$ 
Replacing $q$ by $-q$ in (5.41) and (5.44) gives 
$$\spreadlines{10 pt}\allowdisplaybreaks\align
& H_n^{(1)}(\{G_{2\nu +1}(q)
-(-1)^{\nu}{(2^{2\nu + 2}-1)
\over {4(\nu + 1)}}\cdot |B_{2\nu + 2}|\})  \cr
& \kern -2 em  =   \left\{2^{-(2n^2+3n)}
\prod\limits_{r=1}^{2n}r!\right\}
\cdot \vartheta_3 (0,-q)^{4n(n+1)}.\tag 5.45\cr
\endalign$$ 
Solving for $\vartheta_3 (0,-q)^{4n(n+1)}$ in (5.45) yields
(5.37).
\qed\enddemo

The determinant sum form of the $4n(n+1)$ squares
identity is given by the following theorem.   
\proclaim{Theorem 5.6 } 
Let $n=1,2,3,\cdots$.  Then 
$$\vartheta_3(0,-q)^{4n(n+1)} = 1+\sum\limits_{p=1}^n
(-1)^{n-p}2^{2n^2+3n}\prod\limits_{r=1}^{2n}(r!)^{-1}
\kern-.5em
\sum\limits_{{\emptyset \subset S\subseteq I_n}\atop 
{\Vert S\Vert =p}}\kern-.5em \det (M_{n,S}),
\tag 5.46$$
where $\vartheta_3(0,-q)$ is determined by 
\hbox{\rm(1.1)}, 
and $M_{n,S}$ is the 
$n\times n$ matrix whose $i$-th row is 
$$G_{2i+1},G_{2(i+1)+1},\cdots,
G_{2(i+n-1)+1},\quad 
\text{ if} \ i\in S \quad 
\text{ and}\quad
a_i,a_{i+1},\cdots,a_{i+n-1}, 
\quad \text{ if} \ i\notin S, \tag 5.47$$
where $G_{2i+1}$ and $a_i$ are defined by 
\hbox{\rm(5.39)} and \hbox{\rm(5.40)}, respectively, 
with $B_{2i}$ the Bernoulli numbers defined by 
\hbox{\rm(2.60)}.
\endproclaim
\demo{Proof} Specialize $v_i$ and $w_i$ as follows in
equation (5.3) of Lemma 5.1 and then utilize (2.81) to write
$v_i + w_i$ as a Lambert series.  
$$\spreadlines{6 pt}\allowdisplaybreaks\alignat 2
v_i= &\  (-1)^{i}{(2^{2i+2}-1)\over {4(i+1)}}
\cdot |B_{2i+2}|,&\qquad \text{for}\quad 
i=1,2,3,\cdots,\tag 5.48\cr
w_i= &\  (-1)^{i-1}{z^{2i + 2}
\over 2^{2i+ 3}}\cdot 
(s^2d^2/c^2)_{i}(k^2),&\qquad \text{for}\quad 
i=1,2,3,\cdots,\tag 5.49\cr
v_i+w_i= &\  G_{2i +1}(-q)=
\sum\limits_{r=1}^{\infty}
{r^{2i+1}q^r \over  {1-(-1)^rq^r}},
&\quad \text{for}\quad 
i=1,2,3,\cdots.\tag 5.50\cr
\endalignat$$

Equating (5.42) and (5.44) it is immediate that 
$$\det(w_{i+j-1}) =  \left\{2^{-(2n^2+3n)}
\prod\limits_{r=1}^{2n}r!\right\}
\cdot \vartheta_3 (0,q)^{4n(n+1)}.\tag 5.51$$
From (4.57) we have 
$$\det(v_{i+j-1}) =  (-1)^n2^{-(2n^2+3n)}
\prod\limits_{r=1}^{2n}r!.\tag 5.52$$

Equation (5.46) is now a direct consequence of applying the
determinant evaluations in (5.51) and (5.52), replacing $q$
by $-q$, and then multiplying both sides of the resulting
transformation of (5.3) by 
$$2^{2n^2+3n}\prod\limits_{r=1}^{2n}(r!)^{-1},
\tag 5.53$$ and simplifying.
\qed\enddemo

The analysis of the formulas for $r_{4n^2}(N)$ and 
$r_{4n(n+1)}(N)$ obtained by taking the coefficient of
$q^N$ in Theorems 5.4 and 5.6 is analogous to the formulas
for $r_{16}(n)$ and $r_{24}(n)$ in Theorem 1.7.  The
dominate terms for $r_{4n^2}(N)$ and 
$r_{4n(n+1)}(N)$ arise from the $p=n$ terms in (5.29) and
(5.46), respectively.  The other terms are all of a strictly 
decreasing lower order of magnitude.  That is, the terms for 
$r_{4n^2}(N)$ and $r_{4n(n+1)}(N)$ corresponding to the
$p$-th terms in (5.29) and (5.46) have orders of magnitude 
$N^{(4np-2p^2-1)}$ and $N^{(4np-2p^2+2p-1)}$, 
respectively.  The dominate $p=n$ cases are consistent
with \cite{94, Eqn. (9.20), pp. 122}.  Note that this analysis
does not apply to the $n=1$ case of Theorem 5.4.  All of 
this analysis depends upon Lemma 5.37 at the end of this
section.  

In order to state the next four identities related to 
$\nc(u,k)$ we need the Euler numbers $E_n$ defined by
(2.61).

We begin with the single $H_n^{(1)}$ Hankel determinant
identity in the following theorem.
\proclaim{Theorem 5.7 } Let 
$\vartheta_3 (0,-q)=\vartheta_4 (0,q)$ be determined by 
\hbox{\rm(1.1)}, and let $n=1,2,3,\cdots$.   
We then have

$$\spreadlines{6 pt}\allowdisplaybreaks\align
\kern -2 em
\vartheta_3(0,q)^{2n(n-1)}\vartheta_3(0,-q)^{2n^2}= & 
\ \vartheta_3(0,q)^{2n(n-1)}\vartheta_4(0,q)^{2n^2}\cr
= &\ \left\{(-1)^n2^{2n}
\prod\limits_{r=1}^{n-1}(2r)!^{-2}\right\}\cdot 
\det (g_{r+s-1})_{1\leq r,s\leq n},\tag 5.54\cr
\kern -12.25 em \text{with}\kern 12.25 em  
g_i: = &\  R_{2i-2}-b_i,\tag 5.55\cr 
\endalign$$
where $R_{2i-2}$ and $b_i$ are defined by
$$\spreadlines{6 pt}\allowdisplaybreaks\alignat 2
R_{2i-2}\equiv &R_{2i-2}(q) 
 :=\sum\limits_{r=1}^{\infty}(-1)^{r+1}
{(2r-1)^{2i-2}q^{2r-1} \over  {1+q^{2r-1}}},
&\qquad \text{ for}&\quad \ i=1,2,3,\cdots,
\tag 5.56\cr  
\kern -4.5 em\text{and}\kern 4.5 em
b_i:= &\ (-1)^{i-1}\cdot {\tfrac {1}{4}}\cdot |E_{2i-2}|,
&\qquad \text{ for}&\quad \ i=1,2,3,\cdots,
\tag 5.57\cr  
\endalignat$$
with $E_{2i-2}$ the Euler numbers defined by 
\hbox{\rm(2.61)}.
\endproclaim
\demo{Proof} Our analysis deals with 
$\nc(u,k)$.  Starting with (2.82), applying row
operations and (3.66), appealing to (4.12), and then
utilizing both (5.9) and (5.10) we have the following
computation.   
$$\spreadlines{10 pt}\allowdisplaybreaks\align
& H_n^{(1)}(\{R_{2\nu -2}(q)
-(-1)^{\nu -1}\cdot {\tfrac {1}{4}}\cdot
|E_{2\nu -2}|\})\tag 5.58\cr
& \kern -2 em  = H_n^{(1)}(\{(-1)^\nu {z^{2\nu -1}\over 4}\,\sqrt{1-k^2}
\cdot (nc)_{\nu -1}(k^2)\})\tag 5.59\cr
& \kern -2 em =  \left\{(-1)^n2^{-2n}
\prod\limits_{r=1}^{n-1}(2r)!^2\right\}
\cdot z^{2n^2-n}(\sqrt{1-k^2})^{n^2}\tag 5.60\cr 
& \kern -2 em =  \left\{(-1)^n2^{-2n}
\prod\limits_{r=1}^{n-1}(2r)!^2\right\}
\cdot z^{n(n-1)}(z\sqrt{1-k^2})^{n^2}\tag 5.61\cr 
& \kern -2 em  =  \left\{(-1)^n2^{-2n}
\prod\limits_{r=1}^{n-1}(2r)!^2\right\}\cdot 
\vartheta_3(0,q)^{2n(n-1)}\vartheta_3(0,-q)^{2n^2}.
\tag 5.62\cr
\endalign$$ 

Equating (5.58) with (5.62) and then solving for 
$\vartheta_3(0,q)^{2n(n-1)}\vartheta_3(0,-q)^{2n^2}$
yields (5.54).
\qed\enddemo

Taking $q\mapsto -q$ in the $n=1$ case of Theorem 5.7 
gives Jacobi's $2$-squares identity in 
\cite{117, Eqn. (4.), Section 40} and 
\cite{21, Entry 8(i), pp. 114}.  The $n=2$ case of Theorem
5.7 immediately leads to 
$$\vartheta_3(0,q)^4\vartheta_3(0,-q)^8 =
{\tfrac{1}{4}}\det\vmatrix  4R_0-1 & 4R_2+1 \\
  &    \\
4R_2+1 & 4R_4-5  
\endvmatrix,\tag 5.63$$
where $R_{2i-2}$ is defined by (5.56), and 
$\vartheta_3(0,-q)=\vartheta_4(0,q)$ is determined by
(1.1).  The relation in (5.63) is an elegant formula for the
product of the two identities in \cite{21, Entry 8(ii), pp. 114;
 Examples (i)($q\mapsto -q$), pp. 139}.  It is also analogous
to the Lambert series expansions in Section 40 of
\cite{117}.  

We generally want the powers of the theta functions in this
paper to be quadratic in $n$.  However, it is sometimes
useful to write the left-hand side of (5.54) in the form 
$$\vartheta_3(0,-q^2)^{4n(n-1)}
\vartheta_3(0,-q)^{2n},\tag 5.64$$
by appealing to the relation 
$$\vartheta_3(0,q)\vartheta_3(0,-q)=
\vartheta_3(0,-q^2)^2\tag 5.65$$
from \cite{21, Entry 25(iii), pp. 40}.

The $H_n^{(1)}$ Hankel determinant sum 
identity is given by the following theorem.
\proclaim{Theorem 5.8} 
Let $n=1,2,3,\cdots$.  Then 
$$\spreadlines{6 pt}\allowdisplaybreaks\align
\kern -2 em
\vartheta_3(0,q)^{2n(n-1)}\vartheta_3(0,-q)^{2n^2}= & 
\ \vartheta_3(0,q)^{2n(n-1)}\vartheta_4(0,q)^{2n^2}\cr
= &\ 1+\sum\limits_{p=1}^n
(-1)^p2^{2n}\prod\limits_{r=1}^{n-1}(2r)!^{-2}
\kern-.5em
\sum\limits_{{\emptyset \subset S\subseteq I_n}\atop 
{\Vert S\Vert =p}}\kern-.5em \det (M_{n,S}),
\tag 5.66\cr
\endalign$$
where $\vartheta_3(0,-q)=\vartheta_4(0,q)$ is determined
by \hbox{\rm(1.1)}, and $M_{n,S}$ is the 
$n\times n$ matrix whose $i$-th row is 
$$R_{2i-2},R_{2(i+1)-2},\cdots, R_{2(i+n-1)-2},\quad 
\text{ if} \ i\in S \quad \text{ and}\quad
b_i,b_{i+1},\cdots,b_{i+n-1}, 
\quad \text{ if} \ i\notin S, \tag 5.67$$
where $R_{2i-2}$ and $b_i$ are defined by 
\hbox{\rm(5.56)} and \hbox{\rm(5.57)}, respectively, 
with $E_{2i-2}$ the Euler numbers defined by 
\hbox{\rm(2.61)}.
\endproclaim
\demo{Proof} Specialize $v_i$ and $w_i$ as follows in
equation (5.3) of Lemma 5.1 and then utilize (2.82) to write
$v_i + w_i$ as a Lambert series.  
$$\spreadlines{6 pt}\allowdisplaybreaks\alignat 2
v_i= &\  (-1)^{i-1}\cdot {\tfrac {1}{4}}\cdot
|E_{2i-2}|,&\qquad \text{for}\quad 
i=1,2,3,\cdots,\tag 5.68\cr
w_i= &\  (-1)^i{z^{2i-1}\over 4}\,\sqrt{1-k^2}
\cdot (nc)_{i-1}(k^2),&\qquad \text{for}\quad 
i=1,2,3,\cdots,\tag 5.69\cr
v_i+w_i= &\  R_{2i-2}(q)=
\sum\limits_{r=1}^{\infty}(-1)^{r+1}
{(2r-1)^{2i-2}q^{2r-1} \over  {1+q^{2r-1}}},
&\quad \text{for}\quad 
i=1,2,3,\cdots.\tag 5.70\cr
\endalignat$$

Equating (5.59) and (5.62) it is immediate that 
$$\det(w_{i+j-1}) =  \left\{(-1)^n2^{-2n}
\prod\limits_{r=1}^{n-1}(2r)!^2\right\}\cdot 
\vartheta_3(0,q)^{2n(n-1)}
\vartheta_3(0,-q)^{2n^2}.\tag 5.71$$
From (4.58) we have 
$$\det(v_{i+j-1}) =  2^{-2n}
\prod\limits_{r=1}^{n-1}(2r)!^2.\tag 5.72$$

Equation (5.66) is now a direct consequence of applying the
determinant evaluations in (5.71) and (5.72), and then
multiplying both sides of the resulting transformation of (5.3)
by 
$$(-1)^n2^{2n}
\prod\limits_{r=1}^{n-1}(2r)!^{-2},\tag 5.73$$
and simplifying.
\qed\enddemo

Taking $q\mapsto -q$ in the $n=1$ case of Theorem 5.8 
gives Jacobi's $2$-squares identity in 
\cite{117, Eqn. (4.), Section 40} and 
\cite{21, Entry 8(i), pp. 114}.  The $n=2$ case of Theorem
5.8 immediately leads to 
$$\vartheta_3(0,q)^4\vartheta_3(0,-q)^8 =
1-\left[5R_0+2R_2+R_4\right]+
4\left[R_0R_4-R_2^2\right],\tag 5.74$$
where $R_{2i-2}$ is defined by (5.56), and 
$\vartheta_3(0,-q)=\vartheta_4(0,q)$ is determined by
(1.1).

The single $H_n^{(2)}$ Hankel determinant identity 
is given by the following theorem.
\proclaim{Theorem 5.9 } Let 
$\vartheta_3 (0,-q)=\vartheta_4 (0,q)$ be determined by 
\hbox{\rm(1.1)}, and let $n=1,2,3,\cdots$.   
We then have
$$\spreadlines{6 pt}\allowdisplaybreaks\align
\kern -2 em
&\ \vartheta_3(0,q)^{2n(n+1)}\vartheta_3(0,-q)^{2n^2}=  
\ \vartheta_3(0,q)^{2n(n+1)}\vartheta_4(0,q)^{2n^2}\cr
= &\ \left\{2^{2n}
\prod\limits_{r=1}^{n}(2r-1)!^{-2}\right\}\cdot 
\det (g_{r+s-1})_{1\leq r,s\leq n},\tag 5.75\cr
\kern -6.8 em \text{with}\kern 6.8 em  
g_i: = &\  R_{2i}-b_{i+1},\tag 5.76\cr 
\endalign$$
where $R_{2i}$ and $b_{i+1}$ are determined by 
\hbox{\rm(5.56)} and \hbox{\rm(5.57)}, respectively, 
with $E_{2i}$ the Euler numbers defined by 
\hbox{\rm(2.61)}.
\endproclaim
\demo{Proof} Our analysis still deals with 
$\nc(u,k)$.  Starting with (2.82), applying row
operations and (3.66), appealing to (4.24), and then
utilizing both (5.9) and (5.10) we have the following
computation.   
$$\spreadlines{10 pt}\allowdisplaybreaks\align
& H_n^{(2)}(\{R_{2\nu -2}(q)
-(-1)^{\nu -1}\cdot {\tfrac {1}{4}}\cdot
|E_{2\nu -2}|\})\tag 5.77\cr
& \kern -2 em  = H_n^{(1)}(\{R_{2\nu}(q)
-(-1)^{\nu}\cdot {\tfrac {1}{4}}\cdot
|E_{2\nu}|\})\tag 5.78\cr
& \kern -2 em  = H_n^{(1)}(\{(-1)^{\nu +1}{z^{2\nu+1}
\over 4}\,\sqrt{1-k^2}
\cdot (nc)_{\nu}(k^2)\})\tag 5.79\cr
& \kern -2 em =  \left\{2^{-2n}
\prod\limits_{r=1}^{n}(2r-1)!^2\right\}
\cdot z^{2n^2+n}(\sqrt{1-k^2})^{n^2}\tag 5.80\cr 
& \kern -2 em =  \left\{2^{-2n}
\prod\limits_{r=1}^{n}(2r-1)!^2\right\}
\cdot z^{n(n+1)}(z\sqrt{1-k^2})^{n^2}\tag 5.81\cr 
& \kern -2 em  =  \left\{2^{-2n}
\prod\limits_{r=1}^{n}(2r-1)!^2\right\}\cdot 
\vartheta_3(0,q)^{2n(n+1)}\vartheta_3(0,-q)^{2n^2}.
\tag 5.82\cr
\endalign$$ 

Equating (5.78) with (5.82) and then solving for 
$\vartheta_3(0,q)^{2n(n+1)}\vartheta_3(0,-q)^{2n^2}$
yields (5.75).
\qed\enddemo

The $n=1$ case of Theorem 5.9 is 
$$\vartheta_3(0,q)^{4}\vartheta_3(0,-q)^{2}=\   
1+4R_2=\ 1-4\sum\limits_{r=1}^{\infty}(-1)^{r}
{(2r-1)^{2}q^{2r-1} \over  {1+q^{2r-1}}},\tag 5.83$$
which is Jacobi's degree $6$ Lambert series expansion in 
\cite{117, Eqn. (41.), Section 40}.  Note that 
(5.83) gives an elegant formula for the
product of the two identities in \cite{21, Entry 8(ii), pp. 114;
 Entry 8(v), pp. 114}. 

The $n=2$ case of Theorem 5.9 is 
$$\vartheta_3(0,q)^{12}\vartheta_3(0,-q)^8 =
{\tfrac{1}{(3!)^2}}\det\vmatrix  4R_2+1 & 4R_4-5 \\
  &    \\
4R_4-5 & 4R_6+61  
\endvmatrix,\tag 5.84$$
where $R_{2i}$ is determined by (5.56), and 
$\vartheta_3(0,-q)=\vartheta_4(0,q)$ is determined by
(1.1).  

The corresponding $H_n^{(2)}$ Hankel determinant sum 
identity is given by the following theorem.
\proclaim{Theorem 5.10} 
Let $n=1,2,3,\cdots$.  Then 
$$\spreadlines{6 pt}\allowdisplaybreaks\align
\kern -2 em
\vartheta_3(0,q)^{2n(n+1)}\vartheta_3(0,-q)^{2n^2}=&  
\ \vartheta_3(0,q)^{2n(n+1)}\vartheta_4(0,q)^{2n^2}\cr
= &\ 1+\sum\limits_{p=1}^n
(-1)^{n-p}2^{2n}\prod\limits_{r=1}^{n}(2r-1)!^{-2}
\kern-.5em
\sum\limits_{{\emptyset \subset S\subseteq I_n}\atop 
{\Vert S\Vert =p}}\kern-.5em \det (M_{n,S}),
\tag 5.85\cr
\endalign$$
where $\vartheta_3(0,-q)=\vartheta_4(0,q)$ is determined
by \hbox{\rm(1.1)}, and $M_{n,S}$ is the 
$n\times n$ matrix whose $i$-th row is 
$$R_{2i},R_{2(i+1)},\cdots, R_{2(i+n-1)},\quad 
\text{ if} \ i\in S \quad \text{ and}\quad
b_{i+1},b_{i+2},\cdots,b_{i+n}, 
\quad \text{ if} \ i\notin S, \tag 5.86$$
where $R_{2i}$ and $b_{i+1}$ are determined by 
\hbox{\rm(5.56)} and \hbox{\rm(5.57)}, respectively, 
with $E_{2i}$ the Euler numbers defined by 
\hbox{\rm(2.61)}.
\endproclaim
\demo{Proof} Specialize $v_i$ and $w_i$ as follows in
equation (5.3) of Lemma 5.1 and then utilize (2.82) to write
$v_i + w_i$ as a Lambert series.  
$$\spreadlines{6 pt}\allowdisplaybreaks\alignat 2
v_i= &\  (-1)^{i}\cdot {\tfrac {1}{4}}\cdot
|E_{2i}|,&\qquad \text{for}\quad 
i=1,2,3,\cdots,\tag 5.87\cr
w_i= &\  (-1)^{i+1}{z^{2i+1}\over 4}\,\sqrt{1-k^2}
\cdot (nc)_{i}(k^2),&\qquad \text{for}\quad 
i=1,2,3,\cdots,\tag 5.88\cr
v_i+w_i= &\  R_{2i}(q)=
\sum\limits_{r=1}^{\infty}(-1)^{r+1}
{(2r-1)^{2i}q^{2r-1} \over  {1+q^{2r-1}}},
&\quad \text{for}\quad 
i=1,2,3,\cdots.\tag 5.89\cr
\endalignat$$

Equating (5.79) and (5.82) it is immediate that 
$$\det(w_{i+j-1}) =  \left\{2^{-2n}
\prod\limits_{r=1}^{n}(2r-1)!^2\right\}\cdot 
\vartheta_3(0,q)^{2n(n+1)}
\vartheta_3(0,-q)^{2n^2}.\tag 5.90$$
From (4.59) we have 
$$\det(v_{i+j-1}) =  (-1)^n2^{-2n}
\prod\limits_{r=1}^{n}(2r-1)!^2.\tag 5.91$$

Equation (5.85) is now a direct consequence of applying the
determinant evaluations in (5.90) and (5.91), and then
multiplying both sides of the resulting transformation of (5.3)
by 
$$2^{2n}
\prod\limits_{r=1}^{n}(2r-1)!^{-2},\tag 5.92$$
and simplifying.
\qed\enddemo

Keeping in mind (2.90), the Maclaurin series expansions for 
$\nc(u,k)$ and $\sn(u,k)\dn(u,k)/\cn^2(u,k)$ in (2.66) and
(2.79), and the derivative formula for $\nc(u,k)$, it follows
that Theorem 5.9 is equivalent to the single $H_n^{(1)}$
Hankel determinant theta function identity related to 
$\sn(u,k)\dn(u,k)/\cn^2(u,k)$. 

The next two single $H_n^{(1)}$ Hankel determinant
identities involve the $z=0$ case of the theta function
$\vartheta_2 (z,q)$ in \cite{249, pp. 464}, defined by
(1.2). These identities are related to $\sd(u,k)\cn(u,k)$ and 
$\sn^2(u,k)$, respectively.  They are the first major step in
our proof of the two Kac--Wakimoto conjectured identities for triangular
numbers in \cite{120, pp. 452}.

We have the following theorem. 
\proclaim{Theorem 5.11} Let 
$\vartheta_2 (0,q)$ be defined by 
\hbox{\rm(1.2)}, and let $n=1,2,3,\cdots$.   
We then have
$$\spreadlines{6 pt}\allowdisplaybreaks\align
\kern -4 em\vartheta_2 (0,q)^{4n^2}= &\ 
\left\{4^{n(n+1)}
\prod\limits_{r=1}^{2n-1}(r!)^{-1}\right\}\cdot 
\det (C_{2(r+s-1)-1})_{1\leq r,s\leq n},
\tag 5.93\cr
\kern -7.37 em \text{and}\kern 7.37 em  
\kern -4 em\vartheta_2 (0,q^{1/2})^{4n(n+1)}= &\ 
\left\{2^{n(4n+5)}
\prod\limits_{r=1}^{2n}(r!)^{-1}\right\}\cdot 
\det (D_{2(r+s-1)+1})_{1\leq r,s\leq n},
\tag 5.94\cr
\endalign$$
where $C_{2i-1}$ and $D_{2i+1}$ are defined by ,
respectively, 
$$\spreadlines{6 pt}\allowdisplaybreaks\alignat 2
\kern -2 em C_{2i-1}\equiv &\ C_{2i-1}(q) 
 :=\sum\limits_{r=1}^{\infty}
{(2r-1)^{2i-1}q^{2r-1} \over  {1-q^{2(2r-1)}}},
&\qquad \text{ for}&\quad \ i=1,2,3,\cdots,
\tag 5.95\cr  
\kern -6.3 em\text{and}\kern 6.3 em
\kern -2 em D_{2i+1}\equiv &\ D_{2i+1}(q) 
 :=\sum\limits_{r=1}^{\infty}
{r^{2i+1}q^r \over  {1-q^{2r}}},
&\qquad \text{ for}&\quad \ i=1,2,3,\cdots.
\tag 5.96\cr  
\endalignat$$
\endproclaim 
\demo{Proof} The proof of the first identity deals with 
$\sd(u,k)\cn(u,k)$.  Starting with (2.83), applying row
operations and (3.66), appealing to (4.9), and then
utilizing (5.11) we have the following computation.   
$$\spreadlines{10 pt}\allowdisplaybreaks\align
& H_n^{(1)}(\{C_{2\nu-1}(q)\})\tag 5.97\cr
& \kern -2 em  = H_n^{(1)}(\{(-1)^{\nu-1}{z^{2\nu}k^2 
\over{2^{2\nu+2}}}
\cdot (sc/d)_{\nu}(k^2)\})\tag 5.98\cr
& \kern -2 em =  \left\{4^{-n(n+1)}
\prod\limits_{r=1}^{2n-1} r!\right\}
\cdot (zk)^{2n^2}\tag 5.99\cr 
& \kern -2 em  =  \left\{4^{-n(n+1)}
\prod\limits_{r=1}^{2n-1}r! \right\}
\cdot \vartheta_2 (0,q)^{4n^2}.
\tag 5.100\cr
\endalign$$ 

Equating (5.97) with (5.100) and then solving for 
$\vartheta_2(0,q)^{4n^2}$ yields (5.93).

The proof of the second identity deals with 
$\sn^2(u,k)$.  Starting with (2.84), applying row
operations and (3.66), appealing to (4.7), and then
utilizing (5.12) we have the following computation.   
$$\spreadlines{10 pt}\allowdisplaybreaks\align
& H_n^{(1)}(\{D_{2\nu+1}(q)\})\tag 5.101\cr
& \kern -2 em  = H_n^{(1)}(\{(-1)^{\nu-1}{z^{2\nu+2}k^2 
\over{2^{2\nu+3}}}
\cdot (sn^2)_{\nu}(k^2)\})\tag 5.102\cr
& \kern -2 em =  \left\{2^{-n(4n+5)}
\prod\limits_{r=1}^{2n} r! \right\}
\cdot (4z^2k)^{n(n+1)}\tag 5.103\cr 
& \kern -2 em  =  \left\{2^{-n(4n+5)}
\prod\limits_{r=1}^{2n} r! \right\}
\cdot \vartheta_2 (0,q^{1/2})^{4n(n+1)}.
\tag 5.104\cr
\endalign$$ 

Equating (5.101) with (5.104) and then solving for 
$\vartheta_2 (0,q^{1/2})^{4n(n+1)}$ yields (5.94).
\qed\enddemo

It is sometimes useful to rewrite (5.94) by observing from 
(5.12), (5.11), and (5.9) that 
$$\vartheta_2 (0,q^{1/2})^4 = 4(zk)z = 
4 \vartheta_2 (0,q)^2\vartheta_3(0,q)^2.\tag 5.105$$
Equation (5.94) then becomes 
$$\left[\vartheta_2 (0,q)
\vartheta_3 (0,q)\right]^{2n(n+1)} =
\left\{2^{n(2n+3)}
\prod\limits_{r=1}^{2n}(r!)^{-1}\right\}\cdot 
\det (D_{2(r+s-1)+1})_{1\leq r,s\leq n},
\tag 5.106$$
where $D_{2i+1}$ is defined by (5.96).

The $n=1$ case of equation (5.93) of Theorem 5.11,
written as in (5.99), is given by \cite{117, Eqn. (9.), Section
40}.  Similarly, the $n=1$ case of (5.94), expressed as in 
(5.103), is given by \cite{117, Eqn. (3.), Section 41}.  
Moreover, taking $q\mapsto -q$ or $q\mapsto q^2$ in the 
$n=1$ case of (5.103) yields Eqn. (4.) or Eqn. (5.),
respectively, of \cite{117, Section 41}.

Keeping in mind the Kac--Wakimoto conjectured identities for triangular
numbers in \cite{120, pp. 452}, 
we next rewrite Theorem 5.11 in terms of the sum 
$\vartriangle \kern-.3em (q)$ in (5.17).  Applying the
relation (5.18) to the left-hand sides of (5.93) and (5.94) 
immediately gives the following corollary. 
\proclaim{Corollary 5.12} Let 
$\vartriangle \kern-.3em (q)$ be defined by 
\hbox{\rm(5.17)}, and let $n=1,2,3,\cdots$.   
We then have
$$\spreadlines{6 pt}\allowdisplaybreaks\align
\kern -4 em\vartriangle \kern-.3em (q^2)^{4n^2}= &\ 
\left\{4^{-n(n-1)}q^{-n^2}
\prod\limits_{r=1}^{2n-1}(r!)^{-1}\right\}\cdot 
\det (C_{2(r+s-1)-1})_{1\leq r,s\leq n},\tag 5.107\cr
\kern -7.17 em \text{and}\kern 7.17 em  
\kern -4 em\vartriangle \kern-.3em (q)^{4n(n+1)}= &\ 
\left\{2^{n}q^{-n(n+1)/2}
\prod\limits_{r=1}^{2n}(r!)^{-1}\right\}\cdot 
\det (D_{2(r+s-1)+1})_{1\leq r,s\leq n},\tag 5.108\cr
\endalign$$
where $C_{2i-1}$ and $D_{2i+1}$ are defined by 
\hbox{\rm(5.95)} and \hbox{\rm(5.96)},
respectively. 
\endproclaim 

The $n=1$ cases of (5.107) and (5.108) are the classical
identities of Legendre \cite{139}, \cite{21, Eqns. (ii) and
(iii), pp. 139} given by Theorem 7.7. 

We next have the four single $H_n^{(1)}$ and 
$H_n^{(2)}$ Hankel determinant identities related to 
$\cn(u,k)$ and $\dn(u,k)$.  

We begin with the following $H_n^{(1)}$ theorem.
\proclaim{Theorem 5.13} Let 
$\vartheta_2 (0,q)$ and $\vartheta_3 (0,q)$ be defined by 
\hbox{\rm(1.2)} and \hbox{\rm(1.1)}, respectively. 
Let $n=1,2,3,\cdots$.    We then have
$$\vartheta_2 (0,q)^{2n^2}
\vartheta_3 (0,q)^{2n(n-1)}=\,
\left\{4^{n}
\prod\limits_{r=1}^{n-1}(2r)!^{-2}\right\}\cdot 
\det (T_{2(r+s-1)-2})_{1\leq r,s\leq n},
\tag 5.109$$
where $T_{2i-2}$ is defined by
$$T_{2i-2}\equiv\ T_{2i-2}(q) 
 :=\sum\limits_{r=1}^{\infty}
{(2r-1)^{2i-2}q^{r-\tfrac{1}{2}}\over  {1+q^{2r-1}}},
\qquad \text{ for}\quad i=1,2,3,\cdots.\tag 5.110$$
\endproclaim 
\demo{Proof} Our analysis deals with 
$\cn(u,k)$.  Starting with (2.85), applying row
operations and (3.67), appealing to (4.2), and then
utilizing both (5.9) and (5.11) we have the following
computation.   
$$\spreadlines{10 pt}\allowdisplaybreaks\align
& H_n^{(1)}(\{T_{2\nu -2}(q)\})\tag 5.111\cr
& \kern -2 em  = H_n^{(1)}(\{(-1)^{\nu-1}
{z^{2\nu-1}k\over 4}
\cdot (cn)_{\nu -1}(k^2)\})\tag 5.112\cr
& \kern -2 em =  \left\{4^{-n}
\prod\limits_{r=1}^{n-1}(2r)!^2\right\}
\cdot z^{2n^2-n}k^{n^2}\tag 5.113\cr 
& \kern -2 em =  \left\{4^{-n}
\prod\limits_{r=1}^{n-1}(2r)!^2\right\}
\cdot (zk)^{n^2}z^{n(n-1)}\tag 5.114\cr 
& \kern -2 em  =  \left\{4^{-n}
\prod\limits_{r=1}^{n-1}(2r)!^2\right\}\cdot 
\vartheta_2(0,q)^{2n^2}\vartheta_3(0,q)^{2n(n-1)}.
\tag 5.115\cr
\endalign$$ 

Equating (5.111) with (5.115) and then solving for 
$\vartheta_2(0,q)^{2n^2}\vartheta_3(0,q)^{2n(n-1)}$
yields (5.109).
\qed\enddemo

The $n=1$ case of Theorem 5.13, when written as in 
(5.114), is given by \cite{117, Eqn. (5.), Section 40}. The 
$n=1$ case of (5.109) gives the identity in \cite{21, Ex.
(iv), pp. 139}.

Applying the relation (5.105) to the left-hand-side of
(5.109) immediately implies that 
$$\vartheta_2 (0,q)^{2n}
\vartheta_2 (0,q^{1/2})^{4n(n-1)}=\,
\left\{4^{n^2}
\prod\limits_{r=1}^{n-1}(2r)!^{-2}\right\}\cdot 
\det (T_{2(r+s-1)-2})_{1\leq r,s\leq n},\tag 5.116$$
where $T_{2i-2}$ is defined by (5.110).  The form of the
exponents and arguments in the left-hand-side of (5.116) is
much closer than (5.109) to that of the product  sides of the
eta function identities in Appendix I of \cite{151}. 

Applying (5.18) to the left-hand-side of the $n=2$ case of
(5.116) yields an elegant product formula for 
$\vartriangle \kern-.3em (q^2)^4
\kern-.2em \vartriangle \kern-.3em (q)^8$.  
This expresses the product of the sums in entries (ii)
and (iii) of \cite{21, pp. 139} as a simple $2\times 2$ 
determinant of Lambert series.

The single $H_n^{(2)}$ Hankel determinant identity 
related to $\cn(u,k)$ is given by the following theorem.
\proclaim{Theorem 5.14} Let 
$\vartheta_2 (0,q)$ and $\vartheta_3 (0,q)$ be defined by 
\hbox{\rm(1.2)} and \hbox{\rm(1.1)}, respectively. 
Let $n=1,2,3,\cdots$.    We then have
$$\vartheta_2 (0,q)^{2n^2}
\vartheta_3 (0,q)^{2n(n+1)}=\,
\left\{4^{n}
\prod\limits_{r=1}^{n}(2r-1)!^{-2}\right\}\cdot 
\det (T_{2(r+s-1)})_{1\leq r,s\leq n},
\tag 5.117$$
where $T_{2i}$ is determined by \hbox{\rm(5.110)}.
\endproclaim 
\demo{Proof} Our analysis again deals with 
$\cn(u,k)$.  Starting with (2.85), applying row
operations and (3.66), appealing to (4.19), and then
utilizing both (5.9) and (5.11) we have the following
computation.   
$$\spreadlines{10 pt}\allowdisplaybreaks\align
& H_n^{(2)}(\{T_{2\nu -2}(q)\})\tag 5.118\cr
& \kern -2 em  = H_n^{(1)}(\{T_{2\nu}(q)\})\tag 5.119\cr
& \kern -2 em  = H_n^{(1)}(\{(-1)^{\nu }
{z^{2\nu+1}k\over 4}
\cdot (cn)_{\nu}(k^2)\})\tag 5.120\cr
& \kern -2 em =  \left\{4^{-n}
\prod\limits_{r=1}^{n}(2r-1)!^2\right\}
\cdot z^{2n^2+n}k^{n^2}\tag 5.121\cr 
& \kern -2 em =  \left\{4^{-n}
\prod\limits_{r=1}^{n}(2r-1)!^2\right\}
\cdot (zk)^{n^2}z^{n(n+1)}\tag 5.122\cr 
& \kern -2 em  =  \left\{4^{-n}
\prod\limits_{r=1}^{n}(2r-1)!^2\right\}\cdot 
\vartheta_2 (0,q)^{2n^2}\vartheta_3 (0,q)^{2n(n+1)}.
\tag 5.123\cr
\endalign$$ 

Equating (5.119) with (5.123) and then solving for 
$\vartheta_2 (0,q)^{2n^2}
\vartheta_3 (0,q)^{2n(n+1)}$
yields (5.117).
\qed\enddemo

The $n=1$ case of Theorem 5.14, when written as in
(5.121), is given by \cite{117, Eqn. (40.), Section 40}. 

Keeping in mind (2.88), the Maclaurin series expansions for 
$\cn(u,k)$ and $\sn(u,k)\dn(u,k)$ in (2.62) and
(2.74), and the derivative formula for $\cn(u,k)$, it follows
that Theorem 5.14 is equivalent to the single $H_n^{(1)}$
Hankel determinant theta function identity related to 
$\sn(u,k)\dn(u,k)$. 

In \cite{184, pp. 331} Perron equates the Laplace transform
of $\cn(u,k)$ with its C-fraction expansion and Fourier
series.  This was just a convenient reference.  No applications
of this equality of all three was given.

Applying the relation (5.105) to the left-hand-side of
(5.117) immediately implies that 
$$\vartheta_2 (0,q^{1/2})^{4n^2}
\vartheta_3 (0,q)^{2n}
=\,\left\{4^{n(n+1)}
\prod\limits_{r=1}^{n}(2r-1)!^{-2}\right\}\cdot 
\det (T_{2(r+s-1)})_{1\leq r,s\leq n},\tag 5.124$$
where $T_{2i}$ is determined by (5.110). 

We next utilize (5.18), with $q\mapsto q^{1/2}$, to rewrite
the left-hand-side of (5.124).  Some further simplification
then yields the following corollary.
\proclaim{Corollary 5.15} Let 
$\vartheta_3 (0,q)$ and $\vartriangle \kern-.3em (q)$ be
defined by \hbox{\rm(1.1)} and \hbox{\rm(5.17)},
respectively.  Let $n=1,2,3,\cdots$.    We then have
$$\vartheta_3 (0,q)^{2n}
\vartriangle \kern-.3em (q)^{(2n)^2}=\,
\left\{4^{-n(n-1)}q^{-n(n-1)/2}
\prod\limits_{r=1}^{n}(2r-1)!^{-2}\right\}\cdot 
\det ({\widetilde T}_{2(r+s-1)})_{1\leq r,s\leq n},
\tag 5.125$$
where ${\widetilde T}_{2i}$ is defined by
$${\widetilde T}_{2i}\equiv\ {\widetilde T}_{2i}(q) 
 :=\sum\limits_{r=1}^{\infty}
{(2r-1)^{2i}q^{r-1}\over  {1+q^{2r-1}}},
\qquad \text{ for}\quad i=1,2,3,\cdots.\tag 5.126$$
\endproclaim 

When viewed as a generating function in $q^N$, 
(for $N=0,1,2,3,\cdots$), the left-hand-side of (5.125) has
an elegant combinatorial interpretation.  Motivated by
\cite{6, pp. 506--508}, the coefficient of $q^N$ is the
number of ways of writing $N$ as a sum of $2n$ squares
and $(2n)^2$ triangular numbers.  The resulting expansion
on the right-hand-side of (5.125) is just as simple as those in
Corollary 5.12.

The $n=1$ case of Corollary 5.15 is given by \cite{21, Ex.
(v), pp. 139}.  

We now consider the theta function identities related to 
$\dn(u,k)$.  We start with the following $H_n^{(1)}$ 
theorem.  
\proclaim{Theorem 5.16 } Let 
$\vartheta_2 (0,q)$ and $\vartheta_3 (0,q)$ be defined by 
\hbox{\rm(1.2)} and \hbox{\rm(1.1)}, respectively. 
Let $n=1,2,3,\cdots$.    We then have
$$\spreadlines{6 pt}\allowdisplaybreaks\align 
\kern -2 em\vartheta_2 (0,q)^{2n(n-1)}
\vartheta_3 (0,q)^{2n^2}=\,
&\left\{4^{n^2}
\prod\limits_{r=1}^{n-1}(2r)!^{-2}\right\}\cdot 
\det (N_{2(r+s-1)-2})_{1\leq r,s\leq n}\tag 5.127a\cr
\kern -2 em+&\left\{4^{n^2-1}
\prod\limits_{r=1}^{n-1}(2r)!^{-2}\right\}\cdot 
\det (N_{2(r+s)})_{1\leq r,s\leq n-1},\tag 5.127b\cr
\endalign$$
where $N_{2i-2}$ is defined by 
$$N_{2i-2}\equiv N_{2i-2}(q) 
 :=\sum\limits_{r=1}^{\infty}
{r^{2i-2}q^{r} \over {1+q^{2r}}},
\qquad \text{ for}\quad i=1,2,3,\cdots.\tag 5.128$$
\endproclaim 
\demo{Proof}  Our analysis deals with $\dn(u,k)$. Starting
with (2.86) and (2.87), solving for \newline 
$(-1)^m(z^{2m+1}/2^{2m+2})
\cdot (dn)_{m}(k^2)$ for $m\geq 0$,  applying
row operations, recalling (3.55), and simplifying, leads to the
following identity.
$$\spreadlines{10 pt}\allowdisplaybreaks\align
&H_n^{(1)}(\{N_{2\nu-2}(q)\}) + 
{\tfrac{1}{4}}H_{n-1}^{(3)}(\{N_{2\nu-2}(q)\}) 
\tag 5.129\cr
& \kern -2 em  = H_n^{(0)}(\{(-1)^{\nu }
{z^{2\nu+1}\over 4^{\nu+1}}
\cdot (dn)_{\nu}(k^2)\})\cr
& \kern -2 em  = H_n^{(1)}(\{(-1)^{\nu-1 }
{z^{2\nu-1}\over 4^{\nu }}
\cdot (dn)_{\nu-1}(k^2)\}).\tag 5.130\cr
\endalign$$ 

Next, applying row operations and (3.67) to (5.130),
appealing to (4.2), and then utilizing both (5.9) and (5.11)
we obtain
$$\spreadlines{10 pt}\allowdisplaybreaks\align
&\left\{4^{-n^2}
\prod\limits_{r=1}^{n-1}(2r)!^2\right\}
\cdot z^{2n^2-n}k^{n^2-n}\tag 5.131\cr 
& \kern -2 em =  \left\{4^{-n^2}
\prod\limits_{r=1}^{n-1}(2r)!^2\right\}
\cdot (zk)^{n(n-1)}z^{n^2}\tag 5.132\cr 
& \kern -2 em  =  \left\{4^{-n^2}
\prod\limits_{r=1}^{n-1}(2r)!^2\right\}
\cdot\vartheta_2 (0,q)^{2n(n-1)}
\vartheta_3 (0,q)^{2n^2}.
\tag 5.133\cr
\endalign$$ 

Equating (5.129) with (5.133) and then solving for 
$\vartheta_2 (0,q)^{2n(n-1)}
\vartheta_3 (0,q)^{2n^2}$ yields (5.127).
\qed\enddemo

The $n=1$ case of Theorem 5.16, when written as in
(5.131), is given by \cite{117, Eqn. (4.), Section 40}.  

The $n=1$ case of Theorem 5.16 is equivalent to 
Jacobi's $2$-squares identity in 
\cite{21, Entry 8(i),  pp. 114}.  To see this, note that 
the $0\times 0$ determinant in  (5.127b) is defined to be
$1$, and that an elementary argument \cite{21, pp.
115} involving geometric series and an interchange of
summation gives 
$$\sum\limits_{r=1}^{\infty}
{q^{r} \over {1+q^{2r}}}=
\sum\limits_{r=1}^{\infty}{(-1)^{r-1}
q^{2r-1} \over {1-q^{2r-1}}}.\tag 5.134$$

Multiplying both sides of (5.127) by $4^{n^2}$, it is
immediate from (5.105) that the product side of (5.127)
can be written as $\vartheta_2 (0,q^{1/2})^{4n^2}
\vartheta_2 (0,q)^{-2n}$.  Just as in (5.116), this form is
much closer to \cite{151, Appendix I}.

The single $H_n^{(2)}$ Hankel determinant identity 
related to $\dn(u,k)$ is given by the following theorem.
\proclaim{Theorem 5.17} Let 
$\vartheta_2 (0,q)$ and $\vartheta_3 (0,q)$ be defined by 
\hbox{\rm(1.2)} and \hbox{\rm(1.1)}, respectively. 
Let $n=1,2,3,\cdots$.    We then have
$$\vartheta_2 (0,q)^{2n(n+1)}
\vartheta_3 (0,q)^{2n^2}=\,
\left\{4^{n(n+1)}
\prod\limits_{r=1}^{n}(2r-1)!^{-2}\right\}\cdot 
\det (N_{2(r+s-1)})_{1\leq r,s\leq n},
\tag 5.135$$
where $N_{2i}$ is determined by \hbox{\rm(5.128)}.
\endproclaim 
\demo{Proof} Our analysis again deals with 
$\dn(u,k)$.  Starting with (2.86), applying row
operations and (3.66), appealing to (4.20), and then
utilizing both (5.9) and (5.11) we have the following
computation.   
$$\spreadlines{10 pt}\allowdisplaybreaks\align
& H_n^{(2)}(\{N_{2\nu -2}(q)\})\tag 5.136\cr
& \kern -2 em  = H_n^{(1)}(\{N_{2\nu}(q)\})\tag 5.137\cr
& \kern -2 em  = H_n^{(1)}(\{(-1)^{\nu }
{z^{2\nu+1}\over 4^{\nu+1}}
\cdot (dn)_{\nu}(k^2)\})\tag 5.138\cr
& \kern -2 em =  \left\{4^{-n(n+1)}
\prod\limits_{r=1}^{n}(2r-1)!^2\right\}
\cdot z^{2n^2+n}k^{n^2+n}\tag 5.139\cr 
& \kern -2 em =  \left\{4^{-n(n+1)}
\prod\limits_{r=1}^{n}(2r-1)!^2\right\}
\cdot (zk)^{n(n+1)}z^{n^2}\tag 5.140\cr 
& \kern -2 em  =  \left\{4^{-n(n+1)}
\prod\limits_{r=1}^{n}(2r-1)!^2\right\}\cdot 
\vartheta_2 (0,q)^{2n(n+1)}\vartheta_3 (0,q)^{2n^2}.
\tag 5.141\cr
\endalign$$ 

Equating (5.137) with (5.141) and then solving for 
$\vartheta_2 (0,q)^{2n(n+1)}
\vartheta_3 (0,q)^{2n^2}$
yields (5.135).
\qed\enddemo

The $n=1$ case of Theorem 5.17, when written as in
(5.139), is given by \cite{117, Eqn. (42.), Section 40}. 
See also \cite{21, Entry 17(ii), pp. 138}.  

Keeping in mind (2.89), the Maclaurin series expansions for 
$\dn(u,k)$ and $\sn(u,k)\cn(u,k)$ in (2.63) and
(2.75), and the derivative formula for $\dn(u,k)$, it follows
that Theorem 5.17 is equivalent to the single $H_n^{(1)}$
Hankel determinant theta function identity related to 
$\sn(u,k)\cn(u,k)$. 

Applying the relation (5.105) to the left-hand-side of
(5.135) immediately implies that 
$$\vartheta_2 (0,q^{1/2})^{4n^2}
\vartheta_2 (0,q)^{2n}
=\,\left\{4^{n(2n+1)}
\prod\limits_{r=1}^{n}(2r-1)!^{-2}\right\}\cdot 
\det (N_{2(r+s-1)})_{1\leq r,s\leq n},\tag 5.142$$
where $N_{2i}$ is determined by (5.128). 

We next utilize (5.18), with $q\mapsto q^{1/2}$, to rewrite
the left-hand-side of (5.142).  Some further simplification
then yields the following corollary.
\proclaim{Corollary 5.18} Let 
$\vartriangle \kern-.3em (q)$ be
defined by \hbox{\rm(5.17)}.
 Let $n=1,2,3,\cdots$.    We then have
$$\vartriangle \kern-.3em (q)^{(2n)^2}
\vartriangle \kern-.3em (q^2)^{2n}=\,
\left\{q^{-n(n+1)/2}
\prod\limits_{r=1}^{n}(2r-1)!^{-2}\right\}\cdot 
\det ({N}_{2(r+s-1)})_{1\leq r,s\leq n},
\tag 5.143$$
where ${N}_{2i}$ is determined by \hbox{\rm(5.128)}.
\endproclaim 

Analogous to (5.125), the coefficient of $q^N$ in (5.143) is
the number of ways of writing $N$ as a sum of 
$(2n)^2$ triangular numbers and $2n$ doubles of triangular
numbers.

By considering quotients of our theta function identities
related to a given Jacobi elliptic function and its derivative, or
equivalently, to the pairs of $H_n^{(1)}$ and $H_n^{(2)}$
identities, we are able to write $\vartheta_2(0,q)^{4n}$ 
and $\vartheta_3 (0,q)^{4n}$ as a quotient of 
$n\times n$ determinants of Lambert series.  To obtain
such an identity for $\vartheta_2(0,q)^{4n}$, divide 
(5.135) by (5.127).  (Here, the denominator will involve a
sum of two determinants.)  The quotient of determinants
identities for $\vartheta_3(0,q)^{4n}$ are simpler.  The
first two identities for $\vartheta_3(0,q)^{4n}$ are
obtained by either dividing (5.75) by (5.54), or (5.117) by
(5.109).  The third results by dividing the $q\mapsto -q$
cases of (5.37) by (5.20).  Equating pairs of these three
quotients of determinants gives three determinant
identities.  Finally, equating the $n\mapsto n^2$ case of
each of the three quotient of determinants identities for 
$\vartheta_3(0,q)^{4n}$ with the $q\mapsto -q$ case of
(5.20) yields three more interesting determinantal identities
involving Lambert series.  Each of these last three identities
equates a quotient of $n^2\times n^2$ determinants to
the product of an $n\times n$ determinant and a certain
constant.  

For convenience and future reference we now extract from
the above analysis in this section our generalization to
infinite families of all $21$ of Jacobi's \cite{117, Sections 40,
41, 42} explicitly stated degree $2, 4, 6, 8$ Lambert series
expansions of classical theta functions.  We follow Jacobi in
expressing everything in terms of the elliptic function
parameters $z$, $k$, and $q$.  For three of our identities
we need the additional Lambert series defined by 
$$\spreadlines{6 pt}\allowdisplaybreaks\align 
{\widehat C}_{2m-1}(q):&=\sum\limits_{r=1}^{\infty}
{(2r-1)^{2m-1}(-1)^{r}
q^{r-\tfrac{1}{2}} \over {1+q^{2r-1}}},\tag 5.144\cr
{\widehat T}_{2m-2}(q):&=\sum\limits_{r=1}^{\infty}
{(2r-1)^{2m-2}(-1)^{r}
q^{r-\tfrac{1}{2}} \over {1-q^{2r-1}}}.\tag 5.145\cr
\endalign$$

We have the following theorem.  
\proclaim{Theorem 5.19} Let  $z:= 2\bk(k)/\pi
\equiv 2\bk/\pi$, as in \hbox{\rm(2.1)}, with $k$ the
modulus.  Take $k':=\sqrt{1-k^2}$ and $q$ as in
\hbox{\rm(2.4)}.  Let the Bernoulli numbers $B_n$ and 
Euler numbers $E_n$ be defined by \hbox{\rm(2.60)} and
\hbox{\rm(2.61)},  respectively.  Take $U_{2m-1}(-q)$, 
$G_{2m+1}(-q)$, $R_{2m-2}(q)$, $C_{2m-1}(q)$,
$D_{2m+1}(q)$, $T_{2m-2}(q)$, $N_{2m}(q)$, and
$N_{0}(q)$ to be the Lambert series in 
\hbox{\rm Theorem 2.4}.  Let ${\widehat C}_{2m-1}(q)$ 
and ${\widehat T}_{2m-2}(q)$ be the Lambert series in
\hbox{\rm(5.144)} and \hbox{\rm(5.145)}, respectively.  
Take $H_n^{(m)}(\{c_\nu\})$ to be the $n\times n$
determinants  in \hbox{\rm Definition 3.3}.   Let
$n=1,2,3,\cdots$.  We then have the following expansions.
$$\spreadlines{6 pt}\allowdisplaybreaks\align 
\kern -2 em k^{n^2-n}z^{2n^2-n}=\,
&\left\{4^{n^2}
\prod\limits_{r=1}^{n-1}(2r)!^{-2}\right\}\cdot 
H_n^{(1)}(\{N_{2\nu-2}(q)\})\tag 5.146a\cr
\kern -2 em&+\left\{4^{n^2-1}
\prod\limits_{r=1}^{n-1}(2r)!^{-2}\right\}\cdot 
H_{n-1}^{(3)}(\{N_{2\nu-2}(q)\}),\tag 5.146b\cr
\kern -2 em {(k')}^{n^2-n}z^{2n^2-n}=\,
&\left\{(-1)^{n}4^{n}
\prod\limits_{r=1}^{n-1}(2r)!^{-2}\right\}\cdot\cr
\kern -2 em&H_n^{(1)}(\{R_{2\nu -2}(-q)
-(-1)^{\nu -1}\cdot {\tfrac {1}{4}}\cdot
|E_{2\nu -2}|\}),\tag 5.147\cr
\kern -2 em k^{n^2}z^{2n^2-n}=\,
&\left\{4^{n}
\prod\limits_{r=1}^{n-1}(2r)!^{-2}\right\}\cdot 
H_n^{(1)}(\{T_{2\nu-2}(q)\}),\tag 5.148\cr
\kern -2 em k^{n^2}{(k')}^{n^2-n}z^{2n^2-n}=\,
&\left\{(-1)^{n(n+1)/2}4^{n}
\prod\limits_{r=1}^{n-1}(2r)!^{-2}\right\}\cdot 
H_n^{(1)}(\{{\widehat T}_{2\nu-2}(q)\}),\tag 5.149\cr
\kern -2 em {(k')}^{n^2}z^{2n^2-n}=\,
&\left\{(-1)^{n}4^{n}
\prod\limits_{r=1}^{n-1}(2r)!^{-2}\right\}\cdot\cr
\kern -2 em&H_n^{(1)}(\{R_{2\nu -2}(q)
-(-1)^{\nu -1}\cdot {\tfrac {1}{4}}\cdot
|E_{2\nu -2}|\}),\tag 5.150\cr
\kern -2 em k^{n^2-n}{(k')}^{n^2}z^{2n^2-n}=\,
&\left\{(-1)^{n(n-1)/2}4^{n^2}
\prod\limits_{r=1}^{n-1}(2r)!^{-2}\right\}\cdot 
H_n^{(1)}(\{N_{2\nu-2}(-q)\})\tag 5.151a\cr
\kern -2 em&+\left\{(-1)^{n(n-1)/2}4^{n^2-1}
\prod\limits_{r=1}^{n-1}(2r)!^{-2}\right\}\cdot 
H_{n-1}^{(3)}(\{N_{2\nu-2}(-q)\}),\tag 5.151b\cr
\kern -2 em {(k')}^{n^{2}/2}(1+k')^{n^2-n}&z^{2n^2-n}=\,
\left\{(-1)^{n}2^{n(n+1)}
\prod\limits_{r=1}^{n-1}(2r)!^{-2}\right\}\cdot\cr
\kern -2 em&H_n^{(1)}(\{R_{2\nu -2}(q^2)
-(-1)^{\nu -1}\cdot {\tfrac {1}{4}}\cdot
|E_{2\nu -2}|\}),\tag 5.152\cr
\kern -2 em z^{2n^2}=\,
&\left\{(-1)^n2^{2n^2+n}
\prod\limits_{r=1}^{2n-1}(r!)^{-1}\right\}\cdot\cr
\kern -2 em& H_n^{(1)}(\{U_{2\nu -1}(-q)
-(-1)^{\nu-1}{(2^{2\nu}-1)
\over {4\nu}}\cdot |B_{2\nu}|\}),\tag 5.153\cr
\kern -2 em (kz)^{2n^2}=\,
&\left\{4^{n(n+1)}
\prod\limits_{r=1}^{2n-1}(r!)^{-1}\right\}\cdot 
H_n^{(1)}(\{C_{2\nu-1}(q)\}),\tag 5.154\cr
\kern -2 em (k'z)^{2n^2}=\,
&\left\{(-1)^n2^{2n^2+n}
\prod\limits_{r=1}^{2n-1}(r!)^{-1}\right\}\cdot\cr
\kern -2 em& H_n^{(1)}(\{U_{2\nu -1}(q)
-(-1)^{\nu-1}{(2^{2\nu}-1)
\over {4\nu}}\cdot |B_{2\nu}|\}),\tag 5.155\cr
\kern -2 em (kk')^{n^2}z^{2n^2}=\,
&\left\{(-1)^{n(n+1)/2}4^{n}
\prod\limits_{r=1}^{2n-1}(r!)^{-1}\right\}\cdot 
H_n^{(1)}(\{{\widehat C}_{2\nu-1}(q)\}),\tag 5.156\cr
\kern -2 em (k')^{n^2}z^{2n^2}=\,
&\left\{(-1)^n2^{2n^2+n}
\prod\limits_{r=1}^{2n-1}(r!)^{-1}\right\}\cdot\cr
\kern -2 em& H_n^{(1)}(\{U_{2\nu -1}(q^2)
-(-1)^{\nu-1}{(2^{2\nu}-1)
\over {4\nu}}\cdot |B_{2\nu}|\}),\tag 5.157\cr
\kern -2 em k^{n^2}z^{2n^2}=\,
&\left\{4^{n}
\prod\limits_{r=1}^{2n-1}(r!)^{-1}\right\}\cdot 
H_n^{(1)}(\{C_{2\nu-1}(\sqrt{q})\}),\tag 5.158\cr
\kern -2 em k^{n^2}z^{2n^2+n}=\,
&\left\{4^{n}
\prod\limits_{r=1}^{n}(2r-1)!^{-2}\right\}\cdot 
H_n^{(1)}(\{T_{2\nu}(q)\}),\tag 5.159\cr
\kern -2 em (k')^{n^2}z^{2n^2+n}=\,
&\left\{4^{n}
\prod\limits_{r=1}^{n}(2r-1)!^{-2}\right\}\cdot\cr
\kern -2 em& H_n^{(1)}(\{R_{2\nu}(q)
-(-1)^{\nu}\cdot {\tfrac {1}{4}}\cdot
|E_{2\nu}|\}),\tag 5.160\cr
\kern -2 em k^{n^2+n}z^{2n^2+n}=\,
&\left\{4^{n(n+1)}
\prod\limits_{r=1}^{n}(2r-1)!^{-2}\right\}\cdot
 H_n^{(1)}(\{N_{2\nu}(q)\}),\tag 5.161\cr
\kern -2 em k^{n^2}(k')^{n^2+n}z^{2n^2+n}=\,
&\left\{(-1)^{n(n+1)/2}4^{n}
\prod\limits_{r=1}^{n}(2r-1)!^{-2}\right\}\cdot
 H_n^{(1)}(\{{\widehat T}_{2\nu}(q)\}),\tag 5.162\cr
\kern -2 em (k')^{n^2+n}z^{2n^2+n}=\,
&\left\{4^{n}
\prod\limits_{r=1}^{n}(2r-1)!^{-2}\right\}\cdot\cr
\kern -2 em& H_n^{(1)}(\{R_{2\nu}(-q)
-(-1)^{\nu}\cdot {\tfrac {1}{4}}\cdot
|E_{2\nu}|\}),\tag 5.163\cr
\kern -2 em k^{n^2+n}(k')^{n^2}z^{2n^2+n}=\,
&\left\{(-1)^{n(n+1)/2}4^{n(n+1)}
\prod\limits_{r=1}^{n}(2r-1)!^{-2}\right\}\cdot\cr
\kern -2 em& H_n^{(1)}(\{N_{2\nu}(-q)\}),
\tag 5.164\cr
\kern -2 em k^{n^2+n}z^{2n^2+2n}=\,
&\left\{2^{2n^2+3n}
\prod\limits_{r=1}^{2n}(r!)^{-1}\right\}\cdot
 H_n^{(1)}(\{D_{2\nu+1}(q)\}),\tag 5.165\cr
\kern -2 em k^{n^2+n}(k')^{n^2+n}z^{2n^2+2n}=\,
&\left\{(-1)^{n(n+1)/2}2^{2n^2+3n}
\prod\limits_{r=1}^{2n}(r!)^{-1}\right\}\cdot\cr
\kern -2 em& H_n^{(1)}(\{D_{2\nu+1}(-q)\}),
\tag 5.166\cr
\kern -2 em (kz)^{2n^2+2n}=\,
&\left\{2^{4n^2+5n}
\prod\limits_{r=1}^{2n}(r!)^{-1}\right\}\cdot
 H_n^{(1)}(\{D_{2\nu+1}(q^2)\}),\tag 5.167\cr
\kern -2 em (k'z)^{2n^2+2n}=\,
&\left\{2^{2n^2+3n}
\prod\limits_{r=1}^{2n}(r!)^{-1}\right\}\cdot\cr
\kern -2 em& H_n^{(1)}(\{G_{2\nu +1}(q)
-(-1)^{\nu}{(2^{2\nu + 2}-1)
\over {4(\nu + 1)}}\cdot |B_{2\nu + 2}|\}),
\tag 5.168\cr
\kern -2 em z^{2n^2+2n}=\,
&\left\{2^{2n^2+3n}
\prod\limits_{r=1}^{2n}(r!)^{-1}\right\}\cdot\cr
\kern -2 em& H_n^{(1)}(\{G_{2\nu +1}(-q)
-(-1)^{\nu}{(2^{2\nu + 2}-1)
\over {4(\nu + 1)}}\cdot |B_{2\nu + 2}|\}).
\tag 5.169\cr
\endalign$$
\endproclaim 
\demo{Proof} Ten of the identities in (5.146)--(5.169) are
immediate from the above analysis in this section.  In
particular, (5.146), (5.148), (5.150), (5.153), (5.154),
(5.159), (5.160), (5.161), (5.165), (5.169) are just 
(5.129) and (5.131), (5.111) and (5.113),  (5.58) and
(5.60), (5.24) and (5.26), (5.97) and (5.99), (5.119) and
(5.121), (5.78) and (5.80), (5.137) and (5.139), (5.101)
and (5.103), (5.41) and (5.43), respectively.  

Equation (5.158) results from applying the Gau{\ss}
transformation ($q\mapsto \sqrt{q}$, $kz\mapsto 
2\sqrt{k}z$) to (5.154).  

Equations (5.147), (5.149), (5.151), (5.155), (5.156),
(5.162), (5.163), (5.164), (5.166), (5.168) follow by
applying Jacobi's transformation  ($q\mapsto -q$,
$k'z\mapsto z$, $z\mapsto k'z$,
$kz\mapsto i kz$) to equations (5.150), (5.148),
(5.146), (5.153), (5.158), (5.159), (5.160), (5.161),
(5.165), (5.169), respectively.  

Finally, equations (5.152), (5.157), (5.167) follow by
applying Landen's transformation ($q\mapsto q^2$,
$k'z\mapsto \sqrt{k'}z$, $z\mapsto {\frac{1+k'}{2}}z$,
$kz\mapsto  {\frac{1-k'}{2}}z$) to equations (5.150),
(5.155), and (5.165), respectively. 
\qed\enddemo

We have listed the identities in Theorem 5.19 in the order in
which their $n=1$ cases appear in \cite{117}.  Both (5.146)
and (5.147), (5.148) and (5.149), (5.150) and (5.151),
generalize, respectively, \cite{117, Eqn. (4.) (see also (30.)),
Eqn. (5.) (see also (31.)), Eqn. (6.) (see also (32.)),
Section 40}.  Equations (5.152)--(5.169) generalize (in
order) the following equations in \cite{117}: (7.)--(13.) 
(see also (33.)--(39.)), (40.)--(45.), of Section 40; (3.)--(5.),
of Section 41; (7.), (8.) of Section 42.

The power $4n(n+1)$ in (5.37), and powers $2n^2+n$ and
$2n^2-n$ of $z$ in Theorem 5.19 also occur as the powers
of  $\eta(q)$ in Macdonald's
expansions corresponding to $A_{\ell}\ (\ell\geq 1)$ (take
$\ell = 2n$ here), $B_{\ell}\ (\ell\geq 3)$, $C_{\ell}\
(\ell\geq 2)$, $BC_{\ell}\ (\ell\geq 1)$, and $D_{\ell}\
(\ell\geq 4)$.  These expansions with powers 
$\ell^2+2\ell$, $2\ell^2+\ell$, $2\ell^2+\ell$,
$2\ell^2-\ell$, and $2\ell^2-\ell$, can be found in pages 
134--135 (eq. (6)(a)-(b)), 135 (eq. (6)(a)), 136 (eq. (6)),
137--138 (eq. (6)(c)), and 138 (eq. (6)), respectively, of
\cite{151}. 

Just after equation (45.) in Section 40 of \cite{117}, Jacobi
described how to immediately obtain Lambert series
expansions for $z^3$, $(kz)^3$, and $(k'z)^3$.  The
coefficient of $q^N$ in these expansions quickly leads to
elegant formulas for counting the number of ways that $N$
can be represented as a sum of $6$ squares, and also as a
sum of $6$ triangular numbers.  Later, Smith obtained both
the analytic and combinatorial formulas for $2$, $4$, $6$,
$8$ squares and triangles in \cite{212, Section
127, pp. 306--311}.  See also \cite{213, pp. 206} for a
discussion of the $6$ and $8$ squares work of Jacobi,
Eisenstein, and M. Liouville.  Glaisher derived the elliptic
function and combinatorial formulations of the $6$ squares
and $6$ triangles results in \cite{86, pp. 9--10}.  
More recent derivations of the $6$ squares formulas can be
found in Ramanujan \cite{193, Eqns. (135), (136),
(145)--(147), Table VI. (entry 1), pp. 158--159}, K.
Ananda-Rau \cite{3, pp. 86}, Carlitz \cite{37}, Grosswald
\cite{94, Eqn. (9.19), pp. 121}, Hardy and Wright \cite{102,
pp. 314--315}, Kac and Wakimoto \cite{120, Eqn. (0.3), pp.
416; Ex. 5.2, pp. 444}, and Nathanson \cite{182, pp. 424;
Section 14.5, pp. 436--440}.  Several more recent
derivations of the $6$ triangles formulas appear in 
Ramanujan \cite{194, Eqn. (3.23) (Lambert series version),
pp. 356}, Kac and Wakimoto \cite{120, Ex. 5.2, pp. 444},
Ono, Robins, and Wahl \cite{183, Theorem 4, pp. 81}, Berndt
\cite{24, Entry 6 (Lambert series version); Corollary 6.2
(arithmetical formula)}, Andrews and Berndt \cite{7}, and
Huard, Ou, Spearman, and Williams \cite{108, Theorem 11,
Section 6}.  

Motivated by Jacobi and Glaisher's treatment of the $6$
squares and triangles formulas we next set $n=2$ in
equations (5.146)--(5.151) of Theorem 5.19 and then 
derive elegant new Lambert series expansions for 
$z^6$, $k^2z^6$, $k^4z^6$, and $(kz)^6$, which do not
use cusp forms. We have the following theorem.
\proclaim{Theorem 5.20} Let  $z:= 2\bk(k)/\pi
\equiv 2\bk/\pi$, as in \hbox{\rm(2.1)}, with $k$ the
modulus. Take $q$ as in \hbox{\rm(2.4)}.  
Let $T_{2m-2}(q)$, ${\widehat T}_{2m-2}(q)$, 
$N_{2m-2}(q)$, and $R_{2m-2}(-q)$ be the Lambert series
determined by \hbox{\rm(2.85)}, \hbox{\rm(5.145)},
\hbox{\rm(2.86)} and \hbox{\rm(2.87)}, and
\hbox{\rm(2.82)}, respectively.We then have the following
expansions.
$$\spreadlines{10 pt}\allowdisplaybreaks\align 
(kz)^6= \ & 
{4}\det\vmatrix  T_0(q) & T_2(q)  \\
  &    \\
T_2(q)  & T_4(q)   
\endvmatrix
+
{4}\det\vmatrix  {\widehat T}_0(q) & {\widehat T}_2(q)  \\
  &    \\
{\widehat T}_2(q)  & {\widehat T}_4(q)   
\endvmatrix,\tag 5.170\cr
k^4z^6= \ & 
{4}\det\vmatrix  T_0(q) & T_2(q)  \\
  &    \\
T_2(q)  & T_4(q)   
\endvmatrix,\tag 5.171\cr
k^2z^6= \ & 
4^2N_4(q)+
{4^3}\det\vmatrix  N_0(q) & N_2(q)  \\
  &    \\
N_2(q)  & N_4(q)   
\endvmatrix,\tag 5.172\cr
z^6=\ & 
4^2N_4(q)+
{4^3}\det\vmatrix  N_0(q) & N_2(q)  \\
  &    \\
N_2(q)  & N_4(q)   
\endvmatrix +\cr
 \ & 
\kern 1em{\tfrac{1}{4}}\det\vmatrix  4R_0(-q)-1 &
4R_2(-q)+1 \\
  &    \\
4R_2(-q)+1 & 4R_4(-q)-5 
\endvmatrix.\tag 5.173\cr
\endalign$$
\endproclaim 
\demo{Proof} Equations (5.172) and (5.171) are the $n=2$
cases of (5.146) and (5.148), respectively.  Equation
(5.170) is immediate from subtracting the $n=2$ case of
(5.149) from (5.171).  Equation (5.173) results from
adding the $n=2$ case of (5.147) to (5.172).  
\qed\enddemo

Equations (5.170) and (5.173) can be viewed as additional
two-dimensional generalizations of Jacobi's $2$ triangles and
$2$ squares results in \cite{117, Eqns. (5.) and (4.), Section
40}.  Furthermore, (5.170) and (5.173) are not the same as
the classical results of Liouvllle \cite{144}, \cite{150,
(setting for number theory publications), pp. 227--230; 
(list of number theory publications), pp. 805--812}, 
\cite{56, pp. 306}, and Glaisher \cite{86, pp. 36}, 
\cite{87, pp. 201--202}, or the more recent treatments of
Ramanujan \cite{193, Eqns. (135), (136), (145)--(147),
Table VI. (entry 2), pp. 158--159}, Lomadze \cite{147,
Theorem 6, pp. 264}, Gundlach \cite{95, pp. 196},
Grosswald \cite{94, Eqn. (9.19), pp. 121}, Kac and
Wakimoto \cite{120, $m=2$ case of Ex. 5.3, pp. 444}, and
Ono, Robins, and Wahl \cite{183, Theorem 7, pp. 85}.  Note
that Theorem 7 of \cite{183, pp. 85} and its resulting
congruence mod $256$ also appears at the very end of
Section 61 of Glaisher's article \cite{86, pp. 37}.  

In order to generalize the above classical work of Jacobi,
Glaisher, and Ramanujan from $6$ squares and $6$ triangles
to $20$ squares and $20$ triangles we first set $n=2$ in
the $6$ equations (5.159)--(5.164) of Theorem 5.19 and
then derive the Lambert series expansions for $z^{10}$, 
$k^2z^{10}$, $k^4z^{10}$, $k^6z^{10}$, $k^8z^{10}$,
and $(kz)^{10}$ in the following theorem.
\proclaim{Theorem 5.21} Let  $z:= 2\bk(k)/\pi
\equiv 2\bk/\pi$, as in \hbox{\rm(2.1)}, with $k$ the
modulus. Take $q$ as in \hbox{\rm(2.4)}.  
Let $T_{2m-2}(q)$, ${\widehat T}_{2m-2}(q)$, 
$N_{2m-2}(q)$, and $R_{2m-2}(q)$ be the Lambert series
determined by \hbox{\rm(2.85)}, \hbox{\rm(5.145)},
\hbox{\rm(2.86)} and \hbox{\rm(2.87)}, and
\hbox{\rm(2.82)}, respectively.We then have the following
expansions.
$$\spreadlines{10 pt}\allowdisplaybreaks\align 
(kz)^{10}= \ & 
{\tfrac{2^{10}}{3}}\left\{
\det\vmatrix  N_2(q) & N_4(q)  \\
  &    \\
N_4(q)  & N_6(q)   
\endvmatrix
-
\det\vmatrix  N_2(-q) & N_4(-q)  \\
  &    \\
N_4(-q)  & N_6(-q)   
\endvmatrix\right\}\cr
\ & 
\kern 2em-{\tfrac{2^{3}}{3^2}}\left\{
\det\vmatrix  T_2(q) & T_4(q)  \\
  &    \\
T_4(q)  & T_6(q)   
\endvmatrix
+
\det\vmatrix  {\widehat T}_2(q) & {\widehat T}_4(q)  \\
  &    \\
{\widehat T}_4(q)  & {\widehat T}_6(q)   
\endvmatrix\right\},\tag 5.174\cr
k^8z^{10}= \ & 
{\tfrac{2^{10}}{3^2}}\left\{
2\det\vmatrix  N_2(q) & N_4(q)  \\
  &    \\
N_4(q)  & N_6(q)   
\endvmatrix
-
\det\vmatrix  N_2(-q) & N_4(-q)  \\
  &    \\
N_4(-q)  & N_6(-q)   
\endvmatrix\right\}\cr
\ & 
\kern 2em-{\tfrac{2^{2}}{3^2}}\left\{
\det\vmatrix  T_2(q) & T_4(q)  \\
  &    \\
T_4(q)  & T_6(q)   
\endvmatrix
+
\det\vmatrix  {\widehat T}_2(q) & {\widehat T}_4(q)  \\
  &    \\
{\widehat T}_4(q)  & {\widehat T}_6(q)   
\endvmatrix\right\},\tag 5.175\cr
k^6z^{10}= \ & 
{\tfrac{2^{10}}{3^2}}\det\vmatrix  N_2(q) & N_4(q)  \\
  &    \\
N_4(q)  & N_6(q)   
\endvmatrix,\tag 5.176\cr
k^4z^{10}= \ & 
{\tfrac{2^{2}}{3^2}}\det\vmatrix  T_2(q) & T_4(q)  \\
  &    \\
T_4(q)  & T_6(q)   
\endvmatrix,\tag 5.177\cr
k^2z^{10}=\ & 
{\tfrac{2^{3}}{3^2}}\det\vmatrix  T_2(q) & T_4(q)  \\
  &    \\
T_4(q)  & T_6(q)   
\endvmatrix 
-{\tfrac{2^{10}}{3^2}}\det\vmatrix  N_2(q) & N_4(q)  \\
  &    \\
N_4(q)  & N_6(q)   
\endvmatrix \cr
 \ & 
\kern 2em+{\tfrac{1}{6^2}}\left\{
\det\vmatrix  4R_2(q)+1 & 4R_4(q)-5 \\
  &    \\
4R_4(q)-5 & 4R_6(q)+61
\endvmatrix\right.\cr
 \ & 
\kern 4em-\left. \det\vmatrix  4R_2(-q)+1 & 4R_4(-q)-5 \\
  &    \\
4R_4(-q)-5 & 4R_6(-q)+61
\endvmatrix\right\},\tag 5.178\cr
z^{10}=\ & 
{\tfrac{2^{2}}{3}}\det\vmatrix  T_2(q) & T_4(q)  \\
  &    \\
T_4(q)  & T_6(q)   
\endvmatrix 
-{\tfrac{2^{11}}{3^2}}\det\vmatrix  N_2(q) & N_4(q)  \\
  &    \\
N_4(q)  & N_6(q)   
\endvmatrix \cr
 \ & 
\kern 2em+{\tfrac{1}{18}}\left\{
{\tfrac{3}{2}}\det\vmatrix  4R_2(q)+1 & 4R_4(q)-5 \\
  &    \\
4R_4(q)-5 & 4R_6(q)+61
\endvmatrix\right.\cr
 \ & 
\kern 4em-\left. \det\vmatrix  4R_2(-q)+1 & 4R_4(-q)-5 \\
  &    \\
4R_4(-q)-5 & 4R_6(-q)+61
\endvmatrix\right\}.\tag 5.179\cr
\endalign$$
\endproclaim 
\demo{Proof} Let $\alpha_1$, $\alpha_2$, $\alpha_3$,
$\alpha_4$, $\alpha_5$, $\alpha_6$ denote the $n=2$
cases, respectively, of equations (5.159), (5.160), (5.161),
(5.162), (5.163), and (5.164).  Then, equations (5.176)
and (5.177) are $\alpha_3$ and $\alpha_1$, respectively. 
Equation (5.175) is
$(-\alpha_1+2\alpha_3+\alpha_4+\alpha_6)$.  Equation
(5.174) results from 
$(-2\alpha_1+2\alpha_4+3\alpha_3+3\alpha_6)$. Equation
(5.178) follows from
$(2\alpha_1+\alpha_2-\alpha_3-\alpha_5)$.  Finally,
equation (5.179) is a consequence of
$(3\alpha_1+3\alpha_2-2\alpha_3-2\alpha_5)$.
\qed\enddemo

The formula for $z^{10}$ in (5.179) is not the same as the
classical formulas of Ramanujan
\cite{193, Eqns. (135), (136), (145)--(147), Table VI.
(entry 5), pp. 158--159}, Rankin \cite{201}, and Lomadze
\cite{147, Theorem 7, pp. 266}.  Furthermore, it appears
likely that the Schur function expansion in Theorem 6.7
applied to the determinants on the right-hand-side of
(5.174) will lead to a proof of the $m=2$ case of Conjecture
5.1 in \cite{120, pp. 445}.  The equivalence should be
nontrivial.

We have developed many more consequences of Theorem
5.19.  This work will appear elsewhere.

We next survey our $\chi_n$ determinant identities.  
The derivation of several of these identities requires the
simplification 
$$\sum\limits_{r=1}^{\infty}
{(2r-1)q^{2r-1} \over {1-q^{2r-1}}}=
\sum\limits_{r=1}^{\infty}{rq^{r} 
\over {1+q^{r}}}.\tag 5.180$$
Equation (5.180) is the $x\mapsto \sqrt{q}$ case of the
identity at the top of page 34 of \cite{94}.  Apply the partial
fraction ${q^{r}/{(1+q^{r})}}={q^{r}/{(1-q^{r})}}-
{2q^{2r}/{(1-q^{2r})}}$ termwise to the right-hand-side, and
then cancell the even index terms. This simplification came
up in the context of some Lambert series identities of
Ramanujan.

We also need the observation that the proof of Lemma 5.1
immediately implies that equation (5.3) remains valid if 
Hankel determinants are replaced by $\chi_n$ determinants.
That is, we have the following lemma.
\proclaim{Lemma 5.22} 
Let $v_1,\dots ,v_{2n}$ and $w_1,\dots ,w_{2n}$ be 
indeterminate, and let $n=1,2,3,\cdots$.  Suppose that
$\chi_n (\{w_\nu\})$ and $\chi_n (\{v_\nu\})$ are the
determinants of the $n\times n$ square matrices given by 
\hbox{\rm(3.56)}, with the $n=1$ cases equal to $w_2$
and $v_2$, respectively.  Let $M_{n,S}$ be the $n\times n$
matrix whose $i$-th row is 
$$\spreadlines{8 pt}\allowdisplaybreaks\alignat 2
\kern -4 em
&\kern -1.5 em v_i+w_i,v_{i+1}+w_{i+1},\cdots,
v_{i+n-2}+w_{i+n-2},v_{i+n}+w_{i+n},\kern 3 em
&\kern -4 em\text{if}&\quad \ i\in S,\cr  
\kern -6.0 em\text{and}\kern 6.0 em
&\kern -1.5 em v_i,v_{i+1},\cdots,v_{i+n-2},v_{i+n},
&\kern -4 em\text{if}&\quad \ i\notin S. 
\tag 5.181\cr  
\kern -5.50 em \text{Then}\kern 5.50 em &
\kern -2.0 em \chi_n (\{w_\nu\})=
\sum\limits_{\emptyset \subseteq S
\subseteq I_n}\kern-.5em 
(-1)^{n-\Vert S\Vert }\det (M_{n,S})&\tag 5.182\cr
&\kern 1.9 em = (-1)^n\chi_n (\{v_\nu\})+
\sum\limits_{p=1}^n(-1)^{n-p}
\kern-.5em
\sum\limits_{{\emptyset \subset S\subseteq I_n}
\atop {\Vert S\Vert =p}}\kern-.5em \det (M_{n,S}).
&\tag 5.183\cr
\endalignat$$
\endproclaim

We first have the following theorem.
\proclaim{Theorem 5.23 } Let 
$\vartheta_3 (0,-q)$ be determined by 
\hbox{\rm(1.1)}, and let $n=1,2,3,\cdots$.   
We then have
$$\spreadlines{6 pt}\allowdisplaybreaks\align 
&\vartheta_3(0,-q)^{4n^2}
\left[1+24\sum\limits_{r=1}^{\infty}
{rq^{r} \over {1+q^{r}}}\right]\cr
=\,&\left\{(-1)^{n-1}2^{2n^2+n+1}\tfrac{3}{n(4n^2-1)}
\prod\limits_{r=1}^{2n-1}(r!)^{-1}\right\}\cdot 
\chi_n(\{g_\nu\}),\tag 5.184\cr
\endalign$$
where $\chi_n(\{g_\nu\})$ is the determinant of the 
$n\times n$ matrix whose $i$-th row is  
$$\spreadlines{6 pt}\allowdisplaybreaks\align 
&\ g_i,g_{i+1},\cdots,g_{i+n-2},g_{i+n},
\kern 4 em \text{ for}\quad \ i=1,2,\cdots,n,\tag 5.185\cr
\kern -7.45 em\text{where}\kern 7.45 em 
&\ g_i:= U_{2i-1}-c_i,\tag 5.186\cr
\endalign$$
when $n\geq 2$.  If $n=1$, then $\chi_1(\{g_\nu\})$  
equals $g_2=U_3-c_2$.  
The $U_{2i-1}$ and $c_i$ are defined by 
\hbox{\rm(5.22)} and \hbox{\rm(5.23)}, respectively, 
with $B_{2i}$ the Bernoulli numbers defined by 
\hbox{\rm(2.60)}.
\endproclaim
\demo{Proof} Our analysis deals with 
$\sc(u,k)\dn(u,k)$.  Starting with (2.80), applying row
operations and (3.68), appealing to (4.52), and then
utilizing both (5.9) and (5.15) we have the following
computation.   
$$\spreadlines{10 pt}\allowdisplaybreaks\align
& \chi_n(\{U_{2\nu -1}(-q)
-(-1)^{\nu-1}{(2^{2\nu}-1)
\over {4\nu}}\cdot |B_{2\nu}|\})\tag 5.187\cr
& \kern -2 em  = \chi_n(\{(-1)^{\nu}{z^{2\nu}
\over 2^{2\nu + 1}}\cdot 
(sd/c)_{\nu}(k^2)\})\tag 5.188\cr
& \kern -2 em =  \left\{(-1)^{n-1}2^{-(2n^2+n+1)}
\tfrac{n(4n^2-1)}{3}(1-2k^2)
\prod\limits_{r=1}^{2n-1}r!\right\}
\cdot z^{2n^2+2}\tag 5.189\cr 
& \kern -2 em  =   \left\{(-1)^{n-1}2^{-(2n^2+n+1)}
\tfrac{n(4n^2-1)}{3}
\prod\limits_{r=1}^{2n-1}r!\right\}\cr
& 
\cdot \vartheta_3 (0,q)^{4n^2}
\left[1-24\sum\limits_{r=1}^{\infty}
{(2r-1)q^{2r-1} \over {1+q^{2r-1}}}\right]
.\tag 5.190\cr
\endalign$$ 
Replacing $q$ by $-q$ in (5.187) and (5.190) and then
using (5.180) gives 
$$\spreadlines{10 pt}\allowdisplaybreaks\align
&\chi_n(\{U_{2\nu -1}(q)
-(-1)^{\nu-1}{(2^{2\nu}-1)
\over {4\nu}}\cdot |B_{2\nu}|\}) \tag 5.191a\cr
& \kern -2 em  =   \left\{(-1)^{n-1}2^{-(2n^2+n+1)}
\tfrac{n(4n^2-1)}{3}
\prod\limits_{r=1}^{2n-1}r!\right\}\tag 5.191b\cr
& 
\cdot \vartheta_3 (0,-q)^{4n^2}
\left[1+24\sum\limits_{r=1}^{\infty}
{rq^{r} \over {1+q^{r}}}\right].\tag 5.191c\cr
\endalign$$ 
Solving for (5.191c) in (5.191) yields (5.184).
\qed\enddemo

Lemma 5.22 and Theorem 5.23 lead to the $\chi_n$
determinant sum identity in the following theorem.
\proclaim{Theorem 5.24} 
Let $n=1,2,3,\cdots$.  Then 
$$\spreadlines{6 pt}\allowdisplaybreaks\align 
&\vartheta_3(0,-q)^{4n^2}
\left[1+24\sum\limits_{r=1}^{\infty}
{rq^{r} \over {1+q^{r}}}\right]\cr
=\,&1+\sum\limits_{p=1}^n
(-1)^{p-1}2^{2n^2+n+1}\tfrac{3}{n(4n^2-1)}
\prod\limits_{r=1}^{2n-1}(r!)^{-1}
\kern-.5em
\sum\limits_{{\emptyset \subset S\subseteq I_n}\atop 
{\Vert S\Vert =p}}\kern-.5em \det (M_{n,S}),
\tag 5.192\cr
\endalign$$
where $\vartheta_3(0,-q)$ is determined by 
\hbox{\rm(1.1)}, 
and $M_{n,S}$ is the $n\times n$ matrix whose $i$-th row
is  
$$\spreadlines{6 pt}\allowdisplaybreaks\alignat 2
\kern -4 em
&U_{2i-1},U_{2(i+1)-1},\cdots,
U_{2(i+n-2)-1},U_{2(i+n)-1},
&\kern 4 em\text{if}&\quad \ i\in S,\cr  
\kern -6.35 em\text{and}\kern 6.35 em
&c_i,c_{i+1},\cdots,c_{i+n-2},c_{i+n},
&\kern 4 em\text{if}&\quad \ i\notin S, 
\tag 5.193\cr  
\endalignat$$
when $n\geq 2$.  If $n=1$, then $M_{n,S}$ is 
the $1\times 1$ matrix
$$(U_3),\quad\text{ since} \quad S=\{1\} \quad 
\text{ and}\quad
1\in S .\tag 5.194$$ 
The $U_{2i-1}$ and $c_i$ are defined by 
\hbox{\rm(5.22)} and \hbox{\rm(5.23)}, respectively, 
with $B_{2i}$ the Bernoulli numbers defined by 
\hbox{\rm(2.60)}.
\endproclaim
\demo{Proof} Specialize $v_i$ and $w_i$ in equation
(5.183) of Lemma 5.22 as in (5.31) and (5.32),
respectively. Recall (5.33).

Equating (5.188) and (5.190) it is immediate that 
$$\spreadlines{10 pt}\allowdisplaybreaks\align
\chi_n (\{w_\nu\})=   &\left\{(-1)^{n-1}2^{-(2n^2+n+1)}
\tfrac{n(4n^2-1)}{3}
\prod\limits_{r=1}^{2n-1}r!\right\}\cr
& 
\cdot \vartheta_3 (0,q)^{4n^2}
\left[1-24\sum\limits_{r=1}^{\infty}
{(2r-1)q^{2r-1} \over {1+q^{2r-1}}}\right]
.\tag 5.195\cr
\endalign$$ 
From (4.60) we have 
$$\chi_n (\{v_\nu\})= -\tfrac{n(4n^2-1)}{3}
2^{-(2n^2 + n +1)}
\prod\limits_{r=1}^{2n-1} r!.\tag 5.196$$

Equation (5.192) is now a direct consequence of applying 
the determinant evaluations in (5.195) and (5.196),
replacing $q$ by $-q$, using (5.180), and then multiplying
both sides of the resulting transformation of (5.183) by 
$$(-1)^{n-1}2^{2n^2+n+1}\tfrac{3}{n(4n^2-1)}
\prod\limits_{r=1}^{2n-1}(r!)^{-1},\tag 5.197$$
and simplifying.
\qed\enddemo

The process of ``obtaining a formula by duplication'' from 
\cite{21, pp. 125} applied to the $n=1$ case of Theorem
5.24  gives the identity in 
\cite{258, Eqn. (4) of Table 1(x), pp. 201}.  The $n=2$
case of Theorem 5.24 immediately leads to 
$$\spreadlines{6 pt}\allowdisplaybreaks\align
\vartheta_3(0,-q)^{16}
\left[1+24\sum\limits_{r=1}^{\infty}
{rq^{r} \over {1+q^{r}}}\right]=
1&-\tfrac{8}{15}
\left[17U_1+4U_3-2U_5-4U_7\right]\cr
&-\tfrac{256}{15}
\left[U_1U_7-U_3U_5\right],\tag 5.198\cr
\endalign$$
where $U_{2i-1}$ is defined by (5.22), and 
$\vartheta_3(0,-q)$ is determined by (1.1).

We next have the following theorem.
\proclaim{Theorem 5.25 } Let 
$\vartheta_3 (0,-q)$ be determined by 
\hbox{\rm(1.1)}, and let $n=1,2,3,\cdots$.   
We then have
$$\spreadlines{6 pt}\allowdisplaybreaks\align 
&\vartheta_3(0,-q)^{4n(n+1)}
\left[1+24\sum\limits_{r=1}^{\infty}
{rq^{r} \over {1+q^{r}}}\right]\cr
=\,&\left\{(-1)2^{2n^2+3n}\tfrac{3}{n(n+1)(2n+1)}
\prod\limits_{r=1}^{2n}(r!)^{-1}\right\}\cdot 
\chi_n(\{g_\nu\}),\tag 5.199\cr
\endalign$$
where $\chi_n(\{g_\nu\})$ is the determinant of the 
$n\times n$ matrix whose $i$-th row is  
$$\spreadlines{6 pt}\allowdisplaybreaks\align 
&\ g_i,g_{i+1},\cdots,g_{i+n-2},g_{i+n},
\kern 4 em \text{ for}\quad \ i=1,2,\cdots,n,\tag 5.200\cr
\kern -7.45 em\text{where}\kern 7.45 em 
&\ g_i:= G_{2i+1}-a_i,\tag 5.201\cr
\endalign$$
when $n\geq 2$.  If $n=1$, then $\chi_1(\{g_\nu\})$  
equals $g_2=G_5-a_2$.  
The $G_{2i+1}$ and $a_i$ are defined by 
\hbox{\rm(5.39)} and \hbox{\rm(5.40)}, respectively, 
with $B_{2i+2}$ the Bernoulli numbers defined by 
\hbox{\rm(2.60)}.
\endproclaim
\demo{Proof} Our analysis deals with 
$\sc^2(u,k)\dn^2(u,k)$.  Starting with (2.81), applying row
operations and (3.68), appealing to (4.54), and then
utilizing both (5.9) and (5.15) we have the following
computation.   
$$\spreadlines{10 pt}\allowdisplaybreaks\align
&\chi_n(\{G_{2\nu +1}(-q)
-(-1)^{\nu}{(2^{2\nu + 2}-1)
\over {4(\nu + 1)}}\cdot |B_{2\nu + 2}|\})\tag 5.202\cr
& \kern -2 em  = \chi_n(\{(-1)^{\nu - 1}{z^{2\nu + 2}
\over 2^{2\nu + 3}}\cdot 
(s^2d^2/c^2)_{\nu}(k^2)\})\tag 5.203\cr
& \kern -2 em =  \left\{-2^{-(2n^2+3n)}
\tfrac{n(n+1)(2n+1)}{3}(1-2k^2)
\prod\limits_{r=1}^{2n}r!\right\}
\cdot z^{2n^2+2n+2}\tag 5.204\cr 
& \kern -2 em  =   \left\{-2^{-(2n^2+3n)}
\tfrac{n(n+1)(2n+1)}{3}
\prod\limits_{r=1}^{2n}r!\right\}\cr
& 
\cdot \vartheta_3 (0,q)^{4n(n+1)}
\left[1-24\sum\limits_{r=1}^{\infty}
{(2r-1)q^{2r-1} \over {1+q^{2r-1}}}\right]
.\tag 5.205\cr
\endalign$$ 
Replacing $q$ by $-q$ in (5.202) and (5.205) and then
using (5.180) gives 
$$\spreadlines{10 pt}\allowdisplaybreaks\align
&\chi_n(\{G_{2\nu +1}(q)
-(-1)^{\nu}{(2^{2\nu + 2}-1)
\over {4(\nu + 1)}}\cdot |B_{2\nu + 2}|\})\tag 5.206a\cr
& \kern -2 em  =   \left\{-2^{-(2n^2+3n)}
\tfrac{n(n+1)(2n+1)}{3}
\prod\limits_{r=1}^{2n}r!\right\}\tag 5.206b\cr
& 
\cdot \vartheta_3 (0,-q)^{4n(n+1)}
\left[1+24\sum\limits_{r=1}^{\infty}
{rq^{r} \over {1+q^{r}}}\right].\tag 5.206c\cr
\endalign$$ 
Solving for (5.206c) in (5.206) yields (5.199).
\qed\enddemo

Lemma 5.22 and Theorem 5.25 lead to the $\chi_n$
determinant sum identity in the following theorem.
\proclaim{Theorem 5.26} 
Let $n=1,2,3,\cdots$.  Then 
$$\spreadlines{6 pt}\allowdisplaybreaks\align 
&\vartheta_3(0,-q)^{4n(n+1)}
\left[1+24\sum\limits_{r=1}^{\infty}
{rq^{r} \over {1+q^{r}}}\right]\cr
=\,&1+\sum\limits_{p=1}^n
(-1)^{n-p+1}2^{2n^2+3n}\tfrac{3}{n(n+1)(2n+1)}
\prod\limits_{r=1}^{2n}(r!)^{-1}
\kern-.5em
\sum\limits_{{\emptyset \subset S\subseteq I_n}\atop 
{\Vert S\Vert =p}}\kern-.5em \det (M_{n,S}),
\tag 5.207\cr
\endalign$$
where $\vartheta_3(0,-q)$ is determined by 
\hbox{\rm(1.1)}, 
and $M_{n,S}$ is the $n\times n$ matrix whose $i$-th row
is  
$$\spreadlines{6 pt}\allowdisplaybreaks\alignat 2
\kern -4 em
&G_{2i+1},G_{2(i+1)+1},\cdots,
G_{2(i+n-2)+1},G_{2(i+n)+1},
&\kern 4 em\text{if}&\quad \ i\in S,\cr  
\kern -6.35 em\text{and}\kern 6.35 em
&a_i,a_{i+1},\cdots,a_{i+n-2},a_{i+n},
&\kern 4 em\text{if}&\quad \ i\notin S, 
\tag 5.208\cr  
\endalignat$$
when $n\geq 2$.  If $n=1$, then $M_{n,S}$ is 
the $1\times 1$ matrix
$$(G_5),\quad\text{ since} \quad S=\{1\} \quad 
\text{ and}\quad
1\in S .\tag 5.209$$ 
The $G_{2i+1}$ and $a_i$ are defined by 
\hbox{\rm(5.39)} and \hbox{\rm(5.40)}, respectively, 
with $B_{2i+2}$ the Bernoulli numbers defined by 
\hbox{\rm(2.60)}.
\endproclaim
\demo{Proof} Specialize $v_i$ and $w_i$ in equation
(5.183) of Lemma 5.22 as in (5.48) and (5.49),
respectively. Recall (5.50).

Equating (5.203) and (5.205) it is immediate that 
$$\spreadlines{10 pt}\allowdisplaybreaks\align
\chi_n (\{w_\nu\})=   &\left\{-2^{-(2n^2+3n)}
\tfrac{n(n+1)(2n+1)}{3}
\prod\limits_{r=1}^{2n}r!\right\}\cr
& 
\cdot \vartheta_3 (0,q)^{4n(n+1)}
\left[1-24\sum\limits_{r=1}^{\infty}
{(2r-1)q^{2r-1} \over {1+q^{2r-1}}}\right]
.\tag 5.210\cr
\endalign$$ 
From (4.61) we have 
$$\chi_n (\{v_\nu\})= (-1)^{n-1}
\cdot\tfrac{n(n+1)(2n+1)}{3}2^{-(2n^2 + 3n)}
\prod\limits_{r=1}^{2n} r!.\tag 5.211$$

Equation (5.207) is now a direct consequence of applying 
the determinant evaluations in (5.210) and (5.211),
replacing $q$ by $-q$, using (5.180), and then multiplying
both sides of the resulting transformation of (5.183) by 
$$-2^{2n^2+3n}\tfrac{3}{n(n+1)(2n+1)}
\prod\limits_{r=1}^{2n}(r!)^{-1},\tag 5.212$$
and simplifying.
\qed\enddemo

The process of ``obtaining a formula by duplication'' from 
\cite{21, pp. 125} applied to the $n=1$ case of Theorem
5.26  gives the identity in 
\cite{258, Eqn. (5) of Table 1(ii), pp. 197}.  The $n=2$
case of Theorem 5.26 immediately leads to 
$$\spreadlines{6 pt}\allowdisplaybreaks\align
\vartheta_3(0,-q)^{24}
\left[1+24\sum\limits_{r=1}^{\infty}
{rq^{r} \over {1+q^{r}}}\right]=
1&+\tfrac{8}{45}
\left[124G_3+17G_5-4G_7-2G_9\right]\cr
&-\tfrac{256}{45}
\left[G_3G_9-G_5G_7\right],\tag 5.213\cr
\endalign$$
where $G_{2i+1}$ is defined by (5.39), and 
$\vartheta_3(0,-q)$ is determined by (1.1).

We next have the following theorem.
\proclaim{Theorem 5.27 } Let 
$\vartheta_3 (0,-q)$ be determined by 
\hbox{\rm(1.1)}, and let $n=1,2,3,\cdots$.   
We then have
$$\spreadlines{6 pt}\allowdisplaybreaks\align 
&\vartheta_3(0,q)^{2n(n-1)}\vartheta_3(0,-q)^{2n^2}\cr
&\cdot\left\{2n\left[2+24\sum\limits_{r=1}^{\infty}
{2rq^{2r} \over {1+q^{2r}}}\right]-
\left[1-24\sum\limits_{r=1}^{\infty}
{(2r-1)q^{2r-1} \over {1+q^{2r-1}}}\right]\right\}\cr
=\,&\left\{(-1)^{n-1}2^{2n}\tfrac{3}{n(2n-1)}
\prod\limits_{r=1}^{n-1}(2r)!^{-2}\right\}\cdot 
\chi_n(\{g_\nu\}),\tag 5.214\cr
\endalign$$
where $\chi_n(\{g_\nu\})$ is the determinant of the 
$n\times n$ matrix whose $i$-th row is  
$$\spreadlines{6 pt}\allowdisplaybreaks\align 
&\ g_i,g_{i+1},\cdots,g_{i+n-2},g_{i+n},
\kern 4 em \text{ for}\quad \ i=1,2,\cdots,n,\tag 5.215\cr
\kern -7.45 em\text{where}\kern 7.45 em 
&\ g_i:= R_{2i-2}-b_i,\tag 5.216\cr
\endalign$$
when $n\geq 2$.  If $n=1$, then $\chi_1(\{g_\nu\})$  
equals $g_2=R_2-b_2$.  
The $R_{2i-2}$ and $b_i$ are defined by 
\hbox{\rm(5.56)} and \hbox{\rm(5.57)}, respectively, 
with $E_{2i-2}$ the Euler numbers defined by 
\hbox{\rm(2.61)}.
\endproclaim
\demo{Proof} Our analysis deals with 
$\nc(u,k)$.  Starting with (2.82), applying row
operations and (3.69), appealing to (4.48), and then
utilizing (5.9), (5.10), (5.14), and (5.15) we have the
following computation.   
$$\spreadlines{10 pt}\allowdisplaybreaks\align
&\chi_n(\{R_{2\nu -2}(q)
-(-1)^{\nu -1}\cdot {\tfrac {1}{4}}\cdot
|E_{2\nu -2}|\})\tag 5.217\cr
& \kern -2 em  =\chi_n(\{(-1)^\nu {z^{2\nu -1}
\over 4}\,\sqrt{1-k^2}
\cdot (nc)_{\nu -1}(k^2)\})\tag 5.218\cr
& \kern -2 em =  \left\{(-1)^{n-1}2^{-2n}
\tfrac{n(2n-1)}{3}
\prod\limits_{r=1}^{n-1}(2r)!^2\right\}\cr 
& 
\cdot(1-k^2)^{n^2/2}
\left[2n(2-k^2)-(1-2k^2)\right] 
z^{2n^2-n+2}\tag 5.219\cr
& \kern -2 em  =  \left\{(-1)^{n-1}2^{-2n}
\tfrac{n(2n-1)}{3}
\prod\limits_{r=1}^{n-1}(2r)!^2\right\}\tag 5.220a\cr
& 
\cdot \vartheta_3(0,q)^{2n(n-1)}
\vartheta_3(0,-q)^{2n^2}
\tag 5.220b\cr
& 
\cdot\left\{2n\left[2+24\sum\limits_{r=1}^{\infty}
{2rq^{2r} \over {1+q^{2r}}}\right]-
\left[1-24\sum\limits_{r=1}^{\infty}
{(2r-1)q^{2r-1} \over {1+q^{2r-1}}}\right]\right\}
.\tag 5.220c\cr
\endalign$$

Equating (5.217) with (5.220) and then solving for 
(5.220b)--(5.220c) yields (5.214).
\qed\enddemo

Noting that 
$$\left[2n(2-k^2)-(1-2k^2)\right] = 
\left[(2n+1)+2(n-1)(1-k^2)\right],\tag 5.221$$
and recalling (5.9) and (5.10), it is sometimes useful to
rewrite (5.220c) as 
$$\left[(2n+1)\vartheta_3(0,q)^4+
2(n-1)\vartheta_3(0,-q)^4\right].\tag 5.222$$

Lemma 5.22 and Theorem 5.27 lead to the $\chi_n$
determinant sum identity in the following theorem.
\proclaim{Theorem 5.28} 
Let $n=1,2,3,\cdots$.  Then 
$$\spreadlines{6 pt}\allowdisplaybreaks\align 
&\vartheta_3(0,q)^{2n(n-1)}\vartheta_3(0,-q)^{2n^2}\cr
&\cdot\left\{2n\left[2+24\sum\limits_{r=1}^{\infty}
{2rq^{2r} \over {1+q^{2r}}}\right]-
\left[1-24\sum\limits_{r=1}^{\infty}
{(2r-1)q^{2r-1} \over {1+q^{2r-1}}}\right]\right\}\cr
=\,&(4n-1)+\sum\limits_{p=1}^n
(-1)^{p-1}2^{2n}\tfrac{3}{n(2n-1)}
\prod\limits_{r=1}^{n-1}(2r)!^{-2}
\kern-.5em
\sum\limits_{{\emptyset \subset S\subseteq I_n}\atop 
{\Vert S\Vert =p}}\kern-.5em \det (M_{n,S}),
\tag 5.223\cr
\endalign$$
where $\vartheta_3(0,-q)$ is determined by 
\hbox{\rm(1.1)}, 
and $M_{n,S}$ is the $n\times n$ matrix whose $i$-th row
is  
$$\spreadlines{6 pt}\allowdisplaybreaks\alignat 2
\kern -4 em
&R_{2i-2},R_{2(i+1)-2},\cdots,
R_{2(i+n-2)-2},R_{2(i+n)-2},
&\kern 4 em\text{if}&\quad \ i\in S,\cr  
\kern -6.35 em\text{and}\kern 6.35 em
&b_i,b_{i+1},\cdots,b_{i+n-2},b_{i+n},
&\kern 4 em\text{if}&\quad \ i\notin S, 
\tag 5.224\cr  
\endalignat$$
when $n\geq 2$.  If $n=1$, then $M_{n,S}$ is 
the $1\times 1$ matrix
$$(R_2),\quad\text{ since} \quad S=\{1\} \quad 
\text{ and}\quad
1\in S .\tag 5.225$$ 
The $R_{2i-2}$ and $b_i$ are defined by 
\hbox{\rm(5.56)} and \hbox{\rm(5.57)}, respectively, 
with $E_{2i-2}$ the Euler numbers defined by 
\hbox{\rm(2.61)}.
\endproclaim
\demo{Proof} Specialize $v_i$ and $w_i$ in equation
(5.183) of Lemma 5.22 as in (5.68) and (5.69),
respectively. Recall (5.70).

Equating (5.218) and (5.220) it is immediate that 
$$\spreadlines{10 pt}\allowdisplaybreaks\align
\kern -4.5 em \chi_n (\{w_\nu\})= &
\left\{(-1)^{n-1}2^{-2n}
\tfrac{n(2n-1)}{3}
\prod\limits_{r=1}^{n-1}(2r)!^2\right\}\cr
& 
\kern 1 em 
\cdot \vartheta_3(0,q)^{2n(n-1)}
\vartheta_3(0,-q)^{2n^2}\cr
& 
\kern 1 em 
\cdot\left\{2n\left[2+24\sum\limits_{r=1}^{\infty}
{2rq^{2r} \over {1+q^{2r}}}\right]-
\left[1-24\sum\limits_{r=1}^{\infty}
{(2r-1)q^{2r-1} \over {1+q^{2r-1}}}\right]\right\}
.\tag 5.226\cr
\endalign$$

From (4.62) we have 
$$\chi_n (\{v_\nu\})=-\tfrac{n(2n-1)(4n-1)}{3}2^{-2n}
\prod\limits_{r=1}^{n-1}(2r)!^2.\tag 5.227$$

Equation (5.223) is now a direct consequence of applying 
the determinant evaluations in (5.226) and (5.227), and
then multiplying both sides of the resulting transformation of
(5.183) by 
$$(-1)^{n-1}2^{2n}
\tfrac{3}{n(2n-1)}
\prod\limits_{r=1}^{n-1}(2r)!^{-2},\tag 5.228$$
and simplifying.
\qed\enddemo

The $n=1$ case of Theorem 5.28, rewritten using 
(5.222), gives the identity in 
\cite{258, Eqn. (4) of Table 1(xiv), pp. 203}.  
Similarly, the $n=2$ case of Theorem 5.28 immediately
leads to 
$$\spreadlines{6 pt}\allowdisplaybreaks\align
&\vartheta_3(0,q)^4\vartheta_3(0,-q)^8
\left[5\vartheta_3(0,q)^4+2\vartheta_3(0,-q)^4\right]
\cr 
=\,&7-\tfrac{1}{2}\left[61R_0+5R_2-R_4-R_6\right]-
2\left[R_0R_6-R_2R_4\right],\tag 5.229\cr
\endalign$$
where $R_{2i-2}$ is defined by (5.56), and 
$\vartheta_3(0,-q)$ is determined by (1.1).

We next have the following theorem.
\proclaim{Theorem 5.29 } Let 
$\vartheta_3 (0,-q)$ be determined by 
\hbox{\rm(1.1)}, and let $n=1,2,3,\cdots$.   
We then have
$$\spreadlines{6 pt}\allowdisplaybreaks\align 
&\vartheta_3(0,q)^{2n(n+1)}\vartheta_3(0,-q)^{2n^2}\cr
&\cdot\left\{2n\left[2+24\sum\limits_{r=1}^{\infty}
{2rq^{2r} \over {1+q^{2r}}}\right]+
\left[1-24\sum\limits_{r=1}^{\infty}
{(2r-1)q^{2r-1} \over {1+q^{2r-1}}}\right]\right\}\cr
=\,&\left\{-2^{2n}\tfrac{3}{n(2n+1)}
\prod\limits_{r=1}^{n}(2r-1)!^{-2}\right\}\cdot 
\chi_n(\{g_\nu\}),\tag 5.230\cr
\endalign$$
where $\chi_n(\{g_\nu\})$ is the determinant of the 
$n\times n$ matrix whose $i$-th row is  
$$\spreadlines{6 pt}\allowdisplaybreaks\align 
&\ g_i,g_{i+1},\cdots,g_{i+n-2},g_{i+n},
\kern 4 em \text{ for}\quad \ i=1,2,\cdots,n,\tag 5.231\cr
\kern -7.45 em\text{where}\kern 7.45 em 
&\ g_i:= R_{2i}-b_{i+1},\tag 5.232\cr
\endalign$$
when $n\geq 2$.  If $n=1$, then $\chi_1(\{g_\nu\})$  
equals $g_2=R_4-b_3$.  
The $R_{2i}$ and $b_{i+1}$ are determined by 
\hbox{\rm(5.56)} and \hbox{\rm(5.57)}, respectively, 
with $E_{2i}$ the Euler numbers defined by 
\hbox{\rm(2.61)}.
\endproclaim
\demo{Proof} Our analysis deals with 
$\sc(u,k)\dc(u,k)$.  Starting with (2.90), applying row
operations and (3.68), appealing to (4.49), and then
utilizing (5.9), (5.10), (5.14), and (5.15) we have the
following computation.   
$$\spreadlines{10 pt}\allowdisplaybreaks\align
&\chi_n(\{R_{2\nu}(q)
-(-1)^{\nu}\cdot {\tfrac {1}{4}}\cdot
|E_{2\nu }|\})\tag 5.233\cr
& \kern -2 em  =\chi_n(\{(-1)^{\nu+1}{z^{2\nu+1}
\over 4}\,\sqrt{1-k^2}
\cdot (sd/c^2)_{\nu}(k^2)\})\tag 5.234\cr
& \kern -2 em =  \left\{-2^{-2n}
\tfrac{n(2n+1)}{3}
\prod\limits_{r=1}^{n}(2r-1)!^2\right\}\cr 
& 
\cdot(1-k^2)^{n^2/2}
\left[2n(2-k^2)+(1-2k^2)\right] 
z^{2n^2+n+2}\tag 5.235\cr
& \kern -2 em  =   \left\{-2^{-2n}
\tfrac{n(2n+1)}{3}
\prod\limits_{r=1}^{n}(2r-1)!^2\right\}\tag 5.236a\cr
& 
\cdot \vartheta_3(0,q)^{2n(n+1)}
\vartheta_3(0,-q)^{2n^2}
\tag 5.236b\cr
& 
\cdot\left\{2n\left[2+24\sum\limits_{r=1}^{\infty}
{2rq^{2r} \over {1+q^{2r}}}\right]+
\left[1-24\sum\limits_{r=1}^{\infty}
{(2r-1)q^{2r-1} \over {1+q^{2r-1}}}\right]\right\}
.\tag 5.236c\cr
\endalign$$

Equating (5.233) with (5.236) and then solving for 
(5.236b)--(5.236c) yields (5.230).
\qed\enddemo

Noting that 
$$\left[2n(2-k^2)+(1-2k^2)\right] = 
\left[(2n-1)+2(n+1)(1-k^2)\right],\tag 5.237$$
and recalling (5.9) and (5.10), it is sometimes useful to
rewrite (5.236c) as 
$$\left[(2n-1)\vartheta_3(0,q)^4+
2(n+1)\vartheta_3(0,-q)^4\right].\tag 5.238$$

Lemma 5.22 and Theorem 5.29 lead to the $\chi_n$
determinant sum identity in the following theorem.
\proclaim{Theorem 5.30} 
Let $n=1,2,3,\cdots$.  Then 
$$\spreadlines{6 pt}\allowdisplaybreaks\align 
&\vartheta_3(0,q)^{2n(n+1)}\vartheta_3(0,-q)^{2n^2}\cr
&\cdot\left\{2n\left[2+24\sum\limits_{r=1}^{\infty}
{2rq^{2r} \over {1+q^{2r}}}\right]+
\left[1-24\sum\limits_{r=1}^{\infty}
{(2r-1)q^{2r-1} \over {1+q^{2r-1}}}\right]\right\}\cr
=\,&(4n+1)+\sum\limits_{p=1}^n
(-1)^{n-p+1}2^{2n}\tfrac{3}{n(2n+1)}
\prod\limits_{r=1}^{n}(2r-1)!^{-2}
\kern-.5em
\sum\limits_{{\emptyset \subset S\subseteq I_n}\atop 
{\Vert S\Vert =p}}\kern-.5em \det (M_{n,S}),
\tag 5.239\cr
\endalign$$
where $\vartheta_3(0,-q)$ is determined by 
\hbox{\rm(1.1)}, 
and $M_{n,S}$ is the $n\times n$ matrix whose $i$-th row
is  
$$\spreadlines{6 pt}\allowdisplaybreaks\alignat 2
\kern -4 em
&R_{2i},R_{2(i+1)},\cdots,
R_{2(i+n-2)},R_{2(i+n)},
&\kern 4 em\text{if}&\quad \ i\in S,\cr  
\kern -8.35 em\text{and}\kern 8.35 em
&b_{i+1},b_{i+2},\cdots,b_{i+n-1},b_{i+n+1},
&\kern 4 em\text{if}&\quad \ i\notin S, 
\tag 5.240\cr  
\endalignat$$
when $n\geq 2$.  If $n=1$, then $M_{n,S}$ is 
the $1\times 1$ matrix
$$(R_4),\quad\text{ since} \quad S=\{1\} \quad 
\text{ and}\quad
1\in S .\tag 5.241$$ 
The $R_{2i}$ and $b_{i+1}$ are determined by 
\hbox{\rm(5.56)} and \hbox{\rm(5.57)}, respectively, 
with $E_{2i}$ the Euler numbers defined by 
\hbox{\rm(2.61)}.
\endproclaim
\demo{Proof} Specialize $v_i$ and $w_i$ in equation
(5.183) of Lemma 5.22 as in (5.87) and (5.88),
respectively, with the $(nc)_i(k^2)$ in (5.88) replaced by 
$(sd/c^2)_i(k^2)$. Utilize (2.90) to write $v_i+w_i$ as the
Lambert series $R_{2i}(q)$ in (5.89).

Equating (5.234) and (5.236) it is immediate that 
$$\spreadlines{10 pt}\allowdisplaybreaks\align
\kern -4.5 em \chi_n (\{w_\nu\})= &
\left\{-2^{-2n}\tfrac{n(2n+1)}{3}
\prod\limits_{r=1}^{n}(2r-1)!^2\right\}\cr
& 
\kern 1 em 
\cdot \vartheta_3(0,q)^{2n(n+1)}
\vartheta_3(0,-q)^{2n^2}\cr
& 
\kern 1 em 
\cdot\left\{2n\left[2+24\sum\limits_{r=1}^{\infty}
{2rq^{2r} \over {1+q^{2r}}}\right]+
\left[1-24\sum\limits_{r=1}^{\infty}
{(2r-1)q^{2r-1} \over {1+q^{2r-1}}}\right]\right\}
.\tag 5.242\cr
\endalign$$

From (4.63) we have 
$$\chi_n (\{v_\nu\})=
-(-1)^{n}\tfrac{n(2n+1)(4n+1)}{3}2^{-2n}
\prod\limits_{r=1}^{n}(2r-1)!^2.\tag 5.243$$

Equation (5.239) is now a direct consequence of applying 
the determinant evaluations in (5.242) and (5.243), and
then multiplying both sides of the resulting transformation of
(5.183) by 
$$-2^{2n}\tfrac{3}{n(2n+1)}
\prod\limits_{r=1}^{n}(2r-1)!^{-2},\tag 5.244$$
and simplifying.
\qed\enddemo

The $n=1$ case of Theorem 5.30, rewritten using 
(5.238), gives the identity in 
\cite{258, Eqn. (5) of Table 1(xiv), pp. 203}.  
Similarly, the $n=2$ case of Theorem 5.30 immediately
leads to 
$$\spreadlines{6 pt}\allowdisplaybreaks\align
&\vartheta_3(0,q)^{12}\vartheta_3(0,-q)^8
\left[3\vartheta_3(0,q)^4+6\vartheta_3(0,-q)^4\right]
\cr 
=\,&9+\tfrac{1}{30}\left[1385R_2+61R_4-
5R_6-R_8\right]-
\tfrac{2}{15}\left[R_2R_8-R_4R_6\right],\tag 5.245\cr
\endalign$$
where $R_{2i}$ and $\vartheta_3(0,-q)$ are determined by
(5.56) and (1.1), respectively.

The $1\times 1$ matrices $(U_3)$, $(G_5)$, $(R_2)$, and
$(R_4)$, in (5.194), (5.209), (5.225), and (5.241),
respectively, were also verified directly by first computing
the $n=1$ cases of the right-hand-sides of (5.184), 
(5.199), (5.214), and (5.230), and then comparing
termwise with the $n=1$ cases of the right-hand-sides of
(5.192),  (5.207), (5.223), and (5.239), respectively.

The rest of our $\chi_n$ determinant identities involve the
classical theta function $\vartheta_2(0,q)$ defined by
(1.2).  The entries in the $\chi_n(\{g_\nu\})$ determinant
are Lambert series.  Motivated by Theorems 5.11, 5.13,
5.14, 5.16, and 5.17 we exhibit these entries explicitly.  We
utilize the notation $\det (M_{n})$ and list the Lambert
series that appear in the $i$-th row of 
$\chi_n(\{g_\nu\})$.

We first have the following theorem. 
\proclaim{Theorem 5.31}  Let 
$\vartheta_2 (0,q)$ be defined by 
\hbox{\rm(1.2)}, and let $n=1,2,3,\cdots$.   
We then have
$$\spreadlines{6 pt}\allowdisplaybreaks\align 
&\vartheta_2 (0,q)^{4n^2}
\left[1+24\sum\limits_{r=1}^{\infty}
{rq^{2r} \over {1+q^{2r}}}\right]\cr
=\,&\left\{4^{n(n+1)}\tfrac{3}{n(4n^2-1)}
\prod\limits_{r=1}^{2n-1}(r!)^{-1}\right\}\cdot 
\det (M_{n}),\tag 5.246\cr
\endalign$$
where $M_{n}$ is the $n\times n$ matrix whose $i$-th
row is 
$$C_{2i-1},C_{2(i+1)-1},\cdots,
C_{2(i+n-2)-1},C_{2(i+n)-1},\quad 
\text{ for}\quad i=1,2,\cdots,n,\tag 5.247$$
when $n\geq 2$.  If $n=1$, then $M_{n}$ is 
the $1\times 1$ matrix $(C_3)$.  The $C_{2i-1}$ 
are defined by \hbox{\rm(5.95)}.
\endproclaim 
\demo{Proof} Our analysis deals with 
$\sd(u,k)\cn(u,k)$.  Starting with (2.83), applying row
operations and (3.68), appealing to (4.44), and then
utilizing (5.11) and (5.14) we have the following
computation.   
$$\spreadlines{10 pt}\allowdisplaybreaks\align
&\chi_n(\{C_{2\nu-1}(q)\})\tag 5.248\cr
& \kern -2 em  =\chi_n(\{(-1)^{\nu-1}{z^{2\nu}k^2 
\over{2^{2\nu+2}}}
\cdot (sc/d)_{\nu}(k^2)\})\tag 5.249\cr
& \kern -2 em =  \left\{4^{-n(n+1)}
\tfrac{n(4n^2-1)}{3}
\prod\limits_{r=1}^{2n-1}r!\right\}\cr 
& 
\cdot\tfrac{1}{2}k^{2n^2}(2-k^2)
z^{2n^2+2}\tag 5.250\cr
& \kern -2 em  =  \left\{4^{-n(n+1)}
\tfrac{n(4n^2-1)}{3}
\prod\limits_{r=1}^{2n-1}r!\right\}\tag 5.251a\cr
& 
\cdot \vartheta_2 (0,q)^{4n^2}
\left[1+24\sum\limits_{r=1}^{\infty}
{rq^{2r} \over {1+q^{2r}}}\right].
\tag 5.251b\cr
\endalign$$

Equating (5.248) with (5.251) and then solving for 
(5.251b) yields (5.246).
\qed\enddemo

The $n=1$ case of Theorem 5.31 is equivalent to the
identity in \cite{258, Eqn. (4) of Table 1(vii), pp. 199}.  

We next have the theorem. 
\proclaim{Theorem 5.32 }  Let 
$\vartheta_2 (0,q)$ be defined by 
\hbox{\rm(1.2)}, and let $n=1,2,3,\cdots$.   
We then have
$$\spreadlines{6 pt}\allowdisplaybreaks\align 
&\vartheta_2 (0,q^{1/2})^{4n(n+1)}
\left[1+24\sum\limits_{r=1}^{\infty}
{rq^{r} \over {1+q^{r}}}\right]\cr
=\,&\left\{2^{n(4n+5)}\tfrac{6}{n(n+1)(2n+1)}
\prod\limits_{r=1}^{2n}(r!)^{-1}\right\}\cdot 
\det (M_{n}),\tag 5.252\cr
\endalign$$
where $M_{n}$ is the $n\times n$ matrix whose $i$-th
row is 
$$D_{2i+1},D_{2(i+1)+1},\cdots,
D_{2(i+n-2)+1},D_{2(i+n)+1},\quad 
\text{ for}\quad i=1,2,\cdots,n,\tag 5.253$$
when $n\geq 2$.  If $n=1$, then $M_{n}$ is 
the $1\times 1$ matrix $(D_5)$.  The $D_{2i+1}$ 
are defined by \hbox{\rm(5.96)}.
\endproclaim 
\demo{Proof} Our analysis deals with 
$\sn^2(u,k)$.  Starting with (2.84), applying row
operations and (3.68), appealing to (4.42), and then
utilizing (5.12) and (5.13) we have the following
computation.   
$$\spreadlines{10 pt}\allowdisplaybreaks\align
&\chi_n(\{D_{2\nu+1}(q)\})\tag 5.254\cr
& \kern -2 em  =\chi_n(\{(-1)^{\nu-1}{z^{2\nu+2}k^2 
\over{2^{2\nu+3}}}
\cdot (sn^2)_{\nu}(k^2)\})\tag 5.255\cr
& \kern -2 em =  \left\{2^{-n(4n+5)}
\tfrac{n(n+1)(2n+1)}{6}
\prod\limits_{r=1}^{2n}r!\right\}\cr 
& 
\cdot 4^{n(n+1)}k^{n(n+1)}(1+k^2)
z^{2n^2+2n+2}\tag 5.256\cr
& \kern -2 em  =  \left\{2^{-n(4n+5)}
\tfrac{n(n+1)(2n+1)}{6}
\prod\limits_{r=1}^{2n}r!\right\}\tag 5.257a\cr
& 
\cdot \vartheta_2 (0,q^{1/2})^{4n(n+1)}
\left[1+24\sum\limits_{r=1}^{\infty}
{rq^{r} \over {1+q^{r}}}\right].
\tag 5.257b\cr
\endalign$$

Equating (5.254) with (5.257) and then solving for 
(5.257b) yields (5.252).
\qed\enddemo

The $n=1$ case of Theorem 5.32 is equivalent to the
identity in \cite{258, Eqn. (5) of Table 1(iii), pp. 198}.  

Just as in (5.106), it is sometimes useful to rewrite (5.252)
by utilizing (5.105). In addition, the relation (5.18) applied
to the left-hand sides of (5.246) and (5.252) yields
identities analogous to those in Corollary 5.12.

We next consider the $\chi_n$ determinant identities
related to $\cn(u,k)$, $\sn(u,k)\dn(u,k)$, $\dn(u,k)$, and 
$\sn(u,k)\cn(u,k)$.  

We first have the following theorem.
\proclaim{Theorem 5.33 }  Let 
$\vartheta_2 (0,q)$ and $\vartheta_3 (0,q)$ be defined by 
\hbox{\rm(1.2)} and \hbox{\rm(1.1)}, respectively. 
Let $n=1,2,3,\cdots$.    We then have
$$\spreadlines{6 pt}\allowdisplaybreaks\align 
&\vartheta_2(0,q)^{2n^2}
\vartheta_3(0,q)^{2n(n-1)}\cr
\ \ &\cdot\left\{2n\left[1+24\sum\limits_{r=1}^{\infty}
{rq^{r} \over {1+q^{r}}}\right]+
\left[1-24\sum\limits_{r=1}^{\infty}
{(2r-1)q^{2r-1} \over {1+q^{2r-1}}}\right]
\right\}\cr
=\,&\left\{4^{n}\tfrac{3}{n(2n-1)}
\prod\limits_{r=1}^{n-1}(2r)!^{-2}\right\}\cdot 
\det (M_{n}),\tag 5.258\cr
\endalign$$
where $M_{n}$ is the $n\times n$ matrix whose $i$-th
row is 
$$T_{2i-2},T_{2(i+1)-2},\cdots,
T_{2(i+n-2)-2},T_{2(i+n)-2},\quad 
\text{ for}\quad i=1,2,\cdots,n,\tag 5.259$$
when $n\geq 2$.  If $n=1$, then $M_{n}$ is 
the $1\times 1$ matrix $(T_2)$.  The $T_{2i-2}$ 
are defined by \hbox{\rm(5.110)}.
\endproclaim 
\demo{Proof} Our analysis deals with 
$\cn(u,k)$.  Starting with (2.85), applying row
operations and (3.69), appealing to (4.33), and then
utilizing (5.9), (5.11), (5.13), and (5.15) we have the
following computation.   
$$\spreadlines{10 pt}\allowdisplaybreaks\align
&\chi_n(\{T_{2\nu -2}(q)\})\tag 5.260\cr
& \kern -2 em  =\chi_n(\{(-1)^{\nu-1}
{z^{2\nu-1}k\over 4}
\cdot (cn)_{\nu -1}(k^2)\})\tag 5.261\cr
& \kern -2 em =  \left\{4^{-n}
\tfrac{n(2n-1)}{3}
\prod\limits_{r=1}^{n-1}(2r)!^2\right\}\cr 
& 
\cdot k^{n^2}
\left[2n(1+k^2)+(1-2k^2)\right] 
z^{2n^2-n+2}\tag 5.262\cr
& \kern -2 em  =   \left\{4^{-n}
\tfrac{n(2n-1)}{3}
\prod\limits_{r=1}^{n-1}(2r)!^2\right\}\tag 5.263a\cr
& 
\cdot \vartheta_2(0,q)^{2n^2}
\vartheta_3(0,q)^{2n(n-1)}
\tag 5.263b\cr
& 
\cdot\left\{2n\left[1+24\sum\limits_{r=1}^{\infty}
{rq^{r} \over {1+q^{r}}}\right]+
\left[1-24\sum\limits_{r=1}^{\infty}
{(2r-1)q^{2r-1} \over {1+q^{2r-1}}}\right]
\right\}
.\tag 5.263c\cr
\endalign$$

Equating (5.260) with (5.263) and then solving for 
(5.263b)--(5.263c) yields (5.258).
\qed\enddemo

The $n=1$ case of Theorem 5.33 is equivalent to the
identity in \cite{258, Eqn. (4) of Table 1(xv), pp. 203}.  

Applying the relation (5.105) to the left-hand-side of
(5.258) immediately implies that 
$$\spreadlines{6 pt}\allowdisplaybreaks\align 
&\vartheta_2 (0,q)^{2n}
\vartheta_2 (0,q^{1/2})^{4n(n-1)}\cr
\ \ &\cdot\left\{2n\left[1+24\sum\limits_{r=1}^{\infty}
{rq^{r} \over {1+q^{r}}}\right]+
\left[1-24\sum\limits_{r=1}^{\infty}
{(2r-1)q^{2r-1} \over {1+q^{2r-1}}}\right]
\right\}\cr
=\,&\left\{4^{n^2}\tfrac{3}{n(2n-1)}
\prod\limits_{r=1}^{n-1}(2r)!^{-2}\right\}\cdot 
\det (M_{n}),\tag 5.264\cr
\endalign$$
where $M_{n}$ is given by (5.259).  

Noting that 
$$\left[2n(1+k^2)+(1-2k^2)\right] = 
\left[(4n-1)-2(n-1)(1-k^2)\right],\tag 5.265$$
and recalling (5.9) and (5.10), it is sometimes useful to
rewrite (5.263c) as 
$$\left[(4n-1)\vartheta_3(0,q)^4-
2(n-1)\vartheta_3(0,-q)^4\right].\tag 5.266$$

We next have the following theorem.
\proclaim{Theorem 5.34 }  Let 
$\vartheta_2 (0,q)$ and $\vartheta_3 (0,q)$ be defined by 
\hbox{\rm(1.2)} and \hbox{\rm(1.1)}, respectively. 
Let $n=1,2,3,\cdots$.    We then have
$$\spreadlines{6 pt}\allowdisplaybreaks\align 
&\vartheta_2(0,q)^{2n^2}
\vartheta_3(0,q)^{2n(n+1)}\cr
\ \ &\cdot\left\{2n\left[1+24\sum\limits_{r=1}^{\infty}
{rq^{r} \over {1+q^{r}}}\right]-
\left[1-24\sum\limits_{r=1}^{\infty}
{(2r-1)q^{2r-1} \over {1+q^{2r-1}}}\right]
\right\}\cr
=\,&\left\{4^{n}\tfrac{3}{n(2n+1)}
\prod\limits_{r=1}^{n}(2r-1)!^{-2}\right\}\cdot 
\det (M_{n}),\tag 5.267\cr
\endalign$$
where $M_{n}$ is the $n\times n$ matrix whose $i$-th
row is 
$$T_{2i},T_{2(i+1)},\cdots,
T_{2(i+n-2)},T_{2(i+n)},\quad 
\text{ for}\quad i=1,2,\cdots,n,\tag 5.268$$
when $n\geq 2$.  If $n=1$, then $M_{n}$ is 
the $1\times 1$ matrix $(T_4)$.  The $T_{2i}$ 
are determined by \hbox{\rm(5.110)}.
\endproclaim 
\demo{Proof} Our analysis deals with 
$\sn(u,k)\dn(u,k)$.  Starting with (2.88), applying row
operations and (3.68), appealing to (4.36), and then
utilizing (5.9), (5.11), (5.13), and (5.15) we have the
following computation.   
$$\spreadlines{10 pt}\allowdisplaybreaks\align
&\chi_n(\{T_{2\nu}(q)\})\tag 5.269\cr
& \kern -2 em  =\chi_n(\{(-1)^{\nu+1}
{z^{2\nu+1}k\over 4}
\cdot (sd)_{\nu}(k^2)\})\tag 5.270\cr
& \kern -2 em =  \left\{4^{-n}
\tfrac{n(2n+1)}{3}
\prod\limits_{r=1}^{n}(2r-1)!^2\right\}\cr 
& 
\cdot k^{n^2}
\left[2n(1+k^2)-(1-2k^2)\right] 
z^{2n^2+n+2}\tag 5.271\cr
& \kern -2 em  =   \left\{4^{-n}
\tfrac{n(2n+1)}{3}
\prod\limits_{r=1}^{n}(2r-1)!^2\right\}\tag 5.272a\cr
& 
\cdot \vartheta_2(0,q)^{2n^2}
\vartheta_3(0,q)^{2n(n+1)}
\tag 5.272b\cr
& 
\cdot\left\{2n\left[1+24\sum\limits_{r=1}^{\infty}
{rq^{r} \over {1+q^{r}}}\right]-
\left[1-24\sum\limits_{r=1}^{\infty}
{(2r-1)q^{2r-1} \over {1+q^{2r-1}}}\right]
\right\}
.\tag 5.272c\cr
\endalign$$

Equating (5.269) with (5.272) and then solving for 
(5.272b)--(5.272c) yields (5.267).
\qed\enddemo

The $n=1$ case of Theorem 5.34 is equivalent to the
identity in \cite{258, Eqn. (5) of Table 1(xv), pp. 203}.  

Applying the relation (5.105) to the left-hand-side of
(5.267) immediately implies that 
$$\spreadlines{6 pt}\allowdisplaybreaks\align 
&\vartheta_2 (0,q^{1/2})^{4n^2}
\vartheta_3 (0,q)^{2n}\cr
\ \ &\cdot\left\{2n\left[1+24\sum\limits_{r=1}^{\infty}
{rq^{r} \over {1+q^{r}}}\right]-
\left[1-24\sum\limits_{r=1}^{\infty}
{(2r-1)q^{2r-1} \over {1+q^{2r-1}}}\right]
\right\}\cr
=\,&\left\{4^{n(n+1)}\tfrac{3}{n(2n+1)}
\prod\limits_{r=1}^{n}(2r-1)!^{-2}\right\}\cdot 
\det (M_{n}),\tag 5.273\cr
\endalign$$
where $M_{n}$ is given by (5.268).  

Equation (5.273) can be rewritten as in Corollary 5.15 by
utilizing (5.18), with $q\mapsto q^{1/2}$.  

Noting that 
$$\left[2n(1+k^2)-(1-2k^2)\right] = 
\left[(4n+1)-2(n+1)(1-k^2)\right],\tag 5.274$$
and recalling (5.9) and (5.10), it is sometimes useful to
rewrite (5.272c) as 
$$\left[(4n+1)\vartheta_3(0,q)^4-
2(n+1)\vartheta_3(0,-q)^4\right].\tag 5.275$$

We next have the following theorem.
\proclaim{Theorem 5.35 }  Let 
$\vartheta_2 (0,q)$ and 
$\vartheta_3(0,q)$ be defined by 
\hbox{\rm(1.2)} and
\hbox{\rm(1.1)}, respectively.  Let
$n=1,2,3,\cdots$.    We then have
$$\spreadlines{6 pt}\allowdisplaybreaks\align 
&\vartheta_2(0,q)^{2n(n-1)}
\vartheta_3(0,q)^{2n^2}\cr
\ \ &\cdot\left\{n\left[1+24\sum\limits_{r=1}^{\infty}
{rq^{r} \over {1+q^{r}}}\right]-
\left[1+24\sum\limits_{r=1}^{\infty}
{rq^{2r} \over {1+q^{2r}}}\right]
\right\}\cr
=\,&\left\{4^{n^2}\tfrac{3}{2n(2n-1)}
\prod\limits_{r=1}^{n-1}(2r)!^{-2}\right\}
\cdot\det (M_{n})\tag 5.276a\cr
&+\left\{4^{n^2-1}\tfrac{3}{2n(2n-1)}
\prod\limits_{r=1}^{n-1}(2r)!^{-2}\right\}
\cdot\det (\overline M_{n-1}),\tag 5.276b\cr
\endalign$$
where the matrices $M_{n}$ and $\overline M_{n-1}$ 
are defined as follows:  

The $n\times n$ matrix $M_{n}$ has $i$-th row given by 
$$N_{2i-2},N_{2(i+1)-2},\cdots,
N_{2(i+n-2)-2},N_{2(i+n)-2},\quad 
\text{ for}\quad i=1,2,\cdots,n,\tag 5.277$$
when $n\geq 2$.  If $n=1$, then $M_{n}$ is 
the $1\times 1$ matrix $(N_2)$.  The $N_{2i-2}$ 
are defined by \hbox{\rm(5.128)}. 

The $(n-1)\times (n-1)$ matrix $\overline M_{n-1}$ 
has $i$-th row given by 
$$N_{2i+2},N_{2(i+1)+2},\cdots,
N_{2(i+n-3)+2},N_{2(i+n-1)+2},\quad 
\text{ for}\quad i=1,2,\cdots,n-1,\tag 5.278$$
when $n\geq 3$.  If $n=1$, then $\overline M_{n-1}$ is
defined to be $0$.  If $n=2$, then $\overline M_{n-1}$ is 
the $1\times 1$ matrix $(N_6)$.  The $N_{2i-2}$ 
are defined by \hbox{\rm(5.128)}. 
\endproclaim 
\demo{Proof}  Our analysis deals with $\dn(u,k)$. Starting
with (2.86) and (2.87), solving for \newline 
$(-1)^m(z^{2m+1}/2^{2m+2})
\cdot (dn)_{m}(k^2)$ for $m\geq 0$,  applying
row operations, recalling (3.56), and simplifying, leads to the
following identity.
$$\spreadlines{10 pt}\allowdisplaybreaks\align
&\chi_n(\{N_{2\nu-2}(q)\}) + 
{\tfrac{1}{4}}\chi_{n-1}(\{N_{2\nu+2}(q)\}) 
\tag 5.279\cr
& \kern -2 em  = \chi_n(\{(-1)^{\nu-1 }
{z^{2\nu-1}\over 4^{\nu }}
\cdot (dn)_{\nu-1}(k^2)\}).\tag 5.280\cr
\endalign$$ 

Next, applying row operations and (3.69) to (5.280),
appealing to (4.34), and then utilizing (5.9), (5.11), (5.13),
and (5.14)  we obtain
$$\spreadlines{10 pt}\allowdisplaybreaks\align
& \kern -2 em =  \left\{4^{-n^2}
\tfrac{n(2n-1)}{3}
\prod\limits_{r=1}^{n-1}(2r)!^2\right\}\cr 
& 
\cdot k^{n(n-1)}
\left[2n(1+k^2)-(2-k^2)\right] 
z^{2n^2-n+2}\tag 5.281\cr
& \kern -2 em  =  \left\{4^{-n^2}
\tfrac{n(2n-1)}{3}
\prod\limits_{r=1}^{n-1}(2r)!^2\right\}\tag 5.282a\cr
& 
\cdot \vartheta_2(0,q)^{2n(n-1)}
\vartheta_3(0,q)^{2n^2}
\tag 5.282b\cr
& 
\cdot\left\{2n\left[1+24\sum\limits_{r=1}^{\infty}
{rq^{r} \over {1+q^{r}}}\right]-
2\left[1+24\sum\limits_{r=1}^{\infty}
{rq^{2r} \over {1+q^{2r}}}\right]
\right\}
.\tag 5.282c\cr
\endalign$$

Equating (5.279) with (5.282) and then solving for 
(5.282b)--(5.282c) yields (5.276).
\qed\enddemo

The $n=1$ case of Theorem 5.35 is equivalent to the
identity in \cite{258, Eqn. (4) of Table 1(xi), pp. 201}.  
Here, the $0\times 0$ determinant in (5.276b) 
is defined to be $0$.  

Noting that 
$$\left[2n(1+k^2)-(2-k^2)\right] = 
\left[(4n-1)-(2n+1)(1-k^2)\right],\tag 5.283$$
and recalling (5.9) and (5.10), it is sometimes useful to
rewrite (5.282c) as 
$$\left[(4n-1)\vartheta_3(0,q)^4-
(2n+1)\vartheta_3(0,-q)^4\right].\tag 5.284$$

We next have the following theorem.
\proclaim{Theorem 5.36 }  Let 
$\vartheta_2 (0,q)$ and $\vartheta_3 (0,q)$ be defined by 
\hbox{\rm(1.2)} and \hbox{\rm(1.1)}, respectively. 
Let $n=1,2,3,\cdots$.    We then have
$$\spreadlines{6 pt}\allowdisplaybreaks\align 
&\vartheta_2(0,q)^{2n(n+1)}
\vartheta_3(0,q)^{2n^2}\cr
\ \ &\cdot\left\{n\left[1+24\sum\limits_{r=1}^{\infty}
{rq^{r} \over {1+q^{r}}}\right]+
\left[1+24\sum\limits_{r=1}^{\infty}
{rq^{2r} \over {1+q^{2r}}}\right]
\right\}\cr
=\,&\left\{4^{n(n+1)}\tfrac{6}{n(2n+1)}
\prod\limits_{r=1}^{n}(2r-1)!^{-2}\right\}\cdot 
\det (M_{n}),\tag 5.285\cr
\endalign$$
where $M_{n}$ is the $n\times n$ matrix whose $i$-th
row is 
$$N_{2i},N_{2(i+1)},\cdots,
N_{2(i+n-2)},N_{2(i+n)},\quad 
\text{ for}\quad i=1,2,\cdots,n,\tag 5.286$$
when $n\geq 2$.  If $n=1$, then $M_{n}$ is 
the $1\times 1$ matrix $(N_4)$.  The $N_{2i}$ 
are determined by \hbox{\rm(5.128)}.
\endproclaim 
\demo{Proof} Our analysis deals with 
$\sn(u,k)\cn(u,k)$.  Starting with (2.89), applying row
operations and (3.68), appealing to (4.35), and then
utilizing (5.9), (5.11), (5.13), and (5.14) we have the
following computation.   
$$\spreadlines{10 pt}\allowdisplaybreaks\align
&\chi_n(\{N_{2\nu}(q)\})\tag 5.287\cr
& \kern -2 em  =\chi_n(\{(-1)^{\nu+1}
{z^{2\nu+1}k^2\over 2^{2\nu+2}}
\cdot (sc)_{\nu}(k^2)\})\tag 5.288\cr
& \kern -2 em =  \left\{4^{-(n^2+n+1)}
\tfrac{n(2n+1)}{3}
\prod\limits_{r=1}^{n}(2r-1)!^2\right\}\cr 
& 
\cdot k^{n(n+1)}
\left[2n(1+k^2)+(2-k^2)\right] 
z^{2n^2+n+2}\tag 5.289\cr
& \kern -2 em  =   \left\{4^{-(n^2+n+1)}
\tfrac{n(2n+1)}{3}
\prod\limits_{r=1}^{n}(2r-1)!^2\right\}\tag 5.290a\cr
& 
\cdot \vartheta_2(0,q)^{2n(n+1)}
\vartheta_3(0,q)^{2n^2}
\tag 5.290b\cr
& 
\cdot\left\{2n\left[1+24\sum\limits_{r=1}^{\infty}
{rq^{r} \over {1+q^{r}}}\right]+
2\left[1+24\sum\limits_{r=1}^{\infty}
{rq^{2r} \over {1+q^{2r}}}\right]
\right\}
.\tag 5.290c\cr
\endalign$$

Equating (5.287) with (5.290) and then solving for 
(5.290b)--(5.290c) yields (5.285).
\qed\enddemo

The $n=1$ case of Theorem 5.36 is equivalent to the
identity in \cite{258, Eqn. (5) of Table 1(xi), pp. 201}.  

Applying the relation (5.105) to the left-hand-side of
(5.285) immediately implies that 
$$\spreadlines{6 pt}\allowdisplaybreaks\align 
&\vartheta_2 (0,q^{1/2})^{4n^2}
\vartheta_2 (0,q)^{2n}\cr
\ \ &\cdot\left\{n\left[1+24\sum\limits_{r=1}^{\infty}
{rq^{r} \over {1+q^{r}}}\right]+
\left[1+24\sum\limits_{r=1}^{\infty}
{rq^{2r} \over {1+q^{2r}}}\right]
\right\}\cr
=\,&\left\{4^{n(2n+1)}\tfrac{6}{n(2n+1)}
\prod\limits_{r=1}^{n}(2r-1)!^{-2}\right\}\cdot 
\det (M_{n}),\tag 5.291\cr
\endalign$$
where $M_{n}$ is given by (5.286).  

Equation (5.291) can be rewritten as in Corollary 5.15 by
utilizing (5.18), with $q\mapsto q^{1/2}$.  

Noting that 
$$\left[2n(1+k^2)+(2-k^2)\right] = 
\left[(4n+1)-(2n-1)(1-k^2)\right],\tag 5.292$$
and recalling (5.9) and (5.10), it is sometimes useful to
rewrite (5.290c) as 
$$\left[(4n+1)\vartheta_3(0,q)^4-
(2n-1)\vartheta_3(0,-q)^4\right].\tag 5.293$$

By appealing to the infinite product expansions for the
classical theta functions in \cite{249, pp. 472--473} and 
\cite{6, Eqn. (10.7.8), pp. 510} our expansion formulas in
this section for powers of various products of classical theta
functions is transformed into expansions for the
corresponding infinite products.

We close this section with some more detailed observations
about the formulas for 
$r_{16}(n)$ and $r_{24}(n)$ in (1.26) and (1.28), 
respectively.  We also provide more information about the 
analysis of the formulas for $r_{4n^2}(N)$ and 
$r_{4n(n+1)}(N)$ obtained by taking the coefficient of
$q^N$ in Theorems 5.4 and 5.6.  
This analysis is based upon the elementary estimates for
divisor sums in \cite{94, pp. 122--123; pp. 125} and the
estimate given by the following lemma.
\proclaim{Lemma 5.37 } Let $n=1,2,3,\cdots$ and let 
$b_1,b_2,\ldots,b_n$ be fixed nonnegative integers.  
Let $N\ge n$ be an integer.  Then, there exists positive
constants $c_1$ and $c_2$, independent of $N$, such that
$$c_1N^{b_1+\cdots+b_n +n-1}\leq
\sum\limits_{m_1+\cdots+m_n=N\atop m_i\ge 1}
m_1^{b_1}m_2^{b_2}\cdots m_n^{b_n} 
\leq c_2N^{b_1+\cdots+b_n +n-1}.\tag 5.294$$
\endproclaim 
\demo{Proof} Lemma 5.37 follows by induction from the
$n=2$ case, which in turn is a consequence of the
Euler-Maclaurin sum formula \cite{188, Eqn. (7.2), pp. 14},
a simple case of the beta integral, and several applications
of the binomial theorem.
\qed\enddemo

We first study (1.26) and (1.28).  
Equation (1.28b) is the dominate term for 
$r_{24}(n)$, and (1.28a) is the ``remainder term'' of 
lower order of magnitude.  The elementary analysis in 
\cite{94, pp. 122--123} immediately implies that 
$$(2-\zeta (r))n^r\leq (-1)^n\sigma_r^{\dag}(n)\leq 
\zeta (r)n^r,\tag 5.295$$
where $\sigma_r^{\dag}(n)$ is given by (1.14), 
$\zeta (r)$ is the Riemann $\zeta$-function, and 
$r\geq 2$.  Applying (5.295), and then Lemma 5.37 as
needed, to (1.28) it follows from  computer algebra
\cite{250} computations that  
$$0.0120n^{11}<\alpha_{24}(n)<0.0332n^{11},
\tag 5.296$$
$$\kern -7.77 em\text{and}\kern 7.77 em 
3.52587n^7<\beta_{24}(n)<
36.3288n^3+14.7474n^5+3.58524n^7,\tag 5.297$$
for $n\geq 3$, where 
$$r_{24}(n)=\alpha_{24}(n)+\beta_{24}(n),\tag 5.298$$ 
with $\alpha_{24}(n)$ and $\beta_{24}(n)$ given by
(1.28b) and (1.28a), respectively.  The upper bound in 
(5.297) may be replaced by $\beta_{24}(n)<3.59n^7$, 
when $n\geq 56$.  Adding the bounds in (5.296) and
(5.297) gives upper and lower bounds for $r_{24}(n)$ that
are  consistent with \cite{94, Eqn. (9.20), pp. 122}.  More 
computer algebra suggests that $\alpha_{24}(n)$, and
hence $r_{24}(n)$,  is very close to $0.0231n^{11}$.   
Moreover, for $n$ large, it appears that $\beta_{24}(n)$ 
is very close to either $3.527n^{7}$, $3.529n^{7}$, 
$3.555n^{7}$,  $3.557n^{7}$,  $3.583n^{7}$, 
$3.584n^{7}$, or $3.585n^{7}$, depending on the
congruence class of $n$.  
Noting that $\beta_{24}(1)=48$ and
$\beta_{24}(2)=1104$ we have 
$\beta_{24}(n)>0$, for $n\geq 1$.  Similarly, 
$\alpha_{24}(1)=\alpha_{24}(2)=0$, and  
$\alpha_{24}(n)>0$, for $n\geq 3$.  The power series
form of (1.22) given by Corollary 8.2 implies that 
$\alpha_{24}(n)$ and $\beta_{24}(n)$ are integers for 
$n\geq 1$.   Just do a termwise analysis of (8.2b) which
considers the congruence classes$\mod 3$  of $m_1$ and
$m_2$, and a termwise analysis of (8.2a) with the 
congruence classes$\mod 9$  of $m_1$.     

Equation (1.26b) is the dominate term for $r_{16}(n)$, and
(1.26a) is the ``remainder term'' of  lower order of
magnitude, at least for all sufficiently large
$n$.  We suspect this is true for all $n\geq 1$.  The
elementary analysis in \cite{94, pp. 125} involving
$r_{12}(n)$ immediately implies that 
$$(2-\zeta (r))n^r\leq (-1)^{n+1}\sigmat{r}(n)
\leq\zeta (r)n^r,\tag 5.299$$
where $\sigmat{r}(n)$ is given by (1.27), 
$\zeta (r)$ is the Riemann $\zeta$-function, and 
$r\geq 2$.  We also have 
$$(2-\Cal H_n)n\leq (-1)^{n+1}\sigmat{1}(n)
\leq n\Cal H_n,\tag 5.300$$
where $\Cal H_n$ is the harmonic number given by 
$$\Cal H_n:=\sum\limits_{p=1}^{n}\tfrac{1}{p}.
\tag 5.301$$
Maximizing the standard estimate for $\Cal H_n$ in 
\cite{188, Eqn. (16.1), pp. 28} with $m=2$, and then 
applying (5.299) and (5.300) to (1.26a), we obtain by 
computer algebra that 
$$10.2728n^{5}<\beta_{16}(n)<11.0650n^{5},
\tag 5.302$$
for $n\geq 55$, where 
$$r_{16}(n)=\alpha_{16}(n)+\beta_{16}(n),\tag 5.303$$ 
with $\alpha_{16}(n)$ and $\beta_{16}(n)$ given by
(1.26b) and (1.26a), respectively.  The lower bound in 
(5.302) also holds for $n\geq 1$.  From the $k=16$ 
case of \cite{94, Eqn. (9.20), pp. 122} there are positive 
constants $c_1$ and $c_2$ such that 
$$c_1n^7\leq r_{16}(n)\leq c_2n^7,\tag 5.304$$
for $n\geq 1$.  By combining (5.302), (5.303), and
(5.304)  there are positive constants $c_3$ and $c_4$ 
such that 
$$c_3n^7< \alpha_{16}(n)< c_4n^7,\tag 5.305$$
for all sufficiently large $n$, with $\alpha_{16}(n)$ 
given by (1.26b).  More computer algebra suggests that
$\alpha_{16}(n)$, and hence $r_{16}(n)$,  is close to
$1.87n^{7}$.  Moreover, for $n$ large, it appears that
$\beta_{16}(n)$  is very close to either $10.32n^{5}$,
$10.34n^{5}$, $10.36n^{5}$,  $10.38n^{5}$,  
$10.66n^{5}$, $10.67n^{5}$, $10.71n^{5}$,
$11.00n^{5}$, or $11.04n^{5}$, depending on the
congruence class of $n$.  The power series
form of (1.20) given by Corollary 8.1 implies that 
$\alpha_{16}(n)$ and $\beta_{16}(n)$ are integers for 
$n\geq 1$.   Just do a termwise analysis of (8.1a) and
(8.1b) which considers the congruence classes$\mod 3$  of
$m_1$ and $m_2$.       

The analysis of the terms in the formulas for 
$r_{4n^2}(N)$ and $r_{4n(n+1)}(N)$ described just after
the proof of Theorem 5.6 follows from the elementary
estimates for divisor sums in (5.295), (5.299), and (5.300),
and suitable applications of Lemma 5.37.

\head 6. Schur functions and Lambert series\endhead

In this section we first establish a Schur function
expansion of a general $n\times n$ determinant whose
entries are either constants or Lambert series.  This
expansion applied to Section 5 yields the Schur function 
form of our Hankel determinant identities in
Section 7 and also completes our proof of the
Kac--Wakimoto conjectured identities for triangular
numbers in \cite{120, pp. 452}. At the end of
this section we state an analogous expansion (whose
proof is the same) which leads to the Schur function 
form  of our $\chi_n$ determinant identities in
Section 7.

Throughout this section we use the notation
$I_n:=\{1,2,\ldots,n\}$, $\Vert S\Vert$ is the
cardinality of the set $S\subseteq I_n$, and 
$\det (M)$ is the determinant of the $n\times n$ 
matrix $M$.  
We frequently consider the sets $S$ and $T$, with 
$$S:=\{\ell_1<\ell_2<\cdots<\ell_p\}\qquad 
\text{ and}\qquad
S^c:=\{\ell_{p+1}<\cdots<\ell_n\},\tag 6.1$$
$$T:=\{j_1<j_2<\cdots<j_p\}\qquad 
\text{ and}\qquad
T^c:=\{j_{p+1}<\cdots<j_n\},\tag 6.2$$
where $S^c:=I_n-S$ is the compliment of the set $S$, 
and $p=1,2,\cdots,n$.  We also have 
$$\Sigma (S):=\ell_1+\ell_2+\cdots+\ell_p\qquad
\text{ and}\qquad
\Sigma (T):=j_1+j_2+\cdots+j_p.\tag 6.3$$
 
The Lambert series entries in our determinants are of the
form given by the following definition.  
\definition{Definition 6.1  (A general Lambert series)}
Let $A$, $B$, $C$, $D$, $E$, $F$, $G$, and $u$ be
indeterminant, and assume that $0<|q|<1$.  Then, we
define the Lambert series $L_u$ by
$$\spreadlines{6 pt}\allowdisplaybreaks\align 
L_u\equiv \ & L_u(m_1;A,\dots,G)\equiv 
L_u(m_1;A,B,C,D,E,F,G)\cr
:= \ &\sum\limits_{m_1=1}^{\infty}
{D^uE^{m_1}(Bm_1+C)^u\over 
1+Aq^{Bm_1+C}}\cdot q^{Fm_1+G}\tag 6.4\cr
= \ &(-A)^{-1}q^{G-C}\kern-.65em
\sum\limits_{y_1,m_1\geq 1}
\kern -.6 em(-A)^{y_1}E^{m_1}
(D(Bm_1+C))^u\cr
&\kern 8 em\cdot q^{(F-B)m_1}q^{(Bm_1+C)y_1}.
\tag 6.5\cr
\endalign$$
\enddefinition  

We transform (6.4) into (6.5) by expanding 
$1/(1+Aq^{Bm_1+C})$ as a geometric series and then
interchanging summation.  

Our aim is to obtain a Schur function expansion of the 
$n\times n$ determinant 
$\det(M_{n,S,\{b_r\},\{c_r\},\{a_r\}})$ of the matrix
given by the following definition.
\definition{Definition 6.2  (First Lambert series matrix)}
Let $b_1,b_2,\dots,b_n$ and $c_1,c_2,\dots,c_n$ 
be indeterminant, and let $n=1,2,\cdots$.  Furthermore,
assume that $b_1<b_2<\dots<b_n$ and
$c_1<c_2<\dots<c_n$.  Take 
$\{a_r:r=1,2,\cdots\}$ to be an arbitrary sequence.  Let
$S$ and $S^c$ be the subsets of $I_n$ in (6.1), with
$p=1,2,\cdots,n$.  Let $L_u(r;A,\dots,G)$ be the
Lambert series in Definition 6.1.  Then,
$$M_{n,S,\{b_r\},\{c_r\},\{a_r\}}
\equiv M_{n,S,\{b_1,\dots,b_n\},
\{c_1,\dots,c_n\},\{a_r\}}\tag 6.6$$
is defined to be the $n\times n$ matrix whose 
$i$-th row is
$$\spreadlines{6 pt}\allowdisplaybreaks\alignat 2
\kern -4 em
&L_{c_i+b_1}(m_{\mu};A,\dots,G),
L_{c_i+b_2}(m_{\mu};A,\dots,G),\cdots,
L_{c_i+b_n}(m_{\mu};A,\dots,G),
&\kern 1.6 em\text{if}&\quad \ i=\ell_{\mu}\in S,\cr  
\kern -1.12em\text{and}\kern 1.12 em
&a_i,a_{i+1},\cdots,a_{i+n-1},
&\kern 1.6 em\text{if}&\quad \ i\notin S. 
\tag 6.7\cr  
\endalignat$$
\enddefinition 

In order to state our expansion formula we first need a
few more definitions.
\definition{Definition 6.3  (Divisors)}Let 
$b_1<b_2<\dots<b_n$ and $c_1<c_2<\dots<c_n$ 
be integers.  Then define the divisors $d_b$ and $d_c$
by 
$$\spreadlines{6 pt}\allowdisplaybreaks\align 
d_b:= \ &\text{any common divisor of}
\kern 1.2 em\{b_r-b_1|2\leq r\leq n\},\tag 6.8\cr
\kern -8.00em\text{and}\kern 8.00 em
\kern 1.2 em d_c:= \ &\text{any common divisor of}
\kern 1.2 em\{c_s-c_r|1\leq r<s\leq n\}.\tag 6.9\cr
\endalign$$
\enddefinition 

We often use the greatest common divisor.  

\definition{Definition 6.4  (First Laplace expansion
formula determinant)}Let $\{a_r:r=1,2,\cdots\}$ be
an arbitrary sequence.  Let the sets $S^c$ and $T^c$
be given by (6.1) and (6.2), and let
$p=1,2,\cdots,n$.  Then, 
$$\det(D_{n-p,S^c,T^c})\tag 6.10$$ 
is the determinant of the $(n-p) \times (n-p)$ matrix  
$$D_{n-p,S^c,T^c}:= 
\left[a_{(\ell_{p+r}+j_{p+s}-1)}\right]_{
1\leq r,s\leq n-p}.\tag 6.11$$
\enddefinition

\definition{Definition 6.5  (Schur functions)}
Let $\lambda =(\lambda_1,\lambda_2,
\ldots,\lambda_i,\ldots)$
be a partition of nonnegative integers in decreasing 
order, 
$\lambda_1 \geq \lambda_2 \geq\cdots\geq 
\lambda_i \cdots$, such that only finitely many of the
$\lambda_i$ are nonzero.  The length $\ell(\lambda)$ 
is the number of nonzero parts of $\lambda$.  
Given a partition $\lambda
=(\lambda_1,\lambda_2,\ldots,\lambda_p)$ 
of length $\leq p$,
$$s_\lambda(x)\equiv 
s_\lambda(x_1,x_2,\ldots,x_p):=
{\det (x_i^{\lambda_j+p-j})\over \det (x_i^{p-j})}
\tag 6.12$$
is the Schur function \cite{153} corresponding to the
partition $\lambda$.  (Here, $\det (a_{ij})$ denotes the
determinant of a $p\times p$ matrix with $(i,j)$-th
entry $a_{ij}$).  The Schur function $s_\lambda(x)$ is a
symmetric polynomial in $x_1,x_2,\ldots,x_p$ with
nonnegative integer coefficients.  We typically have 
$p=1,2,\cdots,n$.  
\enddefinition

\definition{Definition 6.6  (Expansion formula partitions)}
Let $b_1,b_2,\dots,b_n$ and $c_1,c_2,\dots,c_n$ 
be indeterminant.  Let $n=1,2,\cdots$ and 
$p=1,2,\cdots,n$.  Furthermore,
assume that $b_1<b_2<\dots<b_n$ and
$c_1<c_2<\dots<c_n$.  Let $d_b$ and $d_c$ be
determined by Definition 6.3.  Take
$\{\ell_1,\ell_2,\dots,\ell_p\}$ and 
$\{j_1,j_2,\dots,j_p\}$ 
as in (6.1) and (6.2).  We then define the partitions 
$\lambda =(\lambda_1,\lambda_2,\ldots,\lambda_p)$ 
and $\nu =(\nu_1,\nu_2,\ldots,\nu_p)$ of length 
$\leq p$ as follows:
$$\spreadlines{6 pt}\allowdisplaybreaks\alignat 2
\kern -4 em
&\lambda_i:={\tfrac{1}{d_c}}c_{\ell_{p-i+1}}-
{\tfrac{1}{d_c}}c_{\ell_1}+i-p,
&\kern 4 em\text{for}&\quad \ i=1,2,\cdots,p,
\tag 6.13\cr  
\kern -7.40 em\text{and}\kern 7.40 em
&\nu_i:={\tfrac{1}{d_b}}b_{j_{p-i+1}}-
{\tfrac{1}{d_b}}b_{j_1}+i-p,
&\kern 4 em\text{for}&\quad \ i=1,2,\cdots,p.
\tag 6.14\cr  
\endalignat$$
\enddefinition 

Keeping in mind Definitions 6.1 to 6.6 we now have the
following theorem.
\proclaim{Theorem 6.7 (Expansion of first Lambert 
series determinant)} Let
$L_u(r;A,\dots,G)$ be the Lambert series in \hbox{\rm
Definition 6.1}.  Let the $n\times n$ matrix
$M_{n,S,\{b_r\},\{c_r\},\{a_r\}}$  of Lambert series 
be given by \hbox{\rm Definition 6.2},  and take 
$n=1,2,\cdots$ and $p=1,2,\cdots,n$.   Let $A$, $B$,
$C$, $D$, $E$, $F$, $G$, and 
$b_1,b_2,\dots,b_n$ and $c_1,c_2,\dots,c_n$ be
indeterminant, and assume that $0<|q|<1$.   
Let $S$, $T$, $S^c$, $T^c$, $\Sigma (S)$, and  
$\Sigma (T)$ be given by \hbox{\rm (6.1)--(6.3)}.
Take $\{a_r:r=1,2,\cdots\}$ to be an arbitrary 
sequence.  Let $d_b$ and $d_c$ be given by 
\hbox{\rm Definition 6.3} and $s_\lambda$ and 
$s_\nu$ be the Schur functions in 
\hbox{\rm Definition 6.5} with partitions $\lambda$ and 
$\nu$ as in \hbox{\rm Definition 6.6}.  Finally, take the 
$(n-p) \times (n-p)$ matrix $D_{n-p,S^c,T^c}$ in 
\hbox{\rm Definition 6.4}.  We then have the expansion
formula
$$\spreadlines{6 pt}\allowdisplaybreaks\align 
&\det(M_{n,S,\{b_r\},\{c_r\},\{a_r\}})\cr
 =\, & (-A)^{-p}q^{p(G-C)}\kern-1.95em
\sum\limits_{{y_1,\ldots,y_p\geq 1}\atop 
{m_1>m_2>\cdots>m_p\geq 1}}
\kern-1.65em(-A)^{y_1+\cdots+y_p}
E^{m_1+\cdots+m_p}\cr
&\kern 2 em\cdot q^{(F-B)(m_1+\cdots+m_p)}
q^{(Bm_1+C)y_1+\cdots+(Bm_p+C)y_p}\cr
&\kern 2 em\cdot 
\prod\limits_{1\leq r<s\leq p}\kern-1.1em
\left((D(Bm_r+C))^{d_c}-(D(Bm_s+C))^{d_c}\right)\cr
&\kern 2 em\cdot 
\prod\limits_{1\leq r<s\leq p}\kern-1.1em
\left((D(Bm_r+C))^{d_b}-(D(Bm_s+C))^{d_b}\right)\cr
&\kern 2 em\cdot 
s_\lambda\left((D(Bm_1+C))^{d_c},\ldots,
(D(Bm_p+C))^{d_c}\right)\cr
&\cdot
\sum\limits_{{\emptyset \subset T\subseteq I_n}
\atop {\Vert T\Vert =p}}\kern-.5em 
(-1)^{\Sigma (S)+\Sigma (T)}\cdot 
\det (D_{n-p,S^c,T^c})\cdot
\left(\prod\limits_{r=1}^{p}D(Bm_r+C)
\right)^{\!\!c_{\ell_1}+b_{j_1}}\cr 
&\kern 3 em\cdot
s_\nu\left((D(Bm_1+C))^{d_b},\ldots,
(D(Bm_p+C))^{d_b}\right).\tag 6.15\cr 
\endalign$$
\endproclaim

We give a derivation proof of Theorem 6.7 that relies
upon the Laplace expansion formula for a determinant,
classical properties of Schur functions, symmetry and
skew-symmetry arguments, and row and column
operations.  

We start with the Laplace expansion formula for a
determinant as given in \cite{118, pp. 396-397}.  
\proclaim{Theorem 6.8 (Laplace expansion formula for a
determinant)} Let $\Lambda$ be an $n\times n$ matrix 
with $n=1,2,\cdots$.  Let $S$, $T$, $S^c$, $T^c$, 
$\Sigma (S)$, and $\Sigma (T)$ be given by 
\hbox{\rm (6.1)--(6.3)}, where $I_n:=\{1,2,\ldots,n\}$, 
$S^c:=I_n-S$ is the compliment of the set $S$, 
and $p=1,2,\cdots,n$.  Let $\Lambda_{S,T}$ denote
the minor of $\Lambda$ obtained from the rows and
columns of $S$ and $T$, respectively.  Similarly, let 
$\Lambda_{S^c,T^c}$ be the minor corresponding to 
$S^c$ and $T^c$.  We then expand the determinant of
$\Lambda$ along a given fixed set 
$S=\{\ell_1<\ell_2<\cdots<\ell_p\}$ of rows as follows:
$$\det(\Lambda) = 
\sum\limits_{{T\subseteq I_n}
\atop {\Vert T\Vert =p}}\kern-.35em 
\Lambda_{S,T}\Lambda_{S^c,T^c}
\cdot (-1)^{\Sigma (S)+\Sigma (T)}.\tag 6.16$$
\endproclaim

If we successively expand 
$\det(M_{n,S,\{b_r\},\{c_r\},\{a_r\}})$ along those
rows of the $a_r$'s that correspond to $i\in S^c$, we
obtain a linear combination of $p\times p$ determinants
of the $L_b$'s. To this end, apply Theorem 6.8, where 
$\Lambda_{S,T}$ is a $p\times p$ minor involving
Lambert series, and $\Lambda_{S^c,T^c}$ depends
upon the constants $\{a_r:r=1,2,\cdots\}$.  We have
the following expansion.  
$$\det(M_{n,S,\{b_r\},\{c_r\},\{a_r\}}) = 
\sum\limits_{{T\subseteq I_n}
\atop {\Vert T\Vert =p}}\kern-.35em
\beta_{S,T}\cdot\det\left(L_{(S,T,r,s)}\right)_{  
1\leq r,s\leq p},\tag 6.17$$
where $S\subseteq I_n$ is fixed, $\beta_{S,T}$ are
constants, $(S,T,r,s)$ is an integer that depends on $S$,
$T$, $r$, and $s$, and $L_{(S,T,r,s)}$ is a Lambert
series from Definition 6.1.  

In the analysis that follows, we need an expansion closely
related to (6.17).  If we replace $L_{(S,T,r,s)}$ in the
right-hand-side of (6.17) by the simpler expression 
$(D(Bm_r+C))^{(S,T,r,s)}$, and then reverse the steps
(via Theorem 6.8) that led to (6.17), we immediately
obtain the following lemma.
\proclaim{Lemma 6.9 (Umbral calculus trick)} Let the
constants $\beta_{S,T}$ and integers $(S,T,r,s)$ be
determined by \hbox{\rm (6.17)}.  Let $S$, $\{b_r\}$,
$\{c_r\}$, $\{a_r\}$ satisfy the conditions in    
\hbox{\rm Definition 6.2}.  Take $n=1,2,\cdots$ and
$p=1,2,\cdots,n$.  We then have the expansion 
$$\det(P_{n,S,\{b_r\},\{c_r\},\{a_r\}}) = 
\sum\limits_{{T\subseteq I_n}
\atop {\Vert T\Vert =p}}\kern-.35em
\beta_{S,T}\cdot\det\left(
(D(Bm_r+C))^{(S,T,r,s)}\right)_{  
1\leq r,s\leq p},\tag 6.18$$
$$\kern -9.95 em\text{where}\kern 9.95 em
P_{n,S,\{b_r\},\{c_r\},\{a_r\}}
\equiv P_{n,S,\{b_1,\dots,b_n\},
\{c_1,\dots,c_n\},\{a_r\}}\tag 6.19$$
is defined to be the $n\times n$ matrix whose 
$i$-th row is
$$\spreadlines{6 pt}\allowdisplaybreaks\alignat 2
\kern -4 em
&(D(Bm_{\mu}+C))^{c_i+b_1},
(D(Bm_{\mu}+C))^{c_i+b_2},\cdots,
(D(Bm_{\mu}+C))^{c_i+b_n},
&\kern 1.6 em\text{if}&\quad \ i=\ell_{\mu}\in S,\cr  
\kern -1.91 em\text{and}\kern 1.91 em
&a_i,a_{i+1},\cdots,a_{i+n-1},
&\kern 1.6 em\text{if}&\quad \ i\notin S. 
\tag 6.20\cr  
\endalignat$$
\endproclaim

It will also be useful to consider the determinant of
(6.19) as a function of $\{m_1,m_2,\dots,m_p\}$.  
That is,
$$f(m_1,\dots,m_p):=
\det(P_{n,S,\{b_r\},\{c_r\},\{a_r\}}).\tag 6.21$$  
Then, if $\rho\in\Cal S_p$ is any permutation of 
$\{1,2,\ldots,p\}$, we have from (6.18) that 
$$f(m_{\rho(1)},\dots,m_{\rho(p)}) = 
\sum\limits_{{T\subseteq I_n}
\atop {\Vert T\Vert =p}}\kern-.35em
\beta_{S,T}\cdot\det\left(
(D(Bm_{\rho(r)}+C))^{(S,T,r,s)}\right)_{  
1\leq r,s\leq p}.\tag 6.22$$

As a first step in transforming (6.17) we write the
determinant on the right-hand-side as
$$\det\left(L_{(S,T,r,s)}\right)_{1\leq r,s\leq p} = 
\sum_{\sigma\in \Cal S_p}\text{sign}(\sigma)
\prod\limits_{r=1}^{p}L_{(S,T,r,\sigma(r))},\tag 6.23$$
where we pick a different element from each row.  

We expand each product on the right-hand-side of (6.23)
into a double multiple sum by appealing to the following
lemma.
\proclaim{Lemma 6.10} Let $b_1,b_2,\dots,b_p$ 
be indeterminant, and let $p=1,2,\cdots$.      
Let $L_u(r;A,\dots,G)$ be the Lambert series in
\hbox{\rm Definition 6.1}.   Let $\rho\in\Cal S_p$ be 
any fixed permutation of $\{1,2,\ldots,p\}$.  Then, 
$$\spreadlines{6 pt}\allowdisplaybreaks\align 
\prod\limits_{r=1}^{p}L_{b_r}\equiv\, & 
\prod\limits_{r=1}^{p}
L_{b_r}(m_r;A,B,C,D,E,F,G)\cr 
 =\, & (-A)^{-p}q^{p(G-C)}\kern-.95em
\sum\limits_{{y_1,\ldots,y_p\geq 1}\atop 
{m_1,\ldots,m_p\geq 1}}
\kern-.65em(-A)^{y_1+\cdots+y_p}
E^{m_1+\cdots+m_p}\cr
&\kern 2 em\cdot q^{(F-B)(m_1+\cdots+m_p)}
q^{(Bm_1+C)y_1+\cdots+(Bm_p+C)y_p}\cr
&\kern 2 em\cdot 
\prod\limits_{r=1}^{p}
(D(Bm_{\rho(r)}+C))^{b_r}.\tag 6.24\cr 
\endalign$$
\endproclaim
\demo{Proof}Apply (6.5), with $y_{\rho(r)}$ and 
$m_{\rho(r)}$, to the $r$-th factor in the left-hand-side
of (6.24), interchange summation, note that 
$(Bm_1+C)y_1+\cdots+(Bm_p+C)y_p = 
B(m_1y_1+\cdots+m_py_p)+C(y_1+\cdots+y_p)$, 
and use symmetry.
\qed\enddemo

Applying Lemma 6.10 with the same $\rho\in\Cal S_p$ 
to each term in the right-hand-side of (6.23), and then
interchanging summation, yields
$$\spreadlines{6 pt}\allowdisplaybreaks\align 
\, &\det\left(L_{(S,T,r,s)}\right)_{1\leq r,s\leq p} \cr 
 =\, & (-A)^{-p}q^{p(G-C)}\kern-.95em
\sum\limits_{{y_1,\ldots,y_p\geq 1}\atop 
{m_1,\ldots,m_p\geq 1}}
\kern-.65em(-A)^{y_1+\cdots+y_p}
E^{m_1+\cdots+m_p}\cr
&\kern 2 em\cdot q^{(F-B)(m_1+\cdots+m_p)}
q^{(Bm_1+C)y_1+\cdots+(Bm_p+C)y_p}\cr
&\kern 2 em\cdot 
\sum_{\sigma\in \Cal S_p}\text{sign}(\sigma)
\prod\limits_{r=1}^{p}
(D(Bm_{\rho(r)}+C))^{(S,T,r,\sigma(r))}\cr
 =\, & (-A)^{-p}q^{p(G-C)}\kern-.95em
\sum\limits_{{y_1,\ldots,y_p\geq 1}\atop 
{m_1,\ldots,m_p\geq 1}}
\kern-.65em(-A)^{y_1+\cdots+y_p}
E^{m_1+\cdots+m_p}\cr
&\kern 2 em\cdot q^{(F-B)(m_1+\cdots+m_p)}
q^{(Bm_1+C)y_1+\cdots+(Bm_p+C)y_p}\cr
&\kern 2 em\cdot 
\left\{\left.\det\left((D(Bm_r+C))^{(S,T,r,s)}\right)_{  
1\leq r,s\leq p}\right|_{m_r
\rightarrow m_{\rho(r)}}\right\}
.\tag 6.25\cr 
\endalign$$
 
\remark{Remark} In equation (6.23) we used the
definition of determinant  in which a different column
element is selected from each row.  This simplified the
above computation.  Had we chosen a different row
element from each column, we would have had to apply
the $\rho\circ\sigma$ case of Lemma 6.10 to the 
$\sigma$ term in (6.23).  
\endremark

It is now not difficult to see that (6.17)--(6.22), (6.25),
and an interchange of summation yields the following
lemma.
\proclaim{Lemma 6.11}Let
$M_{n,S,\{b_r\},\{c_r\},\{a_r\}}$ be determined by 
\hbox{\rm Definition 6.2} and 
$P_{n,S,\{b_r\},\{c_r\},\{a_r\}}$ by 
\hbox{\rm (6.20)}.   Let $n=1,2,\cdots$ and
$p=1,2,\cdots,n$.  Let $\rho\in\Cal S_p$ be any fixed
permutation of $\{1,2,\ldots,p\}$.  Then, 
$$\spreadlines{6 pt}\allowdisplaybreaks\align 
\, &\det(M_{n,S,\{b_r\},\{c_r\},\{a_r\}}) \cr 
 =\, & (-A)^{-p}q^{p(G-C)}\kern-.95em
\sum\limits_{{y_1,\ldots,y_p\geq 1}\atop 
{m_1,\ldots,m_p\geq 1}}
\kern-.65em(-A)^{y_1+\cdots+y_p}
E^{m_1+\cdots+m_p}\cr
&\kern 2 em\cdot q^{(F-B)(m_1+\cdots+m_p)}
q^{(Bm_1+C)y_1+\cdots+(Bm_p+C)y_p}\cr
&\kern 2 em\cdot 
\left\{\left.\det(P_{n,S,\{b_r\},\{c_r\},\{a_r\}})
\right|_{m_r \rightarrow m_{\rho(r)}}\right\}
.\tag 6.26\cr 
\endalign$$
\endproclaim

Before symmetrizing (6.26) with respect to  
$\rho\in\Cal S_p$, we transform the determinant 
\newline 
$\det(P_{n,S,\{b_r\},\{c_r\},\{a_r\}})$ by  some row
operations. 

For each $i=\ell_{r}\in S$, we factor 
$(D(Bm_{r}+C))^{c_{\ell_r}+b_1}$ out of the $i$-th
row of the $n\times n$ matrix 
$P_{n,S,\{b_r\},\{c_r\},\{a_r\}}$.  Keeping in mind 
$d_b$ and $d_c$ from Definition 6.3, and the simple
relation \newline  
$c_{\ell_r} = c_{\ell_1}+d_{c}(c_{\ell_r} - 
c_{\ell_1})/d_c$, we write 
$\det(P_{n,S,\{b_r\},\{c_r\},\{a_r\}})$ as follows:
$$\spreadlines{6 pt}\allowdisplaybreaks\align 
\, &\det(P_{n,S,\{b_r\},\{c_r\},\{a_r\}})\cr 
 =\, &\prod\limits_{r=1}^{p}
(D(Bm_{r}+C))^{c_{\ell_1}+b_1}\cdot
\prod\limits_{r=2}^{p}
\left((D(Bm_{r}+C))^{d_c}\right)^{(c_{\ell_r} - 
c_{\ell_1})/d_c}\cr
&\kern 2 em\cdot 
\det(Q_{n,S,\{b_r\},\{c_r\},\{a_r\}})
,\tag 6.27\cr 
\endalign$$
$$\kern -9.92 em\text{where}\kern 9.92 em
Q_{n,S,\{b_r\},\{c_r\},\{a_r\}}
\equiv Q_{n,S,\{b_1,\dots,b_n\},
\{c_1,\dots,c_n\},\{a_r\}}\tag 6.28$$
is defined to be the $n\times n$ matrix whose 
$i$-th row is
$$\spreadlines{6 pt}\allowdisplaybreaks\alignat 2
\kern -4 em
&1,
\left((D(Bm_{r}+C))^{d_b}\right)^{(b_2-b_1)/d_b},
\cdots,\left((D(Bm_{r}+C))^{d_b}\right)^{
(b_n-b_1)/d_b},
&\kern 1.6 em\text{if}&\quad \ i=\ell_{r}\in S,\cr  
\kern -2.15 em\text{and}\kern 2.15 em
&a_i,a_{i+1},\cdots,a_{i+n-1},
&\kern 1.6 em\text{if}&\quad \ i\notin S. 
\tag 6.29\cr  
\endalignat$$

Letting $\rho\in\Cal S_p$ act on 
$\{m_1,m_2,\dots,m_p\}$ by 
$m_r \rightarrow m_{\rho(r)}$, we observe that  
$$\rho\left(
\det(Q_{n,S,\{b_r\},\{c_r\},\{a_r\}})\right) = 
\text{sign}(\rho)\cdot
\det(Q_{n,S,\{b_r\},\{c_r\},\{a_r\}}),\tag 6.30$$
since permuting $\{m_1,m_2,\dots,m_p\}$ by 
$\rho$ just permutes the rows of 
$Q_{n,S,\{b_r\},\{c_r\},\{a_r\}}$ corresponding to 
$\ i=\ell_{r}\in S$.  

Furthermore, we have 
$$\spreadlines{6 pt}\allowdisplaybreaks\align 
\, &\sum_{\rho\in \Cal S_p}\text{sign}(\rho)
\prod\limits_{r=2}^{p}
\left((D(Bm_{\rho(r)}+C))^{d_c}\right)^{(c_{\ell_r} - 
c_{\ell_1})/d_c}\cr 
=\, &\det\left(\left((D(Bm_{r}+C))^{d_c}\right)^{
(c_{\ell_s} -c_{\ell_1})/d_c}\right)_{1\leq r,s\leq p},
\tag 6.31\cr 
\endalign$$
where in the left-hand-side of (6.31) we used the
definition of determinant in which a different row
element is selected from each column.  Thus, the
determinant in (6.31) is skew-symmetric in 
$\{m_1,m_2,\dots,m_p\}$.  

Finally, by symmetry, we have 
$$\rho\left(
\prod\limits_{r=1}^{p}
(D(Bm_{r}+C))^{c_{\ell_1}+b_1}\right) = 
\prod\limits_{r=1}^{p}
(D(Bm_{r}+C))^{c_{\ell_1}+b_1}.\tag 6.32$$

Since the left-hand-side of (6.26) is independent of 
$\rho\in\Cal S_p$, it now follows from (6.27)--(6.32),
summing both sides of (6.26) over $\rho\in\Cal S_p$,
and interchange of summation that we have the
following lemma.
\proclaim{Lemma 6.12}Let
$M_{n,S,\{b_r\},\{c_r\},\{a_r\}}$ be determined by 
\hbox{\rm Definition 6.2} and 
$Q_{n,S,\{b_r\},\{c_r\},\{a_r\}}$ by 
\hbox{\rm (6.29)}.  Let $d_b$ and $d_c$ be given by 
\hbox{\rm Definition 6.3}.  Let $n=1,2,\cdots$ and
$p=1,2,\cdots,n$.  Then, 
$$\spreadlines{6 pt}\allowdisplaybreaks\align 
\, &\det(M_{n,S,\{b_r\},\{c_r\},\{a_r\}}) \cr 
 =\, &\tfrac{1}{p!} 
(-A)^{-p}q^{p(G-C)}\kern-.95em
\sum\limits_{{y_1,\ldots,y_p\geq 1}\atop 
{m_1,\ldots,m_p\geq 1}}
\kern-.65em
q^{B(m_1y_1+\cdots+m_py_p)}\cr
&\kern 2 em\cdot (-A)^{y_1+\cdots+y_p}
E^{m_1+\cdots+m_p}q^{(F-B)(m_1+\cdots+m_p)}
q^{C(y_1+\cdots+y_p)}\cr
&\kern 2 em\cdot\prod\limits_{r=1}^{p}
(D(Bm_{r}+C))^{c_{\ell_1}+b_1}
\cdot\det(Q_{n,S,\{b_r\},\{c_r\},\{a_r\}})\cr
&\kern 2 em\cdot\det\left(\left((D(Bm_{r}+C))^{d_c}
\right)^{(c_{\ell_s} -c_{\ell_1})/d_c}\right)_{
1\leq r,s\leq p}.\tag 6.33\cr
\endalign$$
\endproclaim

We next rewrite (6.33) by appealing to the following
symmetrization lemma.
\proclaim{Lemma 6.13}Let 
$F(y_1,\dots,y_p;m_1,\dots,m_p)$ be symmetric 
in $\{y_1,\dots,y_p\}$, and in $\{m_1,\dots,m_p\}$. 
Let $G(m_1,\dots,m_p)$ be symmetric  in
$\{m_1,\dots,m_p\}$.  Furthermore, assume that
$F$ equals $0$ if any of $\{m_1,\dots,m_p\}$ are 
equal.  Let $n=1,2,\cdots$ and $p=1,2,\cdots,n$. 
Then,      
$$\spreadlines{6 pt}\allowdisplaybreaks\align 
\, &\sum\limits_{{y_1,\ldots,y_p\geq 1}\atop 
{m_1,\ldots,m_p\geq 1}}
\kern-.65em F(y_1,\dots,y_p;m_1,\dots,m_p)
G(m_1y_1,\dots,m_py_p)\cr 
=\, &p!\kern-1.95em
\sum\limits_{{y_1,\ldots,y_p\geq 1}\atop 
{m_1>m_2>\cdots>m_p\geq 1}}
\kern-1.65em F(y_1,\dots,y_p;m_1,\dots,m_p)
G(m_1y_1,\dots,m_py_p).\tag 6.34\cr 
\endalign$$
\endproclaim
\demo{Proof}Since we can assume that the 
$\{m_1,\dots,m_p\}$ are distinct, we write the
left-hand-side of (6.34) as
$$\sum\limits_{{y_1,\ldots,y_p\geq 1}\atop 
{m_1>m_2>\cdots>m_p\geq 1}}
\sum_{\rho\in \Cal S_p}
F(y_1,\dots,y_p;m_{\rho(1)},\dots,m_{\rho(p)})
G(m_{\rho(1)}y_1,\dots,m_{\rho(p)}y_p).\tag 6.35$$
Next, by an interchange of summation, symmetry of $F$
in $\{m_1,\dots,m_p\}$, and relabelling the 
$\{y_1,\dots,y_p\}$ as 
$\{y_{\rho(1)},\dots,y_{\rho(p)}\}$, we have  
$$\sum_{\rho\in \Cal S_p}
\sum\limits_{{y_{\rho(1)},\ldots,y_{\rho(p)}
\geq 1}\atop  {m_1>m_2>\cdots>m_p\geq 1}}
F(y_{\rho(1)},\dots,y_{\rho(p)};m_1,\dots,m_p)
G(m_{\rho(1)}y_{\rho(1)},
\dots,m_{\rho(p)}y_{\rho(p)}).\tag 6.36$$
The symmetry of $F$ in $\{y_1,\dots,y_p\}$, the
symmetry of $G$ in $\{m_1y_1,\dots,m_py_p\}$, 
the fact that $y_{\rho(1)},\ldots,y_{\rho(p)}\geq 1$ is
equivalent to $y_1,\ldots,y_p\geq 1$, and an
interchange of summation now gives
$$\sum\limits_{{y_1,\ldots,y_p\geq 1}\atop 
{m_1>m_2>\cdots>m_p\geq 1}}
\kern-1.65em F(y_1,\dots,y_p;m_1,\dots,m_p)
G(m_1y_1,\dots,m_py_p)  
\sum_{\rho\in \Cal S_p}1.\tag 6.37$$
Since the innermost sum over $\Cal S_p$ equals $p!$,
we obtain the right-hand-side of (6.34).
\qed\enddemo

It follows immediately from Lemma 6.13 that Lemma
6.12 simplifies to the following lemma.
\proclaim{Lemma 6.14}Let
$M_{n,S,\{b_r\},\{c_r\},\{a_r\}}$ be determined by 
\hbox{\rm Definition 6.2} and 
$Q_{n,S,\{b_r\},\{c_r\},\{a_r\}}$ by 
\hbox{\rm (6.29)}.  Let $d_b$ and $d_c$ be given by 
\hbox{\rm Definition 6.3}.  Let $n=1,2,\cdots$ and
$p=1,2,\cdots,n$.  Then, 
$$\spreadlines{6 pt}\allowdisplaybreaks\align 
\, &\det(M_{n,S,\{b_r\},\{c_r\},\{a_r\}}) \cr 
 =\, & (-A)^{-p}q^{p(G-C)}\kern-1.95em
\sum\limits_{{y_1,\ldots,y_p\geq 1}\atop 
{m_1>m_2>\cdots>m_p\geq 1}}
\kern-1.65em(-A)^{y_1+\cdots+y_p}
E^{m_1+\cdots+m_p}\cr
&\kern 2 em\cdot q^{(F-B)(m_1+\cdots+m_p)}
q^{(Bm_1+C)y_1+\cdots+(Bm_p+C)y_p}\cr
&\kern 2 em\cdot\prod\limits_{r=1}^{p}
(D(Bm_{r}+C))^{c_{\ell_1}+b_1}
\cdot\det(Q_{n,S,\{b_r\},\{c_r\},\{a_r\}})\cr
&\kern 2 em\cdot\det\left(\left((D(Bm_{r}+C))^{d_c}
\right)^{(c_{\ell_s} -c_{\ell_1})/d_c}\right)_{
1\leq r,s\leq p}.\tag 6.38\cr
\endalign$$
\endproclaim
\demo{Proof}Take $G(m_1y_1,\dots,m_py_p):=
q^{B(m_1y_1+\cdots+m_py_p)}$, which is symmetric
in $\{m_1y_1,\dots,m_py_p\}$.  Let 
$F(y_1,\dots,y_p;m_1,\dots,m_p)$ be the rest of the
term in the double multiple sum on the right-hand-side
of (6.33).  The two factors of $F$ which are
determinants are each skew-symmetric in   
$\{m_1,\dots,m_p\}$.  Thus, their product is symmetric
in $\{m_1,\dots,m_p\}$, and vanishes if any of
$\{m_1,\dots,m_p\}$ are equal.  All other factors of
$F$ are clearly symmetric in $\{y_1,\dots,y_p\}$, or in
$\{m_1,\dots,m_p\}$.
\qed\enddemo

In order to finish the proof of Theorem 6.7 we utilize the
Schur functions in Definition 6.5 and the Laplace
expansion formula in Theorem 6.8 to rewrite the two
determinants appearing in each term on the
right-hand-side of (6.38).  

By equation (6.12) of Definition 6.5, column operations,
and the product formula for a Vandermonde  
determinant we have the following lemma.
\proclaim{Lemma 6.15}Let $x_1,\dots,x_p$  
be indeterminant, and $1\leq\mu_1<\mu_2<
\dots<\mu_p$ be positive integers.  Let
$p=1,2,\cdots$.  Then,
$$\det(x_i^{\mu_j-\mu_1})_{1\leq i,j\leq p} 
= (-1)^{p\choose 2}\kern-.8em 
\prod\limits_{1\leq r<s\leq p}\kern-.8em
(x_r-x_s)\cdot  s_\lambda(x_1,\ldots,x_p),\tag 6.39$$
where $s_\lambda(x_1,\ldots,x_p)$ is the Schur
function in \hbox{\rm (6.12)}, with the
partition\newline  
$\lambda =(\lambda_1\geq\lambda_2\geq
\ldots\geq\lambda_p\geq 0)$ given by 
$$\lambda_i = \mu_{p-i+1}-\mu_1+i-p,
\kern 6 em\text{for}\quad \ i=1,2,\cdots,p.
\tag 6.40$$ 
\endproclaim
\demo{Proof}By column operations and a trivial
simplification we first have 
$$\spreadlines{6 pt}\allowdisplaybreaks\align 
\, &\det(x_i^{\mu_j-\mu_1})_{1\leq i,j\leq p} = 
(-1)^{p\choose 2}
\det(x_i^{\mu_{p-j+1}-\mu_1})_{1\leq i,j\leq p}\cr
 =\, &(-1)^{p\choose 2}
\det(x_i^{(\mu_{p-j+1}-\mu_1+j-p)+p-j})_{
1\leq i,j\leq p}.\tag 6.41\cr
\endalign$$
Next, appeal to (6.12) and the product formula for a
Vandermonde determinant.  Finally, since the $\mu_i$
are strictly increasing, it is clear that the $\lambda_i$ in
(6.40) determine the partition 
$\lambda =(\lambda_1\geq\lambda_2\geq
\ldots\geq\lambda_p\geq 0)$.  Just note that
$\lambda_p=0$ and $\lambda_i-\lambda_{i+1} = 
\mu_{p-i+1}-\mu_{p-i}-1\geq 0$, if 
$\ i=1,2,\cdots,p-1$.
\qed\enddemo

We now use Lemma 6.15 to rewrite the second
determinant on the right-hand-side of (6.38).  We
specialize $x_r$ and $\mu_r$ as follows:
$$\spreadlines{6 pt}\allowdisplaybreaks\alignat 2
\kern -4 em
&x_r:=(D(Bm_{r}+C))^{d_c},
&\kern 4 em\text{for}&\quad \ r=1,2,\cdots,p,
\tag 6.42\cr  
\kern -8.95 em\text{and}\kern 8.95 em
&\mu_r:=c_{\ell_r}/d_c,
&\kern 4 em\text{for}&\quad \ r=1,2,\cdots,p.
\tag 6.43\cr  
\endalignat$$
Since $\ell_r$ and $c_r$ are both strictly increasing and
$d_c$ is constant, then the $\mu_r=c_{\ell_r}/d_c$
are also strictly increasing.  We obtain 
$$\spreadlines{6 pt}\allowdisplaybreaks\align 
\, &\det\left(\left((D(Bm_{r}+C))^{d_c}
\right)^{(c_{\ell_s} -c_{\ell_1})/d_c}\right)_{
1\leq r,s\leq p}\cr 
=\, &(-1)^{p\choose 2}\kern-.8em
\prod\limits_{1\leq r<s\leq p}\kern-1.1em
\left((D(Bm_r+C))^{d_c}-(D(Bm_s+C))^{d_c}\right)\cr 
&\kern 2 em\cdot 
s_\lambda\left((D(Bm_1+C))^{d_c},\ldots,
(D(Bm_p+C))^{d_c}\right),\tag 6.44\cr
\endalign$$
where $\lambda_i$ is given by (6.13).  

We next use Theorem 6.8 to expand the first
determinant on the right-hand-side of (6.38) into the
following sum.  
$$\spreadlines{6 pt}\allowdisplaybreaks\align 
&\det(Q_{n,S,\{b_r\},\{c_r\},\{a_r\}})\cr
 =\, & \sum\limits_{{\emptyset \subset T\subseteq I_n}
\atop {\Vert T\Vert =p}}\kern-.5em 
\det(C_{p,S,T})\det(D_{n-p,S^c,T^c})\cdot
(-1)^{\Sigma (S)+\Sigma (T)},\tag 6.45\cr 
\endalign$$
where $S$, $T$, $S^c$, $T^c$, $\Sigma (S)$, and  
$\Sigma (T)$ be given by (6.1)--(6.3), 
$Q_{n,S,\{b_r\},\{c_r\},\{a_r\}}$ is defined by (6.29), 
the $(n-p) \times (n-p)$ matrix $D_{n-p,S^c,T^c}$ is
determined by (6.11), and the $p \times p$ matrix 
$C_{p,S,T}$ is defined by 
$$C_{p,S,T}:= 
\left[\left((D(Bm_{r}+C))^{d_b}
\right)^{(b_{j_s} -b_1)/d_b}\right]_{
1\leq r,s\leq p}.\tag 6.46$$
  
By factoring $\left((D(Bm_{r}+C))^{d_b}
\right)^{(b_{j_1} -b_1)/d_b}$ out of the $r$-th row of 
$C_{p,S,T}$ it is immediate that 
$$\spreadlines{6 pt}\allowdisplaybreaks\align 
\det(C_{p,S,T}) = \, &\prod\limits_{r=1}^{p} 
\left((D(Bm_{r}+C))^{d_b}
\right)^{(b_{j_1} -b_1)/d_b}
\tag 6.47a\cr 
&\kern .5 em\cdot 
\det\left(\left((D(Bm_{r}+C))^{d_b}
\right)^{(b_{j_s} -b_{j_1})/d_b}\right)_{
1\leq r,s\leq p}.\tag 6.47b\cr
\endalign$$

We now use Lemma 6.15 to rewrite the $p \times p$ 
determinant in (6.47b).  We specialize $x_r$ and
$\mu_r$ as follows:  
$$\spreadlines{6 pt}\allowdisplaybreaks\alignat 2
\kern -4 em
&x_r:=(D(Bm_{r}+C))^{d_b},
&\kern 4 em\text{for}&\quad \ r=1,2,\cdots,p,
\tag 6.48\cr  
\kern -8.95 em\text{and}\kern 8.95 em
&\mu_r:=b_{j_r}/d_b,
&\kern 4 em\text{for}&\quad \ r=1,2,\cdots,p.
\tag 6.49\cr  
\endalignat$$
Since $j_r$ and $b_r$ are both strictly increasing and
$d_b$ is constant, then the $\mu_r=b_{j_r}/d_b$
are also strictly increasing.  We obtain 
$$\spreadlines{6 pt}\allowdisplaybreaks\align 
\, &\det\left(\left((D(Bm_{r}+C))^{d_b}
\right)^{(b_{j_s} -b_{j_1})/d_b}\right)_{
1\leq r,s\leq p}\cr 
=\, &(-1)^{p\choose 2}\kern-.8em
\prod\limits_{1\leq r<s\leq p}\kern-1.1em
\left((D(Bm_r+C))^{d_b}-(D(Bm_s+C))^{d_b}\right)\cr 
&\kern 2 em\cdot 
s_\nu\left((D(Bm_1+C))^{d_b},\ldots,
(D(Bm_p+C))^{d_b}\right),\tag 6.50\cr
\endalign$$
where $\nu_i$ is given by (6.14).  
 
By combining (6.45)--(6.50) we have the following
lemma.
\proclaim{Lemma 6.16}Let
$Q_{n,S,\{b_r\},\{c_r\},\{a_r\}}$ be defined by 
\hbox{\rm (6.29)}.  Let $S$, $T$, $S^c$, $T^c$, 
$\Sigma (S)$, $\Sigma (T)$ be given by 
\hbox{\rm (6.1)--(6.3)}, and let $d_b$ and $d_c$ be
given by \hbox{\rm Definition 6.3}.  Let the 
$(n-p) \times (n-p)$ matrix $D_{n-p,S^c,T^c}$ be
determined by \hbox{\rm (6.11)}.  Let the Schur
function $s_\nu$ be defined by \hbox{\rm (6.12)},
with the partition $\nu =(\nu_1\geq\nu_2\geq
\ldots\geq\nu_p\geq 0)$ given by \hbox{\rm
(6.14)}.  Then,
$$\spreadlines{6 pt}\allowdisplaybreaks\align 
\, &\det(Q_{n,S,\{b_r\},\{c_r\},\{a_r\}})\cr 
=\, &(-1)^{p\choose 2}\kern-.8em
\prod\limits_{1\leq r<s\leq p}\kern-1.1em
\left((D(Bm_r+C))^{d_b}-(D(Bm_s+C))^{d_b}\right)\cr 
\, &\cdot
\sum\limits_{{\emptyset \subset T\subseteq I_n}
\atop {\Vert T\Vert =p}}\kern-.5em 
(-1)^{\Sigma (S)+\Sigma (T)}\cdot 
\det (D_{n-p,S^c,T^c})\cdot 
\left(\prod\limits_{r=1}^{p} 
D(Bm_{r}+C)\right)^{b_{j_1} -b_1}\cr
&\kern 4 em\cdot 
s_\nu\left((D(Bm_1+C))^{d_b},\ldots,
(D(Bm_p+C))^{d_b}\right).\tag 6.51\cr
\endalign$$
\endproclaim

The proof of Theorem 6.7 is completed by substituting
(6.44) and (6.51) into the right-hand-side of (6.38) and
then simplifying.

The derivation of the Schur function form of our
$\chi_n$ determinant identities in Section 7 requires 
a slight modification of Theorem 6.7, which arises from 
replacing $a_i,a_{i+1},\cdots,a_{i+n-1}$ in (6.7) by 
$a_i,a_{i+1},\cdots,a_{i+n-2},a_{i+n}$.  We first need
the following two definitions.

\definition{Definition 6.17  (Second Lambert series
matrix)} Let $b_1,b_2,\dots,b_n$ and
$c_1,c_2,\dots,c_n$  be indeterminant, and let
$n=1,2,\cdots$.  Furthermore, assume that
$b_1<b_2<\dots<b_n$ and
$c_1<c_2<\dots<c_n$.  Take 
$\{a_r:r=1,2,\cdots\}$ to be an arbitrary sequence.  Let
$S$ and $S^c$ be the subsets of $I_n$ in (6.1), with
$p=1,2,\cdots,n$.  Let $L_u(r;A,\dots,G)$ be the
Lambert series in Definition 6.1.  Then,
$$\overline M_{n,S,\{b_r\},\{c_r\},\{a_r\}}
\equiv\overline M_{n,S,\{b_1,\dots,b_n\},
\{c_1,\dots,c_n\},\{a_r\}}\tag 6.52$$
is defined to be the $n\times n$ matrix whose 
$i$-th row is
$$\spreadlines{6 pt}\allowdisplaybreaks\alignat 2
\kern -4 em
&L_{c_i+b_1}(m_{\mu};A,\dots,G),
L_{c_i+b_2}(m_{\mu};A,\dots,G),\cdots,
L_{c_i+b_n}(m_{\mu};A,\dots,G),
&\kern 1.6 em\text{if}&\quad \ i=\ell_{\mu}\in S,\cr  
\kern -1.12em\text{and}\kern 1.12 em
&a_i,a_{i+1},\cdots,a_{i+n-2},a_{i+n},
&\kern 1.6 em\text{if}&\quad \ i\notin S. 
\tag 6.53\cr  
\endalignat$$
\enddefinition 

\definition{Definition 6.18  (Second Laplace expansion
formula determinant)}Let $\{a_r:r=1,2,\cdots\}$ be
an arbitrary sequence.  Let the sets $S^c$ and $T^c$
be given by (6.1) and (6.2), and let
$p=1,2,\cdots,n$.  Then, 
$$\det(\overline D_{n-p,S^c,T^c})\tag 6.54$$ 
is the determinant of the $(n-p) \times (n-p)$ matrix  
$$\overline D_{n-p,S^c,T^c}:= 
\left[a_{(\ell_{p+r}+j_{p+s}-1+\chi
(j_{p+s}=n))}\right]_{ 1\leq r,s\leq n-p},\tag 6.55$$
where 
$$\chi (A):=1,\ \text{ if}\ A\ 
\text{ is true, and}\ 0,\ 
\text{ otherwise.}\tag 6.56$$
\enddefinition

The same computations that established Theorem 
6.7 now give the following analogous expansion.
\proclaim{Theorem 6.19 (Expansion of second Lambert 
series determinant)} Let
$L_u(r;A,\dots,G)$ be the Lambert series in \hbox{\rm
Definition 6.1}.  Let the $n\times n$ matrix
$\overline M_{n,S,\{b_r\},\{c_r\},\{a_r\}}$ of Lambert
series  be given by \hbox{\rm Definition 6.17},  and take 
$n=1,2,\cdots$ and $p=1,2,\cdots,n$.   Let $A$, $B$,
$C$, $D$, $E$, $F$, $G$, and 
$b_1,b_2,\dots,b_n$ and $c_1,c_2,\dots,c_n$ be
indeterminant, and assume that $0<|q|<1$.   
Let $S$, $T$, $S^c$, $T^c$, $\Sigma (S)$, and  
$\Sigma (T)$ be given by \hbox{\rm (6.1)--(6.3)}.
Take $\{a_r:r=1,2,\cdots\}$ to be an arbitrary 
sequence.  Let $d_b$ and $d_c$ be given by 
\hbox{\rm Definition 6.3} and $s_\lambda$ and 
$s_\nu$ be the Schur functions in 
\hbox{\rm Definition 6.5} with partitions $\lambda$ and 
$\nu$ as in \hbox{\rm Definition 6.6}.  Finally, take the 
$(n-p) \times (n-p)$ matrix $\overline
D_{n-p,S^c,T^c}$ in \hbox{\rm Definition 6.18}.  We
then have the expansion formula
$$\spreadlines{6 pt}\allowdisplaybreaks\align 
&\det(\overline M_{n,S,\{b_r\},\{c_r\},\{a_r\}})\cr
 =\, & (-A)^{-p}q^{p(G-C)}\kern-1.95em
\sum\limits_{{y_1,\ldots,y_p\geq 1}\atop 
{m_1>m_2>\cdots>m_p\geq 1}}
\kern-1.65em(-A)^{y_1+\cdots+y_p}
E^{m_1+\cdots+m_p}\cr
&\kern 2 em\cdot q^{(F-B)(m_1+\cdots+m_p)}
q^{(Bm_1+C)y_1+\cdots+(Bm_p+C)y_p}\cr
&\kern 2 em\cdot 
\prod\limits_{1\leq r<s\leq p}\kern-1.1em
\left((D(Bm_r+C))^{d_c}-(D(Bm_s+C))^{d_c}\right)\cr
&\kern 2 em\cdot 
\prod\limits_{1\leq r<s\leq p}\kern-1.1em
\left((D(Bm_r+C))^{d_b}-(D(Bm_s+C))^{d_b}\right)\cr
&\kern 2 em\cdot 
s_\lambda\left((D(Bm_1+C))^{d_c},\ldots,
(D(Bm_p+C))^{d_c}\right)\cr
&\cdot
\sum\limits_{{\emptyset \subset T\subseteq I_n}
\atop {\Vert T\Vert =p}}\kern-.5em 
(-1)^{\Sigma (S)+\Sigma (T)}\cdot 
\det (\overline D_{n-p,S^c,T^c})\cdot
\left(\prod\limits_{r=1}^{p}D(Bm_r+C)
\right)^{\!\!c_{\ell_1}+b_{j_1}}\cr 
&\kern 3 em\cdot
s_\nu\left((D(Bm_1+C))^{d_b},\ldots,
(D(Bm_p+C))^{d_b}\right).\tag 6.57\cr 
\endalign$$
\endproclaim

\head 7. The Schur function form of sums of squares 
identities\endhead

In this section we first apply the key determinant expansion
formulas in Theorems 6.7 and  6.17 to most of the main
identities in Section 5 to obtain the Schur function form 
of our infinite families of sums of squares and related
identities.   These include the two Kac--Wakimoto \cite{120, pp. 452}
conjectured identities involving representing a positive
integer by sums of $4n^2$ or $4n(n+1)$ triangular
numbers, respectively. In addition, we obtain in Corollary
7.10 our analog of these Kac--Wakimoto identities which
involves representing a positive integer by sums of $2n$
squares  and $(2n)^2$ triangular numbers.   Motivated by
the analysis in \cite{41, 42}, we next use  Jacobi's
transformation of the theta functions $\vartheta_4$ and
$\vartheta_2$ to derive a direct  connection between our
identities involving $4n^2$ or $4n(n+1)$ squares, and the
identities involving $4n^2$ or $4n(n+1)$ triangular
numbers.    In a different direction we also apply the classical
techniques in \cite{90, pp. 96--101} and \cite{28, pp.
288--305} to several of the theorems in this section to
obtain the corresponding infinite families of lattice sum
transformations. 

We begin with the Schur function form of the
$4n^2$ squares  identity in the following theorem. 
\proclaim{Theorem 7.1 } Let $n=1,2,3,\cdots$. Then
$$\spreadlines{6 pt}\allowdisplaybreaks\align 
&\vartheta_3(0,-q)^{4n^2}\tag 7.1\cr 
=\,& 1+\sum\limits_{p=1}^n
(-1)^p2^{2n^2+n}\prod\limits_{r=1}^{2n-1}(r!)^{-1}
\kern-1.8em
\sum\limits_{{y_1,\ldots,y_p\geq 1}\atop 
{m_1>m_2>\cdots>m_p\geq 1}}
\kern-1.65em(-1)^{y_1+\cdots+y_p}\cr
&\cdot(-1)^{m_1+\cdots+m_p}
q^{m_1y_1+\cdots+m_py_p}\kern-.65em
\prod\limits_{1\leq r<s\leq p}\kern-1.1em
\left(m_r^2-m_s^2\right)^2\cr
&\cdot(m_1m_2\cdots m_p)\kern-.65em 
\sum\limits_{{\emptyset \subset S, T\subseteq I_n}
\atop {\Vert S\Vert=\Vert T\Vert =p}}\kern-.5em 
(-1)^{\Sigma (S)+\Sigma (T)}\cdot 
\det (D_{n-p,S^c,T^c})\cr
&\cdot (m_1m_2\cdots m_p)^{2\ell_1+2j_1-4}
s_{\lambda}(m_1^2,\ldots,m_p^2)\,
s_{\nu}(m_1^2,\ldots,m_p^2),\cr 
\endalign$$
where $\vartheta_3(0,-q)$ is determined by 
\hbox{\rm(1.1)}, 
the sets $S$, $S^c$, $T$, $T^c$ are given by
\hbox{\rm(6.1)-(6.2)}, 
$\Sigma (S)$ and $\Sigma (T)$ by \hbox{\rm(6.3)}, 
the $(n-p) \times (n-p)$ matrix  
$$D_{n-p,S^c,T^c}:= 
\left[c_{(\ell_{p+r}+j_{p+s}-1)}\right]_{
1\leq r,s\leq n-p},\tag 7.2$$
where the $c_i$ are defined by \hbox{\rm(5.23)}, with 
the $B_{2i}$ in \hbox{\rm(2.60)}, and $s_\lambda$ and 
$s_\nu$ are the Schur functions in \hbox{\rm(6.12)},
with the partitions $\lambda$ and $\nu$ given by 
$$\lambda_i:=\ell_{p-i+1}-\ell_1+i-p\quad 
\text{ and}\quad
\nu_i:=j_{p-i+1}-j_1+i-p,\quad 
\text{ for}\ i=1,2,\cdots,p,\tag 7.3$$ 
where the $\ell_i$ and $j_i$ are elements of the sets 
$S$ and $T$, respectively.  
\endproclaim
\demo{Proof} We apply Theorem 6.7 to Theorem 5.4.  
The Lambert series $U_{2i-1}$ in (5.30), determined by
(5.22), is 
$$-L_{2i-1}(r;1,1,0,1,-1,1,0),\tag 7.4$$
as defined in (6.4).  Factor $-1$ out of each row for 
$i\in S$ of the matrix $M_{n,S}$ in (5.29).  This gives
$(-1)^p$.   Next, apply the following case of Theorem 6.7 
to the resulting $\det (M_{n,S})$ in the $p$-th term of
equation (5.29) in Theorem 5.4:
$$\spreadlines{6 pt}\allowdisplaybreaks\align 
A = \ &1,\quad B=1,\quad C=0,\quad D=1,\quad 
E=-1,\quad F=1,\quad G=0,\tag 7.5\cr
a_i  =\ & (-1)^{i-1}{(2^{2i}-1)\over {4i}}\cdot |B_{2i}|,
\quad \text{ for}\quad \ i=1,2,3,\cdots,\tag 7.6\cr
c_i=\ & 2i-1:=2\ell_{\mu}-1,\tag 7.7\cr
b_i=\ & 2(i-1),\quad 
\text{ for}\quad i=1,2,\cdots,n,\tag 7.8\cr
d_b=\ & d_c=2,\tag 7.9\cr
\endalign$$
with $B_{2i}$ the Bernoulli numbers in (2.60). 

It is immediate that
$$\spreadlines{6 pt}\allowdisplaybreaks\align 
&\ G-C = 0,\quad F-B=0,\quad Bm_r+C=m_r,\cr
&\ (D(Bm_r+C))^{d_b} = 
(D(Bm_r+C))^{d_c} = m_r^2.\tag 7.10\cr
\endalign$$

The $\lambda_i$ and $\nu_i$ are given by (7.3).  
Note that $\ell_{p-i+1}$ and $j_{p-i+1}$ are substituted
for $i$ in (7.7) and (7.8), respectively, before computing
the $\lambda_i$ and $\nu_i$ in (6.13) and (6.14).  

Interchanging the inner sum in (5.29) with the double
multiple sum over ${y_1,\ldots,y_p\geq 1}$ and 
${m_1>m_2>\cdots>m_p\geq 1}$ in (6.15), multiplying
and dividing by $(m_1m_2\cdots m_p)$, 
and simplifying, now yields (7.1).
\qed\enddemo

The Schur function form of the $4n(n+1)$ 
squares identity is given by the following theorem. 
\proclaim{Theorem 7.2 } Let
$n=1,2,3,\cdots$.  Then
$$\spreadlines{6 pt}\allowdisplaybreaks\align 
&\vartheta_3(0,-q)^{4n(n+1)}\tag 7.11\cr 
=\,& 1+\sum\limits_{p=1}^n
(-1)^{n-p}2^{2n^2+3n}\prod\limits_{r=1}^{2n}(r!)^{-1}
\kern-1.8em
\sum\limits_{{y_1,\ldots,y_p\geq 1}\atop 
{m_1>m_2>\cdots>m_p\geq 1}}
\kern-1.65em(-1)^{m_1+\cdots+m_p}\cr
&\cdot q^{m_1y_1+\cdots+m_py_p}\,
(m_1m_2\cdots m_p)^3\kern-.65em
\prod\limits_{1\leq r<s\leq p}\kern-1.1em
\left(m_r^2-m_s^2\right)^2\cr
&\cdot
\sum\limits_{{\emptyset \subset S, T\subseteq I_n}
\atop {\Vert S\Vert=\Vert T\Vert =p}}\kern-.5em 
(-1)^{\Sigma (S)+\Sigma (T)}\cdot 
\det (D_{n-p,S^c,T^c})\cr 
&\cdot (m_1m_2\cdots m_p)^{2\ell_1+2j_1-4} 
s_{\lambda}(m_1^2,\ldots,m_p^2)\,
s_{\nu}(m_1^2,\ldots,m_p^2),\cr 
\endalign$$
where the same assumptions hold as in 
\hbox{\rm Theorem 7.1}, except that the 
$(n-p) \times (n-p)$ matrix  
$$D_{n-p,S^c,T^c}:= 
\left[a_{(\ell_{p+r}+j_{p+s}-1)}\right]_{
1\leq r,s\leq n-p},\tag 7.12$$
where the $a_i$ are defined by \hbox{\rm(5.40)}, with 
the $B_{2i+2}$ in \hbox{\rm(2.60)}.
\endproclaim
\demo{Proof} We apply Theorem 6.7 to Theorem 5.6.  
The Lambert series $G_{2i+1}$ in (5.47), determined by
(5.39), is 
$$L_{2i+1}(r;-1,1,0,1,-1,1,0),\tag 7.13$$
as defined in (6.4).  Apply the following case of Theorem
6.7 to the $\det (M_{n,S})$ in the $p$-th term of
equation (5.46) in Theorem 5.6:
$$\spreadlines{6 pt}\allowdisplaybreaks\align 
A = \ &-1,\quad B=1,\quad C=0,\quad D=1,\quad 
E=-1,\quad F=1,\quad G=0,\tag 7.14\cr
a_i  =\ & (-1)^{i}{(2^{2i+2}-1)\over {4(i+1)}}
\cdot |B_{2i+2}|,
\quad \text{ for}\quad \ i=1,2,3,\cdots,\tag 7.15\cr
c_i=\ & 2i+1:=2\ell_{\mu}+1,\tag 7.16\cr
b_i=\ & 2(i-1),\quad 
\text{ for}\quad i=1,2,\cdots,n,\tag 7.17\cr
d_b=\ & d_c=2,\tag 7.18\cr
\endalign$$
with $B_{2i+2}$ the Bernoulli numbers in
(2.60). 

It is immediate that (7.10) holds.
The $\lambda_i$ and $\nu_i$ are given by (7.3).  
Note that $\ell_{p-i+1}$ and $j_{p-i+1}$ are substituted
for $i$ in (7.16) and (7.17), respectively, before computing
the $\lambda_i$ and $\nu_i$ in (6.13) and (6.14).  

Interchanging the inner sum in (5.46) with the double
multiple sum over ${y_1,\ldots,y_p\geq 1}$ and 
${m_1>m_2>\cdots>m_p\geq 1}$ in (6.15), multiplying
and dividing by $(m_1m_2\cdots m_p)^3$, 
and simplifying, now yields (7.11).
\qed\enddemo

We have in Section 8 written down the $n=3$ cases of
Theorems 7.1 and 7.2 which yield explicit formulas for
$36$ and $48$ squares, respectively.   For reference, we
have also written down the $n=2$ cases. 

The next two theorems require the Euler numbers $E_n$
defined by (2.61).  We first have the following theorem.
\proclaim{Theorem 7.3 } Let
$n=1,2,3,\cdots$.  Then
$$\spreadlines{6 pt}\allowdisplaybreaks\align 
&\vartheta_3(0,q)^{2n(n-1)}
\vartheta_3(0,-q)^{2n^2}\tag 7.19\cr 
=\,&1+\sum\limits_{p=1}^n
(-1)^p2^{2n}\prod\limits_{r=1}^{n-1}(2r)!^{-2}
\kern-2.3em
\sum\limits_{{y_1,\ldots,y_p\geq 1}\atop 
{m_1>m_2>\cdots>m_p\geq 1}\sixrm (odd)}
\kern-2.15em(-1)^{y_1+\cdots+y_p}\cr 
&\cdot (-1)^{\tfrac{1}{2}(p+m_1+\cdots+m_p)}
q^{m_1y_1+\cdots+m_py_p}\kern-.65em
\prod\limits_{1\leq r<s\leq p}\kern-1.1em
\left(m_r^2-m_s^2\right)^2\cr
&\cdot
\sum\limits_{{\emptyset \subset S, T\subseteq I_n}
\atop {\Vert S\Vert=\Vert T\Vert =p}}\kern-.5em 
(-1)^{\Sigma (S)+\Sigma (T)}\cdot 
\det (D_{n-p,S^c,T^c})\cr 
&\cdot (m_1m_2\cdots m_p)^{2\ell_1+2j_1-4} 
s_{\lambda}(m_1^2,\ldots,m_p^2)\,
s_{\nu}(m_1^2,\ldots,m_p^2),\cr 
\endalign$$
where the same assumptions hold as in 
\hbox{\rm Theorem 7.1}, except that the 
$(n-p) \times (n-p)$ matrix  
$$D_{n-p,S^c,T^c}:= 
\left[b_{(\ell_{p+r}+j_{p+s}-1)}\right]_{
1\leq r,s\leq n-p},\tag 7.20$$
where the $b_i$ are defined by \hbox{\rm(5.57)}, with 
the $E_{2i-2}$ in \hbox{\rm(2.61)}.
\endproclaim
\demo{Proof} We apply Theorem 6.7 to Theorem 5.8.  
The Lambert series $R_{2i-2}$ in (5.67), determined by
(5.56), is 
$$-L_{2i-2}(r;1,2,-1,1,-1,2,-1),\tag 7.21$$
as defined in (6.4).  Factor $-1$ out of each row for 
$i\in S$ of the matrix $M_{n,S}$ in (5.66).  This gives
$(-1)^p$.   Next, apply the following case of Theorem 6.7 
to the resulting $\det (M_{n,S})$ in the $p$-th term of
equation (5.66) in Theorem 5.8:
$$\spreadlines{6 pt}\allowdisplaybreaks\align 
A = \ &1,\quad B=2,\quad C=-1,\quad D=1,\quad 
E=-1,\quad F=2,\quad G=-1,\tag 7.22\cr
a_i  =\ &(-1)^{i-1}\cdot {\tfrac {1}{4}}\cdot |E_{2i-2}|,
\quad \text{ for}\quad \ i=1,2,3,\cdots,\tag 7.23\cr
c_i=\ & 2i-2:=2\ell_{\mu}-2,\tag 7.24\cr
b_i=\ & 2(i-1),\quad 
\text{ for}\quad i=1,2,\cdots,n,\tag 7.25\cr
d_b=\ & d_c=2,\tag 7.26\cr
\endalign$$
with $E_{2i-2}$ the Euler numbers in (2.61). 

It is immediate that
$$\spreadlines{6 pt}\allowdisplaybreaks\align 
&\ G-C = 0,\quad F-B=0,\quad Bm_r+C=2m_r-1,\cr
&\ (D(Bm_r+C))^{d_b} = 
(D(Bm_r+C))^{d_c} = (2m_r-1)^2.\tag 7.27\cr
\endalign$$

The $\lambda_i$ and $\nu_i$ are given by (7.3).  
Note that $\ell_{p-i+1}$ and $j_{p-i+1}$ are substituted
for $i$ in (7.24) and (7.25), respectively, before computing
the $\lambda_i$ and $\nu_i$ in (6.13) and (6.14).  

Interchanging the inner sum in (5.66) with the double
multiple sum over ${y_1,\ldots,y_p\geq 1}$ and 
${m_1>m_2>\cdots>m_p\geq 1}$ in (6.15), relabelling 
$(2m_r-1)$ as $m_r$ (odd), 
and simplifying, now yields (7.19).
\qed\enddemo

 We next have the following theorem.
\proclaim{Theorem 7.4 } Let
$n=1,2,3,\cdots$.  Then
$$\spreadlines{6 pt}\allowdisplaybreaks\align 
&\vartheta_3(0,q)^{2n(n+1)}
\vartheta_3(0,-q)^{2n^2}\tag 7.28\cr
=\,&1+\sum\limits_{p=1}^n
(-1)^{n-p}2^{2n}\prod\limits_{r=1}^{n}(2r-1)!^{-2}
\kern-2.3em
\sum\limits_{{y_1,\ldots,y_p\geq 1}\atop 
{m_1>m_2>\cdots>m_p\geq 1}\sixrm (odd)}
\kern-2.15em(-1)^{y_1+\cdots+y_p}\cr 
&\cdot (-1)^{\tfrac{1}{2}(p+m_1+\cdots+m_p)}
q^{m_1y_1+\cdots+m_py_p}\kern-.65em
\prod\limits_{1\leq r<s\leq p}\kern-1.1em
\left(m_r^2-m_s^2\right)^2\cr
&\cdot(m_1m_2\cdots m_p)^2\kern-.65em 
\sum\limits_{{\emptyset \subset S, T\subseteq I_n}
\atop {\Vert S\Vert=\Vert T\Vert =p}}\kern-.5em 
(-1)^{\Sigma (S)+\Sigma (T)}\cdot 
\det (D_{n-p,S^c,T^c})\cr 
&\cdot (m_1m_2\cdots m_p)^{2\ell_1+2j_1-4} 
s_{\lambda}(m_1^2,\ldots,m_p^2)\,
s_{\nu}(m_1^2,\ldots,m_p^2),\cr 
\endalign$$
where the same assumptions hold as in 
\hbox{\rm Theorem 7.1}, except that the 
$(n-p) \times (n-p)$ matrix  
$$D_{n-p,S^c,T^c}:= 
\left[b_{(\ell_{p+r}+j_{p+s})}\right]_{
1\leq r,s\leq n-p},\tag 7.29$$
where the $b_i$ are determined by \hbox{\rm(5.57)}, with 
the $E_{2i-2}$ in \hbox{\rm(2.61)}.
\endproclaim
\demo{Proof} We apply Theorem 6.7 to Theorem 5.10.  
The Lambert series $R_{2i}$ in (5.86), determined by
(5.56), is 
$$-L_{2i}(r;1,2,-1,1,-1,2,-1),\tag 7.30$$
as defined in (6.4).  Factor $-1$ out of each row for 
$i\in S$ of the matrix $M_{n,S}$ in (5.85).  This gives
$(-1)^p$.   Next, apply the following case of Theorem 6.7 
to the resulting $\det (M_{n,S})$ in the $p$-th term of
equation (5.85) in Theorem 5.10:
$$\spreadlines{6 pt}\allowdisplaybreaks\align 
A = \ &1,\quad B=2,\quad C=-1,\quad D=1,\quad 
E=-1,\quad F=2,\quad G=-1,\tag 7.31\cr
a_i  =\ &(-1)^{i}\cdot {\tfrac {1}{4}}\cdot |E_{2i}|,
\quad \text{ for}\quad \ i=1,2,3,\cdots,\tag 7.32\cr
c_i=\ & 2i:=2\ell_{\mu},\tag 7.33\cr
b_i=\ & 2(i-1),\quad 
\text{ for}\quad i=1,2,\cdots,n,\tag 7.34\cr
d_b=\ & d_c=2,\tag 7.35\cr
\endalign$$
with $E_{2i}$ the Euler numbers in (2.61). 

It is immediate that (7.27) holds. 
The $\lambda_i$ and $\nu_i$ are given by (7.3).  
Note that $\ell_{p-i+1}$ and $j_{p-i+1}$ are substituted
for $i$ in (7.33) and (7.34), respectively, before computing
the $\lambda_i$ and $\nu_i$ in (6.13) and (6.14).  

Interchanging the inner sum in (5.85) with the double
multiple sum over ${y_1,\ldots,y_p\geq 1}$ and 
${m_1>m_2>\cdots>m_p\geq 1}$ in (6.15), relabelling 
$(2m_r-1)$ as $m_r$ (odd), multiplying
and dividing by $(m_1m_2\cdots m_p)^2$, 
and simplifying, now yields (7.28).
\qed\enddemo

The next two results complete our proof of the 
Kac--Wakimoto conjectured identities for triangular
numbers in \cite{120, pp. 452}.  

We first have the following theorem.
\proclaim{Theorem 7.5} 
Let $\vartheta_2 (0,q)$ be defined by \hbox{\rm(1.2)},
and let $n=1,2,3,\cdots$.    We then have
$$\spreadlines{6 pt}\allowdisplaybreaks\align 
\vartheta_2(0,q)^{4n^2}=\,&
4^{n(n+1)}\prod\limits_{r=1}^{2n-1}(r!)^{-1}
\kern-2.6em\sum\limits_{{y_1,\ldots ,y_n\geq 1
\hbox {\sixrm (odd)}}\atop 
{m_1>m_2>\cdots >m_n\geq 1 \hbox {\sixrm (odd)}}}
\kern-2.6em m_1m_2\cdots m_n\cr
&\cdot q^{m_1y_1+\cdots +m_ny_n}
\kern-.65em
\prod\limits_{1\leq r<s\leq n}\kern-1.1em
\left(m_r^2-m_s^2\right)^2,
\tag 7.36\cr
\kern -13.8 em \text{and}\kern 13.8 em\quad &\cr 
\vartheta_2(0,q^{1/2})^{4n(n+1)}=\,&
2^{n(4n+5)}
\prod\limits_{r=1}^{2n}(r!)^{-1}
\kern-1.65em\sum\limits_{{y_1,\ldots ,y_n\geq 1
\hbox {\sixrm (odd)}}\atop 
{m_1>m_2>\cdots >m_n\geq 1 }}
\kern-1.8em 
\left(m_1m_2\cdots m_n\right)^3\cr
&\cdot q^{m_1y_1+\cdots +m_ny_n}
\kern-.65em
\prod\limits_{1\leq r<s\leq n}\kern-1.1em
\left(m_r^2-m_s^2\right)^2.
\tag 7.37\cr
\endalign$$
\endproclaim 
\demo{Proof} We apply Theorem 6.7 to Theorem 5.11.  
The Lambert series $C_{2i-1}$ in (5.93), determined by
(5.95), is 
$$L_{2i-1}(r;-1,4,-2,\tfrac{1}{2},1,2,-1),\tag 7.38$$
as defined in (6.4).  Apply the following case of Theorem
6.7 to the single  $n \times n$ determinant in the
right-hand-side of equation (5.93) in Theorem 5.11:
$$\spreadlines{6 pt}\allowdisplaybreaks\align 
A = \ &-1,\quad B=4,\quad C=-2,\quad
D=\tfrac{1}{2},\quad  E=1,\quad F=2,\quad G=-1,
\tag 7.39\cr 
c_i=\ & 2i-1:=2\mu-1,\tag 7.40\cr
b_i=\ & 2(i-1),\quad 
\text{ for}\quad i=1,2,\cdots,n,\tag 7.41\cr
d_b=\ & d_c=2.\tag 7.42\cr
\endalign$$

It is immediate that
$$\spreadlines{6 pt}\allowdisplaybreaks\align 
&\ G-C = 1,\quad F-B=-2,\quad Bm_r+C=4m_r-2,\cr
&\ (D(Bm_r+C))^{d_b} = 
(D(Bm_r+C))^{d_c} = (2m_r-1)^2.\tag 7.43\cr
\endalign$$

In the above application of the determinant expansion of
(6.15) we also have 
$$p=n,\quad S=T=I_n,\quad \text{ and}\quad 
\ell_{\mu}=j_{\mu}=\mu,\quad 
\text{ for}\quad \mu=1,2,\cdots,n.\tag 7.44$$
Equations (7.40)--(7.42) and (7.44) imply that 
$$\lambda_i = 0\quad 
\text{ and}\quad
\nu_i = 0,\quad 
\text{ for}\ i=1,2,\cdots,n.\tag 7.45$$ 
Thus, we have 
$$s_{\lambda}(x) = s_{\nu}(x) = 
\det (D_{0,\emptyset ,\emptyset }) = 1.\tag 7.46$$

Keeping in mind (7.39)--(7.46), noting that the inner sum in
(6.15) consists of just one term corresponding to 
$T = I_n$, observing that 
$$1-2m_r + (2m_r-1)2y_r = (2m_r-1)(2y_r-1),\tag 7.47$$
relabelling $(2m_r-1)$ and $(2y_r-1)$ as $m_r$ (odd) 
and $y_r$ (odd), respectively,  
and simplifying, now yields (7.36).

The proof of (7.37) is similar to that of (7.36).  
The Lambert series $D_{2i+1}$ in (5.94), determined by
(5.96), is 
$$L_{2i+1}(r;-1,2,0,\tfrac{1}{2},1,1,0),\tag 7.48$$
as defined in (6.4).  Apply the following case of Theorem
6.7 to the single  $n \times n$ determinant in the
right-hand-side of equation (5.94) in Theorem 5.11:
$$\spreadlines{6 pt}\allowdisplaybreaks\align 
A = \ &-1,\quad B=2,\quad C=0,\quad
D=\tfrac{1}{2},\quad  E=1,\quad F=1,\quad G=0,
\tag 7.49\cr 
c_i=\ & 2i+1:=2\mu+1,\tag 7.50\cr
b_i=\ & 2(i-1),\quad 
\text{ for}\quad i=1,2,\cdots,n,\tag 7.51\cr
d_b=\ & d_c=2.\tag 7.52\cr
\endalign$$

It is immediate that
$$\spreadlines{6 pt}\allowdisplaybreaks\align 
&\ G-C = 0,\quad F-B=-1,\quad Bm_r+C=2m_r,\cr
&\ (D(Bm_r+C))^{d_b} = 
(D(Bm_r+C))^{d_c} = m_r^2.\tag 7.53\cr
\endalign$$

Keeping in mind (7.44)--(7.46) and (7.49)--(7.53), noting
that the inner sum in (6.15) consists of just one term
corresponding to $T = I_n$, observing that 
$-m_r + 2m_ry_r = m_r(2y_r-1)$, 
relabelling $(2y_r-1)$ as $y_r$ (odd), 
and simplifying, now yields (7.37).
\qed\enddemo

The same analysis involving (5.18) that took Theorem 5.11
into Corollary 5.12 transforms Theorem 7.5 into the
following Corollary.
\proclaim{Corollary 7.6} 
 Let $\vartriangle \kern-.3em (q)$ be defined by 
\hbox{\rm(5.17)}, and let $n=1,2,3,\cdots$.   We then have
$$\spreadlines{6 pt}\allowdisplaybreaks\align 
\vartriangle \kern-.3em (q^2)^{4n^2}=\,&
4^{-n(n-1)}q^{-n^2}\prod\limits_{r=1}^{2n-1}(r!)^{-1}
\kern-2.6em\sum\limits_{{y_1,\ldots ,y_n\geq 1
\hbox {\sixrm (odd)}}\atop 
{m_1>m_2>\cdots >m_n\geq 1 \hbox {\sixrm (odd)}}}
\kern-2.6em m_1m_2\cdots m_n\cr
&\cdot q^{m_1y_1+\cdots +m_ny_n}
\kern-.65em
\prod\limits_{1\leq r<s\leq n}\kern-1.1em
\left(m_r^2-m_s^2\right)^2,
\tag 7.54\cr
\kern -11.5 em \text{and}\kern 11.5 em\quad &\cr 
\vartriangle \kern-.3em (q)^{4n(n+1)}=\,&
2^{n}q^{-n(n+1)/2}
\prod\limits_{r=1}^{2n}(r!)^{-1}
\kern-1.65em\sum\limits_{{y_1,\ldots ,y_n\geq 1
\hbox {\sixrm (odd)}}\atop 
{m_1>m_2>\cdots >m_n\geq 1 }}
\kern-1.8em 
\left(m_1m_2\cdots m_n\right)^3\cr
&\cdot q^{m_1y_1+\cdots +m_ny_n}
\kern-.65em
\prod\limits_{1\leq r<s\leq n}\kern-1.1em
\left(m_r^2-m_s^2\right)^2.
\tag 7.55\cr
\endalign$$
\endproclaim 

\remark{Remark} Taking the coefficient of $q^N$ on both
sides of (7.54) and (7.55), respectively, immediately gives
the first and second Kac--Wakimoto conjectured identities for triangular
numbers in \cite{120, pp. 452}.  Our proof of Corollary 7.6 does not require
inclusion/exclusion, and the analysis involving Schur
functions is simpler than in (7.1), (7.11), (7.19), and
(7.28).  Zagier in \cite{253} has also recently independently
proven these two identities in \cite{120, pp. 452} by 
utilizing cusp forms.  In addition, he proved the
more general Conjecture 7.2 of Kac--Wakimoto \cite{120, pp. 451}, and its 
rewritten $m=2$ (first unproven) special case in \cite{120, pp. 451}. 
Information regarding affine superalgebras and Appell's generalized theta
functions can be found in \cite{121}.
\endremark

The $n=1$ cases of (7.54) and (7.55) are the classical
identities of Legendre \cite{139}, \cite{21, Eqns. (ii) and
(iii), pp. 139} given by 
\proclaim{Theorem 7.7 (Legendre)} Let 
$\vartriangle \kern-.3em (q)$ be defined by 
\hbox{\rm(5.17)}. Then  
$$\spreadlines{6 pt}\allowdisplaybreaks\align 
\vartriangle \kern-.3em (q^2)^{4} = &\ 
\sum\limits_{r=1}^{\infty}{(2r-1)q^{2(r-1)} \over 
{1-q^{2(2r-1)}}}=
\sum\limits_{{y_1\geq 1\hbox {\sixrm (odd)}}
\atop{m_1\geq 1 \hbox {\sixrm (odd)}}}
\kern-1.3em m_1\cdot
q^{m_1y_1-1},\tag 7.56\cr
\kern -9.2 em \text{and}\kern 9.2 em\quad 
\vartriangle \kern-.3em (q)^{8} = &\ 
\sum\limits_{r=1}^{\infty}{r^3q^{r-1} \over 
{1-q^{2r}}}\kern 1.63 em =
\sum\limits_{{y_1\geq 1\hbox {\sixrm (odd)}}
\atop{m_1\geq 1}}
\kern-1.3em m_1^3\cdot
q^{m_1y_1-1}.\tag 7.57\cr
\endalign$$
\endproclaim 
The Legendre identities involve the first sum over $r$ 
in (7.56) and (7.57).  The $n=1$ case of the
Kac--Wakimoto conjectured identities for triangular
numbers in \cite{120, pp. 452} takes the 
coefficient of $q^N$ on both sides of (7.56) and (7.57) 
while using the second double sum over $y_1$ and $m_1$.

We next consider the Schur function form of our Hankel
determinant identities related to 
$\cn(u,k)$ and $\dn(u,k)$.  

We begin with the following theorem related to $\cn(u,k)$.
\proclaim{Theorem 7.8 } Let 
$\vartheta_2 (0,q)$ and $\vartheta_3 (0,q)$ be defined by 
\hbox{\rm(1.2)} and \hbox{\rm(1.1)}, respectively. 
Let $n=1,2,3,\cdots$.    We then have
$$\spreadlines{6 pt}\allowdisplaybreaks\align 
&\vartheta_2 (0,q)^{2n^2}
\vartheta_3 (0,q)^{2n(n-1)}\tag 7.58\cr
=\,&4^{n}\prod\limits_{r=1}^{n-1}(2r)!^{-2}
\kern-2.6em\sum\limits_{{y_1,\ldots ,y_n\geq 1
\hbox {\sixrm (odd)}}\atop 
{m_1>m_2>\cdots >m_n\geq 1 \hbox {\sixrm (odd)}}}
\kern-2.6em (-1)^{\tfrac{1}{2}(-n+y_1+\cdots+y_n)}\cr
&\cdot q^{\tfrac{1}{2}(m_1y_1+\cdots +m_ny_n)}
\kern-.65em
\prod\limits_{1\leq r<s\leq n}\kern-1.1em
\left(m_r^2-m_s^2\right)^2.\cr
\endalign$$
\endproclaim 
\demo{Proof} We apply Theorem 6.7 to Theorem 5.13.  
The Lambert series $T_{2i-2}$ in (5.109), determined by
(5.110), is 
$$L_{2i-2}(r;1,2,-1,1,1,1,-\tfrac{1}{2}),\tag 7.59$$
as defined in (6.4).  Apply the following case of Theorem
6.7 to the single  $n \times n$ determinant in the
right-hand-side of equation (5.109) in Theorem 5.13:
$$\spreadlines{6 pt}\allowdisplaybreaks\align 
A = \ &1,\quad B=2,\quad C=-1,\quad
D=1,\quad  E=1,\quad F=1,\quad G=-\tfrac{1}{2},
\tag 7.60\cr 
c_i=\ & 2i-2:=2\mu-2,\tag 7.61\cr
b_i=\ & 2(i-1),\quad 
\text{ for}\quad i=1,2,\cdots,n,\tag 7.62\cr
d_b=\ & d_c=2.\tag 7.63\cr
\endalign$$

It is immediate that
$$\spreadlines{6 pt}\allowdisplaybreaks\align 
&\ G-C = \tfrac{1}{2},\quad F-B=-1,\quad
Bm_r+C=2m_r-1,\cr 
&\ (D(Bm_r+C))^{d_b} = 
(D(Bm_r+C))^{d_c} = (2m_r-1)^2.\tag 7.64\cr
\endalign$$

In the above application of the determinant expansion of
(6.15) we also have (7.44)--(7.46).  
Keeping in mind (7.44)--(7.46) and (7.60)--(7.64), noting
that the inner sum in (6.15) consists of just one term
corresponding to $T = I_n$, observing that (7.47) holds, 
relabelling $(2m_r-1)$ and $(2y_r-1)$ as $m_r$ (odd) 
and $y_r$ (odd), respectively,  
and simplifying, now yields (7.58).
\qed\enddemo

Applying the relation (5.105) to the left-hand-side of
(7.58) immediately implies that 
$$\spreadlines{6 pt}\allowdisplaybreaks\align 
&\vartheta_2 (0,q)^{2n}
\vartheta_2 (0,q^{1/2})^{4n(n-1)}\tag 7.65\cr
=\,&4^{n^2}\prod\limits_{r=1}^{n-1}(2r)!^{-2}
\kern-2.6em\sum\limits_{{y_1,\ldots ,y_n\geq 1
\hbox {\sixrm (odd)}}\atop 
{m_1>m_2>\cdots >m_n\geq 1 \hbox {\sixrm (odd)}}}
\kern-2.6em (-1)^{\tfrac{1}{2}(-n+y_1+\cdots+y_n)}\cr
&\cdot q^{\tfrac{1}{2}(m_1y_1+\cdots +m_ny_n)}
\kern-.65em
\prod\limits_{1\leq r<s\leq n}\kern-1.1em
\left(m_r^2-m_s^2\right)^2.\cr
\endalign$$

Keeping in mind (7.65) and (5.18), the $n=1$ case of
Theorem 7.8 is equivalent to the identity in \cite{21, Ex.
(iv), pp. 139}.  Similarly, the $n=2$ case of Theorem 7.8
leads to 
$$\spreadlines{6 pt}\allowdisplaybreaks\align 
&q\vartriangle \kern-.3em (q^2)^{4}
\cdot q\vartriangle \kern-.3em (q)^{8}=\,
\tfrac{1}{64}\left[T_0T_4-T_2^2\right]
\tag 7.66a\cr
=\,&2^{-6}\kern-2.0em\sum\limits_{{y_1,y_2\geq1}
\atop{m_1>m_2\geq 1 \hbox {\sixrm (odd)}}}
\kern-1.5em (-1)^{y_1+y_2}
\left(m_1^2-m_2^2\right)^2
q^{m_1(y_1-\tfrac{1}{2})+m_2(y_2-\tfrac{1}{2})},
\tag 7.66b\cr
\endalign$$
where $T_{2i-2}$ is defined by (5.110).  

Note that Theorem 7.7 and the product 
$q\vartriangle \kern-.3em (q^2)^{4}
\cdot q\vartriangle \kern-.3em (q)^{8}$ in (7.66a) 
implies that the $4$-fold sum in (7.66b) is the product of
the double sums in (7.56) and (7.57), each multiplied by
$q$.   That is, (7.66b) is the product of the sums in 
\cite{21, Eqns. (ii) and (iii), pp. 139}. The double sums in
(7.56) and (7.57) have no alternating signs, while the sum in
(7.66b) does.    

We next have the following theorem.
\proclaim{Theorem 7.9 }  Let 
$\vartheta_2 (0,q)$ and $\vartheta_3 (0,q)$ be defined by 
\hbox{\rm(1.2)} and \hbox{\rm(1.1)}, respectively. 
Let $n=1,2,3,\cdots$.    We then have
$$\spreadlines{6 pt}\allowdisplaybreaks\align 
&\vartheta_2 (0,q)^{2n^2}
\vartheta_3 (0,q)^{2n(n+1)}\tag 7.67\cr
=\,&4^{n}\prod\limits_{r=1}^{n}(2r-1)!^{-2}
\kern-2.6em\sum\limits_{{y_1,\ldots ,y_n\geq 1
\hbox {\sixrm (odd)}}\atop 
{m_1>m_2>\cdots >m_n\geq 1 \hbox {\sixrm (odd)}}}
\kern-2.6em (-1)^{\tfrac{1}{2}(-n+y_1+\cdots+y_n)}
\cr
&\cdot \left(m_1m_2\cdots m_n\right)^2
q^{\tfrac{1}{2}(m_1y_1+\cdots +m_ny_n)}
\kern-.65em
\prod\limits_{1\leq r<s\leq n}\kern-1.1em
\left(m_r^2-m_s^2\right)^2.\cr
\endalign$$
\endproclaim 
\demo{Proof} We apply Theorem 6.7 to Theorem 5.14.  
The Lambert series $T_{2i}$ in (5.117), determined by
(5.110), is 
$$L_{2i}(r;1,2,-1,1,1,1,-\tfrac{1}{2}),\tag 7.68$$
as defined in (6.4).  Apply the following case of Theorem
6.7 to the single  $n \times n$ determinant in the
right-hand-side of equation (5.117) in Theorem 5.14:
$$\spreadlines{6 pt}\allowdisplaybreaks\align 
A = \ &1,\quad B=2,\quad C=-1,\quad
D=1,\quad  E=1,\quad F=1,\quad G=-\tfrac{1}{2},
\tag 7.69\cr 
c_i=\ & 2i:=2\mu,\tag 7.70\cr
b_i=\ & 2(i-1),\quad 
\text{ for}\quad i=1,2,\cdots,n,\tag 7.71\cr
d_b=\ & d_c=2.\tag 7.72\cr
\endalign$$

It is immediate that (7.64) holds. 
In the above application of the determinant expansion of
(6.15) we also have (7.44)--(7.46).  
Keeping in mind (7.44)--(7.46), (7.64), and (7.69)--(7.72), 
noting that the inner sum in (6.15) consists of just one
term corresponding to $T = I_n$, observing that (7.47)
holds,  relabelling $(2m_r-1)$ and $(2y_r-1)$ as $m_r$
(odd)  and $y_r$ (odd), respectively,  
and simplifying, now yields (7.67).
\qed\enddemo

Applying (5.18) and (5.105) to the left-hand-side of 
(7.67) immediately gives the Schur function form of
Corollary 5.15 in the following corollary.
\proclaim{Corollary 7.10 } Let 
$\vartheta_3 (0,q)$ and $\vartriangle \kern-.3em (q)$ be
defined by \hbox{\rm(1.1)} and \hbox{\rm(5.17)},
respectively.  Let $n=1,2,3,\cdots$.    We then have
$$\spreadlines{6 pt}\allowdisplaybreaks\align 
&\vartheta_3 (0,q)^{2n}
\vartriangle \kern -.3em (q)^{(2n)^2}\tag 7.73\cr
=\,&4^{-n(n-1)}q^{-n^2/2}
\prod\limits_{r=1}^{n}(2r-1)!^{-2}\cr
&\cdot \sum\limits_{{y_1,\ldots ,y_n\geq 1
\hbox {\sixrm (odd)}}\atop 
{m_1>m_2>\cdots >m_n\geq 1 \hbox {\sixrm (odd)}}}
\kern-2.6em (-1)^{\tfrac{1}{2}(-n+y_1+\cdots+y_n)}
\left(m_1m_2\cdots m_n\right)^2\cr
&\cdot q^{\tfrac{1}{2}(m_1y_1+\cdots +m_ny_n)}
\kern-.65em
\prod\limits_{1\leq r<s\leq n}\kern-1.1em
\left(m_r^2-m_s^2\right)^2.\cr
\endalign$$
\endproclaim 

We now consider the Schur function form of the theta
function identities related to $\dn(u,k)$.  We start with the
following theorem.  
\proclaim{Theorem 7.11 } Let 
$\vartheta_2 (0,q)$ and $\vartheta_3 (0,q)$ be defined by 
\hbox{\rm(1.2)} and \hbox{\rm(1.1)}, respectively. 
Let $n=1,2,3,\cdots$.    We then have
$$\spreadlines{6 pt}\allowdisplaybreaks\align 
&\vartheta_2 (0,q)^{2n(n-1)}
\vartheta_3 (0,q)^{2n^2}\cr
=\,&4^{n^2}
\prod\limits_{r=1}^{n-1}(2r)!^{-2} 
\kern-2.0em
\sum\limits_{{y_1,\ldots ,y_n\geq 1
\hbox {\sixrm (odd)}}\atop 
{m_1>m_2>\cdots >m_n\geq 1 }}
\kern-2.0em (-1)^{\tfrac{1}{2}(-n+y_1+\cdots+y_n)}
\cr
&\cdot q^{m_1y_1+\cdots +m_ny_n}
\kern-.65em
\prod\limits_{1\leq r<s\leq n}\kern-1.1em
\left(m_r^2-m_s^2\right)^2
\tag 7.74a\cr
+&4^{n^2-1}
\prod\limits_{r=1}^{n-1}(2r)!^{-2}
\kern-2.3em
\sum\limits_{{y_1,\ldots ,y_{n-1}\geq 1
\hbox {\sixrm (odd)}}\atop 
{m_1>m_2>\cdots >m_{n-1}\geq 1 }}
\kern-2.3em (-1)^{\tfrac{1}{2}
(-n+1+y_1+\cdots+y_{n-1})}\cr
&\cdot q^{m_1y_1+\cdots +m_{n-1}y_{n-1}} 
\left(m_1m_2\cdots m_{n-1}\right)^4
\kern-.65em
\prod\limits_{1\leq r<s\leq {n-1}}\kern-1.1em
\left(m_r^2-m_s^2\right)^2,
\tag 7.74b\cr
\endalign$$
where the term in  \hbox{\rm(7.74b)} is defined to be $1$
if $n=1$.
\endproclaim 
\demo{Proof} We apply Theorem 6.7 to Theorem 5.16.  
The Lambert series $N_{2i-2}$ in (5.127a), determined by
(5.128), is 
$$L_{2i-2}(r;1,2,0,\tfrac{1}{2},1,1,0),\tag 7.75$$
as defined in (6.4).  Apply the following case of Theorem
6.7 to the single  $n \times n$ determinant in the
right-hand-side of equation (5.127a) in Theorem 5.16:
$$\spreadlines{6 pt}\allowdisplaybreaks\align 
A = \ &1,\quad B=2,\quad C=0,\quad
D=\tfrac{1}{2},\quad  E=1,\quad F=1,\quad G=0,
\tag 7.76\cr 
c_i=\ & 2i-2:=2\mu-2,\tag 7.77\cr
b_i=\ & 2(i-1),\quad 
\text{ for}\quad i=1,2,\cdots,n,\tag 7.78\cr
d_b=\ & d_c=2.\tag 7.79\cr
\endalign$$

It is immediate that (7.53) holds.  In the above application 
of the determinant expansion of
(6.15) we also have (7.44)--(7.46).  
Keeping in mind (7.44)--(7.46), (7.53),  and (7.76)--(7.79),
noting that the inner sum in (6.15) consists of just one term
corresponding to $T = I_n$, observing that 
$-m_r + 2m_ry_r = m_r(2y_r-1)$, 
relabelling $(2y_r-1)$ as $y_r$ (odd), 
and simplifying, now yields the term in (7.74a).

The derivation of the term in (7.74b) from (5.127b) is the
same as for (7.74a), except that 
$n\mapsto n-1$, $c_i\mapsto  2i+2$, and 
$c_{\ell_1}+b_{j_1}\mapsto 4$.
\qed\enddemo

The $n=1$ case of Theorem 7.11 is equivalent to 
Jacobi's $2$-squares identity in 
\cite{21, Entry 8(i),  pp. 114}.  The calculation is similar to
that for the $n=1$ case of Theorem 5.16.  Here, just note
that the $n=1$ case of the sum in (7.74b) equals $1$.  

We next have the following theorem.  
\proclaim{Theorem 7.12 } Let 
$\vartheta_2 (0,q)$ and $\vartheta_3 (0,q)$ be defined by 
\hbox{\rm(1.2)} and \hbox{\rm(1.1)}, respectively. 
Let $n=1,2,3,\cdots$.    We then have
$$\spreadlines{6 pt}\allowdisplaybreaks\align 
&\vartheta_2 (0,q)^{2n(n+1)}
\vartheta_3 (0,q)^{2n^2}\tag 7.80\cr
=\,&4^{n(n+1)}
\prod\limits_{r=1}^{n}(2r-1)!^{-2}\cr 
&\cdot 
\sum\limits_{{y_1,\ldots ,y_n\geq 1
\hbox {\sixrm (odd)}}\atop 
{m_1>m_2>\cdots >m_n\geq 1 }}
\kern-2.0em (-1)^{\tfrac{1}{2}(-n+y_1+\cdots+y_n)}
\left(m_1m_2\cdots m_n\right)^2\cr
&\cdot q^{m_1y_1+\cdots +m_ny_n}
\kern-.65em
\prod\limits_{1\leq r<s\leq n}\kern-1.1em
\left(m_r^2-m_s^2\right)^2.\cr
\endalign$$
\endproclaim 
\demo{Proof} We apply Theorem 6.7 to Theorem 5.17.  
The Lambert series $N_{2i}$ in (5.135), determined by
(5.128), is 
$$L_{2i}(r;1,2,0,\tfrac{1}{2},1,1,0),\tag 7.81$$
as defined in (6.4).  Apply the following case of Theorem
6.7 to the single  $n \times n$ determinant in the
right-hand-side of equation (5.135) in Theorem 5.17:
$$\spreadlines{6 pt}\allowdisplaybreaks\align 
A = \ &1,\quad B=2,\quad C=0,\quad
D=\tfrac{1}{2},\quad  E=1,\quad F=1,\quad G=0,
\tag 7.82\cr 
c_i=\ & 2i:=2\mu,\tag 7.83\cr
b_i=\ & 2(i-1),\quad 
\text{ for}\quad i=1,2,\cdots,n,\tag 7.84\cr
d_b=\ & d_c=2.\tag 7.85\cr
\endalign$$

It is immediate that (7.53) holds.  In the above application 
of the determinant expansion of
(6.15) we also have (7.44)--(7.46).  
Keeping in mind (7.44)--(7.46), (7.53),  and (7.82)--(7.85),
noting that the inner sum in (6.15) consists of just one term
corresponding to $T = I_n$, observing that 
$-m_r + 2m_ry_r = m_r(2y_r-1)$, 
relabelling $(2y_r-1)$ as $y_r$ (odd), 
and simplifying, now yields (7.80).
\qed\enddemo

Applying (5.18) and (5.105) to the left-hand-side of 
(7.80) immediately gives the Schur function form of
Corollary 5.18 in the following corollary.
\proclaim{Corollary 7.13 }  Let 
$\vartriangle \kern-.3em (q)$ be
defined by \hbox{\rm(5.17)}.
 Let $n=1,2,3,\cdots$.    We then have
$$\spreadlines{6 pt}\allowdisplaybreaks\align 
&\vartriangle \kern-.3em (q)^{(2n)^2}
\vartriangle \kern-.3em (q^2)^{2n}\tag 7.86\cr
=\,&q^{-n(n+1)/2}
\prod\limits_{r=1}^{n}(2r-1)!^{-2}\cr 
&\cdot 
\sum\limits_{{y_1,\ldots ,y_n\geq 1
\hbox {\sixrm (odd)}}\atop 
{m_1>m_2>\cdots >m_n\geq 1 }}
\kern-2.0em (-1)^{\tfrac{1}{2}(-n+y_1+\cdots+y_n)}
\left(m_1m_2\cdots m_n\right)^2\cr
&\cdot q^{m_1y_1+\cdots +m_ny_n}
\kern-.65em
\prod\limits_{1\leq r<s\leq n}\kern-1.1em
\left(m_r^2-m_s^2\right)^2.\cr
\endalign$$
\endproclaim 

We next survey the Schur function form of our $\chi_n$
determinant identities.  We find it convenient to recall the 
definition of $\chi (A)$ given by 
$$\chi (A):=1,\ \text{ if}\ A\ 
\text{ is true, and}\ 0,\ 
\text{ otherwise.}\tag 7.87$$

We first have the following theorem.
\proclaim{Theorem 7.14} 
Let $n=1,2,3,\cdots$.  Then 
$$\spreadlines{6 pt}\allowdisplaybreaks\align 
&\vartheta_3(0,-q)^{4n^2}
\left[1+24\sum\limits_{r=1}^{\infty}
{rq^{r} \over {1+q^{r}}}\right]
\tag 7.88\cr
=\,&1+\sum\limits_{p=1}^n
(-1)^{p-1}2^{2n^2+n+1}\tfrac{3}{n(4n^2-1)}
\prod\limits_{r=1}^{2n-1}(r!)^{-1}
\kern-1.8em
\sum\limits_{{y_1,\ldots,y_p\geq 1}\atop 
{m_1>m_2>\cdots>m_p\geq 1}}
\kern-1.65em(-1)^{y_1+\cdots+y_p}\cr
&\cdot(-1)^{m_1+\cdots+m_p}
q^{m_1y_1+\cdots+m_py_p}\kern-.65em
\prod\limits_{1\leq r<s\leq p}\kern-1.1em
\left(m_r^2-m_s^2\right)^2\cr
&\cdot(m_1m_2\cdots m_p)\kern-.65em 
\sum\limits_{{\emptyset \subset S, T\subseteq I_n}
\atop {\Vert S\Vert=\Vert T\Vert =p}}\kern-.5em 
(-1)^{\Sigma (S)+\Sigma (T)}\cdot 
\det (\overline D_{n-p,S^c,T^c})\cr
&\cdot (m_1m_2\cdots m_p)^
{2\ell_1+2j_1-4+2\chi (j_1=n)}
s_{\lambda}(m_1^2,\ldots,m_p^2)\,
s_{\nu}(m_1^2,\ldots,m_p^2),\cr 
\endalign$$
where $\vartheta_3(0,-q)$ is determined by 
\hbox{\rm(1.1)}, 
the sets $S$, $S^c$, $T$, $T^c$ are given by
\hbox{\rm(6.1)-(6.2)}, 
$\Sigma (S)$ and $\Sigma (T)$ by \hbox{\rm(6.3)}, 
and the $(n-p) \times (n-p)$ matrix  
$$\overline D_{n-p,S^c,T^c}:= 
\left[c_{(\ell_{p+r}+j_{p+s}-1+\chi (j_{p+s}=n))}
\right]_{1\leq r,s\leq n-p},\tag 7.89$$
where the $c_i$ are determined by \hbox{\rm(5.23)}, with 
the $B_{2i}$ in \hbox{\rm(2.60)}, and $\chi (A)$ 
defined by (7.87). 
Finally, $s_\lambda$ and $s_\nu$ are the Schur functions
in \hbox{\rm(6.12)}, with the partitions $\lambda$ and
$\nu$ given by 
$$\spreadlines{6 pt}\allowdisplaybreaks\alignat 2 
&\lambda_i:=\ell_{p-i+1}-\ell_1+i-p,
&\quad\text{ for}\quad i=1,2,\cdots,p,\tag 7.90\cr 
\kern -6.15 em \text{and,}\kern 6.15 em 
&\nu_1:=0,&\quad\text{ if}\quad p=1,\tag 7.91a\cr
&\nu_1:=j_{p}-j_1+1-p+\chi (j_{p}=n),&\quad 
\text{ if}\quad p>1,\tag 7.91b\cr
&\nu_i:=j_{p-i+1}-j_1+i-p,&\quad 
\text{ if}\quad 2\leq i\leq p,\ \text{ and}, 
\ p>1. \tag 7.91c\cr
\endalignat$$
\endproclaim
\demo{Proof} We apply Theorem 6.19 to Theorem 5.24.  
We utilize the same specializations as in (7.4)--(7.9),
except that 
$$b_i = 2(i-1)+2\chi (i=n),\quad 
\text{ for}\quad i=1,2,\cdots,n,\tag 7.92$$
where $\chi (A)$ is defined by (7.87). 

It is immediate that (7.10) holds.  
The $\lambda_i$ and $\nu_i$ are given by (7.90)--(7.91).  
(Note that the formula for $\nu_i$ in (7.3) is transformed
to (7.91a--c) by adding $1$ to the largest part of $\nu$ 
if the largest element of the set $T$ is $n$).

Interchanging the inner sum in (5.192) with the double
multiple sum over ${y_1,\ldots,y_p\geq 1}$ and 
${m_1>m_2>\cdots>m_p\geq 1}$ in (6.57), multiplying
and dividing by $(m_1m_2\cdots m_p)$, 
and simplifying, now yields (7.88), for $n\geq 2$. 

By (5.194), the $n=1$ case of (7.88) is established by 
showing that this case of the multiple sum in (7.88) is 
$U_3$, as defined by (5.22).  This equality follows
termwise by using (6.5) to compute the $i=2$ case of
(7.4).
\qed\enddemo

We next have the following theorem.
\proclaim{Theorem 7.15} 
Let $n=1,2,3,\cdots$.  Then 
$$\spreadlines{6 pt}\allowdisplaybreaks\align 
&\vartheta_3(0,-q)^{4n(n+1)}
\left[1+24\sum\limits_{r=1}^{\infty}
{rq^{r} \over {1+q^{r}}}\right]
\tag 7.93\cr
=\,&1+\sum\limits_{p=1}^n
(-1)^{n-p+1}2^{2n^2+3n}\tfrac{3}{n(n+1)(2n+1)}
\prod\limits_{r=1}^{2n}(r!)^{-1}
\kern-1.8em
\sum\limits_{{y_1,\ldots,y_p\geq 1}\atop 
{m_1>m_2>\cdots>m_p\geq 1}}
\kern-1.65em(-1)^{m_1+\cdots+m_p}\cr
&\cdot q^{m_1y_1+\cdots+m_py_p}\,
(m_1m_2\cdots m_p)^3\kern-.65em
\prod\limits_{1\leq r<s\leq p}\kern-1.1em
\left(m_r^2-m_s^2\right)^2\cr
&\cdot
\sum\limits_{{\emptyset \subset S, T\subseteq I_n}
\atop {\Vert S\Vert=\Vert T\Vert =p}}\kern-.5em 
(-1)^{\Sigma (S)+\Sigma (T)}\cdot 
\det (\overline D_{n-p,S^c,T^c})\cr 
&\cdot (m_1m_2\cdots m_p)^
{2\ell_1+2j_1-4+2\chi (j_1=n)} 
s_{\lambda}(m_1^2,\ldots,m_p^2)\,
s_{\nu}(m_1^2,\ldots,m_p^2),\cr 
\endalign$$
where the same assumptions hold as in 
\hbox{\rm Theorem 7.14}, except that the 
$(n-p) \times (n-p)$ matrix  
$$\overline D_{n-p,S^c,T^c}:= 
\left[a_{(\ell_{p+r}+j_{p+s}-1+\chi (j_{p+s}=n))}
\right]_{1\leq r,s\leq n-p},\tag 7.94$$
where the $a_i$ are determined by \hbox{\rm(5.40)}, with 
the $B_{2i+2}$ in \hbox{\rm(2.60)}, and $\chi (A)$ 
defined by (7.87). 
\endproclaim
\demo{Proof} We apply Theorem 6.19 to Theorem 5.26.  
We utilize the same specializations as in (7.13)--(7.18),
except that $b_i$ is given by (7.92).  

It is immediate that (7.10) holds.  
The $\lambda_i$ and $\nu_i$ are given by (7.90)--(7.91).  

Interchanging the inner sum in (5.207) with the double
multiple sum over ${y_1,\ldots,y_p\geq 1}$ and 
${m_1>m_2>\cdots>m_p\geq 1}$ in (6.57), multiplying
and dividing by $(m_1m_2\cdots m_p)^3$, 
and simplifying, now yields (7.93), for $n\geq 2$. 

By (5.209), the $n=1$ case of (7.93) is established by 
showing that this case of the multiple sum in (7.93) is 
$G_5$, as defined by (5.39).  This equality follows
termwise by using (6.5) to compute the $i=2$ case of
(7.13).
\qed\enddemo

The next two theorems require the Euler numbers $E_n$
defined by (2.61).  We first have the following theorem.
\proclaim{Theorem 7.16 } Let
$n=1,2,3,\cdots$.  Then
$$\spreadlines{6 pt}\allowdisplaybreaks\align 
&\vartheta_3(0,q)^{2n(n-1)}\vartheta_3(0,-q)^{2n^2}\cr
&\cdot\left\{2n\left[2+24\sum\limits_{r=1}^{\infty}
{2rq^{2r} \over {1+q^{2r}}}\right]-
\left[1-24\sum\limits_{r=1}^{\infty}
{(2r-1)q^{2r-1} \over {1+q^{2r-1}}}\right]\right\}
\tag 7.95\cr
=\,&(4n-1)+\sum\limits_{p=1}^n
(-1)^{p-1}2^{2n}\tfrac{3}{n(2n-1)}
\prod\limits_{r=1}^{n-1}(2r)!^{-2}\kern-2.3em
\sum\limits_{{y_1,\ldots,y_p\geq 1}\atop 
{m_1>m_2>\cdots>m_p\geq 1}\sixrm (odd)}
\kern-2.15em(-1)^{y_1+\cdots+y_p}\cr 
&\cdot (-1)^{\tfrac{1}{2}(p+m_1+\cdots+m_p)}
q^{m_1y_1+\cdots+m_py_p}\kern-.65em
\prod\limits_{1\leq r<s\leq p}\kern-1.1em
\left(m_r^2-m_s^2\right)^2\cr
&\cdot
\sum\limits_{{\emptyset \subset S, T\subseteq I_n}
\atop {\Vert S\Vert=\Vert T\Vert =p}}\kern-.5em 
(-1)^{\Sigma (S)+\Sigma (T)}\cdot 
\det (\overline D_{n-p,S^c,T^c})\cr 
&\cdot (m_1m_2\cdots m_p)^{2\ell_1+2j_1-4
+2\chi (j_1=n)} 
s_{\lambda}(m_1^2,\ldots,m_p^2)\,
s_{\nu}(m_1^2,\ldots,m_p^2),\cr 
\endalign$$
where the same assumptions hold as in 
\hbox{\rm Theorem 7.14}, except that the 
$(n-p) \times (n-p)$ matrix  
$$\overline D_{n-p,S^c,T^c}:= 
\left[b_{(\ell_{p+r}+j_{p+s}-1
+\chi (j_{p+s}=n))}\right]_{1\leq r,s\leq n-p},\tag 7.96$$
where the $b_i$ are defined by \hbox{\rm(5.57)}, with 
the $E_{2i-2}$ in \hbox{\rm(2.61)}, and $\chi (A)$ 
defined by (7.87). 
\endproclaim
\demo{Proof} We apply Theorem 6.19 to Theorem 5.28.  
We utilize the same specializations as in (7.21)--(7.26),
except that $b_i$ is given by (7.92).  

It is immediate that (7.27) holds.  
The $\lambda_i$ and $\nu_i$ are given by (7.90)--(7.91).  

Interchanging the inner sum in (5.223) with the double
multiple sum over ${y_1,\ldots,y_p\geq 1}$ and 
${m_1>m_2>\cdots>m_p\geq 1}$ in (6.57), relabelling 
$(2m_r-1)$ as $m_r$ (odd), 
and simplifying, now yields (7.95), for $n\geq 2$. 

By (5.225), the $n=1$ case of (7.95) is established by 
showing that this case of the multiple sum in (7.95) is 
$R_2$, as defined by (5.56).  This equality follows
termwise by using (6.5) to compute the $i=2$ case of
(7.21).
\qed\enddemo

 We next have the following theorem.
\proclaim{Theorem 7.17 } Let
$n=1,2,3,\cdots$.  Then
$$\spreadlines{6 pt}\allowdisplaybreaks\align 
&\vartheta_3(0,q)^{2n(n+1)}\vartheta_3(0,-q)^{2n^2}\cr
&\cdot\left\{2n\left[2+24\sum\limits_{r=1}^{\infty}
{2rq^{2r} \over {1+q^{2r}}}\right]+
\left[1-24\sum\limits_{r=1}^{\infty}
{(2r-1)q^{2r-1} \over {1+q^{2r-1}}}\right]\right\}
\tag 7.97\cr
=\,&(4n+1)+\sum\limits_{p=1}^n
(-1)^{n-p+1}2^{2n}\tfrac{3}{n(2n+1)}
\prod\limits_{r=1}^{n}(2r-1)!^{-2}
\kern-2.3em
\sum\limits_{{y_1,\ldots,y_p\geq 1}\atop 
{m_1>m_2>\cdots>m_p\geq 1}\sixrm (odd)}
\kern-2.15em(-1)^{y_1+\cdots+y_p}\cr 
&\cdot (-1)^{\tfrac{1}{2}(p+m_1+\cdots+m_p)}
q^{m_1y_1+\cdots+m_py_p}\kern-.65em
\prod\limits_{1\leq r<s\leq p}\kern-1.1em
\left(m_r^2-m_s^2\right)^2\cr
&\cdot(m_1m_2\cdots m_p)^2\kern-.65em 
\sum\limits_{{\emptyset \subset S, T\subseteq I_n}
\atop {\Vert S\Vert=\Vert T\Vert =p}}\kern-.5em 
(-1)^{\Sigma (S)+\Sigma (T)}\cdot 
\det (\overline D_{n-p,S^c,T^c})\cr 
&\cdot (m_1m_2\cdots m_p)^{2\ell_1+2j_1-4
+2\chi (j_1=n)} 
s_{\lambda}(m_1^2,\ldots,m_p^2)\,
s_{\nu}(m_1^2,\ldots,m_p^2),\cr 
\endalign$$
where the same assumptions hold as in 
\hbox{\rm Theorem 7.14}, except that the 
$(n-p) \times (n-p)$ matrix  
$$\overline D_{n-p,S^c,T^c}:= 
\left[b_{(\ell_{p+r}+j_{p+s}
+\chi (j_{p+s}=n))}\right]_{1\leq r,s\leq n-p},\tag 7.98$$
where the $b_i$ are determined by \hbox{\rm(5.57)}, with 
the $E_{2i-2}$ in \hbox{\rm(2.61)}, and $\chi (A)$ 
defined by (7.87). 
\endproclaim
\demo{Proof} We apply Theorem 6.19 to Theorem 5.30.  
We utilize the same specializations as in (7.30)--(7.35),
except that $b_i$ is given by (7.92).  

It is immediate that (7.27) holds.  
The $\lambda_i$ and $\nu_i$ are given by (7.90)--(7.91).  

Interchanging the inner sum in (5.239) with the double
multiple sum over ${y_1,\ldots,y_p\geq 1}$ and 
${m_1>m_2>\cdots>m_p\geq 1}$ in (6.57), relabelling 
$(2m_r-1)$ as $m_r$ (odd), multiplying
and dividing by $(m_1m_2\cdots m_p)^2$, 
and simplifying, now yields (7.97), for $n\geq 2$. 

By (5.241), the $n=1$ case of (7.97) is established by 
showing that this case of the multiple sum in (7.97) is 
$R_4$, as defined by (5.56).  This equality follows
termwise by using (6.5) to compute the $i=2$ case of
(7.30).
\qed\enddemo

The analog of Theorem 7.5 is given by the following
theorem.
\proclaim{Theorem 7.18 }  Let 
$\vartheta_2 (0,q)$ be defined by 
\hbox{\rm(1.2)}, and let $n=1,2,3,\cdots$.   
We then have
$$\spreadlines{6 pt}\allowdisplaybreaks\align 
&\vartheta_2 (0,q)^{4n^2}
\left[1+24\sum\limits_{r=1}^{\infty}
{rq^{2r} \over {1+q^{2r}}}\right]\cr
=\,&4^{n(n+1)}\tfrac{3}{n(4n^2-1)}
\prod\limits_{r=1}^{2n-1}(r!)^{-1}
\kern-2.6em\sum\limits_{{y_1,\ldots ,y_n\geq 1
\hbox {\sixrm (odd)}}\atop 
{m_1>m_2>\cdots >m_n\geq 1 \hbox {\sixrm (odd)}}}
\kern-2.6em m_1m_2\cdots m_n\cr
&\cdot \left(m_1^2+\cdots +m_n^2\right)
q^{m_1y_1+\cdots +m_ny_n}
\kern-.65em
\prod\limits_{1\leq r<s\leq n}\kern-1.1em
\left(m_r^2-m_s^2\right)^2.
\tag 7.99\cr
\kern -8.35 em \text{and}\kern 8.35 em\quad &\cr
&\vartheta_2 (0,q^{1/2})^{4n(n+1)}
\left[1+24\sum\limits_{r=1}^{\infty}
{rq^{r} \over {1+q^{r}}}\right]\cr
=\,&2^{n(4n+5)}\tfrac{6}{n(n+1)(2n+1)}
\prod\limits_{r=1}^{2n}(r!)^{-1}
\kern-1.65em\sum\limits_{{y_1,\ldots ,y_n\geq 1
\hbox {\sixrm (odd)}}\atop 
{m_1>m_2>\cdots >m_n\geq 1 }}
\kern-1.8em 
\left(m_1m_2\cdots m_n\right)^3\cr
&\cdot \left(m_1^2+\cdots +m_n^2\right)
q^{m_1y_1+\cdots +m_ny_n}
\kern-.65em
\prod\limits_{1\leq r<s\leq n}\kern-1.1em
\left(m_r^2-m_s^2\right)^2.
\tag 7.100\cr
\endalign$$
\endproclaim 
\demo{Proof} We apply Theorem 6.19 to Theorems 5.31
and 5.32.   We take $n\geq 2$ and first consider (7.99). 
We utilize the same specializations as in (7.38)--(7.42),
except that $b_i$ is given by (7.92).  It is immediate that
(7.43) holds.  

In the above application of the determinant expansion of
(6.57) we also have (7.44).   
Equations (7.40), (7.42), (7.44), and (7.92) imply that 
$$\lambda_i = 0,\quad 
\text{ for}\ i=1,2,\cdots,n,\tag 7.101$$ 
$$\nu_1 = 1,\quad 
\text{ and}\quad
\nu_i = 0,\quad 
\text{ for}\ i=2,3,\cdots,n.\tag 7.102$$ 
Thus, we have 
$$s_{\lambda}(x) = 
\det (\overline D_{0,\emptyset ,\emptyset }) 
= 1,\tag 7.103$$
$$s_{\nu}(x_1,\ldots,x_n) = x_1+\cdots +x_n
.\tag 7.104$$

Keeping in mind (7.39), (7.40), (7.42), (7.43),
(7.44), (7.92), (7.101)--(7.104), 
noting that the inner sum in (6.57) consists of just one
term corresponding to $T = I_n$, observing (7.47), 
relabelling $(2m_r-1)$ and $(2y_r-1)$ as $m_r$ (odd) 
and $y_r$ (odd), respectively, and simplifying, now yields
(7.99), for $n\geq 2$.

By Theorem 5.31, the $n=1$ case of (7.99) is established
by showing that this case of the multiple sum in (7.99) is 
$C_3$, as defined by (5.95).  This equality follows
termwise by using (6.5) to compute the $i=2$ case of
(7.38).

The proof of (7.100) is similar to that of (7.99). We take 
$n\geq 2$ and apply Theorem 6.19 to Theorem 5.32.  
We utilize the same specializations as in (7.48)--(7.52),
except that $b_i$ is given by (7.92).  It is immediate that
(7.53) holds.  

Keeping in mind (7.44), (7.49), (7.50), (7.52),
(7.53), (7.92), (7.101)--(7.104), 
noting that the inner sum in (6.57) consists of just one
term corresponding to $T = I_n$, observing that 
$-m_r + 2m_ry_r = m_r(2y_r-1)$, 
relabelling $(2y_r-1)$ as $y_r$ (odd), 
and simplifying, now yields (7.100), for $n\geq 2$.

By Theorem 5.32, the $n=1$ case of (7.100) is established
by showing that this case of the multiple sum in (7.100) is 
$D_5$, as defined by (5.96).  This equality follows
termwise by using (6.5) to compute the $i=2$ case of
(7.48).
\qed\enddemo

The same analysis involving (5.18) that took Theorem 5.11
into Corollary 5.12 transforms Theorem 7.18 into an analog
of Corollary 7.6.  

We next consider the Schur function form of our $\chi_n$
determinant identities related to $\cn(u,k)$,
$\sn(u,k)\dn(u,k)$, $\dn(u,k)$, and 
$\sn(u,k)\cn(u,k)$.  

We first have the following theorem.
\proclaim{Theorem 7.19 } Let 
$\vartheta_2 (0,q)$ and $\vartheta_3 (0,q)$ be defined by 
\hbox{\rm(1.2)} and \hbox{\rm(1.1)}, respectively. 
Let $n=1,2,3,\cdots$.    We then have
$$\spreadlines{6 pt}\allowdisplaybreaks\align 
&\vartheta_2(0,q)^{2n^2}
\vartheta_3(0,q)^{2n(n-1)}\cr
\ \ &\cdot\left\{2n\left[1+24\sum\limits_{r=1}^{\infty}
{rq^{r} \over {1+q^{r}}}\right]+
\left[1-24\sum\limits_{r=1}^{\infty}
{(2r-1)q^{2r-1} \over {1+q^{2r-1}}}\right]
\right\}\tag 7.105\cr
=\,&4^{n}\tfrac{3}{n(2n-1)}
\prod\limits_{r=1}^{n-1}(2r)!^{-2}
\kern-2.6em\sum\limits_{{y_1,\ldots,y_n\geq 1
\hbox {\sixrm (odd)}}\atop 
{m_1>m_2>\cdots >m_n\geq 1 \hbox {\sixrm (odd)}}}
\kern-2.6em (-1)^{\tfrac{1}{2}(-n+y_1+\cdots+y_n)}\cr
&\cdot  \left(m_1^2+\cdots +m_n^2\right)
q^{\tfrac{1}{2}(m_1y_1+\cdots +m_ny_n)}
\kern-.65em
\prod\limits_{1\leq r<s\leq n}\kern-1.1em
\left(m_r^2-m_s^2\right)^2.\cr
\endalign$$
\endproclaim 
\demo{Proof} We take $n\geq 2$ and apply Theorem 6.19
to Theorem 5.33.   We utilize the same specializations as in
(7.59)--(7.63), except that $b_i$ is given by (7.92).  It is
immediate that (7.64) holds.  In this application 
of the determinant expansion of (6.57) we also have
(7.44), and (7.101)--(7.104).  

Keeping in mind (7.44), (7.60), (7.61), (7.63),
(7.64), (7.92), (7.101)--(7.104), 
noting that the inner sum in (6.57) consists of just one
term corresponding to $T = I_n$, observing that (7.47) 
holds, relabelling $(2m_r-1)$ and $(2y_r-1)$ as $m_r$
(odd)  and $y_r$ (odd), respectively, and simplifying, now
yields (7.105), for $n\geq 2$.

By Theorem 5.33, the $n=1$ case of (7.105) is established
by showing that this case of the multiple sum in (7.105) is 
$T_2$, as defined by (5.110).  This equality follows
termwise by using (6.5) to compute the $i=2$ case of
(7.59). 
\qed\enddemo

We next have the following theorem.
\proclaim{Theorem 7.20 }  Let 
$\vartheta_2 (0,q)$ and $\vartheta_3 (0,q)$ be defined by 
\hbox{\rm(1.2)} and \hbox{\rm(1.1)}, respectively. 
Let $n=1,2,3,\cdots$.    We then have
$$\spreadlines{6 pt}\allowdisplaybreaks\align 
&\vartheta_2(0,q)^{2n^2}
\vartheta_3(0,q)^{2n(n+1)}\cr
\ \ &\cdot\left\{2n\left[1+24\sum\limits_{r=1}^{\infty}
{rq^{r} \over {1+q^{r}}}\right]-
\left[1-24\sum\limits_{r=1}^{\infty}
{(2r-1)q^{2r-1} \over {1+q^{2r-1}}}\right]
\right\}\tag 7.106\cr
=\,&4^{n}\tfrac{3}{n(2n+1)}
\prod\limits_{r=1}^{n}(2r-1)!^{-2}
\kern-2.6em\sum\limits_{{y_1,\ldots ,y_n\geq 1
\hbox {\sixrm (odd)}}\atop 
{m_1>m_2>\cdots >m_n\geq 1 \hbox {\sixrm (odd)}}}
\kern-2.6em (-1)^{\tfrac{1}{2}(-n+y_1+\cdots+y_n)}
\cr
&\cdot \left(m_1^2+\cdots +m_n^2\right)
\left(m_1m_2\cdots m_n\right)^2\cr
&\cdot q^{\tfrac{1}{2}(m_1y_1+\cdots +m_ny_n)}
\kern-.65em
\prod\limits_{1\leq r<s\leq n}\kern-1.1em
\left(m_r^2-m_s^2\right)^2.\cr
\endalign$$
\endproclaim 
\demo{Proof} We take $n\geq 2$ and apply Theorem 6.19
to Theorem 5.34.   We utilize the same specializations as in
(7.68)--(7.72), except that $b_i$ is given by (7.92).  It is
immediate that (7.64) holds.  In this application 
of the determinant expansion of (6.57) we also have
(7.44), and (7.101)--(7.104).  

Keeping in mind (7.44), (7.64), (7.69), (7.70),
(7.72), (7.92), (7.101)--(7.104), 
noting that the inner sum in (6.57) consists of just one
term corresponding to $T = I_n$, observing that (7.47) 
holds, relabelling $(2m_r-1)$ and $(2y_r-1)$ as $m_r$
(odd)  and $y_r$ (odd), respectively, and simplifying, now
yields (7.106), for $n\geq 2$.

By Theorem 5.34, the $n=1$ case of (7.106) is established
by showing that this case of the multiple sum in (7.106) is 
$T_4$, as defined by (5.110).  This equality follows
termwise by using (6.5) to compute the $i=2$ case of
(7.68). 
\qed\enddemo

We next have the following theorem.
\proclaim{Theorem 7.21 } Let 
$\vartheta_2 (0,q)$ and $\vartheta_3 (0,q)$ be defined by 
\hbox{\rm(1.2)} and \hbox{\rm(1.1)}, respectively. 
Let $n=1,2,3,\cdots$.    We then have
$$\spreadlines{6 pt}\allowdisplaybreaks\align 
&\vartheta_2(0,q)^{2n(n-1)}
\vartheta_3(0,q)^{2n^2}\cr
\ \ &\cdot\left\{n\left[1+24\sum\limits_{r=1}^{\infty}
{rq^{r} \over {1+q^{r}}}\right]-
\left[1+24\sum\limits_{r=1}^{\infty}
{rq^{2r} \over {1+q^{2r}}}\right]
\right\}\cr
=\,&4^{n^2}\tfrac{3}{2n(2n-1)}
{\prod\limits_{r=1}^{n-1}(2r)!^{-2}} 
\kern-2.0em
\sum\limits_{{y_1,\ldots ,y_n\geq 1
\hbox {\sixrm (odd)}}\atop 
{m_1>m_2>\cdots >m_n\geq 1 }}
\kern-2.0em (-1)^{\tfrac{1}{2}(-n+y_1+\cdots+y_n)}
\cr
&\cdot \left(m_1^2+\cdots +m_n^2\right)
q^{m_1y_1+\cdots +m_ny_n}
\kern-.65em
\prod\limits_{1\leq r<s\leq n}\kern-1.1em
\left(m_r^2-m_s^2\right)^2
\tag 7.107a\cr
+&4^{n^2-1}\tfrac{3}{2n(2n-1)}
\prod\limits_{r=1}^{n-1}(2r)!^{-2}
\kern-2.3em
\sum\limits_{{y_1,\ldots ,y_{n-1}\geq 1
\hbox {\sixrm (odd)}}\atop 
{m_1>m_2>\cdots >m_{n-1}\geq 1 }}
\kern-2.3em (-1)^{\tfrac{1}{2}
(-n+1+y_1+\cdots+y_{n-1})}\cr
&\cdot  \left(m_1^2+\cdots +m_{n-1}^2\right)
\left(m_1m_2\cdots m_{n-1}\right)^4\cr
&\cdot q^{m_1y_1+\cdots +m_{n-1}y_{n-1}} 
\kern-.65em
\prod\limits_{1\leq r<s\leq {n-1}}\kern-1.1em
\left(m_r^2-m_s^2\right)^2,
\tag 7.107b\cr
\endalign$$
where the term in  \hbox{\rm(7.107b)} is defined to be $0$
if $n=1$.
\endproclaim 
\demo{Proof} We apply Theorem 6.19 to Theorem 5.35.  
We first take $n\geq 2$ and apply Theorem 6.19 
to the single  $n \times n$ determinant in the
right-hand-side of equation (5.276a) in Theorem 5.35.  
We utilize the same specializations as in
(7.75)--(7.79), except that $b_i$ is given by (7.92).  It is
immediate that (7.53) holds.  In this application 
of the determinant expansion of (6.57) we also have
(7.44), and (7.101)--(7.104).  

Keeping in mind (7.44), (7.53), (7.76), (7.77),
(7.79), (7.92), (7.101)--(7.104), 
noting that the inner sum in (6.57) consists of just one
term corresponding to $T = I_n$, observing that 
$-m_r + 2m_ry_r = m_r(2y_r-1)$, 
relabelling $(2y_r-1)$ as $y_r$ (odd), 
and simplifying, now yields the term in (7.107a), for 
$n\geq 2$.

By Theorem 5.35, the $n=1$ case of (7.107a) is 
established by showing that this case of the multiple sum in
(7.107a) is $N_2$, as defined by (5.128).  This equality
follows termwise by using (6.5) to compute the $i=2$ case
of (7.75). 

For $n\geq 3$, the derivation of the term in (7.107b) from
(5.276b) is the same as for (7.107a), except that 
$n\mapsto n-1$, $c_i\mapsto  2i+2$, and 
$c_{\ell_1}+b_{j_1}\mapsto 4$.

By Theorem 5.35, the $n=2$ case of (7.107b) is 
established by showing that this case of the multiple sum in
(7.107b) is $N_6$, as defined by (5.128).  This equality
follows termwise by using (6.5) to compute the $i=4$ case
of (7.75). 

If $n=1$, the term in (7.107b) is defined to be $0$.
\qed\enddemo

We next have the following theorem.  
\proclaim{Theorem 7.22 } Let 
$\vartheta_2 (0,q)$ and $\vartheta_3 (0,q)$ be defined by 
\hbox{\rm(1.2)} and \hbox{\rm(1.1)}, respectively. 
Let $n=1,2,3,\cdots$.    We then have
$$\spreadlines{6 pt}\allowdisplaybreaks\align 
&\vartheta_2(0,q)^{2n(n+1)}
\vartheta_3(0,q)^{2n^2}\cr
\ \ &\cdot\left\{n\left[1+24\sum\limits_{r=1}^{\infty}
{rq^{r} \over {1+q^{r}}}\right]+
\left[1+24\sum\limits_{r=1}^{\infty}
{rq^{2r} \over {1+q^{2r}}}\right]
\right\}\tag 7.108\cr
=\,&4^{n(n+1)}\tfrac{6}{n(2n+1)}
\prod\limits_{r=1}^{n}(2r-1)!^{-2}\kern-2.1em
\sum\limits_{{y_1,\ldots ,y_n\geq 1
\hbox {\sixrm (odd)}}\atop 
{m_1>m_2>\cdots >m_n\geq 1 }}
\kern-2.0em (-1)^{\tfrac{1}{2}(-n+y_1+\cdots+y_n)}\cr 
&\cdot \left(m_1^2+\cdots +m_n^2\right)
\left(m_1m_2\cdots m_n\right)^2\cr
&\cdot q^{m_1y_1+\cdots +m_ny_n}
\kern-.65em
\prod\limits_{1\leq r<s\leq n}\kern-1.1em
\left(m_r^2-m_s^2\right)^2.\cr
\endalign$$
\endproclaim 
\demo{Proof} We take $n\geq 2$ and apply Theorem 6.19
to Theorem 5.36.   We utilize the same specializations as in
(7.81)--(7.85), except that $b_i$ is given by (7.92).  It is
immediate that (7.53) holds.  In this application 
of the determinant expansion of (6.57) we also have
(7.44), and (7.101)--(7.104).  

Keeping in mind (7.44), (7.53), (7.82), (7.83),
(7.85), (7.92), (7.101)--(7.104), 
noting that the inner sum in (6.57) consists of just one
term corresponding to $T = I_n$, observing that 
$-m_r + 2m_ry_r = m_r(2y_r-1)$, 
relabelling $(2y_r-1)$ as $y_r$ (odd), 
and simplifying, now yields (7.108), for $n\geq 2$.

By Theorem 5.36, the $n=1$ case of (7.108) is established
by showing that this case of the multiple sum in (7.108) is 
$N_4$, as defined by (5.128).  This equality follows
termwise by using (6.5) to compute the $i=2$ case of
(7.81). 
\qed\enddemo

\remark{Remark} Note that the terms in the multiple sums in
(7.99), (7.100), (7.105), (7.106), (7.107a), and (7.108)
are simply $(m_1^2+\cdots +m_n^2)$ times the 
corresponding terms in (7.36), (7.37), (7.58), (7.67),
(7.74a), and (7.80), respectively. 
\endremark

Chan's paper \cite{41}, which utilized Dedekind's 
$\eta$-function transformation formula and 
Hecke's correspondence  between Fourier series and
Dirichlet series to prove the equivalence of some partition
identities of Ramanujan, motivates a direct relationship
between Theorems 7.1 and 7.2 and the simpler 
identities in Theorem 7.5.   This relationship is a
consequence of 
$$\spreadlines{6 pt}\allowdisplaybreaks\align 
&\kern -3 em\vartheta_3(0,-q) = (-i\tau)^{-1/2}
\vartheta_2(0,q_1),\tag 7.109\cr
\kern -13.75 em \text{where}\kern 13.75 em 
\quad &\kern -3 em q:=e^{i\pi\tau}
\quad\text{ and} \qquad 
q_1:=e^{-i\pi/\tau}.\tag 7.110\cr
\endalign$$ 
To see (7.109), take $z=0$ in Jacobi's transformation 
$$\vartheta_4(\tau z|\tau)=(-i\tau)^{-1/2}
e^{-i\tau z^2/\pi}\vartheta_2(z|\tfrac{-1}{\tau})
\tag 7.111$$
in \cite{135, Eqn. (1.7.15), pp. 17} of the theta functions
$\vartheta_4$ and $\vartheta_2$.  Equation (7.109) 
implies that 
$$\spreadlines{6 pt}\allowdisplaybreaks\align 
\kern -3 em\vartheta_3(0,-q)^{4n^2} = &\  
(-1)^n\tau^{-2n^2}
\vartheta_2(0,q_1)^{4n^2},
\tag 7.112\cr
\kern -11.6 em \text{and}\kern 11.6 em 
\quad \kern -3 em\vartheta_3(0,-q)^{4n(n+1)} = &\   
\tau^{-2n(n+1)}
\vartheta_2(0,q_1)^{4n(n+1)},
\tag 7.113\cr
\endalign$$ 
where $q$ and $q_1$ are given by (7.110).  

Applying Theorem 7.5 to the right-hand sides of 
(7.112) and (7.113) now leads to the following corollary.
\proclaim{Corollary 7.23} Let $\vartheta_3(0,q)$ be 
defined by 
\hbox{\rm(1.1)}. Let $q$ and $q_1$ be given by 
\hbox{\rm(7.110)}.  Let 
$n=1,2,3,\cdots$.   We then have   
$$\spreadlines{6 pt}\allowdisplaybreaks\align 
\kern -2 em\vartheta_3(0,-q)^{4n^2}=\,&
(-1)^n4^{n(n+1)}\tau^{-2n^2}
\prod\limits_{r=1}^{2n-1}(r!)^{-1}
\kern-2.6em\sum\limits_{{y_1,\ldots ,y_n\geq 1
\hbox {\sixrm (odd)}}\atop 
{m_1>m_2>\cdots >m_n\geq 1 \hbox {\sixrm (odd)}}}
\kern-2.6em m_1m_2\cdots m_n\cr
&\cdot q_1^{m_1y_1+\cdots +m_ny_n}
\kern-.65em
\prod\limits_{1\leq r<s\leq n}\kern-1.1em
\left(m_r^2-m_s^2\right)^2,
\tag 7.114\cr
\kern -10.3 em \text{and}\kern 10.3 em\quad &\cr 
\kern -2 em\vartheta_3(0,-q)^{4n(n+1)} =\,&
2^{n(4n+5)}\tau^{-2n(n+1)}
\prod\limits_{r=1}^{2n}(r!)^{-1}
\kern-1.65em\sum\limits_{{y_1,\ldots ,y_n\geq 1
\hbox {\sixrm (odd)}}\atop 
{m_1>m_2>\cdots >m_n\geq 1 }}
\kern-1.8em 
\left(m_1m_2\cdots m_n\right)^3\cr
&\cdot q_1^{2(m_1y_1+\cdots +m_ny_n)}
\kern-.65em
\prod\limits_{1\leq r<s\leq n}\kern-1.1em
\left(m_r^2-m_s^2\right)^2.
\tag 7.115\cr
\endalign$$
\endproclaim 

To obtain an identity just involving multiple sums, 
replace the products $\vartheta_3(0,-q)^{4n^2}$ and 
$\vartheta_3(0,-q)^{4n(n+1)}$ in (7.114) and (7.115) by
the  sums in (7.1) and (7.11). 

A Lambert series version of Corollary 7.23 follows by 
applying Theorems 5.4, 5.6, and 5.11 to the products in 
(7.112) and (7.113).  The analysis in \cite{41, 42} 
suggests that an extension of the  Hecke correspondence 
between Fourier series and Dirichlet series, applied to this   
Lambert series version of Corollary 7.23, may lead to 
functional equations and evaluations for the corresponding 
multiple $L$-series.   

The classical techniques in \cite{90, pp. 96--101} and 
\cite{28, pp. 288--305} applied to Theorems 7.1 through 
7.22 lead to corresponding infinite families of lattice sum
transformations.  We illustrate this procedure by applying
the Mellin transform to Theorems 7.1, 7.2, 7.5, and
the two special cases in Corollaries 8.1 and 8.2.  These
computations involve the specialized theta functions 
$\vartheta_2(q):=\vartheta_2(0,q)$, 
$\vartheta_3(q):=\vartheta_3(0,q)$,
$\vartheta_4(q):=\vartheta_3(0,-q)$, 
  and the Mellin transform $Mf(s)$ given by 
$$Mf(s)\equiv M_s(f):=\int_{0}^{\infty}
\kern-.6em f(t)t^{s-1}\,dt.\tag 7.116$$
We make frequent use of 
$$M_s(e^{-nt})=n^{-s}\cdot\Gamma (s),\tag 7.117$$
where $\Gamma (s)$ is the gamma function.  

We first study the $m$-dimensional cubic lattice sums 
$b_m(2s)$ defined by 
$$b_m(2s):=\sum\limits_
{{(i_1,\dots,i_m)\neq (0,\dots,0)}
\atop {-\infty <i_1,\dots,i_m <\infty}}\kern-.25em  
{{(-1)^{i_1+\cdots+i_m}}\over 
{(i_1^2+\cdots+i_m^2)^s}},\tag 7.118$$
where $\Re(s)>0$. The convergence of these general 
sums is discussed in \cite{28, Ex.1 and 2, pp. 290--291}.  
Sums of the form $b_3(2s)$ occur naturally in chemistry.  
For example, $b_3(1)$ may be viewed as the potential or
Coulomb sum at the origin of a cubic lattice with
alternating unit charges at all nonzero lattice points.   
This is essentially an idealization of a rocksalt crystal.  The
quantity $b_3(1)$ is known as {\sl Madelung's constant}
for NaCl. Different crystals lead to different lattice sums.  

Applying the Mellin transform termwise to 
both sides of the $q=e^{-t}$ case of (7.1) yields  
the following theorem.
\proclaim{Theorem 7.24 } Let $n=1,2,3,\cdots$. 
We then have the formal identity  
$$\spreadlines{6 pt}\allowdisplaybreaks\align 
&b_{4n^2}(2s)\tag 7.119\cr 
=\,&\sum\limits_{p=1}^n
(-1)^p2^{2n^2+n}\prod\limits_{r=1}^{2n-1}(r!)^{-1}
\kern-1.8em
\sum\limits_{{y_1,\ldots,y_p\geq 1}\atop 
{m_1>m_2>\cdots>m_p\geq 1}}
\kern-1.65em(-1)^{y_1+\cdots+y_p}\cr
&\cdot(-1)^{m_1+\cdots+m_p}
(m_1y_1+\cdots+m_py_p)^{-s}\kern-.65em
\prod\limits_{1\leq r<s\leq p}\kern-1.1em
\left(m_r^2-m_s^2\right)^2\cr
&\cdot(m_1m_2\cdots m_p)\kern-.65em 
\sum\limits_{{\emptyset \subset S, T\subseteq I_n}
\atop {\Vert S\Vert=\Vert T\Vert =p}}\kern-.5em 
(-1)^{\Sigma (S)+\Sigma (T)}\cdot 
\det (D_{n-p,S^c,T^c})\cr
&\cdot (m_1m_2\cdots m_p)^{2\ell_1+2j_1-4}
s_{\lambda}(m_1^2,\ldots,m_p^2)\,
s_{\nu}(m_1^2,\ldots,m_p^2),\cr 
\endalign$$
where the same assumptions hold as in 
\hbox{\rm Theorem 7.1}, and $b_{4n^2}(2s)$ is determined
by \hbox{\rm (7.118)}.
\endproclaim

Similarly, applying the Mellin transform termwise to 
both sides of the $q=e^{-t}$ case of (7.11) yields  
the following theorem.
\proclaim{Theorem 7.25 } Let $n=1,2,3,\cdots$. 
We then have the formal identity
$$\spreadlines{6 pt}\allowdisplaybreaks\align 
&b_{4n(n+1)}(2s)\tag 7.120\cr 
=\,&\sum\limits_{p=1}^n
(-1)^{n-p}2^{2n^2+3n}\prod\limits_{r=1}^{2n}(r!)^{-1}
\kern-1.8em
\sum\limits_{{y_1,\ldots,y_p\geq 1}\atop 
{m_1>m_2>\cdots>m_p\geq 1}}
\kern-1.65em(-1)^{m_1+\cdots+m_p}\cr
&\cdot (m_1y_1+\cdots+m_py_p)^{-s}\,
(m_1m_2\cdots m_p)^3\kern-.65em
\prod\limits_{1\leq r<s\leq p}\kern-1.1em
\left(m_r^2-m_s^2\right)^2\cr
&\cdot
\sum\limits_{{\emptyset \subset S, T\subseteq I_n}
\atop {\Vert S\Vert=\Vert T\Vert =p}}\kern-.5em 
(-1)^{\Sigma (S)+\Sigma (T)}\cdot 
\det (D_{n-p,S^c,T^c})\cr 
&\cdot (m_1m_2\cdots m_p)^{2\ell_1+2j_1-4} 
s_{\lambda}(m_1^2,\ldots,m_p^2)\,
s_{\nu}(m_1^2,\ldots,m_p^2),\cr 
\endalign$$
where the same assumptions hold as in 
\hbox{\rm Theorem 7.2}, and 
$b_{4n(n+1)}(2s)$ is determined
by \hbox{\rm (7.118)}.
\endproclaim

It appears that the convergence conditions $\Re(s)>4n-3$ 
and $\Re(s)>4n-1$ are sufficient for Theorems 7.24 and 
7.25, respectively.
  
We next study the $m$-dimensional cubic lattice sums 
$c_m(2s)$ defined by 
$$c_m(2s):=\sum\limits_
{-\infty <i_1,\dots,i_m <\infty}\kern-1.35em  
\left[(i_1+1/2)^2+\cdots+(i_m+1/2)^2\right]^{-s},
\tag 7.121$$ 
where $\Re(s)>m/2$.  This convergence 
condition is determined by utilizing the estimates for 
$r_m(N)$ in \cite{94, Eqn. (9.20), pp. 122} to compare 
(7.121) with $\zeta (s-\tfrac{m}{2}+1)$.    

Applying the Mellin transform termwise to both sides of the
$q=e^{-t}$ cases of (7.36) and (7.37) yields 
the following theorem.
\proclaim{Theorem 7.26}Let $n=1,2,3,\cdots$, and let 
$c_{4n^2}(2s)$ and $c_{4n(n+1)}(2s)$ be determined
by \hbox{\rm (7.121)}. We then have the formal identities  
$$\spreadlines{6 pt}\allowdisplaybreaks\align 
\kern -2 em c_{4n^2}(2s)=\,&
4^{n(n+1)}\cdot\prod\limits_{r=1}^{2n-1}(r!)^{-1}
\kern-2.6em\sum\limits_{{y_1,\ldots ,y_n\geq 1
\hbox {\sixrm (odd)}}\atop 
{m_1>m_2>\cdots >m_n\geq 1 \hbox {\sixrm (odd)}}}
\kern-2.6em m_1m_2\cdots m_n\cr
&\cdot (m_1y_1+\cdots+m_ny_n)^{-s}
\kern-.65em
\prod\limits_{1\leq r<s\leq n}\kern-1.1em
\left(m_r^2-m_s^2\right)^2,
\tag 7.122\cr
\kern -10.75 em \text{and}\kern 10.75 em\quad &\cr 
\kern -2 em c_{4n(n+1)}(2s)=\,&
2^{-s}2^{4n^2+5n}\cdot
\prod\limits_{r=1}^{2n}(r!)^{-1}
\kern-1.65em\sum\limits_{{y_1,\ldots ,y_n\geq 1
\hbox {\sixrm (odd)}}\atop 
{m_1>m_2>\cdots >m_n\geq 1 }}
\kern-1.8em 
\left(m_1m_2\cdots m_n\right)^3\cr
&\cdot (m_1y_1+\cdots+m_ny_n)^{-s}
\kern-.65em
\prod\limits_{1\leq r<s\leq n}\kern-1.1em
\left(m_r^2-m_s^2\right)^2.
\tag 7.123\cr
\endalign$$
\endproclaim 

It appears that the convergence conditions $\Re(s)>2n^2$ 
and $\Re(s)>2n(n+1)$ are sufficient for (7.122) and 
(7.123), respectively.  Furthermore, the right-hand sides of 
(7.122) and (7.123) may be an analytic continuation of the
left-hand sides to $\Re(s)>4n-2$ and $\Re(s)>4n$, 
respectively.  

We now make use of Corollaries 8.1 and 8.2 to write down
simplified versions of the $n=2$ cases of Theorems 7.24
and 7.25.  

We have the following corollaries.
\proclaim{Corollary 7.27 } Let $b_{16}(2s)$ be 
determined by \hbox{\rm(7.118)}.   
We then have the formal identity
$$\spreadlines{6 pt}\allowdisplaybreaks\align 
b_{16}(2s)
=\, & -\tfrac{2^5}{3}
\sum\limits_{y_1\geq 1}
(-1)^{y_1}y_1^{-s}\cdot 
\sum\limits_{m_1\geq 1}
(-1)^{m_1}
(1+m_1^2+m_1^4)m_1^{1-s}\cr
&+\tfrac{2^8}{3}\kern-1.15em 
\sum\limits_{{y_1,y_2\geq 1}\atop 
{m_1>m_2\geq 1}}\kern-.85em 
(-1)^{y_1+y_2+m_1+m_2}(m_1m_2)
(m_1^2-m_2^2)^2(m_1y_1+m_2y_2)^{-s}.
\tag 7.124\cr
\endalign$$
\endproclaim

\proclaim{Corollary 7.28 } Let $b_{24}(2s)$ be 
determined by \hbox{\rm(7.118)}. 
We then have the formal identity
$$\spreadlines{6 pt}\allowdisplaybreaks\align 
b_{24}(2s) =& +\tfrac{2^4}{3^2}\cdot
\zeta (s)\cdot 
\sum\limits_{m_1\geq 1}
(-1)^{m_1}
(17+8m_1^2+2m_1^4)m_1^{3-s}\cr
&+\tfrac{2^9}{3^2}\kern-1.15em 
\sum\limits_{{y_1,y_2\geq 1}\atop 
{m_1>m_2\geq 1}}\kern-.85em 
(-1)^{m_1+m_2}(m_1m_2)^3
(m_1^2-m_2^2)^2(m_1y_1+m_2y_2)^{-s}.
\tag 7.125\cr
\endalign$$
\endproclaim
It appears that the convergence conditions $\Re(s)>5$ 
and $\Re(s)>7$ are sufficient for Corollaries 7.27 and 
7.28, respectively.

The analysis involving Theorems 7.24--7.26 may have
applications to the work in \cite{32}.

\head 8. The $36$ and $48$ squares identities
\endhead

In this section we write down the $n=2$ and $n=3$ cases 
of Theorems 7.1 and 7.2.  These results provide 
explicit multiple power series formulas for
$16$, $24$, $36$, and $48$ squares.  

The first two corollaries give the $n=2$ cases of Theorems 7.1
and 7.2, respectively.
\proclaim{Corollary 8.1 } Let $\vartheta_3(0,-q)$ be determined
by \hbox{\rm(1.1)}.  Then 
$$\spreadlines{6 pt}\allowdisplaybreaks\align 
\vartheta_3(0,-q)^{16} =1& -\tfrac{2^5}{3}\kern-1em 
\sum\limits_{y_1,m_1\geq 1}
\kern-.55em(-1)^{y_1+m_1}m_1
(1+m_1^2+m_1^4)q^{m_1y_1}\tag 8.1a\cr
&+\tfrac{2^8}{3}\kern-1.15em 
\sum\limits_{{y_1,y_2\geq 1}\atop 
{m_1>m_2\geq 1}}\kern-.85em 
(-1)^{y_1+y_2+m_1+m_2}(m_1m_2)
(m_1^2-m_2^2)^2q^{m_1y_1+m_2y_2}.
\tag 8.1b\cr
\endalign$$
\endproclaim

\proclaim{Corollary 8.2 } Let $\vartheta_3(0,-q)$ be determined
by \hbox{\rm(1.1)}.  Then 
$$\spreadlines{6 pt}\allowdisplaybreaks\align 
\vartheta_3(0,-q)^{24} =1& +\tfrac{2^4}{3^2}\kern-1em 
\sum\limits_{y_1,m_1\geq 1}
\kern-.55em(-1)^{m_1}m_1^3
(17+8m_1^2+2m_1^4)q^{m_1y_1}\tag 8.2a\cr
&+\tfrac{2^9}{3^2}\kern-1.15em 
\sum\limits_{{y_1,y_2\geq 1}\atop 
{m_1>m_2\geq 1}}\kern-.85em 
(-1)^{m_1+m_2}(m_1m_2)^3
(m_1^2-m_2^2)^2q^{m_1y_1+m_2y_2}.
\tag 8.2b\cr
\endalign$$
\endproclaim

The next two corollaries give the $n=3$ cases of Theorems 7.1
and 7.2, respectively.
\proclaim{Corollary 8.3 } Let $\vartheta_3(0,-q)$ be determined
by \hbox{\rm(1.1)}.  Then 
$$\spreadlines{6 pt}\allowdisplaybreaks\align 
\vartheta_3(0,-q)^{36} =1&
-\tfrac{2^3}{3^2\cdot 5}\kern-1em 
\sum\limits_{y_1,m_1\geq 1}
\kern-.55em(-1)^{y_1+m_1}q^{m_1y_1}m_1\cr
&\kern 4em
\cdot\left[69+120m_1^2+172m_1^4+40m_1^6+
4m_1^8\right]\cr
&+\tfrac{2^9}{3^3\cdot 5}\kern-1.15em 
\sum\limits_{{y_1,y_2\geq 1}\atop 
{m_1>m_2\geq 1}}\kern-.85em 
(-1)^{y_1+y_2+m_1+m_2}
q^{m_1y_1+m_2y_2}(m_1m_2)(m_1^2-m_2^2)^2\cr
&\kern 4em
\cdot\left[62+17m_1^2+17m_2^2+2m_1^4+
2m_2^4+8m_1^2m_2^2\right.\cr
&\kern 6.05em
\left.+2m_1^4m_2^2+2m_1^2m_2^4+
2m_1^4m_2^4\right]\cr
&-\tfrac{2^{13}}{3^3\cdot 5}\kern-1.8em 
\sum\limits_{{y_1,y_2,y_3\geq 1}\atop 
{m_1>m_2>m_3\geq 1}}\kern-1.6em 
(-1)^{(y_1+y_2+y_3)+(m_1+m_2+m_3)}
q^{m_1y_1+m_2y_2+m_3y_3}\cr
&\kern 4em\cdot (m_1m_2m_3)\kern-.6em
\prod\limits_{1\leq r<s\leq 3}\kern-.9em
\left(m_r^2-m_s^2\right)^2.
\tag 8.3\cr
\endalign$$
\endproclaim

\proclaim{Corollary 8.4 } Let $\vartheta_3(0,-q)$ be determined
by \hbox{\rm(1.1)}.  Then 
$$\spreadlines{6 pt}\allowdisplaybreaks\align 
\vartheta_3(0,-q)^{48} =1&
+\tfrac{2^5}{3^3\cdot 5^2}\kern-1em 
\sum\limits_{y_1,m_1\geq 1}
\kern-.55em(-1)^{m_1}q^{m_1y_1}m_1^3\cr
&\kern 4em
\cdot\left[902+760m_1^2+321m_1^4+40m_1^6+
2m_1^8\right]\cr
&+\tfrac{2^{10}}{3^5\cdot 5^2}\kern-1.15em 
\sum\limits_{{y_1,y_2\geq 1}\atop 
{m_1>m_2\geq 1}}\kern-.85em 
(-1)^{m_1+m_2}
q^{m_1y_1+m_2y_2}(m_1m_2)^3(m_1^2-m_2^2)^2\cr
&\kern 4em
\cdot\left[1382+248m_1^2+248m_2^2+17m_1^4+
17m_2^4+68m_1^2m_2^2\right.\cr
&\kern 7.05em
\left.+8m_1^4m_2^2+8m_1^2m_2^4+
2m_1^4m_2^4\right]\cr
&+\tfrac{2^{15}}{3^5\cdot 5^2}\kern-1.8em 
\sum\limits_{{y_1,y_2,y_3\geq 1}\atop 
{m_1>m_2>m_3\geq 1}}\kern-1.6em 
(-1)^{m_1+m_2+m_3}
q^{m_1y_1+m_2y_2+m_3y_3}\cr
&\kern 4em\cdot (m_1m_2m_3)^3\kern-.6em
\prod\limits_{1\leq r<s\leq 3}\kern-.9em
\left(m_r^2-m_s^2\right)^2.
\tag 8.4\cr
\endalign$$
\endproclaim

\noindent{\bf Acknowledgment:}~The author would like to thank George Andrews,
Richard Askey, Bruce Berndt, Mourad Ismail, and Christian 
Krattenthaler for some of the references to the literature.



\widestnumber\no{9999} 
\Refs

\ref\key 1
\by N. H. Abel\paper Recherches sur les fonctions elliptiques
\jour J. Reine Angew. Math. \vol 2\yr 1827 
\pages 101--181
\moreref reprinted in \OE uvres Compl\`etes {\bf T1}
\inbook\publ Grondahl and Son\publaddr Christiania 
\yr 1881\pages 263--388
\moreref 
\publ reprinted by Johnson Reprint Corporation
\publaddr New York, 1965
\endref  

\ref\key 2
\by W. A. Al-Salam and L. Carlitz
\paper Some determinants of Bernoulli, Euler, and related 
numbers
\jour  Portugal. Math. (2)\vol 18\yr 1959
\pages 91--99 
\endref

\ref\key 3
\by K. Ananda-Rau\paper On the representation of a number as the
sum of an even number of squares\jour J. Madras Univ. Sect. B
\vol 24\yr 1954\pages 61--89
\endref

\ref\key 4
\by G. E. Andrews\paper Applications of basic
hypergeometric functions\jour SIAM Rev.\vol 16\yr 1974\pages 441--484
\endref

\ref\key 5
\by G. E. Andrews
\book $q$-Series: Their development and application in analysis, 
number theory, combinatorics, physics and computer algebra, {\rm NSF
CBMS Regional Conference Series}\vol 66 
\yr 1986
\endref

\ref\key 6
\by G. E. Andrews, R. Askey, and R. Roy
\paper Special Functions
\inbook Encyclopedia of Mathematics and Its
Applications\vol 71\ed G.-C. Rota\publ Cambridge University
Press\publaddr Cambridge\yr 1999
\endref

\ref\key 7
\by G. E. Andrews and B. C. Berndt
\book Ramanujan's Lost Notebook
\bookinfo Part I
\publ Springer-Verlag
\publaddr New York \yr 
\moreref in preparation
\endref

\ref\key 8
\by T. M. Apostol
\paper Modular Functions and Dirichlet Series in Number Theory
\inbook vol. 41 of Graduate Texts in Mathematics
\publ Springer-Verlag\publaddr New York 
\yr 1976
\endref

\ref\key 9
\by R. Askey and M. E. H. Ismail
\paper Recurrence relations, continued fractions and orthogonal
polynomials
\jour Mem. Amer. Math. Soc.\vol 300\yr 1984\pages 108 pages
\endref

\ref\key 10
\by H. Au-Yang and J. H. H. Perk
\paper Critical correlations in a
$Z$-invariant inhomogeneous Ising model
\jour Phys. A\vol 144\yr 1987\pages 44--104
\endref

\ref\key 11
\by I. G. Bashmakova
\paper Diophantus and Diophantine equations
\inbook The Dolciani Mathematical Expositions
\vol 20\publ Mathematical Association of America
\publaddr Washington, DC\yr 1997, xvi+90 pp
\transl Translated from the 1972 Russian original by Abe
Shenitzer and updated by Joseph Silverman
\endref

\ref\key  12
\by I. G. Bashmakova and G. S. Smirnova
\paper The birth of literal algebra
\jour Amer. Math. Monthly
\vol 106\yr 1999
\pages 57--66
\transl Translated from the Russian and edited by Abe Shenitzer 
\endref

\ref\key 13
\by E. F. Beckenbach, W. Seidel, O. Sz\'asz
\paper Recurrent determinants of Legendre and of ultraspherical 
polynomials 
\jour  Duke Math. J.\vol 18\yr 1951 
\pages 1--10
\endref

\ref\key 14
\by E. T. Bell
\paper On the number of representations of $2n$ as a sum of $2r$ 
squares
\jour Bull. Amer. Math. Soc.\vol 26\yr 1919\pages 19--25
\endref

\ref\key 15
\by E. T. Bell
\paper Theta expansions useful in arithmetic 
\jour The Messenger of Mathematics (New Series)
\vol 53\yr 1924\pages 166--176
\endref

\ref\key 16
\by E. T. Bell
\paper On the power series for elliptic functions 
\jour Trans. Amer. Math. Soc.\vol 36\yr 1934\pages 841--852
\endref

\ref\key 17
\by E. T. Bell
\paper The arithmetical function $M(n,f,g)$ and its associates
connected with elliptic power series 
\jour Amer. J. Math.\vol 58\yr 1936\pages 759--768
\endref

\ref\key 18
\by E. T. Bell
\paper Polynomial approximations for elliptic functions 
\jour Trans. Amer. Math. Soc.\vol 44\yr 1938\pages 47--57
\endref

\ref\key 19 
\by C. Berg and G. Valent
\paper The Nevanlinna parameterization for some indeterminate
Stieltjes moment problems associated with birth and death processes
\jour Methods Appl. Anal.\vol 1 
\yr 1994\pages 169--209
\endref

\ref\key 20
\by B. C. Berndt\book Ramanujan's Notebooks
\bookinfo Part II
\publ Springer-Verlag
\publaddr New York \yr 1989
\endref

\ref\key 21
\by B. C. Berndt\book Ramanujan's Notebooks
\bookinfo Part III
\publ Springer-Verlag
\publaddr New York \yr 1991 
\endref

\ref\key 22
\by B. C. Berndt\book Ramanujan's Notebooks
\bookinfo Part V
\publ Springer-Verlag
\publaddr New York \yr 1998 
\endref

\ref\key 23
\by B. C. Berndt\paper Ramanujan's theory of theta-functions
\inbook Theta Functions From the Classical to the Modern
\ed M. Ram Murty
\bookinfo vol. 1 of CRM Proceedings \& Lecture Notes 
\publ American Mathematical Society\publaddr Providence, RI  
\yr 1993\pages 1--63
\endref

\ref\key 24
\by B. C. Berndt
\paper Fragments by Ramanujan on Lambert series 
\inbook in Number Theory and Its Applications (Kyoto, 1997)
\eds K\. Gy\"ory and S\. Kanemitsu
\bookinfo vol. 2 of Dev. Math.
\publ Kluwer Academic Publishers
\publaddr Dordrecht
\yr 1999\pages 35--49
\endref

\ref\key 25
\by M. Bhaskaran
\paper A plausible reconstruction of Ramanujan's proof of his
formula for $\vartheta^{4s}(q)$
\inbook in Ananda Rau memorial volume
\bookinfo
Publications of the Ramanujan Institute;no. 1.
\publ Ramanujan Institute\publaddr Madras\yr 1969
\pages 25--33
\endref 

\ref\key 26
\by M. N. Bleicher and M. I. Knopp
\paper Lattice points in a sphere
\jour Acta Arith. \vol 10\yr 1964/1965\pages 369--376
\endref

\ref\key 27
\by F. van der Blij[i]\paper The function $\tau (n)$ of S. Ramanujan
\jour  Math. Student\vol 18\yr 1950
\pages 83--99
\endref

\ref\key 28
\by J. M. and P. B. Borwein\book Pi and the AGM
\publ John Wiley \& Sons, Inc.
\publaddr New York \yr 1987
\endref

\ref\key 29
\by D. M. Bressoud\book Proofs and confirmations. The story of 
the alternating sign matrix conjecture
\publ MAA Spectrum. Mathematical Association of America,
Washington, DC; Cambridge University Press
\publaddr Cambridge\yr 1999
\moreref see pp. 245--256  
\endref

\ref\key 30
\by D. M. Bressoud and J. Propp 
\paper How the alternating sign matrix conjecture was solved
\jour Notices Amer. Math. Soc.
\vol 46\yr 1999\pages 637--646
\endref

\ref\key 31
\by C. Brezinski
\paper History of Continued Fractions and Pad\'e Approximants
\inbook vol. 12 of Springer Series in Computational Mathematics
\publ Springer-Verlag\publaddr New York 
\yr 1991
\endref

\ref\key 32
\by D. J. Broadhurst
\paper On the enumeration of irreducible $k$-fold Euler sums and
their roles in knot theory and field theory
\jour J. Math. Phys., to appear
\endref

\ref\key 33
\by V. Bulygin\paper Sur une application des 
fonctions elliptiques au probl\`eme de repr\'esentation des nombres
entiers par une somme de carr\'es \jour Bull. Acad. Imp. Sci. St.
Petersbourg Ser. VI\vol 8\yr 1914\pages 389--404
\moreref \by B. Boulyguine\paper Sur la repr\'esentation
d'un nombre entier par une somme de carr\'es
\jour Comptes Rendus Paris\vol 158\yr 1914\pages 328--330
\endref

\ref\key 34
\by V. Bulygin (B. Boulyguine)
\paper Sur la repr\'esentation
d'un nombre entier par une somme de carr\'es
\jour Comptes Rendus Paris\vol 161\yr 1915\pages 28--30
\endref

\ref\key 35
\by J. L. Burchnall 
\paper An algebraic property of the classical polynomials 
\jour Proc. London Math. Soc. (3)
\vol 1\yr 1951\pages 232--240
\endref

\ref\key 36
\by L. Carlitz\paper Hankel determinants and Bernoulli numbers
\jour  T\^ohoku Math. J. (2)\vol 5\yr 1954 
\pages 272--276
\endref

\ref\key 37
\by L. Carlitz\paper Note on sums of 4 and 6 squares 
\jour Proc. Amer. Math. Soc.\vol 8\yr 1957
\pages 120--124
\endref

\ref\key 38
\by L. Carlitz\paper Some orthogonal polynomials related to
elliptic functions 
\jour  Duke Math. J.\vol 27\yr 1960 
\pages 443--459
\endref

\ref\key 39
\by L. Carlitz\paper Bulygin's method for sums of squares. 
The arithmetical theory of quadratic forms, I
\jour (Proc. Conf., Louisiana State Univ., Baton Rouge, La., 1972;
dedicated to Louis Joel Mordell). J. Number Theory\vol 5\yr 1973 
\pages 405--412
\endref

\ref\key 40
\by R. Chalkley\paper A persymmetric determinant
\jour J. Math. Anal. Appl.\vol 187\yr 1994
\pages 107--117
\endref

\ref\key 41
\by H. H. Chan\paper On the equivalence of Ramanujan's partition 
identities and a connection with the Rogers-Ramanujan continued
fraction
\jour J. Math. Anal. Appl.\vol 198\yr 1996
\pages 111--120
\endref

\ref\key 42
\by H. H. Chan
\moreref Private Communication, August 1996
\endref

\ref\key 43
\by K. Chandrasekharan
\paper Elliptic Functions
\inbook vol. 281 of Grundlehren Math. Wiss.
\publ Springer-Verlag\publaddr Berlin
\yr 1985
\endref

\ref\key 44
\by T. S. Chihara
\paper An Introduction to Orthogonal Polynomials
\inbook vol. 13 of Mathematics and Its Applications
\publ Gordon and Breach\publaddr New York 
\yr 1978
\endref

\ref\key 45
\by S. H. Choi and D. Gouyou--Beauchamps
\paper Enum\'eration de tableaux de Young semi-standard
\inbook Series Formelles et Combinatoire Algebrique:Act\'es du
Colloque
\eds M. Delest, G. Jacob, P. Leroux
\publ Universit\'e Bordeaux I
\yr 2--4 May, 1991\pages 229--243
\endref

\ref\key 46
\by D. V. Chudnovsky and G. V. Chudnovsky
\paper Computational problems in arithmetic of linear 
differential equations.  Some diophantine applications
\inbook Number Theory New York 1985--88
\eds D. \& G. Chudnovsky, H. Cohn, M. Nathanson 
\bookinfo
Lecture Notes in Math.\vol 1383
\publ Springer-Verlag\publaddr New York\yr 1989
\pages 12--49
\endref 

\ref\key 47
\by D. V. Chudnovsky and G. V. Chudnovsky
\paper Hypergeometric and modular function identities, and new
rational approximations to and continued fraction expansions of
classical constants and functions
\inbook A Tribute to Emil Grosswald:number theory and related
analysis
\eds M. Knopp and M. Sheingorn
\bookinfo vol. 143 of Contemporary Mathematics 
\publ American Mathematical Society\publaddr Providence, RI  
\yr 1993\pages 117--162
\endref

\ref\key 48
\by L. Comtet\book Advanced Combinatorics
\publ D. Reidel Pub. Co.
\publaddr Dordrecht-Holland/Boston-U.S.A.\yr 1974
\endref

\ref\key 49
\by E. Conrad\paper A note on certain continued fraction expansions
of Laplace transforms of Dumont's bimodular Jacobi elliptic
functions\paperinfo preprint
\endref

\ref\key 50
\by E. Conrad\paper A Handbook of Jacobi Elliptic Functions
\jour Class notes \yr (1996),
\nofrills\finalinfo preprint
\endref

\ref\key 51
\by J. H. Conway and N. J. A. Sloane
\book Sphere Packings, Lattices and Groups
\bookinfo 3rd ed. [With additional contributions by E. Bannai, R.
E. Borcherds, J. Leech, S. P. Norton, A. M. Odlyzko, R. A. Parker,
L. Queen and B. B. Venkov]. Grundlehren der mathematischen
Wissenschaften
\vol 290
\publ Springer-Verlag\publaddr New York\yr 1999
\endref 

\ref\key 52
\by T. L. Curtright and C. B. Thorn
\paper Symmetry patterns in the mass
spectra of dual string models
\jour Nuclear Phys. B\vol 274\yr 1986
\pages 520--558
\endref

\ref\key 53
\by T. W. Cusick\paper Identities involving powers of persymmetric
determinants
\jour Proc. Cambridge Philos. Soc.\vol 65\yr 1969
\pages 371--376
\endref

\ref\key 54
\by H. Datta\paper On the theory of continued fractions
\jour  Proc. Edinburgh Math. Soc.\vol 34\yr 1916
\pages 109--132
\endref

\ref\key 55
\by P. Delsarte
\paper Nombres de Bell et polyn\^omes de Charlier
\jour  C. R. Acad. Sc. Paris (series A)
\vol 287\yr 1978
\pages 271--273
\endref

\ref\key 56
\by L. E. Dickson\book History of the Theory of Numbers
\bookinfo vol. 2
\publ Chelsea
\publaddr New York \yr 1966 
\endref

\ref\key 57
\by A. C. Dixon
\paper On the doubly periodic functions arising out of the curve 
$x^3 + y^3 - 3\alpha xy = 1$
\jour The Quarterly Journal of Pure and Applied Mathematics
\vol 24\yr 1890\pages 167--233
\endref

\ref\key 58
\by D. Dumont
\paper Une approche combinatoire des fonctions elliptiques de Jacobi
\jour Adv. in Math. \vol 41\yr 1981 
\pages 1--39
\endref

\ref\key 59
\by D. Dumont
\paper Pics de cycle et d\'eriv\'ees partielles
\jour S\'eminaire Lotharingien de Combinatoire
\vol 13\yr 1985 \pages B13a (19 pp)
\endref

\ref\key 60
\by D. Dumont
\paper Le param\'etrage de la courbe d'\'equation 
$x^3 + y^3 =1$ (Une introduction \'el\'ementaire aux fonctions
elliptiques)
\paperinfo preprint (May, 1988)
\endref

\ref\key 61
\by F. J. Dyson\paper Missed opportunities\jour Bull.
Amer. Math. Soc.\vol 78\yr 1972\pages 635--653
\endref 

\ref\key  62
\by R. Ehrenborg 
\paper The Hankel determinant of exponential polynomials
\jour Amer. Math. Monthly
\vol 107\yr 2000
\pages 557--560
\endref

\ref\key 63
\by A. Erd\'elyi (with A. Magnus, F. Oberhettinger, and
F. Tricomi)
\paper Higher Transcendental Functions
\inbook Bateman Manuscript Project\vol II
\ed A. Erd\'elyi \publ McGraw-Hill Book Co.
\publaddr New York\yr 1953
\moreref
\publ reissued by Robert E. Krieger Pub. Co. 
\publaddr Malabar, Florida, 1981, 1985
\endref

\ref\key 64
\by T. Estermann\paper On the representations of a number as a sum
of squares\jour Acta Arith. \vol 2\yr 1936\pages 47--79
\endref

\ref\key 65
\by L. Euler
\paper De fractionibus continuis dissertatio 
\jour Comm. Acad. Sci. Imp. St. P\'etersbourg
\vol 9\yr 1737\pages 98--137
\moreref 
reprinted in Works. 1911--. Leonhardi Euleri Opera Omnia 
(F. Rudio, A. Krazer, P. Stackel, eds.)
\inbook Ser. I\vol 14 (C. Boehm and G. Faber, eds.)  
\publ B. G. Teubner
\publaddr Lipsiae 
\yr 1925\pages 187--215
\moreref see also  
\paper An essay on continued fractions
\jour Math. Systems Theory
\vol 18\yr 1985\pages 295--328
\transl Translated from the Latin by Myra F. Wyman and 
Bostwick F. Wyman
\endref  

\ref\key 66
\by L. Euler
\paper De fractionibus continuis observationes 
\jour Comm. Acad. Sci. Imp. St. P\'etersbourg
\vol 11\yr 1739\pages 32--81
\moreref 
reprinted in Works. 1911--. Leonhardi Euleri Opera Omnia 
(F. Rudio, A. Krazer, P. Stackel, eds.)
\inbook Ser. I\vol 14 (C. Boehm and G. Faber, eds.)  
\publ B. G. Teubner
\publaddr Lipsiae 
\yr 1925\pages 291--349
\endref  

\ref\key 67
\by L. Euler
\book Introductio in Analysin Infinitorum
\bookinfo vol. I 
\publ Marcum-Michaelem Bousquet
\publaddr Lausanne \yr 1748
\moreref 
reprinted in Works. 1911--. Leonhardi Euleri Opera Omnia 
(F. Rudio, A. Krazer, P. Stackel, eds.)
\inbook Ser. I\vol 8 (A. Krazer and F. Rudio, eds.)  
\publ B. G. Teubner
\publaddr Lipsiae 
\yr 1922\pages 1--392, (see bibliographie on page b*) 
\moreref see also  
\book Introduction to Analysis of the Infinite:Book I
\publ Springer-Verlag
\publaddr New York\yr 1988
\transl Translated from the Latin by John D. Blanton
\endref 

\ref\key 68
\by L. Euler
\paper De transformatione serierum in fractiones continuas:ubi
simul haec theoria non mediocriter amplificatur
\inbook Opuscula Analytica, t. ii
\publaddr Petropoli: Typis Academiae Imperialis Scientiarum, 
(1783--1785)
\yr 1785\pages 138--177
\moreref 
reprinted in Works. 1911--. Leonhardi Euleri Opera Omnia 
(F. Rudio, A. Krazer, A. Speiser, L. G. du Pasquier, eds.)
\inbook Ser. I\vol 15 (G. Faber, ed.)  
\publ B. G. Teubner
\publaddr Lipsiae 
\yr 1927\pages 661--700 
\endref  

\ref\key 69
\by P. Flajolet
\paper Combinatorial aspects of continued fractions
\jour Discrete Math.\vol 32\yr 1980 
\pages 125--161
\endref

\ref\key 70
\by P. Flajolet
\paper On congruences and continued fractions for some classical
combinatorial quantities
\jour Discrete Math.\vol 41\yr 1982 
\pages 145--153
\endref

\ref\key 71
\by P. Flajolet and J. Fran\c con
\paper Elliptic functions, continued fractions and doubled
permutations
\jour  European J. Combin.\vol 10\yr 1989 
\pages 235--241
\endref

\ref\key 72
\by F. G. Frobenius
\paper \"Uber Relationen zwischen den N\"aherungsbr\"uchen von
Potenzreihen
\jour J. Reine Angew. Math. \vol 90\yr 1881 
\pages 1--17
\moreref 
reprinted in Frobenius' Gesammelte Abhandlungen  
(J.-P. Serre, ed.)
\inbook\vol 2
\publ Springer-Verlag\publaddr Berlin
\yr 1968\pages 47--63 
\endref  

\ref\key 73
\by F. G. Frobenius and L. Stickelberger
\paper Zur Theorie der elliptischen Functionen
\jour J. Reine Angew. Math. \vol 83\yr 1877 
\pages 175--179
\moreref 
reprinted in Frobenius' Gesammelte Abhandlungen  
(J.-P. Serre, ed.)
\inbook\vol 1
\publ Springer-Verlag\publaddr Berlin
\yr 1968\pages 335--339 
\endref  

\ref\key 74
\by F. G. Frobenius and L. Stickelberger
\paper \"Uber die Addition und Multiplication der elliptischen
Functionen
\jour J. Reine Angew. Math. \vol 88\yr 1880 
\pages 146--184
\moreref 
reprinted in Frobenius' Gesammelte Abhandlungen  
(J.-P. Serre, ed.)
\inbook\vol 1
\publ Springer-Verlag\publaddr Berlin
\yr 1968\pages 612--650 
\endref  

\ref\key 75
\by M. Fulmek and C. Krattenthaler
\paper The number of rhombus
tilings of a symmetric hexagon which contain a fixed rhombus on the
symmetry axis, II
\jour European J. Combin.\vol 21\yr 2000
\pages 601--640
\endref

\ref\key 76
\by H. Garland\paper Dedekind's $\eta$-function and the cohomology
of infinite dimensional Lie algebras
\jour Proc. Nat. Acad. Sci., U.S.A.\vol 72\yr 1975
\pages 2493--2495
\endref

\ref\key 77
\by H. Garland and J. Lepowsky\paper Lie algebra homology and the
Macdonald-Kac formulas
\jour Invent. Math.\vol 34\yr 1976
\pages 37--76
\endref

\ref\key 78
\by F. Garvan
\moreref Private Communication, March 1997
\endref

\ref\key 79
\by G. Gasper and M. Rahman\paper Basic 
Hypergeometric Series\inbook Encyclopedia of Mathematics and Its
Applications\vol 35\ed G.-C. Rota\publ Cambridge University
Press\publaddr Cambridge\yr 1990
\endref

\ref\key 80
\by J. Geronimus\paper On some persymmetric determinants 
\jour  Proc. Roy. Soc. Edinburgh\vol 50\yr 1930 
\pages 304--309
\endref

\ref\key 81
\by J. Geronimus \paper On some persymmetric determinants formed by
the polynomials of M. Appell
\jour  J. London Math. Soc.\vol 6\yr 1931
\pages 55--59
\endref

\ref\key 82
\by I. Gessel and G. Viennot
\paper Binomial determinants, paths, and hook length formulae
\jour Adv. in Math.
\vol 58\yr 1985\pages 300--321
\endref 

\ref\key 83
\by F. Gesztesy and R. Weikard
\paper Elliptic algebro-geometric solutions of the KdV and AKNS 
hierarchies - an analytic approach
\jour Bull. Amer. Math. Soc. (N.S.)
\vol 35\yr 1998 \pages 271--317
\endref

\ref\key 84
\by J. W. L. Glaisher
\paper On the square of the series in which the coefficients are
the sums of the divisors of the exponents
\jour Mess. Math., New Series\vol 14\yr 1884-5 
\pages 156--163
\moreref 
reprinted in J. W. L. Glaisher
\inbook Mathematical Papers, Chiefly Connected with the $q$-series
in Elliptic Functions (1883--1885)
\publ Cambridge\publaddr  W. Metcalfe and Son, Trinity Street
\yr 1885\pages 371--379 
\endref  

\ref\key 85
\by J. W. L. Glaisher\paper On the numbers of representations 
of a number as a sum of $2r$ squares, where $2r$ does not exceed 
eighteen\jour  Proc. London Math. Soc. (2)\vol 5\yr 1907
\pages 479--490 
\endref

\ref\key 86
\by J. W. L. Glaisher\paper On the representations 
of a number as the sum of two, four, six, eight, ten, and twelve
squares\jour  Quart. J. Pure and Appl. Math.  Oxford\vol 38\yr
1907\pages 1--62 
\endref

\ref\key 87
\by J. W. L. Glaisher\paper On the representations of a
number as the sum of fourteen and sixteen squares\jour 
Quart. J. Pure and Appl. Math.  Oxford\vol 38\yr 1907\pages 178--236
\endref

\ref\key 88
\by J. W. L. Glaisher\paper On the representations of a
number as the sum of eighteen squares\jour Quart. J. Pure and 
Appl. Math. Oxford\vol 38\yr 1907\pages 289--351
\endref

\ref\key 89
\by M. L. Glasser
\moreref Private Communication, April 1996
\endref

\ref\key 90
\by M. L. Glasser and I. J. Zucker
\paper Lattice Sums
\inbook vol. 5 of Theoretical Chemistry:Advances and Perspectives
\eds H. Eyring and D. Henderson
\publ Academic Press\publaddr New York 
\yr 1980\pages 67--139
\endref

\ref\key 91
\by H. W. Gould\paper Explicit formulas for Bernoulli numbers
\jour  Amer. Math. Monthly\vol 79\yr 1972
\pages 44--51
\endref

\ref\key 92
\by I. P. Goulden and D. M. Jackson\book Combinatorial Enumeration
\publ John Wiley \& Sons
\publaddr New York\yr 1983 
\endref

\ref\key 93
\by I. S. Gradshteyn and I. M. Ryzhik
\book Table of Integrals, Series, and Products
\bookinfo 4th ed. \publ Academic Press
\publaddr San Diego \yr 1980 
\transl Translated from the Russian by Scripta Technica, Inc., and
edited by A. Jeffrey
\endref

\ref\key 94
\by E. Grosswald\book Representations of integers 
as sums of squares\publ Springer-Verlag\publaddr New York\yr
1985
\endref

\ref\key 95
\by K.-B. Gundlach\paper On the representation of a number as a sum
of squares\jour Glasgow Math. J.\vol 19\yr 1978\pages 173--197
\endref

\ref\key 96
\by R. A. Gustafson\paper The Macdonald identities for affine root
systems of classical type and hypergeometric series very well-poised
on semi-simple Lie algebras\inbook Ramanujan International Symposium
on Analysis (Dec. 26th to 28th, 1987, Pune, India)\ed N. K. Thakare
\yr 1989\pages 187--224
\endref

\ref\key 97
\by G.-N. Han, A. Randrianarivony and J. Zeng
\paper Un autre $q$-analogue des nombres d'Euler
\jour S\'eminaire Lotharingien de Combinatoire
\vol 42\yr 1999 \pages B42e (22 pp)
\endref

\ref\key 98
\by G.-N. Han and J. Zeng
\paper $q$-Polyn\^omes de Ghandi et statistique de Denert
\jour Discrete Math.\vol 205\yr 1999\pages 119--143
\endref

\ref\key 99 
\by G. H. Hardy\paper On the representation of a number as the sum
of any number of squares, and in particular of five or seven
\jour Proc. Nat. Acad. Sci., U.S.A.\vol 4\yr 1918\pages 189--193
\endref

\ref\key 100
\by G. H. Hardy\paper On the representation of a number as the sum
of any number of squares, and in particular of five
\jour Trans. Amer. Math. Soc.\vol 21\yr 1920\pages 255--284
\endref

\ref\key 101
\by G. H. Hardy\book Ramanujan
\publ Cambridge University Press
\publaddr Cambridge\yr 1940 
\moreref
\publ reprinted by Chelsea
\publaddr New York, 1978
\moreref 
\publ reprinted by AMS Chelsea
\publaddr Providence, RI, 1999
\moreref 
Now distributed by The American Mathematical Society
\publaddr Providence, RI
\endref  


\ref\key 102
\by G. H. Hardy and E. M. Wright\book An Introduction 
to the Theory of Numbers\bookinfo 5th ed.\publ Oxford University
Press\publaddr Oxford\yr 1979
\endref

\ref\key 103
\by J. B. H. Heilermann
\paper De transformatione serierum in fractiones continuas
\paperinfo Dr. Phil. Dissertation
\publ Royal Academy of M\"unster, 1845
\endref

\ref\key 104
\by J. B. H. Heilermann \paper \"Uber die Verwandlung der Reihen in
Kettenbr\"uche 
\jour J. Reine Angew. Math. \vol 33\yr 1846 
\pages 174--188
\endref

\ref\key 105
\by H. Helfgott and I. M. Gessel
\paper Enumeration of tilings of diamonds and hexagons with
defects
\jour  Electron. J. Combin.\vol 6\yr 1999
\pages \#R16, 26~pp
\endref

\ref\key 106
\by E. Hendriksen and H. Van Rossum
\paper Orthogonal Moments
\jour  Rocky Mountain J. Math.\vol 21\yr 1991
\pages 319--330
\endref

\ref\key 107
\by L. K. Hua\book Introduction to Number Theory
\publ Springer-Verlag
\publaddr New York \yr 1982
\endref

\ref\key 108
\by J. G. Huard, Z. M. Ou, B. K. Spearman, K. S. Williams
\paper Elementary evaluation of certain convolution sums involving
divisor functions
\paperinfo preprint (2-9-00)
\endref

\ref\key 109
\by M. E. H. Ismail, J. Letessier, G. Valent, and J. Wimp
\paper Two families of associated Wilson polynomials
\jour Canad. J. Math.\vol 42\yr 1990\pages 659--695
\endref

\ref\key 110
\by M. E. H. Ismail and D. R. Masson
\paper Generalized orthogonality and continued fractions
\jour J. Approx. Theory\vol 83\yr 1995\pages 1--40
\endref

\ref\key 111
\by M. E. H. Ismail and D. R. Masson
\paper Some continued fractions related to elliptic functions
\inbook Continued Fractions: from analytic number theory to
constructive approximation (Columbia, MO, 1998)
\eds B. C. Berndt and F. Gesztesy
\bookinfo vol. 236 of Contemporary Mathematics 
\publ American Mathematical Society\publaddr Providence, RI  
\yr 1999\pages 149--166
\endref

\ref\key 112
\by M. E. H. Ismail and M. Rahman
\paper The associated Askey-Wilson polynomials
\jour Trans. Amer. Math. Soc.\vol 328\yr 1991\pages 201--237
\endref

\ref\key 113
\by M. E. H. Ismail and D. Stanton
\paper Classical orthogonal polynomials as moments
\jour Canad. J. Math.\vol 49\yr 1997\pages 520--542
\endref

\ref\key 114
\by M. E. H. Ismail and D. Stanton
\paper More orthogonal polynomials as moments
\inbook Mathematical Essays in Honor of Gian-Carlo Rota 
(Cambridge, MA, 1996)
\eds B.~E.~Sagan and R.~P.~Stanley
\bookinfo vol. 161 of Progress in Mathematics 
\publ Birkh\"auser Boston, Inc.\publaddr Boston, MA  
\yr 1998\pages 377--396
\endref

\ref\key 115
\by M. E. H. Ismail and G. Valent
\paper On a family of orthogonal polynomials related to elliptic
functions
\jour Illinois J. Math. \vol 42\yr 1998\pages 294--312
\endref

\ref\key 116
\by M. E. H. Ismail, G. Valent and G. Yoon
\paper Some orthogonal polynomials related to elliptic functions
\jour J. Approx. Theory, to appear
\endref

\ref\key 117
\by C. G. J. Jacobi\paper Fundamenta Nova Theoriae
Functionum Ellipticarum\jour Regiomonti. Sumptibus
fratrum Borntr\"ager, 1829
\moreref reprinted in Jacobi's Gesammelte Werke
\inbook\vol 1\publ Reimer\publaddr Berlin 
\yr 1881--1891\pages 49--239 
\moreref 
\publ reprinted by Chelsea
\publaddr New York, 1969
\moreref 
Now distributed by The American Mathematical Society
\publaddr Providence, RI
\endref  

\ref\key 118
\by N. Jacobson\book Basic Algebra I
\publ W. H. Freeman and Co.
\publaddr San Francisco, CA\yr 1974 
\endref

\ref\key 119
\by W. B. Jones and W. J. Thron\paper Continued Fractions:Analytic
Theory and Applications  
\inbook Encyclopedia of Mathematics and Its Applications\vol 11
\ed G.-C. Rota\publ Addison-Wesley\publaddr London
\yr 1980
\moreref
\publ Now distributed by Cambridge University Press  
\publaddr Cambridge
\endref

\ref\key 120
\by V. G. Kac and M. Wakimoto
\paper Integrable highest weight modules over affine superalgebras
and number theory
\inbook Lie theory and geometry. In honor of Bertram Kostant
\eds Jean-Luc Brylinski, Ranee Brylinski, Victor Guillemin and 
Victor Kac
\bookinfo vol. 123 of Progress in Mathematics 
\publ Birkh\"auser Boston, Inc.\publaddr Boston, MA  
\yr 1994\pages 415--456
\endref

\ref\key 121
\by V. G. Kac and M. Wakimoto
\paper Integrable highest weight modules over affine superalgebras
and Appell's function
\paperinfo preprint arXiv:math-ph/0006007
\endref

\ref\key 122
\by S. Karlin and G. Szeg{\H o}
\paper On certain determinants whose elements are orthogonal 
polynomials
\jour J. Analyse Math.\vol 8\yr 1961\pages 1--157
\moreref\inbook reprinted in Gabor Szeg{\H o}:Collected Papers Vol. 3 
(R. Askey, ed.)
\publ Birkh\"auser Boston, Inc.\publaddr Boston, MA
\yr 1982\pages 603--762 
\endref

\ref\key 123
\by Marvin I. Knopp
\paper On powers of the theta-function greater than the eighth
\jour Acta Arith. \vol 46\yr 1986\pages 271--283
\endref

\ref\key 124
\by R. Koekoek and R. F. Swarttouw
\book The Askey--scheme of hypergeometric orthogonal polynomials and 
its $q$-analogue
\publ TU Delft
\publaddr The Netherlands\yr 1998;~~on the WWW:
{\tt ftp://ftp.twi.tudelft.nl/TWI/publications/
tech-reports/1998/DUT-TWI-98-17.ps.gz}
\endref

\ref\key 125
\by C. Krattenthaler
\paper Advanced determinant calculus
\jour S\'eminaire Lotharingien de Combinatoire
\vol 42\yr 1999 \pages B42q (67 pp)
\endref

\ref\key 126
\by E. Kr\"atzel\paper \"Uber die Anzahl der Darstellungen von
nat\"urlichen Zahlen als Summe von $4k$ Quadraten
\jour Wiss. Z. Friedrich-Schiller-Univ. Jena \vol 10 
\yr 1960/61 \pages 33--37
\endref 

\ref\key 127
\by E. Kr\"atzel\paper \"Uber die Anzahl der Darstellungen von
nat\"urlichen Zahlen als Summe von $4k+2$ Quadraten
\jour Wiss. Z. Friedrich-Schiller-Univ. Jena \vol 11 
\yr 1962 \pages 115--120
\endref 

\ref\key 128
\by D. B. Lahiri\paper On a type of series involving the partition
function with applications to certain congruence relations
\jour  Bull. Calcutta Math. Soc.\vol 38\yr 1946
\pages 125--132
\endref

\ref\key 129
\by D. B. Lahiri\paper On Ramanujan's function $\tau (n)$ and the
divisor function $\sigma_k(n)$ - I
\jour  Bull. Calcutta Math. Soc.\vol 38\yr 1946
\pages 193--206
\endref

\ref\key 130
\by D. B. Lahiri\paper On Ramanujan's function $\tau (n)$ and the
divisor function $\sigma_k(n)$ - II
\jour  Bull. Calcutta Math. Soc.\vol 39\yr 1947
\pages 33--52
\endref

\ref\key 131
\by D. B. Lahiri\paper Identities connecting the partition, divisor
and Ramanujan's functions
\jour Proc. Nat. Inst. Sci. India\vol 34A\yr 1968
\pages 96--103
\endref

\ref\key 132
\by D. B. Lahiri\paper Some arithmetical identities for Ramanujan's
and divisor functions
\jour Bull. Austral. Math. Soc.\vol 1\yr 1969
\pages 307--314
\endref

\ref\key 133
\by A. Lascoux
\paper Inversion des matrices de Hankel
\jour Linear Algebra Appl.
\vol 129\yr 1990
\pages 77--102
\endref

\ref\key 134
\by D. F. Lawden
\paper Elliptic Functions and Applications
\inbook vol. 80 of Applied Mathematical Sciences
\publ Springer-Verlag\publaddr New York 
\yr 1989
\endref

\ref\key 135 
\by B. Leclerc 
\paper On identities satisfied by minors of a matrix 
\jour Adv. in Math.\vol 100\yr 1993
\pages 101--132
\endref

\ref\key 136 
\by B. Leclerc 
\paper Powers of staircase Schur functions and symmetric analogues
of Bessel polynomials 
\jour Discrete Math.\vol 153\yr 1996\pages 213--227
\endref

\ref\key 137
\by B. Leclerc
\moreref Private Communication, July 1997
\endref

\ref\key 138
\by B. Leclerc
\paper On certain formulas of Karlin and Szeg{\H o}
\jour S\'eminaire Lotharingien de Combinatoire
\vol 41\yr 1998 \pages B41d (21 pp)
\endref

\ref\key 139
\by A. M. Legendre\paper Trait\'e des Fonctions Elliptiques et
des Int\'egrales Euleriennes
\inbook t. III
\publ Huzard-Courcier \publaddr Paris 
\yr 1828\pages 133--134
\endref

\ref\key 140
\by D. H. Lehmer\paper Some functions of Ramanujan
\jour  Math. Student\vol 27\yr 1959
\pages 105--116
\endref

\ref\key 141
\by J. Lepowsky\paper Generalized Verma modules, 
loop space cohomology and Macdonald-type identities\jour Ann. Sci.
\'Ecole Norm. Sup. (4)\vol 12\yr 1979\pages 169--234
\endref

\ref\key 142
\by J. Lepowsky\paper Affine Lie algebras and combinatorial 
identities\inbook Lie Algebras and Related Topics\nofrills\bookinfo
(Rutgers Univ. Press., New Brunswick, N. J., 1981), Lecture Notes in
Math.\vol 933\publ Springer-Verlag\publaddr Berlin\yr 1982
\pages 130--156
\endref 

\ref\key 143
\by G. M. Lilly and S. C. Milne\paper The $C_{\ell}$ Bailey
Transform and Bailey Lemma\jour Constr. Approx.\vol 9\yr 1993
\pages 473--500
\endref 

\ref\key 144
\by J. Liouville
\paper Extrait d'une lettre \`a M. Besge
\jour  J. Math. Pures Appl. (2)\vol 9\yr 1864
\pages 296--298
\endref

\ref\key 145
\by D. E. Littlewood
\book The Theory of Group Characters and Matrix Representations of
Groups
\bookinfo 2nd ed. \publ Oxford University Press
\publaddr Oxford \yr 1958 
\endref

\ref\key 146
\by Z-G Liu
\paper The representation of integers as sums of twenty four squares
\paperinfo preprint (1-31-2000)
\endref

\ref\key 147
\by G. A. Lomadze\paper On the representation of numbers by sums of
squares\jour Akad. Nauk Gruzin. SSR Trudy Tbiliss. Mat. Inst. 
Razmadze\vol 16\yr 1948\pages 231--275\lang in Russian; Georgian
summary 
\endref 

\ref\key 148
\by J. S. Lomont and J. D. Brillhart
\book Elliptic Polynomials
\publ Chapman \& Hall/CRC Press
\publaddr Boca Raton, FL
\yr 2000
\endref

\ref\key 149
\by L. Lorentzen and H. Waadeland
\paper Continued Fractions With Applications
\inbook vol. 3 of Studies in Computational Mathematics
\publ North-Holland\publaddr Amsterdam 
\yr 1992
\endref

\ref\key 150
\by J. L\"utzen
\book Joseph Liouville 1809--1882:master of pure and applied
mathematics
\bookinfo
 in Studies in the History of Mathematics and Physical Sciences
\vol 15
\publ Springer-Verlag\publaddr New York\yr 1990
\endref 

\ref\key 151
\by I. G. Macdonald\paper Affine root systems and 
Dedekind's $\eta$-function\jour Invent. Math.\vol 15\yr 1972\pages
91--143
\endref

\ref\key 152
\by I. G. Macdonald\paper Some conjectures for root systems
\jour SIAM J. Math. Anal.\vol 13\yr 1982
\pages 988--1007
\endref

\ref\key 153
\by I. G. Macdonald\book Symmetric Functions and Hall 
Polynomials\bookinfo 2nd ed. \publ Oxford University Press
\publaddr Oxford \yr 1995 
\endref

\ref\key 154
\by G. B. Mathews\paper On the representation of a number as a sum
of squares\jour Proc. London Math. Soc.\vol 27\yr 1895-6\pages 55--60
\endref

\ref\key 155
\by H. McKean and V. Moll
\book Elliptic Curves:function theory, geometry, arithmetic
\publ Cambridge University Press
\publaddr Cambridge\yr 1997
\endref

\ref\key 156
\by M. L. Mehta
\book Elements of Matrix Theory
\publ Hindustan Publishing Corp.
\publaddr Delhi\yr 1977
\endref

\ref\key 157
\by M. L. Mehta
\book Matrix Theory:selected topics and useful results
\publ Les Editions de Physique
\publaddr Les Ulis, France\yr 1989
\moreref
\publ see Appendix A.5
\moreref
(In India, sold and distributed by Hindustan Publishing Corp.)
\endref

\ref\key 158
\by S. C. Milne\paper An elementary proof of the Macdonald 
identities for $A_\ell^{(1)}$\jour Adv. in Math.\vol 57\yr 1985
\pages 34--70
\endref

\ref\key 159
\by S. C. Milne\paper Basic hypergeometric series very well-poised
in $U(n)$\jour J. Math. Anal. Appl.\vol 122\yr 1987\pages 223--256
\endref

\ref\key 160
\by S. C. Milne\paper Classical partition functions and the
$U(n+1)$ Rogers-Selberg identity\jour Discrete
Math.\vol 99\yr 1992\pages 199--246
\endref

\ref\key 161
\by S. C. Milne\paper The $C_{\ell}$ Rogers-Selberg Identity
\jour  SIAM J. Math. Anal.\vol 25\yr 1994
\pages 571--595
\endref

\ref\key 162
\by S. C. Milne\paper New infinite families of exact sums of
squares formulas, Jacobi elliptic functions, and Ramanujan's tau
function
\jour Proc. Nat. Acad. Sci., U.S.A.\vol 93\yr 1996
\pages 15004--15008
\endref

\ref\key 163
\by S. C. Milne\paper Balanced ${}_3\phi_2$ summation theorems
for $U(n)$ basic hypergeometric series\jour Adv. in Math.
\vol 131\yr 1997\pages 93--187
\endref 

\ref\key 164
\by S. C. Milne\paper Hankel determinants of Eisenstein series
\inbook in Symbolic Computation, Number Theory, Special Functions, Physics and
Combinatorics (Gainesville, 1999)
\eds F\.~G\. Garvan and M\. Ismail
\bookinfo vol. 4 of Dev. Math.
\publ Kluwer Academic Publishers
\publaddr Dordrecht
\yr (2001, to appear)
\moreref
preprint arXiv:math.NT/0009130
\endref

\ref\key 165
\by S. C. Milne\paper A new formula for Ramanujan's tau function
\paperinfo in preparation
\endref

\ref\key 166
\by S. C. Milne\paper Continued fractions, Hankel determinants, 
and further identities for powers of classical theta functions
\paperinfo in preparation
\endref

\ref\key 167
\by S. C. Milne\paper Sums of squares, Schur functions, and
multiple basic hypergeometric series
\paperinfo in preparation
\endref

\ref\key 168
\by S. C. Milne and G. M. Lilly\paper The $A_\ell$ and $C_\ell$
Bailey transform and lemma\jour Bull. Amer. Math. Soc. (N.S.)
\vol 26\yr 1992 \pages 258--263
\endref

\ref\key 169
\by S. C. Milne and G. M. Lilly\paper Consequences of the $A_\ell$
and $C_\ell$ Bailey Transform and Bailey Lemma\jour Discrete
Math.\vol 139\yr 1995\pages 319--346
\endref

\ref\key 170
\by S. C. Mitra\paper On the expansion of the Weierstrassian and
Jacobian elliptic functions in powers of the argument
\jour  Bull. Calcutta Math. Soc.\vol 17\yr 1926
\pages 159--172
\endref

\ref\key 171
\by L. J. Mordell\paper On Mr Ramanujan's empirical expansions of
modular functions\jour  Proc. Cambridge Philos. Soc.\vol 19
\yr 1917\pages 117--124 
\endref

\ref\key 172
\by L. J. Mordell\paper On the representation of numbers as the sum
of $2r$ squares\jour  Quart. J. Pure and Appl. Math. Oxford\vol 48
\yr 1917\pages 93--104 
\endref

\ref\key 173
\by L. J. Mordell\paper On the representations of a number as a sum
of an odd number of squares\jour  Trans. Cambridge Philos. Soc.
\vol 22\yr 1919\pages 361--372 
\endref

\ref\key 174
\by T. Muir
\paper New general formulae for the transformation of infinite
series into continued fractions
\jour  Trans. Roy. Soc. Edinburgh\vol 27
\yr 1872--1876; see Part IV. 1875--76 
\pages 467--471
\endref

\ref\key 175
\by T. Muir
\paper On the transformation of Gauss' hypergeometric series into a
continued fraction
\jour Proc. London Math. Soc.\vol 7\yr 1876
\pages 112--119
\endref

\ref\key 176
\by T. Muir
\paper On Eisenstein's continued fractions
\jour  Trans. Roy. Soc. Edinburgh\vol 28
\yr 1876--1878; see Part I. 1876--1877 
\pages 135--143
\endref

\ref\key 177
\by T. Muir\book The theory of determinants in the historical order
of development
\bookinfo Vol.I (1906), Vol. II (1911), Vol. III (1920), Vol. IV
(1923)
\publ Macmillan and Co., Ltd.
\publaddr London\yr 1906 
\endref

\ref\key 178
\by T. Muir\paper The theory of persymmetric determinants in the
historical order of development up to 1860
\jour  Proc. Roy. Soc. Edinburgh\vol 30\yr 1910
\pages 407--431
\endref

\ref\key 179
\by T. Muir\paper The theory of persymmetric determinants from 1894
to 1919
\jour  Proc. Roy. Soc. Edinburgh\vol 47\yr 1926-27
\pages 11--33
\endref

\ref\key 180
\by T. Muir\book Contributions to the history of determinants
1900--1920
\publ Blackie \& Son, Ltd.
\publaddr London and Glasgow\yr 1930 
\endref

\ref\key 181
\by T. Muir\book A Treatise on the Theory of Determinants
\publ Dover Publications
\publaddr New York\yr 1960 
\endref

\ref\key 182
\by M. B. Nathanson
\paper Elementary Methods in Number Theory
\inbook vol. 195 of Graduate Texts in Mathematics
\publ Springer-Verlag\publaddr New York 
\yr 2000
\endref

\ref\key 183
\by K. Ono, S. Robins, P. T. Wahl 
\paper On the representation of integers as sums of triangular
numbers
\jour  Aequationes Math.\vol 50\yr 1995
\pages 73--94
\endref

\ref\key 184
\by O. Perron\book Die Lehre von den Kettenbr\"uchen
\bookinfo 2nd ed. 
\publ B. G. Teubner
\publaddr Leipzig and Berlin \yr 1929 
\moreref 
\publ reprinted by Chelsea
\publaddr New York, 1950
\endref

\ref\key 185
\by Von K. Petr
\paper \"Uber die Anzahl der Darstellungen einer Zahl als Summe von
zehn und zw\"olf Quadraten
\jour Archiv Math. Phys., (3)\vol 11\yr 1907
\pages 83--85
\endref

\ref\key 186
\by B. van der Pol\paper The representation of numbers as sums of
eight, sixteen and twenty-four squares\jour Nederl. Akad. 
Wetensch. Proc. Ser. A \vol 57\moreref\jour = Nederl. Akad. 
Wetensch. Indag. Math.\vol 16\yr 1954\pages 349--361 
\endref 

\ref\key 187
\by G. Prasad
\book An Introduction to the Theory of Elliptic Functions and
Higher Transcendentals
\publ University of Calcutta
\yr 1928
\endref

\ref\key 188
\by H. Rademacher
\paper Topics in Analytic Number Theory
\inbook vol. 169 of Grundlehren Math. Wiss.
\publ Springer-Verlag\publaddr New York 
\yr 1973
\endref

\ref\key 189
\by Ch. Radoux
\paper Calcul effectif de certains d\'eterminants de Hankel
\jour  Bull. Soc. Math. Belg. S\'er B
\vol 31\yr 1979
\pages 49--55
\endref

\ref\key 190
\by Ch. Radoux
\paper D\'eterminant de Hankel construit sur les polyn\^omes 
de H\'ermite
\jour  Ann. Soc. Sci. Bruxelles S\'er.~I \vol 104\yr 1990
\pages 59--61
\endref

\ref\key 191
\by Ch. Radoux
\paper D\'eterminant de Hankel construit sur des polyn\^omes li\'es
aux nombres de d\'erangements
\jour  European J. Combin.\vol 12\yr 1991
\pages 327--329
\endref

\ref\key 192
\by Ch. Radoux
\paper D\'eterminants de Hankel et th\'eor\`eme de Sylvester
\inbook Actes de la 28e session du S\'eminaire Lotharingien de
Combinatoire, publication de l'I.R.M.A. no~498/S--28
\publaddr Strasbourg
\yr 1992\pages 115--122
\endref

\ref\key 193
\by S. Ramanujan\paper On certain arithmetical functions\jour Trans.
Cambridge Philos. Soc. \vol 22\yr 1916\pages 159--184
\moreref\inbook reprinted in Collected Papers of Srinivasa Ramanujan
\publ Chelsea\publaddr New York
\yr 1962\pages 136--162 
\moreref 
\publ reprinted by AMS Chelsea
\publaddr Providence, RI, 2000
\moreref 
Now distributed by The American Mathematical Society
\publaddr Providence, RI
\endref  

\ref\key 194
\by S. Ramanujan
\book The Lost Notebook and Other Unpublished Papers
\publ Narosa
\publaddr New Delhi
\yr 1988
\endref

\ref\key 195
\by A. Randrianarivony
\book Fractions continues, combinatoire et
extensions de nombres classiques
\bookinfo Ph.D. Thesis\publ  Univ. Louis Pasteur
\publaddr Strasbourg,  France\yr 1994
\endref

\ref\key 196
\by A. Randrianarivony
\paper Fractions continues, $q$-nombres de Catalan et 
$q$-polyn\^omes de Genocchi
\jour  European J. Combin. 
\vol 18 \yr 1997 \pages 75--92
\endref

\ref\key 197
\by A. Randrianarivony
\paper $q,p$-analogue des nombres de Catalan
\jour  Discrete Math. 
\vol 178  \yr 1998  \pages 199--211
\endref

\ref\key 198
\by A. Randrianarivony and J. A. Zeng
\paper Extension of Euler numbers
and records of up-down permutations
\jour  J. Combin. Theory Ser. A  
\vol 68  \yr 1994  \pages 86--99
\endref

\ref\key 199
\by A. Randrianarivony and J. Zeng
\paper A family of polynomials 
interpolating several classical series of numbers
\jour  Adv. in Appl. Math.  
\vol 17   \yr 1996   \pages 1--26
\endref

\ref\key 200
\by R. A. Rankin\paper On the representations of a number as 
a sum of squares and certain related identities
\jour Proc. Cambridge Philos. Soc.\vol 41\yr 1945 
\pages 1--11
\endref

\ref\key 201
\by R. A. Rankin\paper On the representation of a number as the sum
of any number of squares, and in particular of twenty
\jour Acta Arith. \vol 7\yr 1962\pages 399--407
\endref

\ref\key 202
\by R. A. Rankin\paper Sums of squares and cusp forms
\jour Amer. J. Math.\vol 87\yr 1965\pages 857--860
\endref

\ref\key 203
\by R. A. Rankin\book Modular Forms and Functions
\publ Cambridge University Press
\publaddr Cambridge \yr 1977 
\endref

\ref\key 204
\by D. Redmond\book Number Theory:An Introduction
\publ Marcel Dekker
\publaddr New York\yr 1996 
\endref

\ref\key  205
\by D. P. Robbins\paper Solution to problem {\bf 10387*}
\jour Amer. Math. Monthly
\vol 104\yr 1997
\pages 366--367  
\endref

\ref\key 206
\by L. J. Rogers\paper On the representation of certain
asymptotic series as convergent continued fractions\jour Proc. 
London Math. Soc (2)\vol 4\yr 1907\pages 72--89
\endref

\ref\key 207
\by A. Schett
\paper Properties of the Taylor series expansion coefficients of the
Jacobian elliptic functions
\jour Math. Comp. \vol 30\yr 1976 
\pages 143--147, with microfiche supplement
\moreref (See also: ``{\it Corrigendum}''.
\vol 31\yr 1977\pages 330)
\endref

\ref\key 208
\by A. Schett
\paper Recurrence formula of the Taylor series expansion
coefficients of the Jacobian elliptic functions
\jour Math. Comp. \vol 31\yr 1977
\pages 1003--1005, with microfiche supplement
\endref

\ref\key 209
\by W. Seidel\paper Note on a persymmetric determinant
\jour Quart. J. Math., Oxford Ser. (2)\vol 4\yr 1953
\pages 150--151
\endref

\ref\key 210
\by J.-P. Serre
\paper A Course in Arithmetic
\inbook vol. 7 of Graduate Texts in Mathematics
\publ Springer-Verlag\publaddr New York 
\yr 1973
\endref

\ref\key 211
\by W. Sierpinski
\paper Wz\'or analityczny na pewna funkcje liczbowa (Une formule 
analytique pour une fonction num\'erique)
\jour Wiadomo\'sci Matematyczne Warszawa
\vol 11\yr 1907\pages 225--231\lang in Polish
\endref

\ref\key 212
\by H. J. S. Smith\book Report on the Theory of Numbers
\publ Part VI, [Report of the British Association for 1865, pp.
322--375] 
\publaddr \yr 1894 
\moreref\inbook reprinted in The Collected Mathematical Papers of 
H. J. S. Smith, vol. 1
\publ (Ed. J. W. L. Glaisher)\publaddr 
\yr 1894\pages 306--311 
\moreref
\publ reprinted by Chelsea
\publaddr New York, 1965
\endref

\ref\key 213
\by H. J. S. Smith
\paper On the orders and genera of quadratic forms containing more
than 3 indeterminates
\jour Proc. Roy. Soc. London
\vol 16\yr 1867\pages 197--208
\moreref\inbook reprinted in The Collected Mathematical Papers of 
H. J. S. Smith, vol. 1
\publ (Ed. J. W. L. Glaisher)\publaddr 
\yr 1894\pages 510--523 
\moreref
\publ reprinted by Chelsea
\publaddr New York, 1965
\endref

\ref\key 214
\by R. P. Stanley\book Enumerative Combinatorics
\bookinfo Vol. I
\publ Wadsworth \& Brooks Cole
\publaddr Belmont, CA \yr 1986 
\endref

\ref\key 215
\by M. A. Stern
\paper Theorie der Kettenbr\"uche und ihre Anwendung
\jour J. Reine Angew. Math. \vol 10\yr 1833
\pages 1--22; 154--166; 241--274; 364--376
\endref  

\ref\key 216
\by M. A. Stern
\paper Theorie der Kettenbr\"uche und ihre Anwendung
\jour J. Reine Angew. Math. \vol 11\yr 1834
\pages 33--66; 142--168; 277--306; 311--350
\endref  

\ref\key 217
\by T. J. Stieltjes
\paper Sur la r\'eduction en fraction continue d'une s\'erie
proc\'edant suivant les puissances descendantes d'une variable
\jour Ann. Fac. Sci. Toulouse\vol 3\yr 1889\pages H. 1--17
\moreref reprinted in \OE uvres Compl\`etes {\bf T2}
\inbook\publ P. Noordhoff\publaddr Groningen 
\yr 1918\pages 184--200 
\moreref see also \OE uvres Compl\`etes = Collected Papers, 
Volume II, (Ed. G. van Dijk)
\publ Springer-Verlag\publaddr Berlin 
\yr 1993\pages 188--204 
\endref  

\ref\key 218
\by T. J. Stieltjes
\paper Sur quelques int\'egrales d\'efinies et leur
d\'eveloppement en fractions continues
\jour Quart. J. Math.\vol 24\yr 1890\pages 370--382
\moreref reprinted in \OE uvres Compl\`etes {\bf T2}
\inbook\publ P. Noordhoff\publaddr Groningen 
\yr 1918\pages 378--391 
\moreref see also \OE uvres Compl\`etes = Collected Papers, 
Volume II, (Ed. G. van Dijk)
\publ Springer-Verlag\publaddr Berlin 
\yr 1993\pages 382--395 
\endref  

\ref\key 219
\by T. J. Stieltjes
\paper Recherches sur les fractions continues
\jour Ann. Fac. Sci. Toulouse\vol 8\yr 1894\pages J. 1--122
\moreref \vol 9\yr 1895\pages A. 1--47
\moreref reprinted in \OE uvres Compl\`etes {\bf T2}
\inbook\publ P. Noordhoff\publaddr Groningen 
\yr 1918\pages 402--566;(see pages 549--554)
\moreref see also \OE uvres Compl\`etes = Collected Papers, 
Volume II, (Ed. G. van Dijk)
\publ Springer-Verlag\publaddr Berlin 
\yr 1993\pages 406--570 
\moreref (see also pp. 609--745 for an English translation.  
Note especially pp. 728--733)
\endref

\ref\key 220
\by O. Sz\'asz
\paper \"Uber Hermitesche Formen mit rekurrierender Determinante
und \"uber rationale Polynome
\jour  Math. Z.\vol 11\yr 1921
\pages 24--57
\endref

\ref\key 221
\by G. Szeg{\H o}
\paper On an inequality of Tur\'an concerning Legendre polynomials 
\jour Bull. Amer. Math. Soc.\vol 54\yr 1948\pages 401--405
\moreref\inbook reprinted in Gabor Szeg{\H o}:Collected Papers Vol. 3, 
1945--1972 (R. Askey, ed.)
\publ Birkh\"auser Boston, Inc.\publaddr Boston, MA
\yr 1982\pages 69--73; 74--75
\endref

\ref\key 222
\by O. Taussky\paper Sums of Squares 
\jour  Amer. Math. Monthly\vol 77\yr 1970 
\pages 805--830
\endref

\ref\key 223
\by J. Touchard
\paper Sur un probl\'eme de configurations et sur les fractions
continues
\jour Canad. J. Math.\vol 4\yr 1952\pages 2--25
\endref

\ref\key 224
\by P. Tur\'an\paper On the zeros of the polynomials of Legendre
\jour \v Casopis pro P\v estov\'ani Matematiky a Fysiky
\vol 75\yr 1950 
\pages 113--122
\endref

\ref\key 225
\by H. W. Turnbull
\book The Theory of Determinants, Matrices, and Invariants
\publ Blackie and Son ltd.
\publaddr London \yr 1928 
\moreref 
\publ reprinted by Dover Publications
\publaddr New York, 1960
\endref

\ref\key 226
\by J. V. Uspensky
\paper Sur la repr\'esentation des nombres 
par les sommes des carr\'es
\jour Communications de la Soci\'et\'e 
math\'ematique de Kharkow s\'erie 2
\vol 14\yr 1913\pages 31--64\lang in Russian
\endref 

\ref\key 227
\by J. V. Uspensky\paper Note sur le nombre de repr\'esentations
des nombres par une somme d'un nombre pair de carr\'es
\jour Bulletin de l'Acad\'emie des Sciences de l'URSS. Leningrad. 
(Izvestija Akademii Nauk Sojuza Sovetskich Respublik. Leningrad.)
Serie 6
\vol 19\yr 1925\pages 647--662\lang in French
\endref 

\ref\key 228
\by J. V. Uspensky\paper On Jacobi's arithmetical theorems
concerning the simultaneous representation of numbers by two
different quadratic forms
\jour Trans. Amer. Math. Soc.
\vol 30\yr 1928\pages 385--404
\endref 

\ref\key 229
\by J. V. Uspensky and M. A. Heaslet
\book Elementary Number Theory
\publ McGraw-Hill Book Co.\publaddr New York 
\yr 1939
\endref

\ref\key 230
\by G. Valent\paper Asymptotic analysis of some associated
orthogonal polynomials connected with elliptic functions
\jour  SIAM J. Math. Anal.\vol 25\yr 1994
\pages 749--775
\endref

\ref\key 231
\by G. Valent\paper Associated Stieltjes-Carlitz polynomials and a
generalization of Heun's differential equation
\jour J. Comput. Appl. Math.\vol 57\yr 1995
\pages 293--307
\endref

\ref\key 232
\by G. Valent and W. Van Assche
\paper The impact of Stieltjes' work on continued fractions and
orthogonal polynomials: additional material 
\jour J. Comput. Appl. Math.\vol 65\yr 1995
\pages 419--447
\moreref this volume was devoted to the Proceedings of the
International Conference on Orthogonality, Moment Problems and
Continued Fractions (Delft, 1994)
\endref

\ref\key 233
\by W. Van Assche
\paper Asymptotics for orthogonal polynomials and three-term
recurrences
\inbook Orthogonal Polynomials: Theory and Practice
\ed P. Nevai
\bookinfo vol. 294 of NATO-ASI Series C: Mathematical and Physical
Sciences
\publ Kluwer Academic Publishers\publaddr Dordrecht 
\yr 1990\pages 435--462
\endref

\ref\key 234
\by W. Van Assche
\paper The impact of Stieltjes work on continued fractions and
orthogonal polynomials
\inbook vol. I of T.J. Stieltjes: \OE uvres Compl\`etes = Collected
Papers
\ed G. van Dijk
\publ Springer-Verlag\publaddr Berlin 
\yr 1993\pages 5--37
\endref

\ref\key 235
\by P. R. Vein 
\paper Persymmetric determinants. I. The derivatives of
determinants with Appell function elements
\jour Linear and Multilinear Algebra
\vol 11 \yr 1982\pages 253--265
\endref

\ref\key 236
\by P. R. Vein 
\paper Persymmetric determinants. II. Families of
distinct submatrices with nondistinct determinants
\jour Linear and Multilinear Algebra 
\vol 11 \yr 1982\pages 267--276 
\endref

\ref\key 237
\by P. R. Vein 
\paper Persymmetric determinants. III. A basic
determinant
\jour Linear and Multilinear Algebra 
\vol 11 \yr 1982\pages 305--315
\endref

\ref\key 238
\by P. R. Vein 
\paper Persymmetric determinants. IV. An alternative
form of the Yamazaki--Hori determinantal solution of the 
Ernst equation
\jour Linear and Multilinear Algebra 
\vol 12 \yr 1982/83\pages 329--339 
\endref

\ref\key 239
\by P. R. Vein 
\paper Persymmetric determinants. V. Families of
overlapping coaxial equivalent determinants
\jour Linear and Multilinear Algebra 
\vol 14 \yr 1983\pages 131--141
\endref

\ref\key 240
\by P. R. Vein and P. Dale 
\paper Determinants, their derivatives and nonlinear differential
equations
\jour J. Math. Anal. Appl.
\vol 74 \yr 1980\pages 599--634
\endref

\ref\key 241
\by B. A. Venkov
\book Elementary Number Theory
\publ Wolters-Noordhoff Publishing
\publaddr Groningen\yr 1970
\transl Translated from the Russian and edited by 
Helen Alderson (Popova)
\endref 

\ref\key 242
\by R. Vermes
\paper Hankel determinants formed from successive derivatives 
\jour  Duke Math. J.\vol 37\yr 1970 
\pages 255--259
\endref

\ref\key 243
\by G. Viennot\paper Une interpr\'etation combinatoire des
coefficients des d\'eveloppements en s\'erie enti\`ere des fonctions
elliptiques de Jacobi
\jour  J. Combin. Theory Ser. A\vol 29\yr 1980
\pages 121--133
\endref

\ref\key 244
\by G. Viennot
\paper Une th\'eorie combinatoire des polyn\^omes orthogonaux
g\'en\'eraux
\jour  Lecture Notes. publication de l'UQAM, Montr\'eal
\yr 1983
\endref

\ref\key 245
\by G. Viennot
\paper A combinatorial interpretation of the quotient-difference
algorithm
\inbook technical report No. 8611
\publ Universit\'e de Bordeaux I
\yr 1986
\endref

\ref\key 246
\by A. Z. Walfisz\paper On the representation of numbers by sums of
squares:asymptotic formulas
\jour  Uspehi Mat. Nauk (N.S.) (6)\vol 52\yr 1952
\pages 91--178\lang in Russian
\transl\nofrills English transl. in
\jour Amer. Math. Soc. Transl. (2)\vol 3\yr 1956\pages 163--248
\endref

\ref\key 247
\by H. S. Wall
\book Analytic Theory of Continued Fractions
\publ D. Van Nostrand
\publaddr New York\yr 1948 
\moreref 
\publ reprinted by Chelsea
\publaddr New York, 1973
\endref

\ref\key 248
\by J. B. Walton\paper Theta series in the Gaussian field
\jour  Duke Math. J.\vol 16\yr 1949
\pages 479--491
\endref

\ref\key 249
\by E. T. Whittaker and G. N. Watson
\book A Course of Modern Analysis\bookinfo 4th ed. 
\publ Cambridge University Press
\publaddr Cambridge \yr 1927 
\endref

\ref\key 250
\by S. Wolfram
\book The Mathematica Book
\bookinfo 4th ed. 
\publ Wolfram Media/Cambridge University Press
\publaddr Cambridge\yr 1999
\endref

\ref\key 251
\by S. Wrigge\paper Calculation of the Taylor series expansion
coefficients of the Jacobian elliptic function $\sn(x,k)$
\jour  Math. Comp.\vol 36\yr 1981
\pages 555--564
\endref

\ref\key 252
\by S. Wrigge\paper A note on the Taylor series expansion
coefficients of the Jacobian elliptic function $\sn(x,k)$
\jour  Math. Comp.\vol 37\yr 1981
\pages 495--497
\endref

\ref\key 253
\by D. Zagier
\paper Proof of a conjecture of Kac and Wakimoto
\paperinfo to appear in MRL (Mathematical Research Letters)
\endref

\ref\key 254
\by D. Zeilberger
\paper Proof of the alternating sign matrix conjecture
\jour  Electron. J. Combin.
\vol 3 
\yr 1996
\pages \#R13, 84~pp
\endref

\ref\key 255
\by D. Zeilberger
\paper Proof of the refined alternating sign matrix conjecture
\jour  New York J. Math.\vol 2\yr 1996\pages 59--68
\endref

\ref\key 256
\by J. Zeng\paper \'Enum\'erations de permutations et $J$-fractions
Continues
\jour  European J. Combin.\vol 14\yr 1993
\pages 373--382
\endref

\ref\key 257
\by J. Zeng
\paper Sur quelques propri\'etes de sym\'etrie des nombres de
Genocchi
\jour  Discrete Math.\vol 153\yr 1996
\pages 319--333
\endref

\ref\key 258
\by I. J. Zucker
\paper The summation of series of hyperbolic functions
\jour SIAM J. Math. Anal. \vol 10\yr 1979 
\pages 192--206
\endref

\endRefs
\enddocument